# On Tournament Design

Leo Fried

Advised by Professor Eric Maskin

Undergraduate Thesis in Mathematics at Harvard College

March 18th, 2024

# Contents

# 1 Tournament Design

## 1.1 Introduction

Although tournaments have been in use for as long as humans have played sports, their formal study is relatively underdeveloped, with many key questions in the field remaining open. In this thesis we aim to address five big open questions in the field, beginning by focusing on the *bracket*: a tournament format in which teams are eliminated upon a loss, games are played until only one team remains, and the matchups between game-winners are determined in advance of any games being played.

**Is there a succinct notation for describing brackets?** The space of brackets is quite large, and fully drawing out a given bracket can be time-consuming, difficult to quickly interpret, and nearly impossible to read. We introduce the *bracket signature*, a new system for compressing brackets into a single succinct lists of digits, allowing for them to be easily communicated and important properties of brackets to be quickly verified.

**Which brackets are fair?** We construct the notion of a *proper bracket*, with a host of desirable fairness properties, noting that nearly all brackets in use by leagues around the world are proper. We then develop insight into the order in which games are played and teams are eliminated in a proper bracket. Finally, we prove the *fundamental theorem of brackets*: there is exactly one proper bracket with each bracket signature.

**Which brackets are accurate?** Edwards's [10] answered this question in 1991 using a measure of accuracy he called *orderedness*, but at the time the previous two questions were still open. With properness and signatures now defined, as well as the fundamental theorem proved, we present a much simpler and more direct proof of the answer and establish a few other generalizable results along the way. Unfortunately, Edwards's Theorem is quite pessimistic about the number of ordered brackets – even the standard eight-team bracket is not ordered – leading to the next question.

**Are there other bracket-like formats that are ordered?** Hwang [12] published a proof that *reseeded brackets* are ordered, but we show that his proof was incorrect. After an analysis of a few other techniques, we conjecture that any bracket-like format that is ordered for any number of teams lacks several other key properties.

**How can brackets select runners-up?** We propose the notion of a *multibracket*, a generalization of the bracket that unifies under a single umbrella a wide variety of bracket-like formats that have been used to select runners-up (including, but not limited to, *consolation brackets*, *semibrackets*, *linear multibrackets*, *swiss systems*, *Page-McIntyre systems*, and *double-elimination*). We use this framework to prove several key results about the number of such formats, their efficiency, and their accuracy, and discuss how a tournament designer might select which of these formats to use.

Throughout this work, we analyze 16 tournament formats in use by leagues across the globe, as well as 44 constructed for analysis. We prove 42 theorems, of which 29 are novel (including one that disproves Hwang's theorem), and define 30 new terms. We hope this thesis will aid leagues in designing more effective tournaments, as well as serve as a jumping-off point for future work in the field of tournament design.

## 1.2 History of the Field

The formal theory of tournament design was born out of the study of paired comparisons, a field that began in 1927 with Thurstone's *A Law of Comparative Judgement* [18]. Thurstone was a psychologist investigating how individuals rank a collection of objects on some axis (weight, beauty, excellence, etc), while only being able to examine two of the objects at a time. The similarity to the problems of tournament design is clear, though Thurstone did not draw the connection.

In 1963, David wrote *The Method of Paired Comparisons* [9], aiming to gather all of the theory that had been developed about paired comparisons thus far, as well as several contributions of his own, into a single monograph. At the time of publication, the field was still viewed through Thurstone's psychological lens, rather than the lens we will use involving teams competing in a game or sport. Where we will say "two teams play a game," David said "a single judge must choose between two objects." Still, the formalizations of the problems are equivalent.

In the years following David's work, the field of tournament design came into its own, with various authors examining the fairness and accuracy of a wide range of tournament formats, most commonly either round robins or brackets. Much of the work at the time, however, was an analysis of what happens when a specific set of teams (or a narrow class of sets of teams) takes part in a specific tournament design, rather than anything more general.

In 1991, Edwards submitted his doctoral thesis, *The Combinatorial Theory of Single-Elimination Tournaments* [10], the single most complete analysis of brackets that has been published to date. Edwards counted and cataloged the full space of brackets, defined *orderedness*, which we will soon see is a natural and desirable property, and completely determined which brackets are ordered. Edwards's Theorem, after which Section 2.4 is named, was first proved in that thesis.

Since then, the field has become much more statistical, with most of the analysis being done by way of Monte Carlo simulations. Dabney's *Tourney Geek* [7], for example, evaluates various tournament designs based on several different statistical measurements of fairness that are estimated via simulation.

In this thesis, we return to the type of study conducted by Edwards: proofs of claims about the outcomes of various tournament designs, rather than statistical results. We will work from first principles, beginning with the definition of a game and a tournament format, constructing various classes of formats, and then examining those formats and the properties they might have. Like most studies in the field of tournament design, we are *game-ambivalent*. We abstract away the underlying game or sport: our results apply to football as well as they will to chess as well as they will to competitive rock-paper-scissors.

In this way, the theory of tournament design and the theory of sorting algorithms are quite similar: the types of questions posed in the fields are nearly identical as well. In both cases, the designer is given a list of objects (teams), may make an arbitrary number of comparisons (games), and then must output a sorting (champion). There are, however, a number of differences that separate the fields.

The first difference is that of noise. The sorting theorist works with the guarantee that if two objects are compared more than once, the comparison will give the same result every time. For this reason, the sorting theorist often finds it wasteful to compare the same pair of objects more than once. But the tournament theorist's job is much harder, as team performance is noisier. When two teams play, there is no guarantee that the better team will win, and when they play more than once, there is no guarantee that the same team will win every time.

The second difference is that of accuracy. An algorithm submitted by the sorting theorist is required to correctly sort any list of objects, otherwise it is not a sorting algorithm. The tournament theorist is under no such constraints: the noise makes such an algorithm impossible. Thus algorithms like "randomly select a winner" and "play lots of games and then declare as champion the team with the fewest wins" are valid tournament designs, even if they are (probably) not particularly good ones.

The third difference is that of priors. While the sorting theorist typically begins their algorithms with no priors on the set of objects, the tournament theorist is often given a "seeding" of teams, identifying which teams have been judged as better. This seeding can be varyingly accurate: in some cases, the tournament theorist begins their algorithms with very strong priors, while in others, the seeding provides minimal information.

The fourth difference is that of fairness. The sorting theorist is working with a set of lifeless objects whose feelings will not be hurt based on the algorithm, freeing the sorting theorist to focus only on the task of accurately sorting the objects. The tournament theorist, on the other hand, must appeal to the sense of fairness held by the competitors: in many cases, fairness is a more important consideration than accuracy.

The final difference is that of viewership. The sorting theorist works in private, comparing objects and gathering data until a sort can be published. The tournament theorist, on the other hand, works in front of an audience, who are looking not just for an accurate tournament, but for an exciting one: the NCAA College Basketball Tournament is a classic example, as we will soon see, of a tournament that is not very accurate but none the less very exciting for viewers.

Still, there is a lot of overlap between the two fields. Knuth's *The Art of Computer Programming: Sorting and Searching* [14] often used the language of teams and games when presenting various sorting algorithms. We borrow from the field of sorting in turn: in particular the concept of a *sorting network*.

Sorting networks, first patented by Armstrong, Nelson, and O'Connor [2], are sorting algorithms with the additional property that, after a comparison is made between $a$ and $b$, the rest of the algorithm is identical no matter the result, except that $a$ and $b$ are swapped. Knuth's text contains a section about the space and properties of sorting networks.

This thesis will primarily examine networked tournament formats, that is, tournament formats with this networking property. These formats are a particularly nice set of formats to study, both because the networking property turns out to be a powerful one, and because many formats used in the postseason of leagues are networked, giving our study applications to many tournaments across many sports.

## 1.3 Tournament Formats

Before we begin our study, we set the stage by defining the key terms in the field of tournament design. Let $\mathcal{T} = [t_1, ..., t_n]$ be a list of teams.

**Definition 1.3.1: Gameplay Function** **(Unattributed)**

A *gameplay function* $g$ on $\mathcal{T}$ is a nondeterministic function $g : \mathcal{T} \times \mathcal{T} \to \mathcal{T}$ with the following properties:

(a) $\mathbb{P}[g(t_i, t_j) = t_i] + \mathbb{P}[g(t_i, t_j) = t_j] = 1$.

(b) $\mathbb{P}[g(t_i, t_j) = t_i] = \mathbb{P}[g(t_j, t_i) = t_i]$.

A gameplay function represents a process in which two teams compete in a game, with one of them emerging as the winner. This model simplifies away effects like home-field advantage or teams improving over the course of a tournament: a gameplay function is fully described by a single probability for each pair of teams in the list.

**Definition 1.3.2: Playing, Winning, Losing, and Beating** **(Unattributed)**

When $g$ is queried on input $(t_i, t_j)$ we say that $t_i$ and $t_j$ *played a game*. We say that the team that got output by $g$ *won*, that the team that did not *lost*, and that the winning team *beat* the losing team.

**Definition 1.3.3: $p_{ij}$** **(Unattributed)**

$p_{ij} = \mathbb{P}[t_i \text{ beats } t_j]$.

The information in a gameplay function can be encoded into a *matchup table.*

**Definition 1.3.4: Matchup Table** **(Unattributed)**

The *matchup table* implied by a gameplay function $g$ on a list of teams $\mathcal{T}$ of length $n$ is an $n$-by-$n$ matrix $\mathcal{M}$ such that $\mathcal{M}_{ij} = p_{ij}$.

For example, let $\mathcal{T}$ = [Favorites, Rock, Paper, Scissors, Conceders], and $g$ be such that the Conceders concede every game they play; the Favorites are 70 percent favorites against Rock, Paper, and Scissors; and Rock, Paper, and Scissors matchup with each other according to the normal rules of rock-paper-scissors. Then the matchup table would look like so:

**Figure 1.3.5: The Matchup Table for $(\mathcal{T}, g)$**

| | Favorites | Rock | Paper | Scissors | Conceders |
|---|---|---|---|---|---|
| Favorites | 0.5 | 0.7 | 0.7 | 0.7 | 1.0 |
| Rock | 0.3 | 0.5 | 0.0 | 1.0 | 1.0 |
| Paper | 0.3 | 1.0 | 0.5 | 0.0 | 1.0 |
| Scissors | 0.3 | 0.0 | 1.0 | 0.5 | 1.0 |
| Conceders | 0.0 | 0.0 | 0.0 | 0.0 | 0.5 |

**Theorem 1.3.6** **(Unattributed)**

If $\mathcal{M}$ is the matchup table for some gameplay function on $\mathcal{T}$, then $\mathcal{M} + \mathcal{M}^T$ is the matrix of all ones.

*Proof.* $(\mathcal{M} + \mathcal{M}^T)_{ij} = \mathcal{M}_{ij} + \mathcal{M}_{ji} = p_{ij} + p_{ji} = 1.$ □

Theorem 1.3.6 implies that matchup tables are defined by the entries below the diagonal, so to reduce busyness we will often display only those entries.

**Figure 1.3.7: The Matchup Table for $(\mathcal{T}, g)$**

| | Favorites | Rock | Paper | Scissors | Conceders |
|---|---|---|---|---|---|
| Favorites | | | | | |
| Rock | 0.3 | | | | |
| Paper | 0.3 | 1.0 | | | |
| Scissors | 0.3 | 0.0 | 1.0 | | |
| Conceders | 0.0 | 0.0 | 0.0 | 0.0 | |

Finally, we define the tournament format.

**Definition 1.3.8: Tournament Format** **(Unattributed)**

A *tournament format* is an algorithm that takes as input a list of teams $\mathcal{T}$ and a gameplay function $g$ and outputs a ranking (potentially including ties) on $\mathcal{T}$.

We also introduce a piece of shorthand to help make notation more concise.

**Definition 1.3.9: $\mathbb{W}_{\mathcal{A}}(t, \mathcal{T})$** **(Unattributed)**

$\mathbb{W}_{\mathcal{A}}(t, \mathcal{T})$ is the probability that team $t \in \mathcal{T}$ wins tournament format $\mathcal{A}$ when it is run on the list of teams $\mathcal{T}$.

Finally, we will focus our study on the subset of tournament formats that fulfill the *network condition*, first patented as a condition for sorting algorithms by Armstrong, Nelson, and O'Connor [2].

**Definition 1.3.10: Deterministic Tournament Format** **(Unattributed)**

A tournament format is *deterministic* if it employs no randomness other than the randomness inherent in the gameplay function $g$.

This definition does not require that a deterministic tournament format always declare the same champion when presented with the same list of teams, only that it declare the same champion when presented with the same list of teams and the game results are all the same.

**Definition 1.3.11: Networked Tournament Format**
**(Armstrong, Nelson, and O'Connor, 1957)**

A tournament format is *networked* if it is deterministic, and after each game between $t_i$ and $t_j$, the rest of the format is identical no matter which team won, except that $t_i$ and $t_j$ are swapped.

# 2 Brackets

## 2.1 Bracket Signatures

**Definition 2.1.1: Bracket** (Unattributed)

A *bracket* is a networked format in which

(a) Teams don't play any games after their first loss, and

(b) Games are played until only one team has no losses, and that team is crowned champion.

We can draw a bracket as a tree-like structure in the following way.

**Figure 2.1.2: 2024 College Football Playoff**

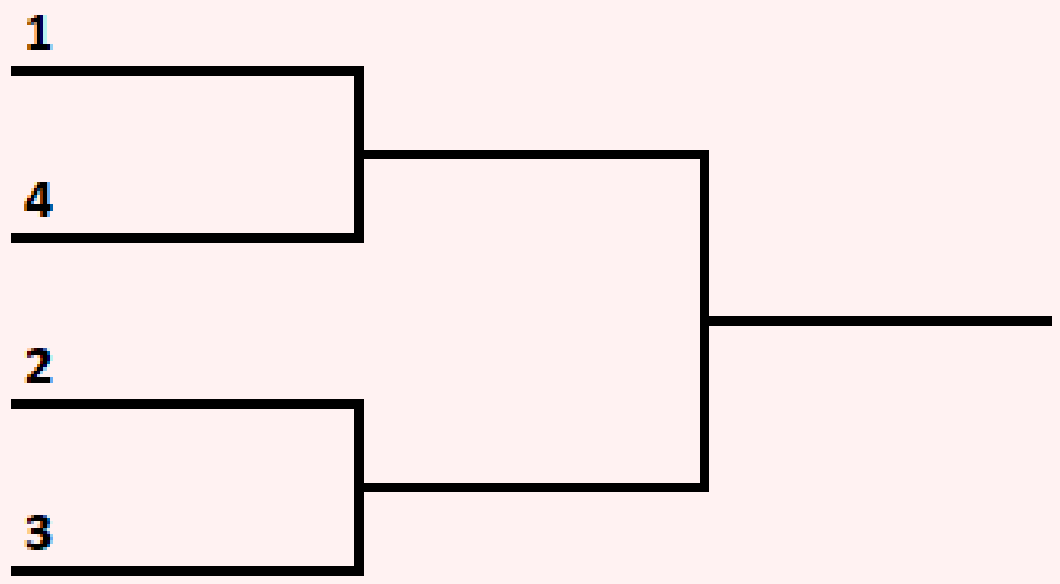

The numbers 1, 2, 3, and 4 indicate where $t_1, t_2, t_3$, and $t_4$ in $\mathcal{T}$ are placed to start. In the actual 2024 College Football Playoff [27], the list of teams $\mathcal{T}$ was [Michigan, Washington, Texas, Alabama], so the bracket was filled in like so.

**Figure 2.1.3: 2024 College Football Playoff After Team Placement**

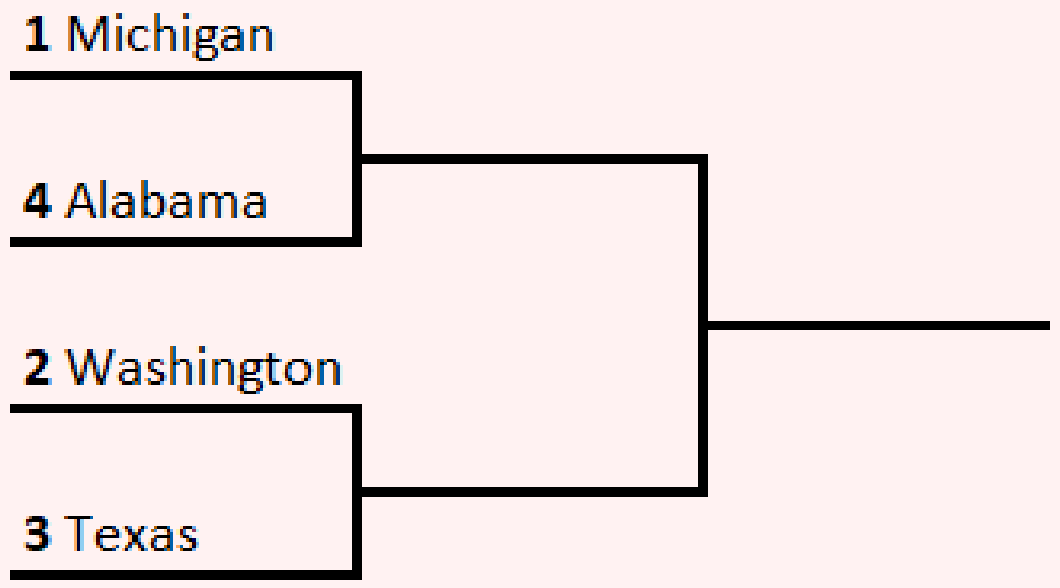

As games are played, we write the name of the winning teams on the corresponding lines. This bracket tells us that Michigan played Alabama, and Washington played Texas. Michigan and Washington won their respective games, and then Michigan beat Washington, winning the tournament.

**Figure 2.1.4: 2024 College Football Playoff After Completion**

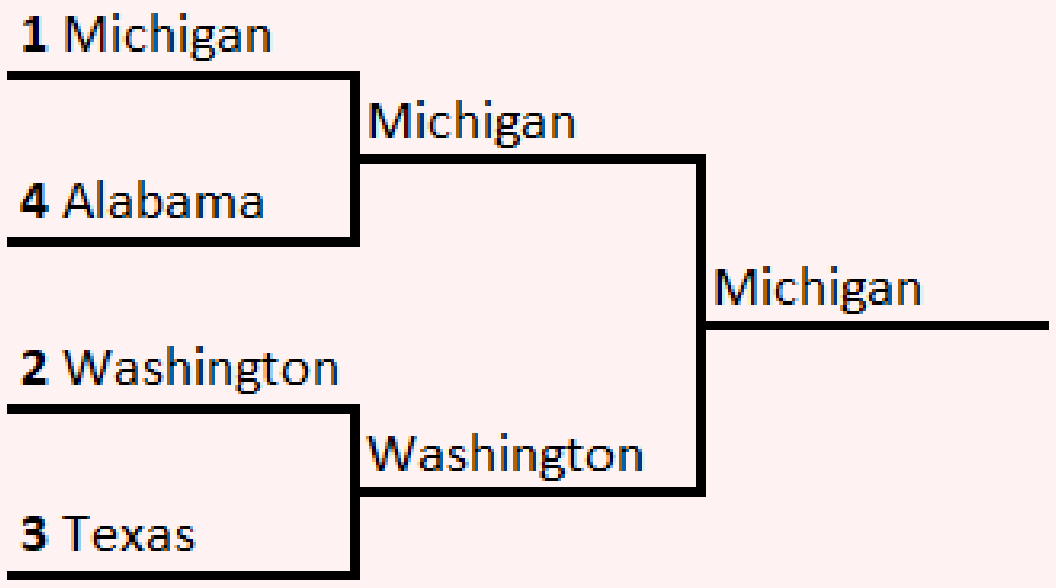

Rearranging the way the bracket is pictured, if it doesn't affect any of the matchups, does not create a new bracket. For example, Figure 2.1.5 is just another way to draw the same bracket.

**Figure 2.1.5: Alternative Drawing of the 2024 College Football Playoffs**

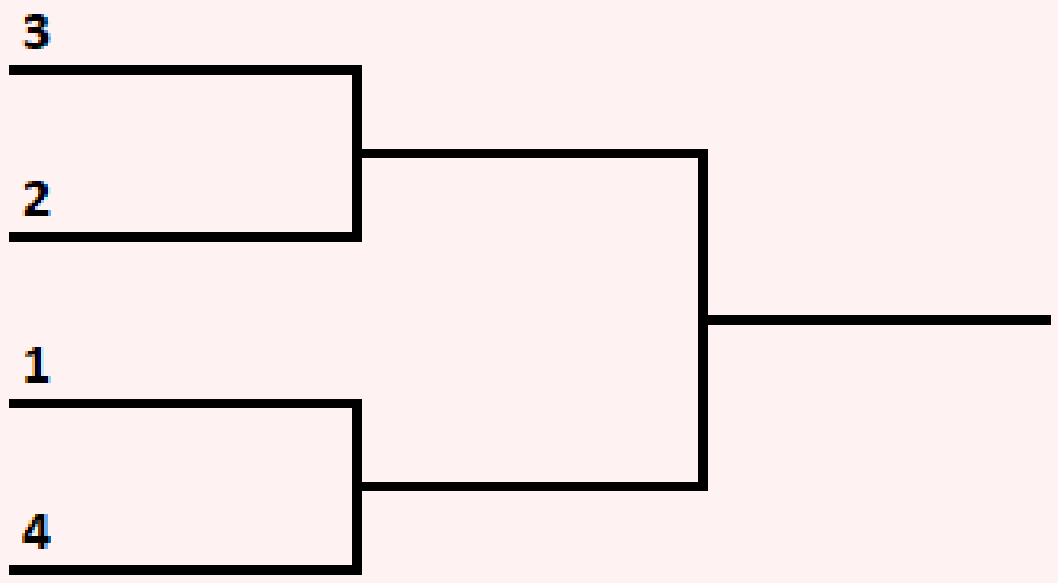

There are a few important pieces of vocabulary when it comes to brackets:

**Definition 2.1.6: Starting Line** **(Unattributed)**

A *starting line* is a line in a bracket where a team is placed before it has played any games.

**Definition 2.1.7: Round** **(Unattributed)**

A *round* is a set of games such that the winners of each of those games have the same number of games remaining to win the tournament.

**Definition 2.1.8: Bye** **(Unattributed)**

A team has a *bye* in round $r$ if it plays no games in round $r$ or before.

The 2024 College Football Playoffs had four starting lines, one for each of its participating teams, and was played over two rounds: The first round consisted of the games Michigan vs

Alabama and Washington vs Texas, and the second round was just the single Michigan vs Washington game. The 2024 College Football Playoffs had no byes.

With the terminology established, we begin by investigating the *shape* of brackets.

**Definition 2.1.9: Shape** **(Unattributed)**

The *shape* of a bracket is the tree that underlies it.

The following two brackets have the same shape.

**Figure 2.1.10: Two Brackets with the Same Shape**

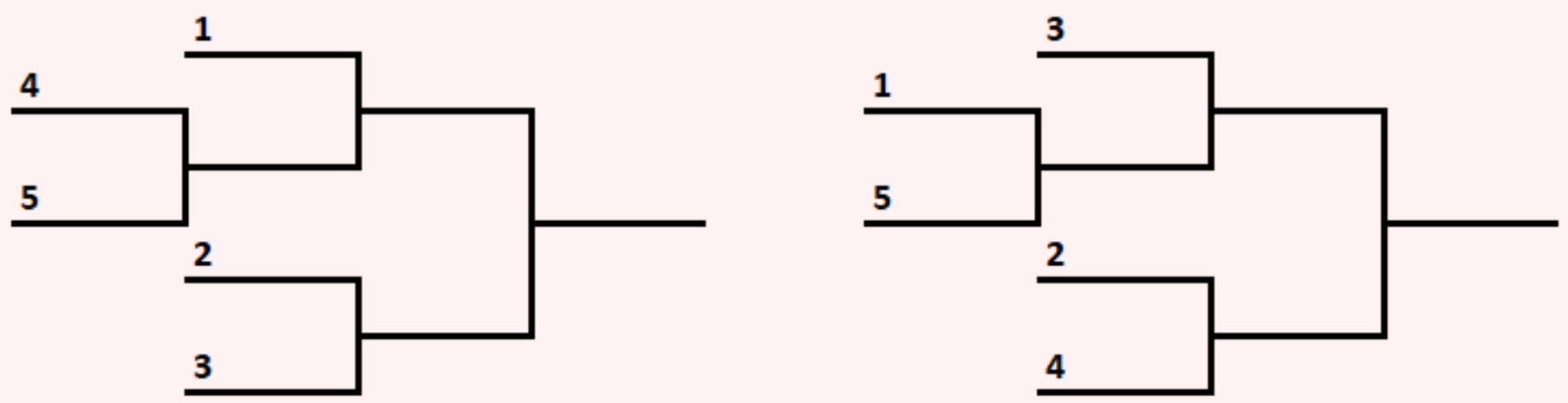

One way to describe the shape of a bracket is its signature.

**Definition 2.1.11: Bracket Signature** **(Fried, 2024)**

The *signature* of an $r$-round bracket $\mathcal{A}$ is the list $[[\mathbf{a_0}; \ldots; \mathbf{a_r}]]$ where $a_i$ is the number of teams that get $i$ byes.

The signature of a bracket is defined by its shape: the two brackets in Figure 2.1.10 have the same shape, so they also have the same signature.

The signatures of the brackets discussed in this section are shown in Figure 2.1.12. It's worth verifying the signatures we've seen so far and trying to draw brackets with the signatures we haven't yet before moving on.

**Figure 2.1.12: The Signatures of Some Brackets**

| Bracket | Signature |
|---|---|
| 2024 College Football Playoff | $[[\mathbf{4}; \mathbf{0}; \mathbf{0}]]$ |
| The brackets in Figure 2.1.10 | $[[\mathbf{2}; \mathbf{3}; \mathbf{0}; \mathbf{0}]]$ |
| The brackets in Figure 2.1.13 | $[[\mathbf{4}; \mathbf{2}; \mathbf{0}; \mathbf{0}]]$ |
| 2023 WCC Women's Basketball Tournament | $[[\mathbf{4}; \mathbf{2}; \mathbf{2}; \mathbf{2}; \mathbf{0}; \mathbf{0}]]$ |

Two brackets with the same shape must have the same signature, but the converse is not true: two brackets with different shapes can have the same signature. For example, both bracket shapes depicted in Figure 2.1.13 have the signature $[[\mathbf{4}; \mathbf{2}; \mathbf{0}; \mathbf{0}]]$.

**Figure 2.1.13: Two Shapes with the Signature $[[4; 2; 0; 0]]$**

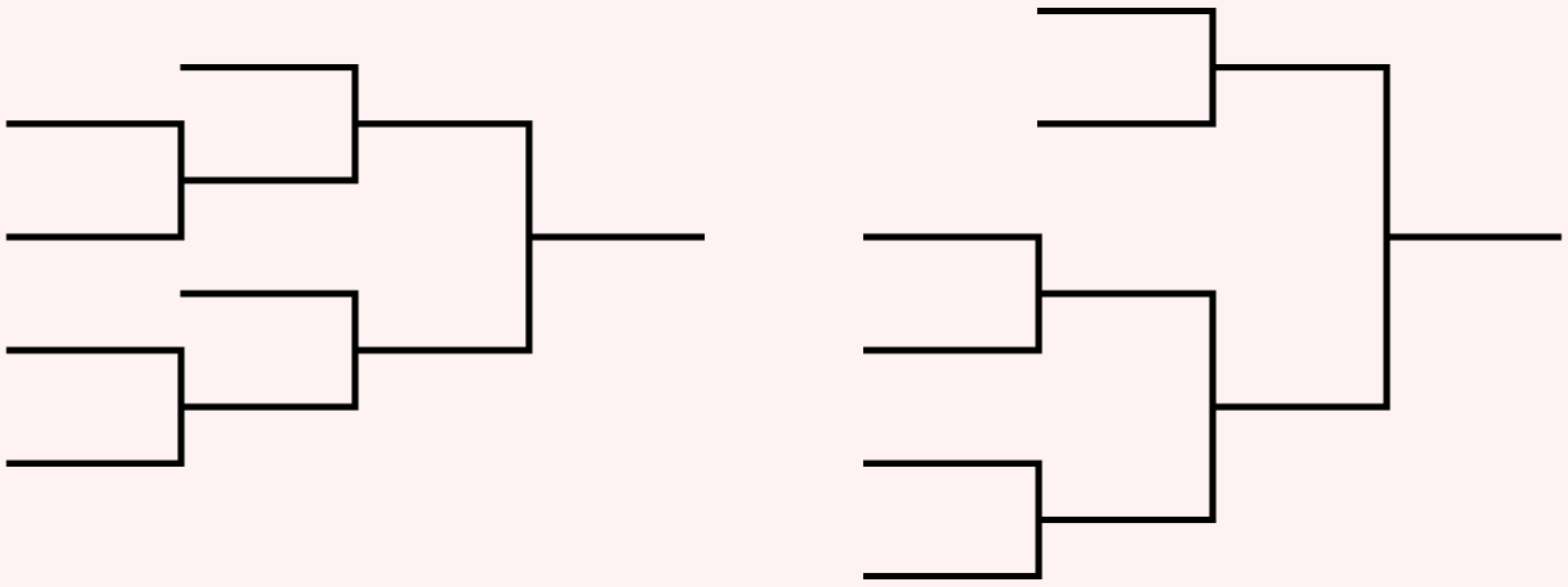

Despite this, bracket signatures are a useful way to talk about the shape of a bracket. Communicating a bracket's signature is a lot easier than communicating its shape, and much of the important information (such as how many games each team must win in order to win the tournament) is contained in the signature.

Bracket signatures have one more important property.

**Theorem 2.1.14** **(Fried, 2024)**

Let $\mathcal{A} = [[\mathbf{a_0}; ...; \mathbf{a_r}]]$ be a list of natural numbers. Then $\mathcal{A}$ is a bracket signature if and only if

$$\sum_{i=0}^{r} a_i \cdot \left(\frac{1}{2}\right)^{r-i} = 1.$$

*Proof.* Let $\mathcal{A}$ be the signature for some bracket. Assume that every game in the bracket is a coin flip, and consider each team's probability of winning the tournament. A team that has $i$ byes must win $r - i$ games to win the tournament, and so will do so with probability $\left(\frac{1}{2}\right)^{r-i}$. For each $i \in \{0, ..., r\}$, there are $a_i$ teams with $i$ byes, so

$$\sum_{i=0}^{r} a_i \cdot \left(\frac{1}{2}\right)^{r-i}$$

is the total probability of the teams winning the tournament, which is just 1.

We prove the other direction by induction on $r$. If $r = 0$, then the only list with the desired property is $[[\mathbf{1}]]$, which is the signature for the unique one-team bracket. For any other $r$, first note that $a_0$ must be even: if it were odd, then

$$\begin{aligned}
\sum_{i=0}^{r} a_i \cdot \left(\frac{1}{2}\right)^{r-i} &= \frac{1}{2^r} \cdot \sum_{i=0}^{r} a_i \cdot 2^i \\
&= \frac{1}{2^r} \cdot \left(a_0 + 2\sum_{i=1}^{r} a_i \cdot 2^{i-1}\right) \\
&= k/2^r && \text{for some odd } k \\
&\neq 1.
\end{aligned}$$

Now, consider the signature $\mathcal{B} = [[\mathbf{a_1 + a_0/2}; \mathbf{a_2}; \ldots; \mathbf{a_r}]]$. By induction, there exists a bracket with signature $\mathcal{B}$. But if we take that bracket and replace $a_0/2$ of the starting lines with no byes with a game whose winner gets placed on that line, we get a new bracket with signature $\mathcal{A}$. □

In the next few sections, we will use the language and properties of bracket signatures to describe the brackets that we work with. For now though, let's return to the 2024 College Football Playoff. The bracket used in the 2024 College Football Playoff has a special property that not all brackets have: it is *balanced*.

**Definition 2.1.15: Balanced Bracket** **(Unattributed)**

A bracket is *balanced* if none of the participating teams have byes.

The 2023 West Coast Conference Women's Basketball Tournament [21], on the other hand, is unbalanced.

**Figure 2.1.16: 2023 WCC Women's Basketball Tournament**

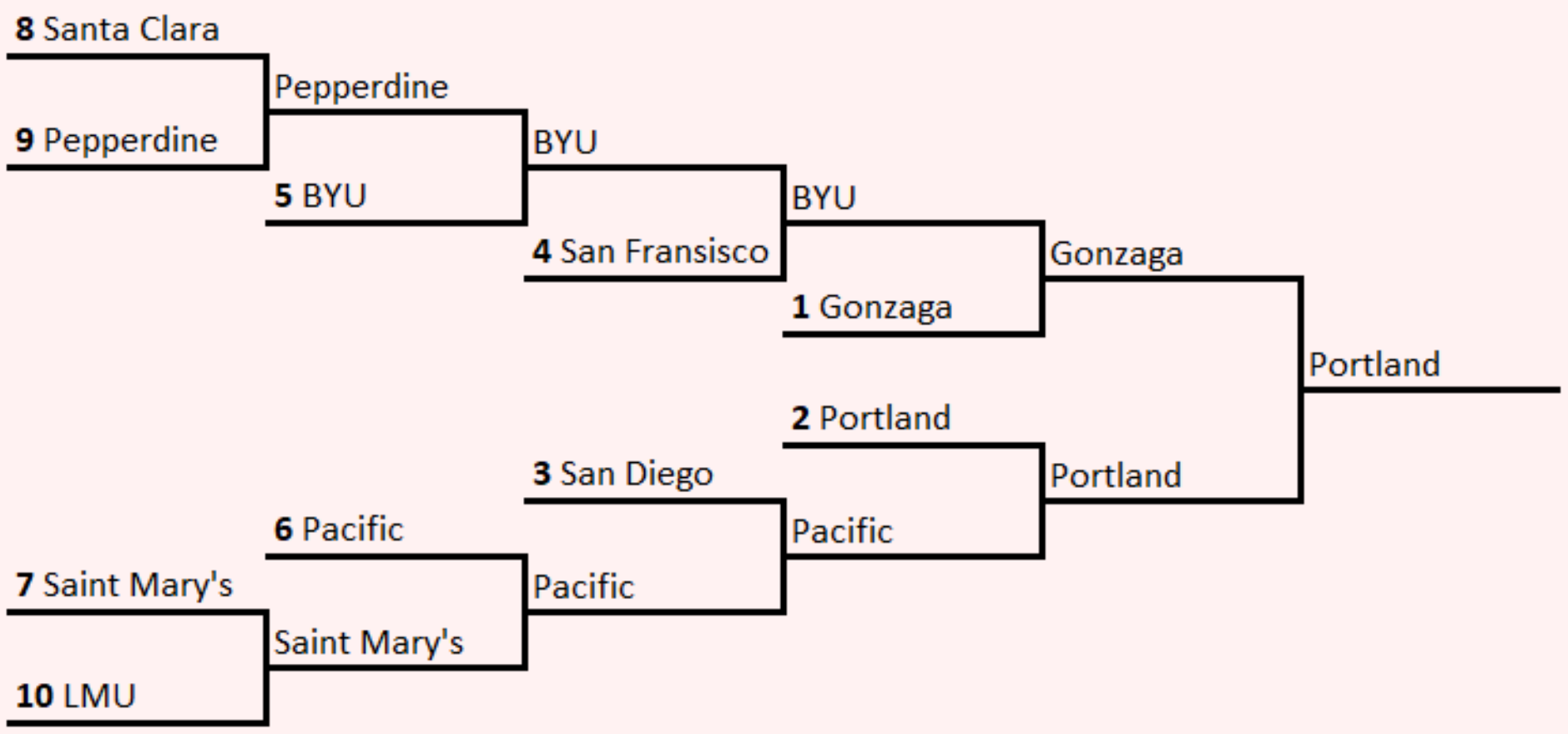

Gonzaga and Portland each have three byes and so only need to win two games to win the tournament, while Santa Clara, Pepperdine, Saint Mary's, and LMU need to win five. Unsurprisingly, this format conveys a massive advantage to Gonzaga and Portland, but this was intentional: those two teams were being rewarded for doing the best during the regular season.

In many cases, however, it is undesirable to grant advantages to certain teams over others. One might hope, for any $n$, to be able to construct a balanced bracket for $n$ teams, but unfortunately this is rarely possible.

**Theorem 2.1.17** **(Unattributed)**

There exists an $n$-team balanced bracket if and only if $n$ is a power of two.

*Proof.* A bracket is balanced if no teams have byes, which is true exactly when its signature is of the form $[[\mathbf{n};\mathbf{0};\ldots;\mathbf{0}]]$, where $n$ is the number of teams in the bracket. By Theorem 2.1.14, such a list is a bracket signature exactly when $n = 2^r$ where $r$ is the number of zeros in the list. Thus there exists an $n$-team balanced bracket if and only if $n$ is a power of two. □

Given this, brackets are not a great option when we want to avoid giving some teams advantages over others unless we have a power of two teams. They are a fantastic tool, however, if doling out advantages is the goal, perhaps after some teams did better during the regular season and ought to be rewarded with an easier path in the bracket.

## 2.2 Proper Brackets

**Definition 2.2.1: Seeding** **(Unattributed)**

The *seeding* of an $n$-team bracket is the arrangement of the numbers 1 through $n$ on the starting lines of a bracket.

Together, the shape and seeding fully specify a bracket.

**Definition 2.2.2: $i$-seed** **(Unattributed)**

In a list of teams $\mathcal{T} = [t_1, ..., t_n]$, we refer to $t_i$ as the $i$-seed.

**Definition 2.2.3: Higher and Lower Seeds** **(Unattributed)**

Somewhat confusingly, convention is that smaller numbers are the *higher seeds*, and bigger numbers are the *lower seeds*.

Seeding is typically used to reward better and more deserving teams. Consider the eight-team bracket used in the 2005 National Basketball Association Eastern Conference Playoffs [38]. At the end of the regular season, the top eight teams in the Eastern Conference were ranked and placed into the bracket which played out as shown below.

**Figure 2.2.4: 2005 NBA Eastern Conference Playoffs**

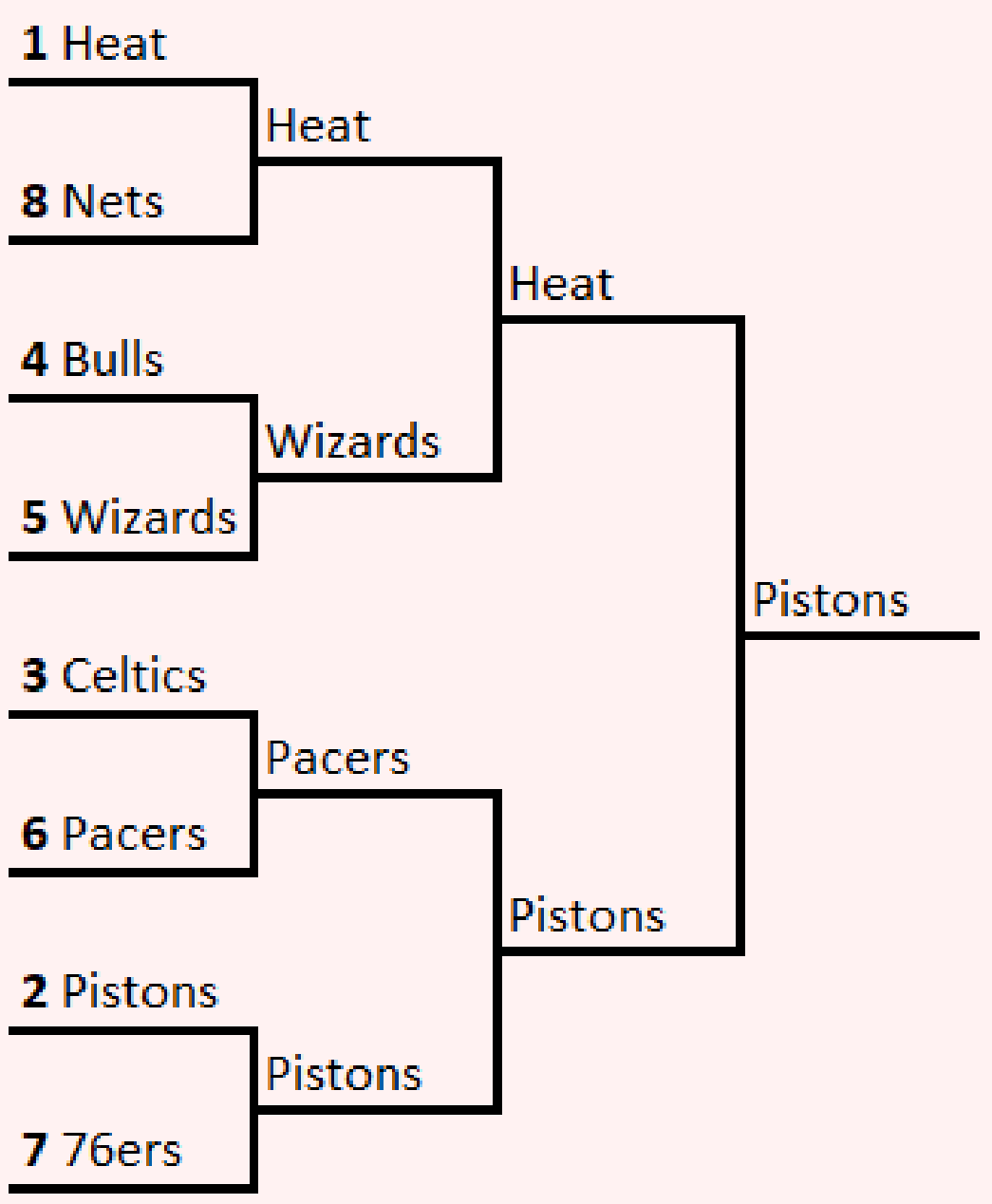

Despite this bracket being balanced, the higher seeds are still at advantage: they have an easier set of opponents. Compare the 1-seeded Heat, whose first two rounds are versus the 8-seeded Nets and then the 5-seeded Wizards, to the 6-seeded Pacers, whose first two rounds are versus the 3-seeded Celtics and then the 2-seeded Pistons. The Heat's schedule is far easier: despite their needing to win the same number of games as the Pacers, the Heat are set up to play lower seeds on their way there than the Pacers are.

Thus, we've identified two ways in which brackets can give an advantage to certain teams: by giving them more byes, and by giving them easier (expected) opponents. Not every seeding of a bracket does this: for example, consider the following alternative seeding for the 2005 NBA Eastern Conference Playoffs.

**Figure 2.2.5: Alternative Seeding of the 2005 NBA Eastern Conference Playoffs**

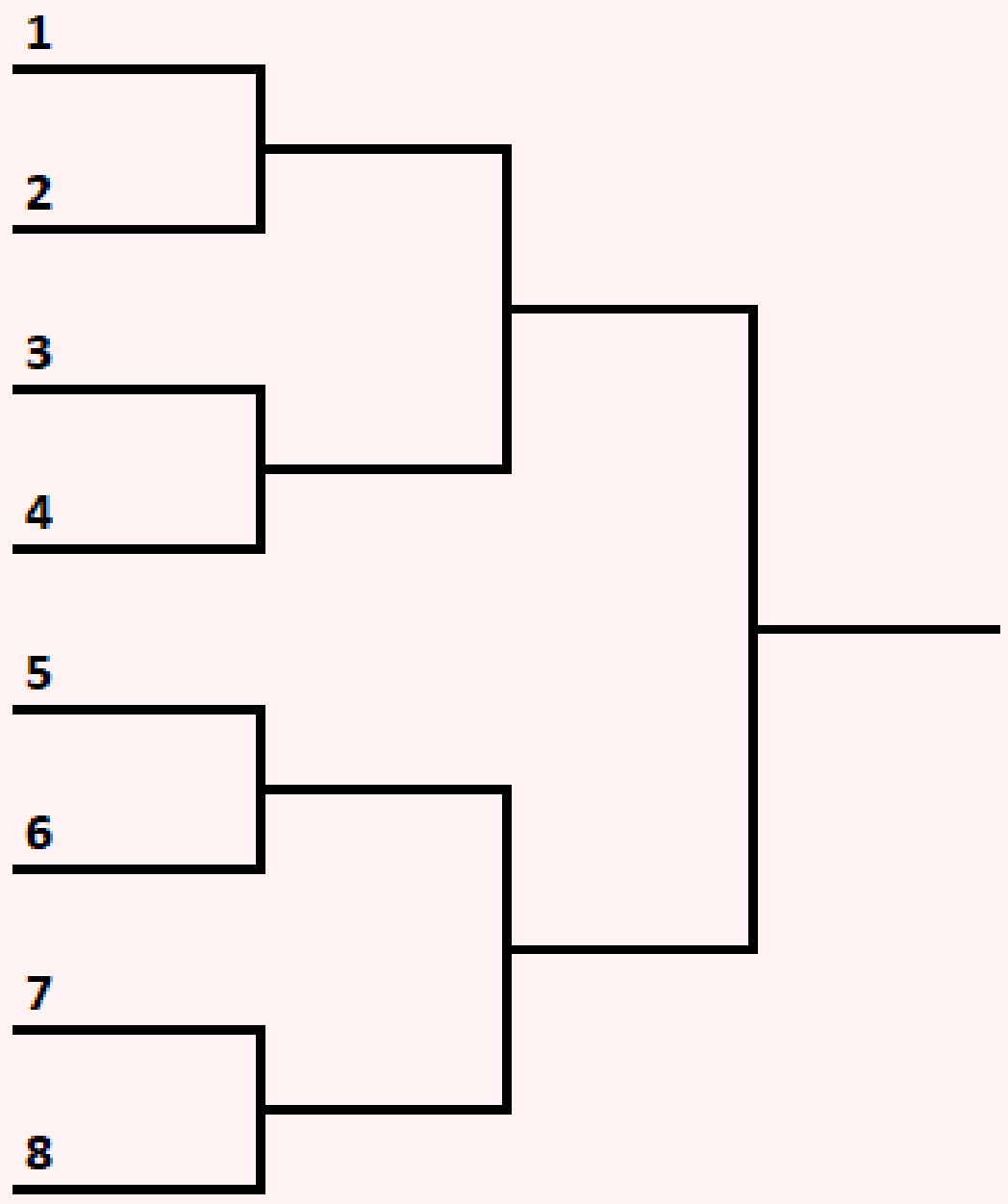

This seeding does a very poor job of rewarding the higher-seeded teams: the 1- and 2-seeds are matched up in the first round, while the easiest road is given to the 7-seed, who plays the 8-seed in the first round and then the 5-seed or 6-seed in the second. Since the whole point of seeding is to give the higher-seeded teams an advantage, we introduce the concept of a *proper seeding.*

**Definition 2.2.6: Chalk** **(Unattributed)**

A tournament went *chalk* if the higher-seeded team won every game during the tournament.

**Definition 2.2.7: Proper Seeding** **(Fried, 2024)**

A seeding of a bracket is *proper* if, as long as the bracket goes chalk, in every round it is better to be a higher-seeded team than a lower-seeded one, where:

(a) It is better to have a bye than to play a game.

(b) It is better to play a lower seed than to play a higher seed.

**Definition 2.2.8: Proper Bracket** **(Fried, 2024)**

A bracket is *proper* if its seeding is proper.

It is clear that the actual 2005 NBA Eastern Conference Playoffs was properly seeded, while our alternative seeding was not.

We now quickly derive a few lemmas about proper brackets.

**Lemma 2.2.9** **(Fried, 2024)**

In a proper bracket, if $m$ teams have a bye in a given round, those teams must be seeds 1 through $m$.

*Proof.* If they did not, the seeding would be in violation of condition (a). □

**Definition 2.2.10: Dramatic Bracket** **(Fried, 2024)**

A bracket is *dramatic* if, as long as the bracket goes chalk, in every round, the $m$ remaining teams are the top $m$ seeds.

**Lemma 2.2.11** **(Fried, 2024)**

Proper brackets are dramatic.

*Proof.* We will prove the contrapositive. Let $\mathcal{A}$ be a bracket that is not dramatic, so for some $i < j$, after some round, $t_i$ has been eliminated but $t_j$ is still alive. Let $k$ be the seed of the team that $t_i$ lost to. Because the bracket went chalk, $k < i$. Now consider what $t_j$ did in that round. If they had a bye, then the bracket violates condition (a). Assume instead they played $t_\ell$. They beat $t_\ell$, so $j < \ell$, giving,

$$k < i < j < \ell.$$

In the round that $t_i$ was eliminated, $t_i$ played $t_k$, while $t_j$ played $t_\ell$, violating condition (b). Thus, $\mathcal{A}$ is not proper. □

**Lemma 2.2.12** **(Fried, 2024)**

In a proper bracket, if in a given round $m$ teams have a bye and $k$ games are being played, then if the bracket goes chalk, those matchups will be seed $m+i$ vs seed $(m+2k+1)-i$ for $i \in \{1, ..., k\}$.

*Proof.* In the given round, there are $m+2k$ teams remaining. Theorem 2.2.11 tells us that (if the bracket goes chalk) those teams must be seeds 1 through $m+2k$. Theorem 2.2.9 tells us that seeds 1 through $m$ must have a bye, so the teams playing must be seeds $m+1$ through $m+2k$. Then condition (b) tells us that the matchups must be exactly $m+i$ vs seed $(m+2k+1)-i$ for $i \in \{1, ..., k\}$. □

We can use Lemmas 2.2.9 through 2.2.12 to properly seed various bracket shapes. For example, consider the following seven-team shape.

**Figure 2.2.13: A Seven-Team Bracket Shape**

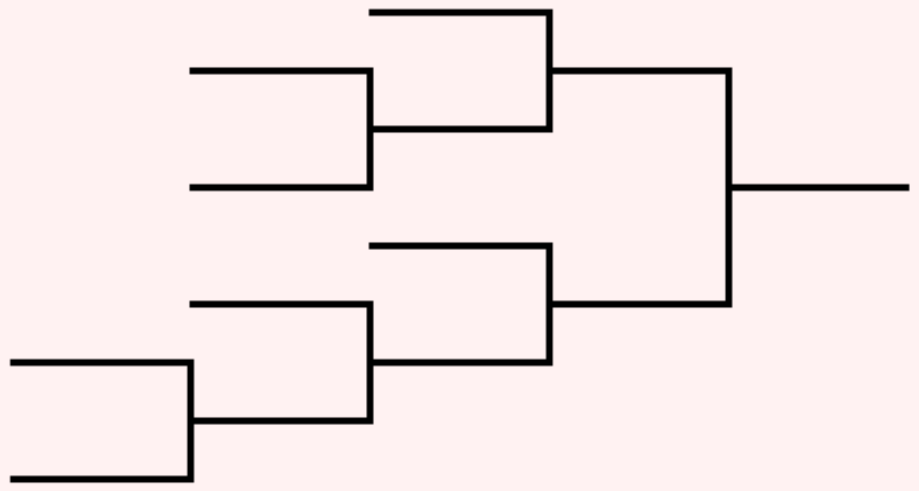

Lemma 2.2.9 tells us that the first-round matchup must be between the 6-seed and the 7-seed. Lemma 2.2.12 tells us that if the bracket goes chalk, the second-round matchups must be 3v6 and 4v5, so the 3-seed plays the winner of the first-round matchup. Finally, we can apply Lemma 2.2.12 again to the semifinals to find that the 1-seed should play the winner of the 4v5 matchup, while the 2-seed should play the winner of the 3v(6v7) matchup. In total, our proper seeding looks like this.

**Figure 2.2.14: A Seven-Team Bracket, Properly Seeded**

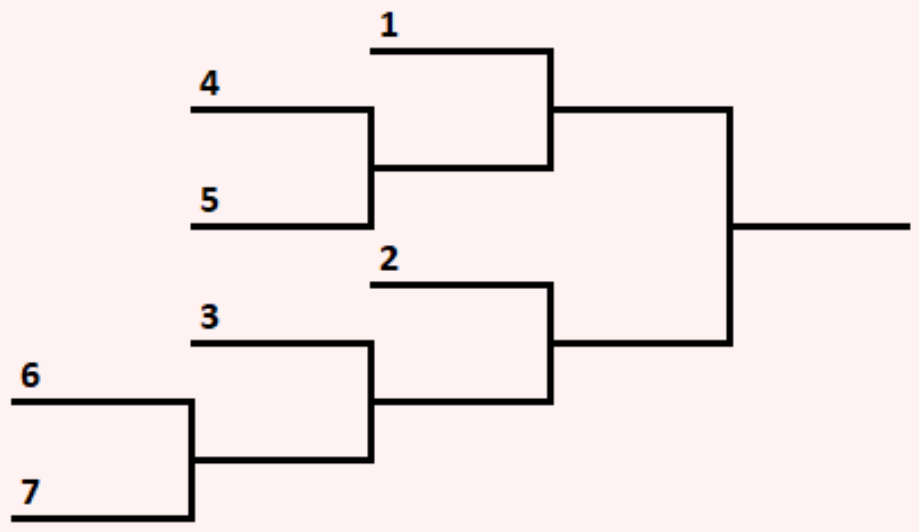

We can also quickly simulate the bracket going chalk to verify Lemma 2.2.11.

Lemmas 2.2.9 through 2.2.12 are quite powerful. It is not a coincidence that we managed to specify exactly what a proper seeding of the above bracket must look like with no room for variation: soon we will prove that the proper seeding for a particular bracket shape is unique.

But not every shape admits even this one proper seeding. Consider the following six-team shape.

**Figure 2.2.15: A Six-Team Bracket Shape**

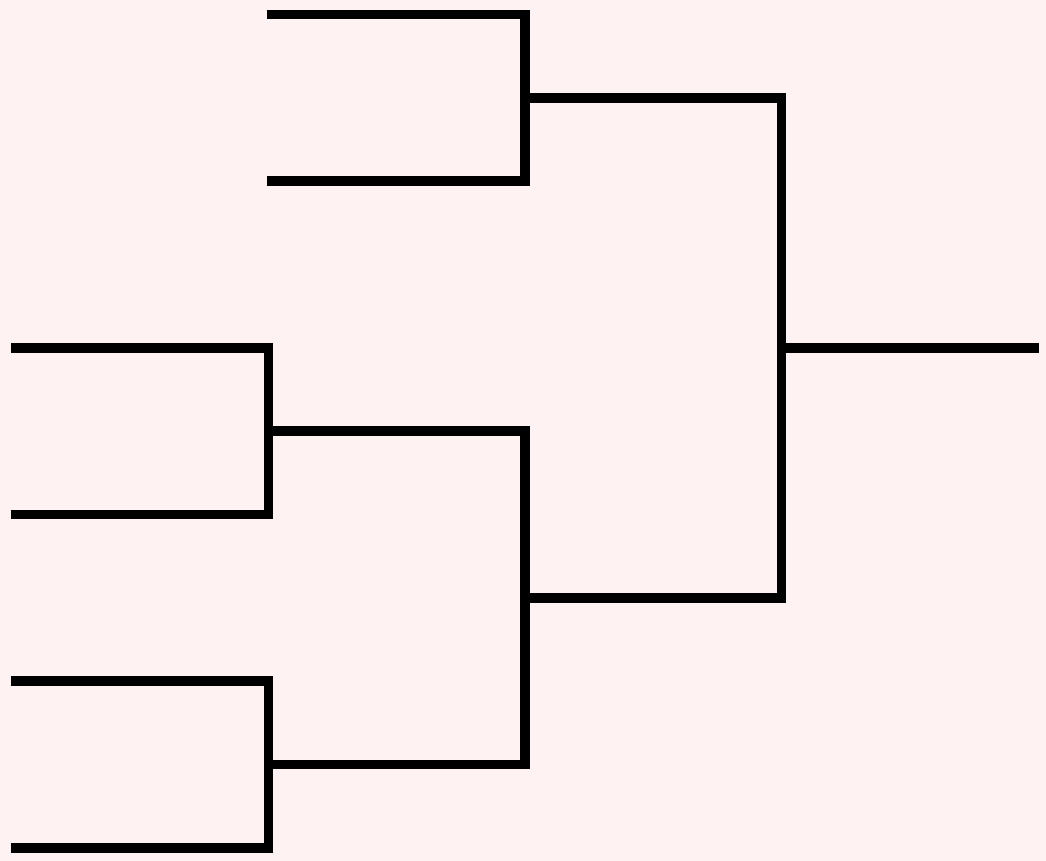

This shape admits no proper seedings. Lemma 2.2.9 requires that the two teams getting byes be the 1- and 2-seed, but this violates Lemma 2.2.12 which requires that in the second round the 1- and 2-seeds do not play each other. So how can we think about which shapes admit proper seedings?

**Theorem 2.2.16: The Fundamental Theorem of Brackets (Fried, 2024)**

There is exactly one proper bracket with each bracket signature.

*Proof.* Let $\mathcal{A} = [[\mathbf{a_0}; \ldots; \mathbf{a_r}]]$ be an $r$-round bracket signature. We proceed by induction on $r$. If $r = 0$, then the only possible bracket signature is $[[\mathbf{1}]]$, and it points to the unique one-team bracket, which is indeed proper.

For any other $r$, the first-round matchups of a proper bracket with signature $\mathcal{A}$ are defined by Lemma 2.2.12. Then if those matchups go chalk, we are left with a proper bracket of signature $[[\mathbf{a_0/2 + a_1}; \mathbf{a_2}; \ldots; \mathbf{a_r}]]$, which induction tells us admits exactly one proper bracket.

Thus both the first-round matchups and the rest of the bracket are determined, and by combining them we get a proper bracket with signature $\mathcal{A}$, so there is exactly one proper bracket with signature $\mathcal{A}$. □

The fundamental theorem of brackets means that we can refer to the proper bracket $\mathcal{A} = [[\mathbf{a_0}; ...; \mathbf{a_r}]]$ in a well-defined way, as long as

$$\sum_{i=0}^{r} a_i \cdot \left(\frac{1}{2}\right)^{r-i} = 1.$$

In practice, virtually every sports league that uses a traditional bracket uses a proper one: while different leagues take very different approaches to how many byes to give teams (compare the 2023 West Coast Conference Women's Basketball Tournament with the 2005 NBA Eastern Conference Playoffs), they are almost all proper. This makes bracket signatures a convenient labeling system for the set of brackets that we might reasonably encounter. They also are a powerful tool for specifying new brackets: if you are interested in (say) an eleven-team bracket where four teams get no byes, four teams get one bye, one team gets two byes and two teams get three byes, we can describe the proper bracket with those specs as $[[\mathbf{4}; \mathbf{4}; \mathbf{1}; \mathbf{2}; \mathbf{0}; \mathbf{0}]]$ and use Lemmas 2.2.9 through 2.2.12 to draw it with ease.

**Figure 2.2.17: The Proper Bracket of Signature $[[4; 4; 1; 2; 0; 0]]$**

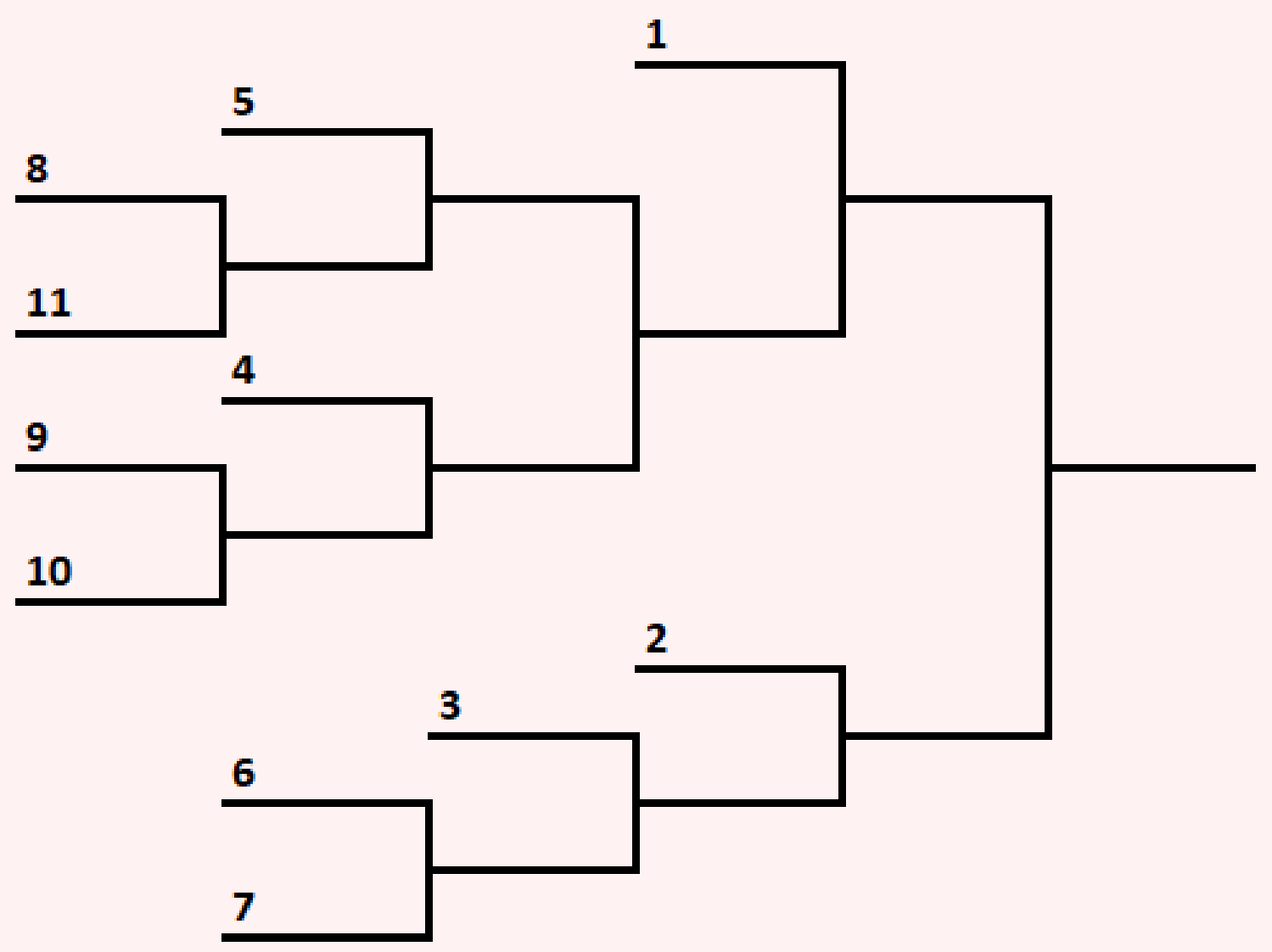

Due to these properties, we will almost exclusively discuss proper bracket from here on out: when we refer to the bracket $\mathcal{A}$ for some signature $\mathcal{A}$, we mean the proper bracket with signature $\mathcal{A}$.

## 2.3 Ordered Brackets

Consider the proper bracket $[[\mathbf{16};\mathbf{0};\mathbf{0};\mathbf{0};\mathbf{0}]]$, which was used in the 2013 NCAA Men's Basketball Tournament South Region [32], and is shown below.

**Figure 2.3.1: 2013 NCAA Men's Basketball Tournament South Region**

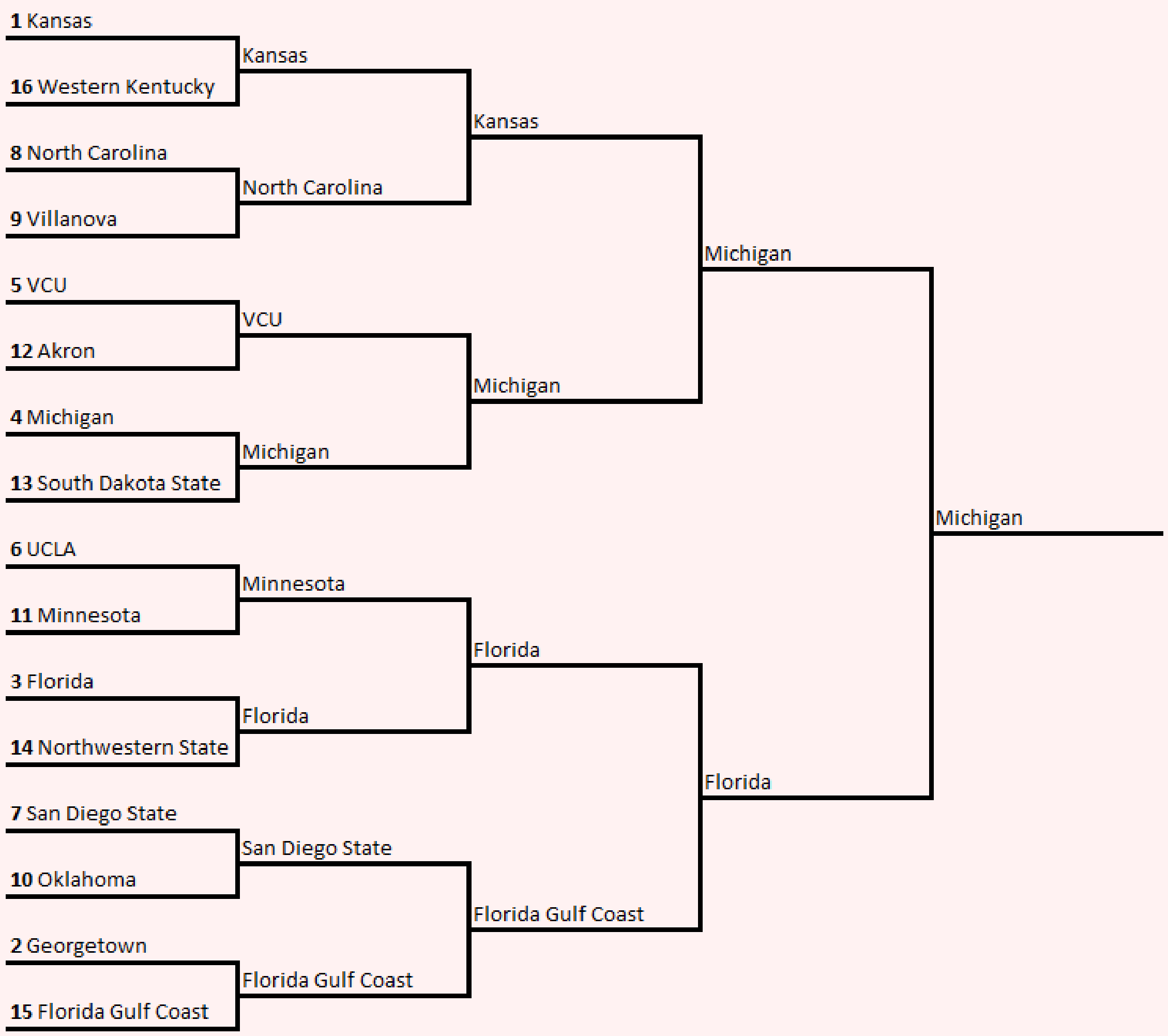

The definition of a proper seeding ensures that as long as the bracket goes chalk, it will always be better to be a higher seed than a lower seed. But what if it doesn't go chalk?

One counterintuitive fact about the NCAA Basketball Tournament is that it is probably better to be a 10-seed than a 9-seed. Why? Let's look at what seeds the $t_9$ and $t_{10}$ are likely to face in the first two rounds.

**Figure 2.3.2: NCAA Basketball Tournament 9- and 10-seed Schedules**

| Seed | First Round | Second Round |
|---|---|---|
| 9 | 8 | 1 |
| 10 | 7 | 2 |

The 9-seed has an easier first-round matchup, while the 10-seed has an easier second-round matchup. However, this isn't quite symmetrical. Because the teams tend to be drawn from a roughly normal distribution, the expected difference in skill between the 1- and 2-seeds is far greater than the expected difference between the 7- and 8-seeds, implying that the 10-seed does, in fact, have an easier route than the 9-seed.

Silver [17] investigated this matter in full, finding that in the NCAA Basketball Tournament, starting lines 10 through 15 give teams better odds of winning the region than starting lines 8 and 9. Of course, this does not mean that the 11-seed (say) has a better chance of winning a given region than the 8-seed does, as the 8-seed is a better team than the 11-seed. But it does mean that the 8-seed would love to swap places with the 11-seed, and that doing so would increase their odds to win the region.

This is not a great state of affairs: the whole point of seeding is confer an advantage to higher-seeded teams, and the proper bracket $[[\mathbf{16}; \mathbf{0}; \mathbf{0}; \mathbf{0}; \mathbf{0}]]$ is failing to do that. Not to mention that giving lower-seeded teams an easier route than higher-seeded ones can incentivize teams to lose during the regular season in order to try to get a lower but more advantageous seed.

To fix this, we need a stronger notion of what makes a bracket effective than properness. The issue with proper seedings is the false assumption that higher-seeded teams will always beat lower-seeded teams. A more nuanced assumption, initially proposed by David [9], might look like this.

**Definition 2.3.3: Strongly Stochastically Transitive (David, 1963)**

A list of teams $\mathcal{T}$ is *strongly stochastically transitive* if for each $i, j, k$ such that $j < k$,

$$\mathbb{P}[t_i \text{ beats } t_j] \leq \mathbb{P}[t_i \text{ beats } t_k].$$

A list of teams being strongly stochastically transitive (SST) captures the intuition that each team ought to do better against lower-seeded teams than against higher-seeded teams. We give a few quick implications of this definition.

**Corollary 2.3.4** **(Unattributed)**

(1) If $\mathcal{T}$ is SST, then for each $i < j$,

$$\mathbb{P}[t_i \text{ beats } t_j] \geq 0.5.$$

(2) If $\mathcal{T}$ is SST, then for each $i, j, k, \ell$ such that $i < j$ and $k < \ell$,

$$\mathbb{P}[t_i \text{ beats } t_\ell] \geq \mathbb{P}[t_j \text{ beats } t_k].$$

(3) If $\mathcal{T}$ is SST, then the matchup table $\mathcal{M}$ is monotonically increasing along each row and monotonically decreasing along each column.

Note that not every set of teams can be seeded to be SST. Consider, for example, the game of rock-paper-scissors. Rock loses to paper which loses to scissors which loses to rock, so no ordering of these "teams" will be SST. For our purposes, however, SST will work well enough.

Our new, nuanced alternative of a proper bracket is an *ordered bracket*. The concept of orderedness was first used by Chung and Hwang [6] and Horen and Riezman [11], but Edwards [10] was the one to formalize and name it.

**Definition 2.3.5: Monotonic** **(Unattributed)**

A tournament format $\mathcal{A}$ is monotonic with respect to a list of teams $\mathcal{T}$ if, for all $i < j$, $\mathbb{W}_{\mathcal{A}}(t_i, \mathcal{T}) \geq \mathbb{W}_{\mathcal{A}}(t_j, \mathcal{T})$.

**Definition 2.3.6: Ordered** **(Edwards, 1991)**

An $n$-team tournament format $\mathcal{A}$ is *ordered* if it is monotonic with respect to every SST list of $n$ teams.

In an informal sense, a tournament format being ordered is the strongest thing we can want without knowing more about why the tournament is being played. Depending on the situation, we might be interested in a format that almost always declares the most-skilled team as the winner, or in a format that gives each team roughly the same chance of winning, or anywhere in between. But certainly, better teams should win more, which is what the ordered condition requires.

In particular, a bracket being ordered is a stronger claim than it being proper.

**Theorem 2.3.7** **(Fried, 2024)**

Every ordered bracket is proper.

*Proof.* We show the contrapositive. Let $\mathcal{A}$ be an $r$-round non-proper bracket.

Assume first that $\mathcal{A}$ violates condition (a). Let $t_i$ and $t_j$ be teams such that $i < j$, but $t_i$ plays its first game in round $r_i$ while $t_j$ plays its first game in round $r_j$ for $r_i < r_j$. Let $\mathcal{T}$ be a list of teams such that $p_{ij} = 0.5$ for all $i, j$. Then,

$$\mathbb{W}_{\mathcal{A}}(t_i, \mathcal{T}) = 0.5^{r-r_i+1} < 0.5^{r-r_j+1} = \mathbb{W}_{\mathcal{A}}(t_j, \mathcal{T}).$$

Thus $\mathcal{A}$ is not monotonic with respect to $\mathcal{T}$, so it is not ordered.

Now assume $\mathcal{A}$ violated condition (b) for the first time in the $s$th round, and let $t_\ell$ be the lowest-seeded team such that there exists a $t_i, t_j$, and $t_k$ where if $\mathcal{A}$ goes chalk, then in round $s$, $t_i$ will play $t_j$ and $t_k$ will play $t_\ell$, but $i < k$ and $j < \ell$ (thus breaking condition (b)). Because $t_\ell$ is the lowest such seed, we also have $k < \ell$.

Let $\mathcal{T}$ be the SST set of teams where all games between teams seeded $\ell - 1$ or better is a coin flip, but all games involving at least one team seeded $\ell$ or worse is always won by the higher seeded team. Then

$$\mathbb{W}_{\mathcal{A}}(t_i, \mathcal{T}) = 0.5^{r-s+1} > 0.5^{r-s} = \mathbb{W}_{\mathcal{A}}(t_k, \mathcal{T}).$$

Thus $\mathcal{A}$ is not monotonic with respect to $\mathcal{T}$, so it is not ordered.

Therefore all ordered brackets are proper. □

With Theorem 2.3.7, we can use the language of bracket signatures to describe ordered brackets without worrying that two ordered brackets might share a signature. Now we examine three particularly important examples of ordered brackets.

We begin with the unique one-team bracket.

**Figure 2.3.8:** [[1]]

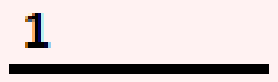

**Theorem 2.3.9** **(Unattributed)**

[[**1**]] is ordered.

*Proof.* Since there is only team, the ordered bracket condition is vacuously true. □

Next we look at the unique two-team bracket.

**Figure 2.3.10:** [[2; 0]]

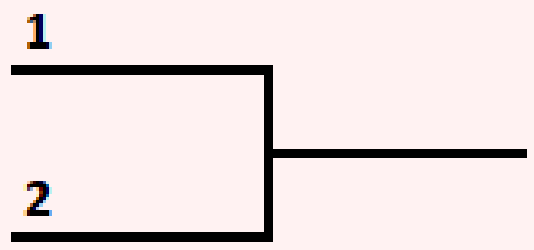

**Theorem 2.3.11** **(Unattributed)**

[[**2**; **0**]] is ordered.

*Proof.* Let $\mathcal{A} = [[\mathbf{2}; \mathbf{0}]]$. Then,

$$\mathbb{W}_{\mathcal{A}}(t_1, \mathcal{T}) = \mathbb{P}[t_1 \text{ beats } t_2] \geq 0.5 \geq \mathbb{P}[t_2 \text{ beats } t_1] = \mathbb{W}_{\mathcal{A}}(t_2, \mathcal{T})$$

so $\mathcal{A}$ is ordered. □

And thirdly, we show that the balanced four-team bracket is ordered, first proved by Horen and Riezman [11].

**Figure 2.3.12:** [[4; 0; 0]]

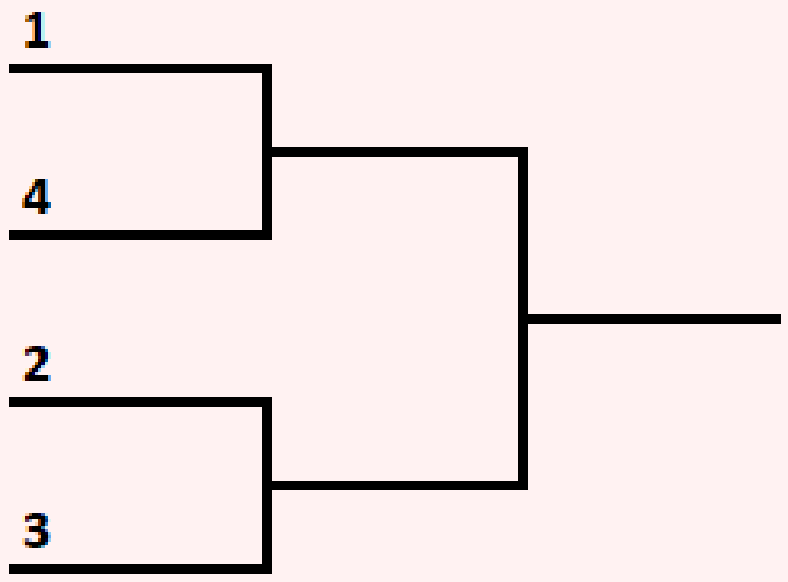

**Theorem 2.3.13** **(Horen and Riezman, 1985)**

$[[\mathbf{4};\mathbf{0};\mathbf{0}]]$ is ordered.

*Proof.* Let $\mathcal{A} = [[\mathbf{4};\mathbf{0};\mathbf{0}]]$. Then,

$$\begin{aligned}
\mathbb{W}_{\mathcal{A}}(t_1, \mathcal{T}) &= p_{14} \cdot (p_{23}p_{12} + p_{32}p_{13}) \\
&= p_{14}p_{23}p_{12} + p_{14}p_{32}p_{13} \\
&\geq p_{14}p_{23}p_{21} + p_{24}p_{41}p_{23} \\
&= p_{23} \cdot (p_{14}p_{21} + p_{41}p_{24}) \\
&= \mathbb{W}_{\mathcal{A}}(t_2, \mathcal{T})
\end{aligned}$$

$$\begin{aligned}
\mathbb{W}_{\mathcal{A}}(t_2, \mathcal{T}) &= p_{23} \cdot (p_{14}p_{21} + p_{41}p_{24}) \\
&\geq p_{32} \cdot (p_{14}p_{31} + p_{41}p_{34}) \\
&= \mathbb{W}_{\mathcal{A}}(t_3, \mathcal{T})
\end{aligned}$$

$$\begin{aligned}
\mathbb{W}_{\mathcal{A}}(t_3, \mathcal{T}) &= p_{32} \cdot (p_{14}p_{31} + p_{41}p_{34}) \\
&= p_{32}p_{14}p_{31} + p_{32}p_{41}p_{34} \\
&\geq p_{42}p_{23}p_{41} + p_{32}p_{41}p_{43} \\
&= p_{41} \cdot (p_{23}p_{42} + p_{32}p_{43}) \\
&= \mathbb{W}_{\mathcal{A}}(t_4, \mathcal{T})
\end{aligned}$$

Thus $\mathcal{A}$ is ordered. □

However, not every proper bracket is ordered. One particularly important example of a non-ordered proper bracket is $[[\mathbf{4}; \mathbf{2}; \mathbf{0}; \mathbf{0}]]$.

**Figure 2.3.14:** $[[4; 2; 0; 0]]$

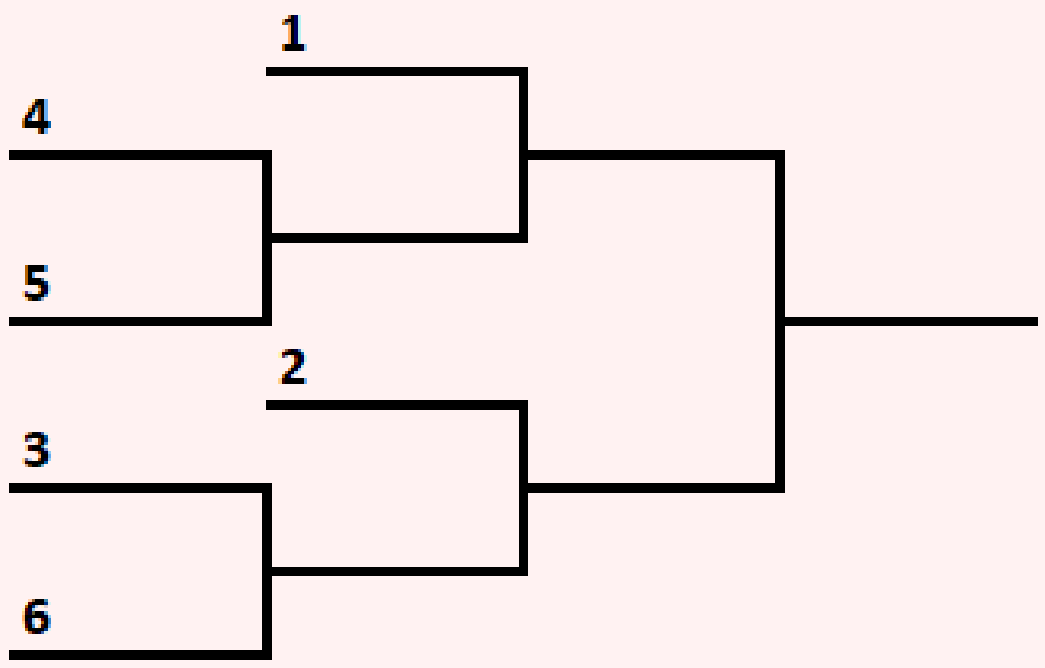

**Theorem 2.3.15** **(Edwards, 1991)**

$[[\mathbf{4}; \mathbf{2}; \mathbf{0}; \mathbf{0}]]$ is not ordered.

*Proof.* Let $\mathcal{A} = [[\mathbf{4}; \mathbf{2}; \mathbf{0}; \mathbf{0}]]$, and let $\mathcal{T}$ have the following matchup table.

| | $t_1$ | $t_2$ | $t_3$ | $t_4$ | $t_5$ | $t_6$ |
|---|---|---|---|---|---|---|
| $t_1$ | | | | | | |
| $t_2$ | 0.5 | | | | | |
| $t_3$ | 0.5 | 0.5 | | | | |
| $t_4$ | 0.5 | 0.5 | 0.5 | | | |
| $t_5$ | 0 | 0.5 | 0.5 | 0.5 | | |
| $t_6$ | 0 | 0.5 | 0.5 | 0.5 | 0.5 | |

Then $\mathbb{W}_{\mathcal{A}}(t_5, \mathcal{T}) = 0$ but $\mathbb{W}_{\mathcal{A}}(t_6, \mathcal{T}) > 0$, so $\mathcal{A}$ is not monotonic with respect to $\mathcal{T}$ and thus not ordered. □

(Note that in this particular example, one could argue that $t_5$ isn't actually better than $t_6$, as their odds of beating each other team is the same, and thus it is not a big deal that $t_6$ is more likely to win the tournament. However, $\mathcal{A}$ is not monotonic with respect to any matchup table where $\mathbb{P}[t_i \text{ beats } t_j] = 0$ if and only if $i \in \{5, 6\}$ and $j = 1$, even ones where $t_5$ is clearly the superior team.)

In the next section, we move on from describing particular ordered and non-ordered brackets in favor of a more general result.

## 2.4 Edwards's Theorem

We now attempt to completely classify the set of ordered brackets. Edwards's [10] original proof, as well as a more recent proof by Alegri and Dimitrov [1], accomplished this without access to the machinery of bracket signatures or proper brackets. We present a quicker proof that makes use of the fundamental theorem of brackets and develops two nice lemmas along the way.

We begin with the stapling lemma, which allows us to combine two smaller ordered brackets into a larger ordered one by having the winner of one of the brackets be treated as the lowest seed in the other. This is depicted in Figure 2.4.1.

**Figure 2.4.1: Setup of the Stapling Lemma with $\mathcal{A} = [[2;1;0]]$, $\mathcal{B} = [[4;0;0]]$, and $\mathcal{C} = [[2;1;3;0;0]]$**

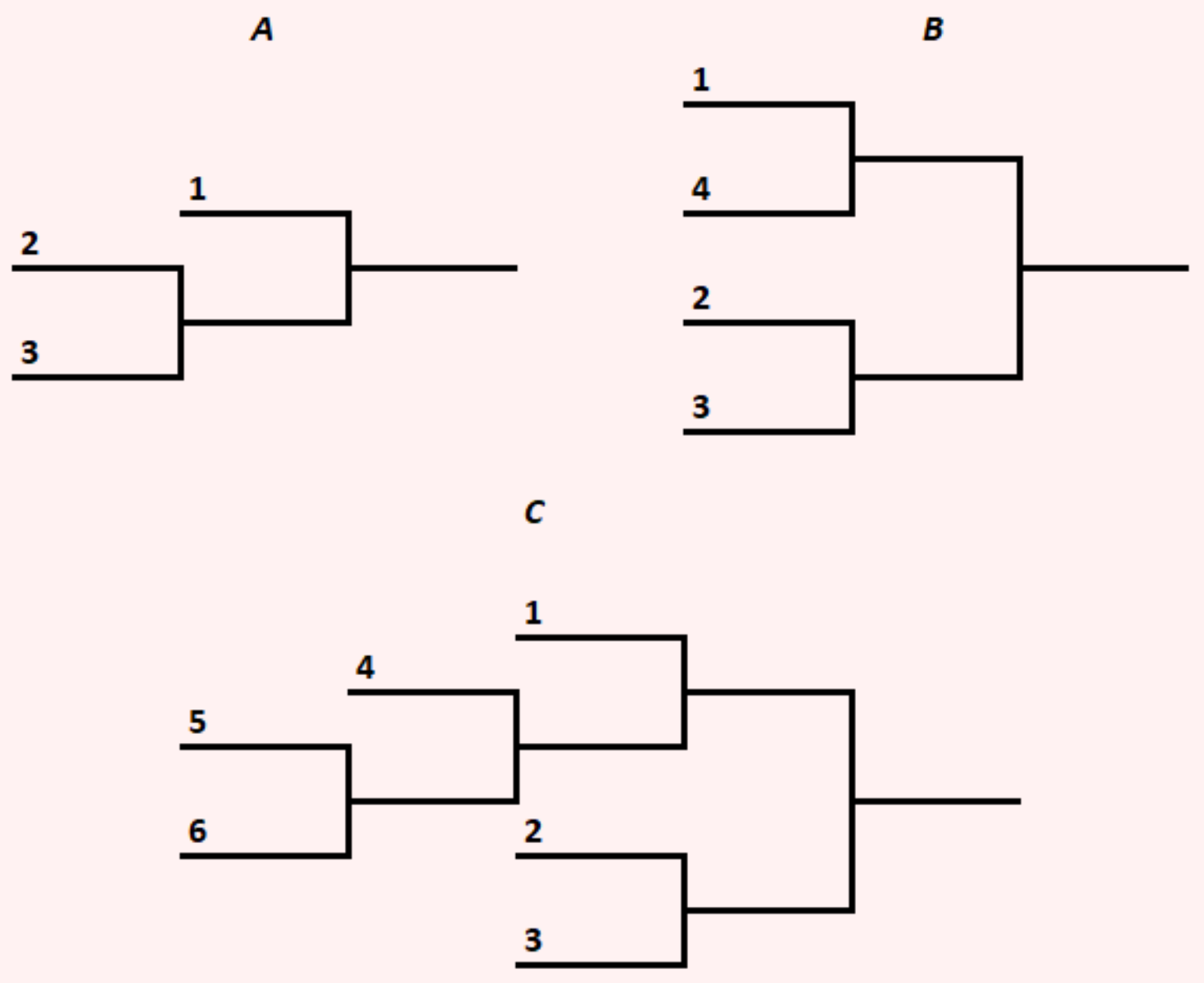

**Lemma 2.4.2: The Stapling Lemma** **(Fried, 2024)**

If $\mathcal{A} = [[\mathbf{a_0};...;\mathbf{a_r}]]$ and $\mathcal{B} = [[\mathbf{b_0};...;\mathbf{b_s}]]$ are ordered brackets, then $\mathcal{C} = [[\mathbf{a_0};...;\mathbf{a_r}+\mathbf{b_0}-\mathbf{1};...;\mathbf{b_s}]]$ is an ordered bracket as well.

*Proof.* Let $\mathcal{A}, \mathcal{B}$, and $\mathcal{C}$ be as specified. Let $\mathcal{T}$ be an SST list of $n+m-1$ teams, and let $\mathcal{R}, \mathcal{S} \subset \mathcal{T}$ be the lowest $n$ and the highest $m-1$ seeds of $\mathcal{T}$ respectively. We divide proving that $\mathcal{C}$ is ordered into proving three sub-statements:

(a) For $i < j < m$, $\mathbb{W}_{\mathcal{C}}(t_i, \mathcal{T}) \geq \mathbb{W}_{\mathcal{C}}(t_j, \mathcal{T})$

(b) $\mathbb{W}_{\mathcal{C}}(t_{m-1}, \mathcal{T}) \geq \mathbb{W}_{\mathcal{C}}(t_m, \mathcal{T})$

(c) For $m \leq i < j$, $\mathbb{W}_{\mathcal{C}}(t_i, \mathcal{T}) \geq \mathbb{W}_{\mathcal{C}}(t_j, \mathcal{T})$

Together, these show that $\mathcal{C}$ is ordered.

We begin with the first sub-statement. Let $i < j < m$. Then,

$$\begin{aligned}
\mathbb{W}_{\mathcal{C}}(t_i, \mathcal{T}) &= \sum_{k=m}^{n+m-1} \mathbb{W}_{\mathcal{A}}(t_k, \mathcal{R}) \cdot \mathbb{W}_{\mathcal{B}}(t_i, \mathcal{S} \cup \{t_k\}) \\
&\geq \sum_{k=m}^{n+m-1} \mathbb{W}_{\mathcal{A}}(t_k, \mathcal{R}) \cdot \mathbb{W}_{\mathcal{B}}(t_j, \mathcal{S} \cup \{t_k\}) \\
&= \mathbb{W}_{\mathcal{C}}(t_j, \mathcal{T})
\end{aligned}$$

The first and last equalities follow from the structure of $\mathcal{C}$, and the inequality follows from $\mathcal{B}$ being ordered.

Now the second sub-statement.

$$\begin{aligned}
\mathbb{W}_{\mathcal{C}}(t_{m-1}, \mathcal{T}) &= \sum_{k=m}^{n+m-1} \mathbb{W}_{\mathcal{A}}(t_k, \mathcal{R}) \cdot \mathbb{W}_{\mathcal{B}}(t_{m-1}, \mathcal{S} \cup \{t_k\}) \\
&\geq \mathbb{W}_{\mathcal{A}}(t_m, \mathcal{R}) \cdot \mathbb{W}_{\mathcal{B}}(t_{m-1}, \mathcal{S} \cup \{t_m\}) \\
&\geq \mathbb{W}_{\mathcal{A}}(t_m, \mathcal{R}) \cdot \mathbb{W}_{\mathcal{B}}(t_m, \mathcal{S} \cup \{t_m\}) \\
&= \mathbb{W}_{\mathcal{C}}(t_m, \mathcal{T})
\end{aligned}$$

The equalities follow from the structure of $\mathcal{C}$, the first inequality follows from probabilities being non-negative, and the second inequality follows from $\mathcal{B}$ being ordered.

Finally, we show the third sub-statement. Let $m \leq i < j$. Then,

$$\begin{aligned}
\mathbb{W}_{\mathcal{C}}(t_i, \mathcal{T}) &= \mathbb{W}_{\mathcal{A}}(t_i, \mathcal{R}) \cdot \mathbb{W}_{\mathcal{B}}(t_i, \mathcal{S} \cup \{t_i\}) \\
&\geq \mathbb{W}_{\mathcal{A}}(t_j, \mathcal{R}) \cdot \mathbb{W}_{\mathcal{B}}(t_i, \mathcal{S} \cup \{t_i\}) \\
&\geq \mathbb{W}_{\mathcal{A}}(t_j, \mathcal{R}) \cdot \mathbb{W}_{\mathcal{B}}(t_j, \mathcal{S} \cup \{t_j\}) \\
&= \mathbb{W}_{\mathcal{C}}(t_j, \mathcal{T})
\end{aligned}$$

The equalities follow from the structure of $\mathcal{C}$, the first inequality from $\mathcal{A}$ being ordered, and the second inequality from the teams being SST.

We have shown all three sub-statements, so $\mathcal{C}$ is ordered. □

Now, if we begin with the set of brackets $\{[[\mathbf{1}]], [[\mathbf{2}; \mathbf{0}]], [[\mathbf{4}; \mathbf{0}; \mathbf{0}]]\}$ and then repeatedly apply the stapling lemma, we can construct a set of brackets that we know are ordered. In other words,

**Corollary 2.4.3** **(Edwards, 1991)**

Proper brackets whose signature is formed by the following process are ordered:

1. Start with the list $[[\mathbf{0}]]$ (note that this is not yet a bracket signature).
2. As many times as desired, prepend the list with $[[\mathbf{1}]]$ or $[[\mathbf{3}; \mathbf{0}]]$.
3. Then, add 1 to the first element in the list, turning it into a bracket signature.

Corollary 2.4.3 uses the tools that we have developed so far to identify a set of ordered brackets. Somewhat surprisingly, this set is complete: any bracket not reachable using the process in Corollary 2.4.3 is not ordered. To prove this we first need to show the containment lemma.

**Definition 2.4.4: Containment** **(Fried, 2024)**

Let $\mathcal{A}$ and $\mathcal{B}$ be bracket signatures. $\mathcal{A}$ *contains* $\mathcal{B}$ if there exists some $i$ such that

(a) At least as many games are played in the $(i+1)$th round of $\mathcal{A}$ as in the first round of $\mathcal{B}$, and

(b) For $1 < j \leq r$ where $r$ is the number of rounds in $\mathcal{B}$, there are exactly as many games played in the $(i+j)$th round of $\mathcal{A}$ as in the $j$th round of $\mathcal{B}$.

Intuitively, $\mathcal{A}$ containing $\mathcal{B}$ means that if $\mathcal{A}$ went chalk, and games within each round were played in order of largest seed-gap to smallest seed-gap, then at some point, there would be a bracket of shape $\mathcal{B}$ used to determine the last team in the rest of bracket $\mathcal{A}$. Figure 2.4.5 shows $\mathcal{A} = [[\mathbf{2}; \mathbf{5}; \mathbf{1}; \mathbf{0}; \mathbf{3}; \mathbf{0}; \mathbf{0}]]$ containing $\mathcal{B} = [[\mathbf{4}; \mathbf{2}; \mathbf{0}; \mathbf{0}]]$. After the 10v11 game and the 5v(10v11) game, there is a bracket of shape $\mathcal{B}$ (the solid lines) that will be played to identify the last team to play in the rest of the bracket.

**Figure 2.4.5: Setup of the Containment Lemma with $\mathcal{A} = [[2; 5; 1; 0; 3; 0; 0]]$ and $\mathcal{B} = [[4; 2; 0; 0]]$.**

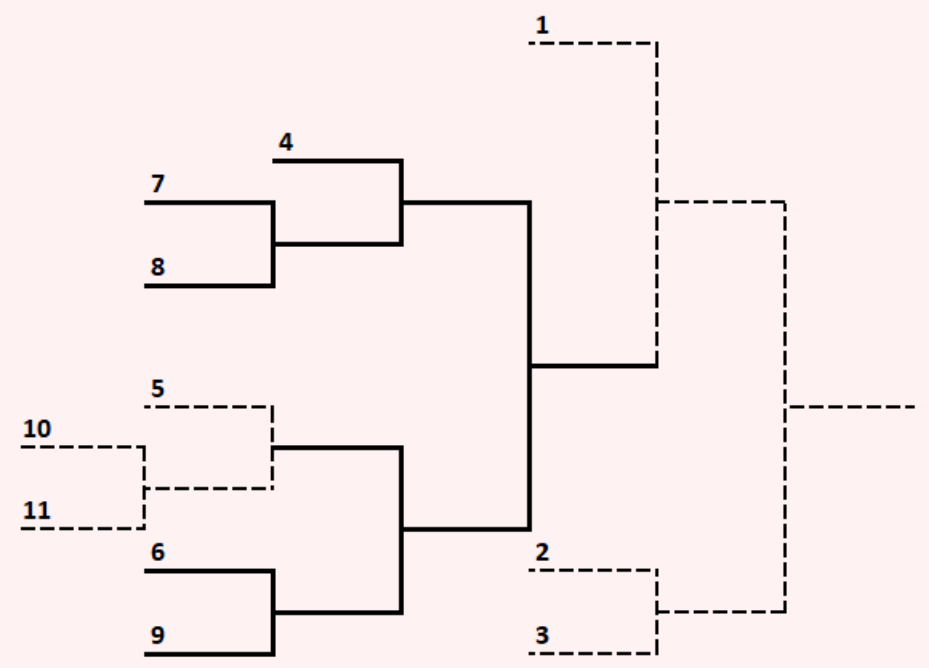

**Lemma 2.4.6: The Containment Lemma** **(Fried, 2024)**

If $\mathcal{A}$ contains $\mathcal{B}$, and $\mathcal{B}$ is not ordered, then neither is $\mathcal{A}$.

*Proof.* Let $\mathcal{A}$ be a bracket signature with $r$ rounds and $n$ teams, and let $\mathcal{B}$ have $s$ rounds and $m$ teams, such that $\mathcal{A}$ contains $\mathcal{B}$ and $\mathcal{B}$ is not ordered. Let $k$ be the number of teams in $\mathcal{A}$ that get at least $s+i$ byes (where $i$ is from the definition of *contains*).

$\mathcal{B}$ is not ordered, so let $\mathcal{M}$ be a matchup table that violates the orderedness condition, where none of the win probabilities are 0. (If we have an $\mathcal{M}$ that includes 0s, we can replace them with $\epsilon$. For small enough $\epsilon$, $\mathcal{M}$ will still violate the condition.) Let $p$ be the minimum probability in $\mathcal{M}$. Let $\mathcal{P}$ be a matchup table in which the lower-seeded team wins with probability $p$, and let $\mathcal{Z}$ be a matchup table in which the lower-seeded team wins with probability 0.

Now, consider the following block matchup table on a list of $n$ teams $\mathcal{T}$.

| | $t_1$ - $t_k$ | $t_{k+1}$ - $t_{k+m}$ | $t_{k+m+1}$ - $t_n$ |
|---|---|---|---|
| $t_1$ - $t_k$ | $\mathcal{P}$ | $\mathcal{P}$ | $\mathcal{Z}$ |
| $t_{k+1}$ - $t_{k+m}$ | $\mathcal{P}$ | $\mathcal{M}$ | $\mathcal{Z}$ |
| $t_{k+m+1}$ - $t_n$ | $\mathcal{Z}$ | $\mathcal{Z}$ | $\mathcal{Z}$ |

Let $\mathcal{S} \subset \mathcal{T}$ be the sublist of teams seeded between $k+1$ and $k+m$. Then, for $t_j \in \mathcal{S}$,

$$\mathbb{W}_{\mathcal{A}}(t_j, \mathcal{T}) = \mathbb{W}_{\mathcal{B}}(t_j, \mathcal{S}) \cdot p^{r-s-i},$$

since $t_j$ wins any games it might have to play in rounds $i$ or before automatically, any games after $s+i$ with probability $p$, and any games in between according to $\mathcal{M}$.

However, $\mathcal{M}$ (and thus $\mathcal{S}$) violates the orderedness condition for $\mathcal{B}$, and so $\mathcal{T}$ does for $\mathcal{A}$. □

With the containment lemma shown, we can proceed to the main theorem.

**Theorem 2.4.7: Edwards's Theorem** **(Edwards, 1991)**

The only ordered brackets are those described by Corollary 2.4.3.

*Proof.* Let $\mathcal{A}$ be a proper bracket not described by Corollary 2.4.3. The corollary describes all proper brackets in which each round either has only game, or has two games but is immediately followed by a round with only one game. Thus $\mathcal{A}$ must include at least two successive rounds with two or more games each.

Such a chain will be followed by a round with a single game, and so the final round in the chain must have only two games. Thus, $\mathcal{A}$ includes a sequence of three rounds, the first of which has at least two games, the second of which has exactly two games, and the third of which has one game.

Therefore, $\mathcal{A}$ contains $[[\mathbf{4}; \mathbf{2}; \mathbf{0}; \mathbf{0}]]$. But we know that $[[\mathbf{4}; \mathbf{2}; \mathbf{0}; \mathbf{0}]]$ is not ordered, and so by the containment lemma, neither is $\mathcal{A}$. $\square$

While quite powerful, what Edwards's Theorem says about the space of ordered brackets is quite disappointing. At most three teams can be introduced in each round of an ordered bracket, so the length of the shortest ordered bracket on $n$ teams grows linearly with $n$ (rather than logarithmically, as is the case for the shortest proper bracket). If we want a bracket on many teams to be ordered, we risk forcing lower-seeded teams to play a large number of games, and we only permit the top-seeded teams to play a few. For example, the shortest ordered bracket that the 2021 NCAA Basketball South Region could have used is $[[\mathbf{4}; \mathbf{0}; \mathbf{3}; \mathbf{0}; \mathbf{3}; \mathbf{0}; \mathbf{3}; \mathbf{0}; \mathbf{3}; \mathbf{0}; \mathbf{0}]]$, which is played over a whopping ten rounds.

**Figure 2.4.8: The Shortest Sixteen-Team Ordered Bracket**

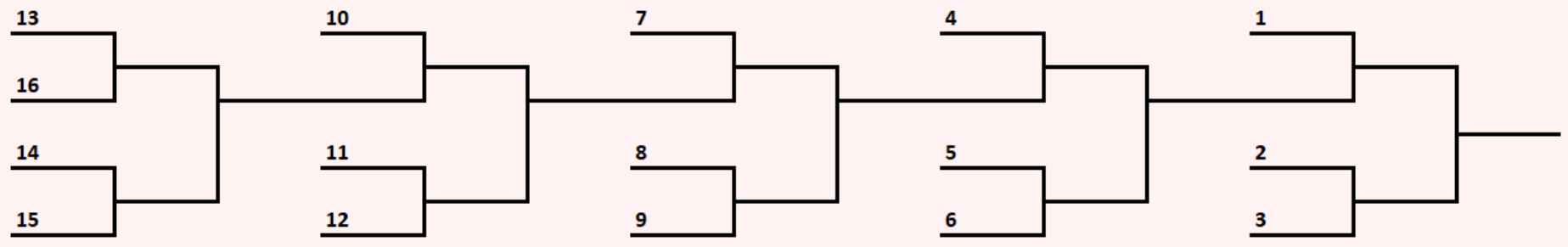

Because of this, few leagues use ordered brackets, and those that do usually have so few teams that every proper bracket is ordered.

## 2.5 Reseeded Brackets

As we discussed last section, Edwards's Theorem tells us that the number of rounds required to construct on ordered bracket grows linearly with the number of teams involved. This is quite frustrating, as part of the power of brackets is the ability to crown a champion in a number of rounds logarithmic in the number of teams participating. As an attempt to combat the problem that Edwards's Theorem presents, we expand the our gaze and consider formats that are similar to brackets, but not necessarily networked. Can we recover some ordered bracket-like formats that way?

**Definition 2.5.1: Knockout Tournament** **(Unattributed)**

A *knockout tournament* is a tournament in that is played over a series of rounds subject to the following constraints:

(a) Each team plays at most one game in each round.

(b) If a team loses in a round, they don't play any games in later rounds.

(c) If a team wins in a round, they play a game in the next round.

(d) Exactly one team finishes undefeated, and that team is crowned champion.

Clearly brackets are just networked knockout tournaments, but there are many knockout tournaments that aren't networked. The definition of a knockout tournament is designed to allow for the notions of signatures and properness to still apply.

**Definition 2.5.2: Knockout Tournament Signature** **(Fried, 2024)**

The *signature* of an $r$-round knockout tournament $\mathcal{A}$ is the list $[[\mathbf{a_0}; ...; \mathbf{a_r}]]$ where $a_i$ is the number of teams that get $i$ byes.

**Definition 2.5.3: Proper Knockout Tournament** **(Fried, 2024)**

A knockout tournament is *proper* if, as long as the tournament goes chalk, in every round it is better to be a higher-seeded team than a lower-seeded one, where:

(a) It is better to have a bye than to play a game.

(b) It is better to play a lower seed than to play a higher seed.

Ultimately, the reason that proper brackets are not, in general, ordered, is that lower-seeded teams are treated, if they win, as the team that they beat for the rest of the tournament. Consider again the proper bracket analyzed by Silver: $[[\mathbf{16}; \mathbf{0}; \mathbf{0}; \mathbf{0}; \mathbf{0}]]$. If an 11-seed wins in the first round, they take on the schedule of a 6-seed for the rest of the tournament, while if the 9-seed wins, they take on the schedule of an 8-seed. Given that a 6-seed has an easier schedule than an 8-seed, it's not hard to see why it might be preferable to be an 11-seed rather than a 9-seed.

But knockout tournaments are under no such restrictions. A knockout tournament could simply pair the highest- and lowest-remaining seeds in every round, potentially avoiding the issues we faced in the last two sections. These formats are called *reseeded brackets*.

**Definition 2.5.4: Reseeded Bracket** (Hwang, 1982)

A *reseededed bracket* is a knockout tournament in which, after each round, the highest-seeded team playing that round is matched up with the lowest-seeded team playing that round, second-highest vs second-lowest, etc.

Note that technically reseeded brackets are not networked and thus not brackets at all, just knockout tournaments. However, because reseeded brackets act so similarly to traditional brackets, and because colloquially they are referred to as brackets, we opt to continue using the word "bracket" to describe them.

In 2024, both National Football League conferences [37] used a reseeded bracket with signature $[[\mathbf{6};\mathbf{1};\mathbf{0};\mathbf{0}]]^R$. (The superscript $R$ indicates this is reseeded bracket.) If the first round of the bracket goes chalk, then it looks just like a normal bracket.

**Figure 2.5.5: 2024 National Football League AFC Playoffs**

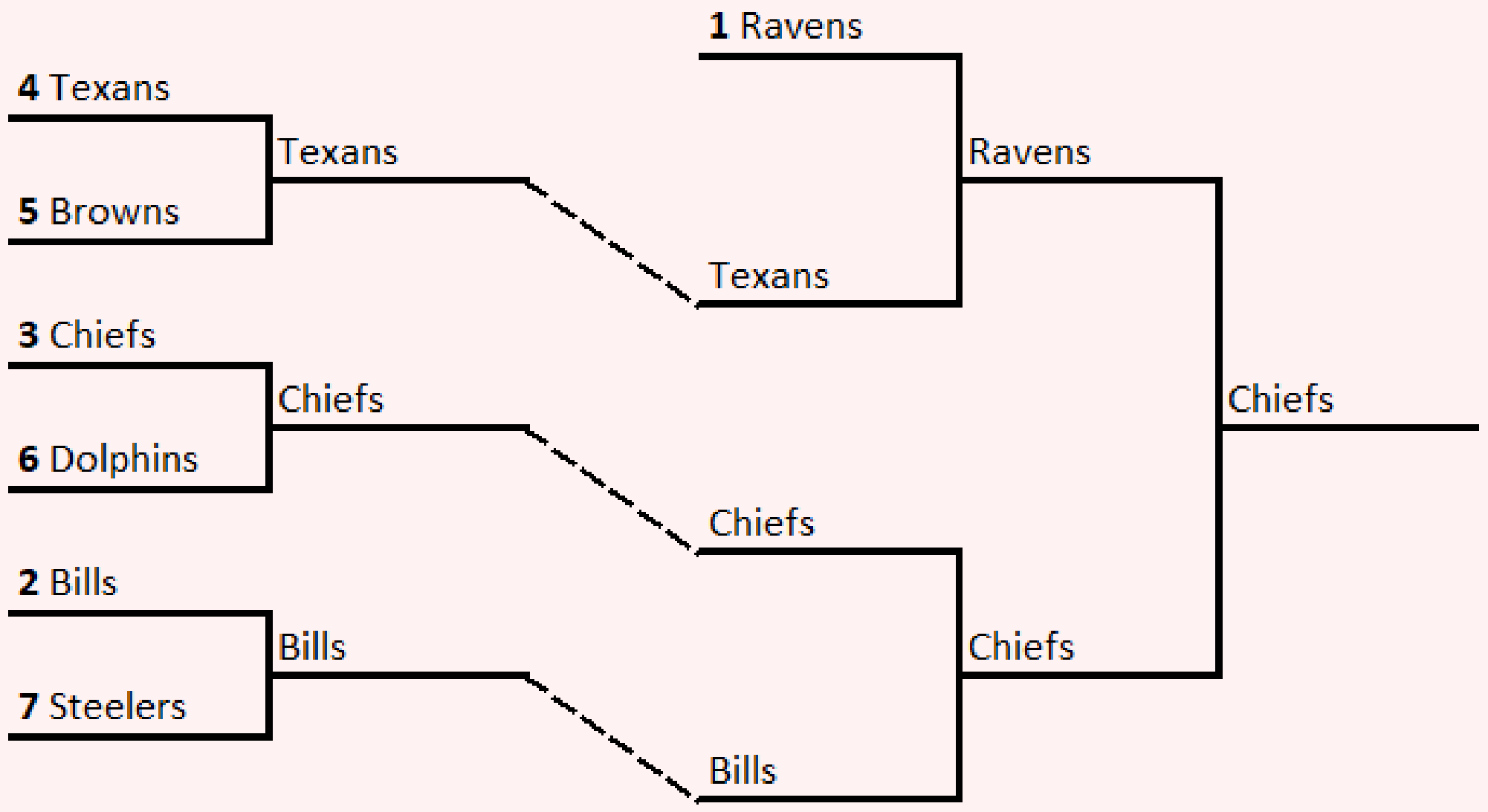

But if there are first-round upsets, then the bracket is rearranged to ensure that it is still better to be a higher seed rather than a lower seed.

**Figure 2.5.6: 2024 National Football League NFC Playoffs**

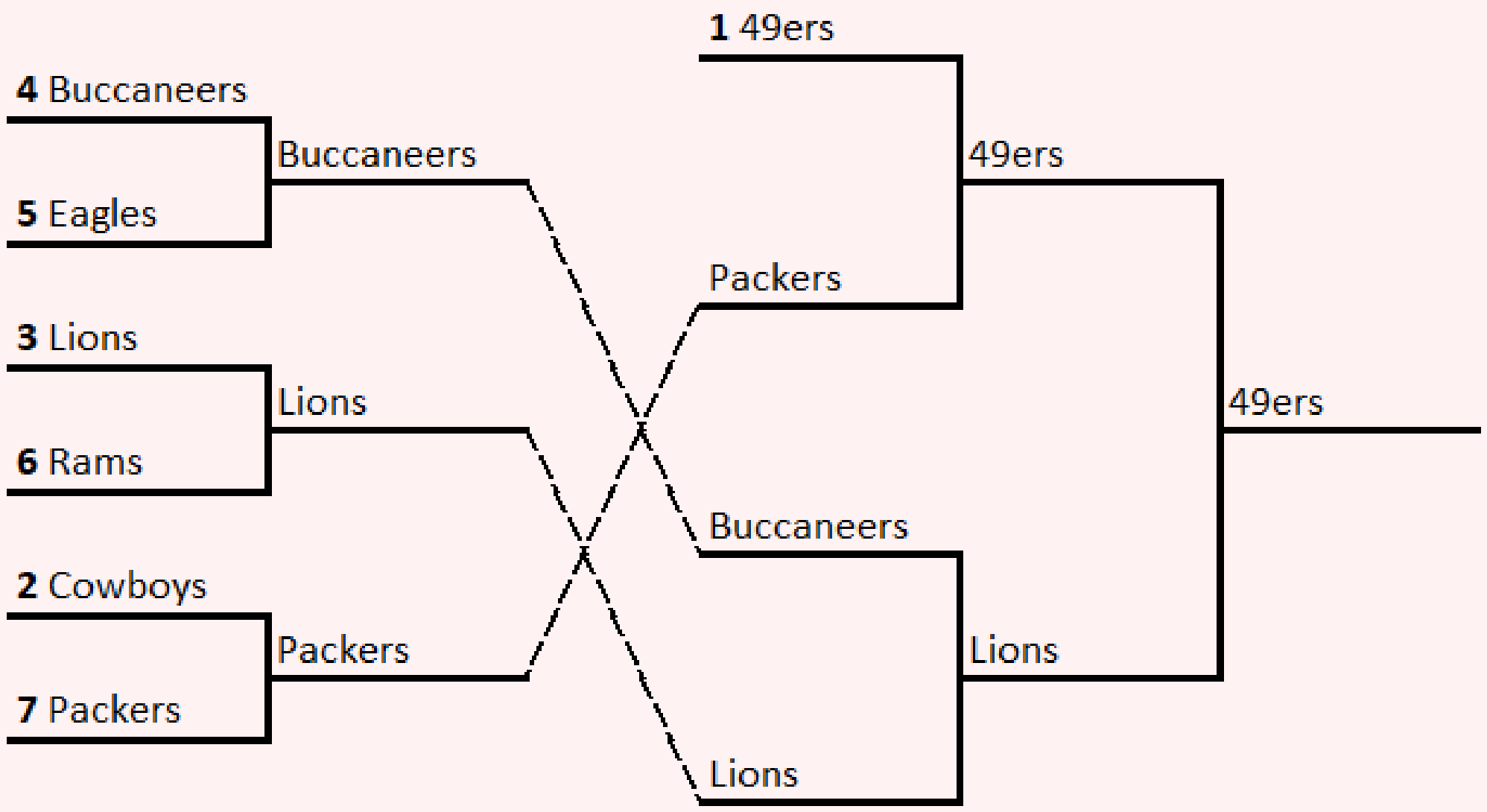

In the NFC, the 7-seed Packers upset the 2-seed Cowboys. Had a conventional bracket been used, the semifinal matchups would have been 1-seed vs 4-seed and 3-seed vs 7-seed: the 3-seed would have had an easier draw than the 1-seed, while the 7-seed would have had an easier draw than the 5-seed. Reseeding fixes this by matching the 7-seed Packers with top-seed 49ers, and the 3-seeded Lions with the 4-seeded Buccaneers.

Reseeding is a powerful technique. For one, the fundamental theorem still applies to reseeded brackets, allowing us to refer to reseeded brackets by their signatures as well.

**Theorem 2.5.7** **(Fried, 2024)**

There is exactly one proper reseeded bracket with each bracket signature.

*Proof.* The definition of properness ensures that there is only one way byes can be distributed in a proper reseeded bracket. Additionally, because reseeded brackets have no additional parameters beyond which seeds get how many byes, there is no more than one reseeded bracket with each signature that could be proper. Finally, that bracket is indeed proper: if the bracket goes to chalk, the matchups will be the exact same as in the proper traditional bracket of the same signature. □

But what about orderedness? It's intuitive to think that all proper reseeded are ordered, almost by definition: in fact, Hwang [12] published a proof of this for balanced reseeded brackets.

**Conjecture 2.5.8** **(Hwang, 1982)**

All balanced reseeded brackets are ordered.

We show for the first time that Hwang's proof is incorrect: neither the stronger claim that all proper reseeded brackets are ordered, nor Hwang's weaker claim are true. Our classification of the ordered reseeded brackets takes the same route as our proof of Edwards's Theorem: we first examine the orderedness of certain important brackets, and then we use the stapling and containment lemmas to specify the complete set of ordered reseeded brackets.

The proofs of the stapling and containment lemmas for reseeded brackets, as well as the fact that all ordered reseeded brackets are proper, are so similar to the corresponding proofs for traditional brackets that we just state them without proof.

**Theorem 2.5.9** **(Fried, 2024)**

All ordered reseeded brackets are proper.

**Lemma 2.5.10** **(Fried, 2024)**

If $\mathcal{A} = [[\mathbf{a_0}; ...; \mathbf{a_r}]]^R$ and $\mathcal{B} = [[\mathbf{b_0}; ...; \mathbf{b_s}]]^R$ are ordered reseeded brackets, then $\mathcal{C} = [[\mathbf{a_0}; ...; \mathbf{a_r} + \mathbf{b_0} - \mathbf{1}; ...; \mathbf{b_s}]]^R$ is an ordered reseeded bracket as well.

**Lemma 2.5.11** **(Fried, 2024)**

If $\mathcal{A}$ and $\mathcal{B}$ are reseeded brackets, $\mathcal{A}$ contains $\mathcal{B}$, and $\mathcal{B}$ is not ordered, then neither is $\mathcal{A}$.

We now examine particular brackets.

**Theorem 2.5.12** **(Fried, 2024)**

$[[\mathbf{1}]]^R$, $[[\mathbf{2}; \mathbf{0}]]^R$, and $[[\mathbf{4}; \mathbf{0}; \mathbf{0}]]^R$ are ordered.

*Proof.* No reseeding is done in a bracket of two or fewer rounds. Thus because the traditional proper brackets of these signatures are ordered, the reseeded brackets are as well. □

Our primary example of a reseeded bracket that is ordered despite the traditional bracket of the same signature not being ordered is $[[\mathbf{4}; \mathbf{2}; \mathbf{0}; \mathbf{0}]]^R$.

**Theorem 2.5.13** **(Fried, 2024)**

$[[\mathbf{4};\mathbf{2};\mathbf{0};\mathbf{0}]]^R$ is ordered.

*Proof.* Let $\mathcal{A} = [[\mathbf{4};\mathbf{2};\mathbf{0};\mathbf{0}]]^R$ and let $\mathcal{B} = [[\mathbf{4};\mathbf{0};\mathbf{0}]]^R = [[\mathbf{4};\mathbf{0};\mathbf{0}]]$. Then,

$$\begin{aligned}
\mathbb{W}_{\mathcal{A}}(t_1,\mathcal{T}) &= p_{36}p_{45}\mathbb{W}_{\mathcal{B}}(t_1,\{t_1,t_2,t_3,t_4\}) + p_{36}p_{54}\mathbb{W}_{\mathcal{B}}(t_1,\{t_1,t_2,t_3,t_5\}) + \\
&\quad p_{63}p_{45}\mathbb{W}_{\mathcal{B}}(t_1,\{t_1,t_2,t_4,t_6\}) + p_{63}p_{54}\mathbb{W}_{\mathcal{B}}(t_1,\{t_1,t_2,t_5,t_6\}) \\
&\geq p_{36}p_{45}\mathbb{W}_{\mathcal{B}}(t_2,\{t_1,t_2,t_3,t_4\}) + p_{36}p_{54}\mathbb{W}_{\mathcal{B}}(t_2,\{t_1,t_2,t_3,t_5\}) + \\
&\quad p_{63}p_{45}\mathbb{W}_{\mathcal{B}}(t_2,\{t_1,t_2,t_4,t_6\}) + p_{63}p_{54}\mathbb{W}_{\mathcal{B}}(t_2,\{t_1,t_2,t_5,t_6\}) \\
&= \mathbb{W}_{\mathcal{A}}(t_2,\mathcal{T})
\end{aligned}$$

$$\begin{aligned}
\mathbb{W}_{\mathcal{A}}(t_2,\mathcal{T}) &= p_{36}p_{45}\mathbb{W}_{\mathcal{B}}(t_1,\{t_1,t_2,t_3,t_4\}) + p_{36}p_{54}\mathbb{W}_{\mathcal{B}}(t_1,\{t_1,t_2,t_3,t_5\}) + \\
&\quad p_{63}p_{45}\mathbb{W}_{\mathcal{B}}(t_1,\{t_1,t_2,t_4,t_6\}) + p_{63}p_{54}\mathbb{W}_{\mathcal{B}}(t_1,\{t_1,t_2,t_5,t_6\}) \\
&\geq p_{36}p_{45}\mathbb{W}_{\mathcal{B}}(t_3,\{t_1,t_2,t_3,t_4\}) + p_{36}p_{54}\mathbb{W}_{\mathcal{B}}(t_3,\{t_1,t_2,t_3,t_5\}) \\
&= \mathbb{W}_{\mathcal{A}}(t_3,\mathcal{T})
\end{aligned}$$

Letting,

$$\begin{aligned}
a &= p_{36}p_{54}p_{32}p_{31} \\
b &= p_{36}p_{54}p_{32}p_{35} \\
c &= p_{63}p_{45}p_{42}p_{41} \\
d &= p_{63}p_{45}p_{42}p_{46} \\
e &= p_{36}p_{45}
\end{aligned}$$

we find,

$$\begin{aligned}
\mathbb{W}_{\mathcal{A}}(t_3,\mathcal{T}) &= p_{36}p_{54}p_{32}(p_{15}p_{31} + p_{51}p_{35}) + p_{36}p_{45}\mathbb{W}_{\mathcal{B}}(t_4,\{t_1,t_2,t_3,t_4\}) \\
&= p_{15}a + p_{51}b + e\mathbb{W}_{\mathcal{B}}(t_3,\{t_1,t_2,t_3,t_4\}) \\
&= p_{15}a + (p_{51} - p_{61})b + p_{61}b + e\mathbb{W}_{\mathcal{B}}(t_4,\{t_1,t_2,t_3,t_4\}) \\
&\geq p_{15}c + (p_{51} - p_{61})c + p_{61}d + e\mathbb{W}_{\mathcal{B}}(t_3,\{t_1,t_2,t_3,t_4\}) \\
&= p_{16}c + p_{61}d + e\mathbb{W}_{\mathcal{B}}(t_4,\{t_1,t_2,t_3,t_4\}) \\
&= p_{63}p_{45}p_{42}(p_{16}p_{41} + p_{61}p_{46}) + p_{45}p_{36}\mathbb{W}_{\mathcal{B}}(t_4,\{t_1,t_2,t_3,t_4\}) \\
&= \mathbb{W}_{\mathcal{A}}(t_4,\mathcal{T})
\end{aligned}$$

$$\begin{aligned}
\mathbb{W}_{\mathcal{A}}(t_4,\mathcal{T}) &= p_{36}p_{45}\mathbb{W}_{\mathcal{B}}(t_4,\{t_1,t_2,t_3,t_4\}) + p_{63}p_{45}\mathbb{W}_{\mathcal{B}}(t_4,\{t_1,t_2,t_4,t_6\}) \\
&\geq p_{36}p_{54}\mathbb{W}_{\mathcal{B}}(t_5,\{t_1,t_2,t_3,t_4\}) + p_{63}p_{54}\mathbb{W}_{\mathcal{B}}(t_5,\{t_1,t_2,t_4,t_6\}) \\
&= \mathbb{W}_{\mathcal{A}}(t_5,\mathcal{T})
\end{aligned}$$

Letting,

$$a = p_{36}p_{54}p_{51}p_{52}$$
$$b = p_{36}p_{54}p_{51}p_{53}$$
$$c = p_{63}p_{45}p_{61}p_{62}$$
$$d = p_{63}p_{45}p_{61}p_{64}$$
$$e = p_{63}p_{54}$$

we find,

$$\begin{aligned}
\mathbb{W}_{\mathcal{A}}(t_5, \mathcal{T}) &= p_{36}p_{54}p_{51}(p_{23}p_{52} + p_{32}p_{53}) + p_{54}p_{63}\mathbb{W}_{\mathcal{B}}(t_5, \{t_1, t_2, t_5, t_6\}) \\
&= p_{23}a + p_{32}b + e\mathbb{W}_{\mathcal{B}}(t_5, \{t_1, t_2, t_5, t_6\}) \\
&= p_{23}a + (p_{32} - p_{42})b + p_{42}b + e\mathbb{W}_{\mathcal{B}}(t_5, \{t_1, t_2, t_5, t_6\}) \\
&\geq p_{23}c + (p_{32} - p_{42})c + p_{42}d + e\mathbb{W}_{\mathcal{B}}(t_6, \{t_1, t_2, t_5, t_6\}) \\
&= p_{23}c + p_{32}d + e\mathbb{W}_{\mathcal{B}}(t_6, \{t_1, t_2, t_5, t_6\}) \\
&= p_{63}p_{45}p_{61}(p_{24}p_{62} + p_{42}p_{64}) + p_{63}p_{54}\mathbb{W}_{\mathcal{B}}(t_5, \{t_1, t_2, t_5, t_6\}) \\
&= \mathbb{W}_{\mathcal{A}}(t_6, \mathcal{T})
\end{aligned}$$

Thus $\mathcal{A}$ is ordered. □

Unfortunately, that is where the power of reseeding to convert non-ordered signatures into ordered ones ends. The following two signatures are not ordered.

**Theorem 2.5.14** **(Fried, 2024)**

$[[\mathbf{6}; \mathbf{1}; \mathbf{0}; \mathbf{0}]]^R$ is not ordered.

*Proof.* Let $\mathcal{A} = [[\mathbf{6}; \mathbf{1}; \mathbf{0}; \mathbf{0}]]^R$, and let $\mathcal{T}$ have the following matchup table.

| | $t_1$ | $t_2$ | $t_3$ | $t_4$ | $t_5$ | $t_6$ | $t_7$ |
|---|---|---|---|---|---|---|---|
| $t_1$ | | | | | | | |
| $t_2$ | $p$ | | | | | | |
| $t_3$ | $p$ | $p$ | | | | | |
| $t_4$ | $p$ | $p$ | 0.5 | | | | |
| $t_5$ | $p$ | $p$ | 0.5 | 0.5 | | | |
| $t_6$ | $p$ | $p$ | $p$ | 0.5 | 0.5 | | |
| $t_7$ | $p$ | $p$ | $p$ | 0.5 | 0.5 | 0.5 | |

For $t_6$ to win the format, three probability $p$ upsets must occur: $t_6$ beating $t_3$ in the first round, $t_6$ beating $t_1$ in the second round, and someone beating $t_2$. Thus,

$$\mathbb{W}_{\mathcal{A}}(t_6, \mathcal{T}) = O(p^3).$$

But for $t_7$ to win the format, only two probability $p$ upsets are necessary: $t_7$ beating $t_2$ in the first round and $t_7$ beating $t_1$ in the second round, as the winner of $t_4$ vs $t_5$ might beat $t_3$ in the semifinals. Thus,

$$\mathbb{W}_{\mathcal{A}}(t_7, \mathcal{T}) = 0.25p^2 + O(p^3).$$

So for small enough $p$, $\mathbb{W}_{\mathcal{A}}(t_6, \mathcal{T}) < \mathbb{W}_{\mathcal{A}}(t_7, \mathcal{T})$, so $\mathcal{A}$ is not monotonic with respect to $\mathcal{T}$ and thus not ordered. □

**Theorem 2.5.15** **(Fried, 2024)**

$[[\mathbf{4}; \mathbf{2}; \mathbf{2}; \mathbf{0}; \mathbf{0}]]^R$ is not ordered.

*Proof.* Let $\mathcal{A} = [[\mathbf{4}; \mathbf{2}; \mathbf{2}; \mathbf{0}; \mathbf{0}]]^R$, and let $\mathcal{T}$ have the following matchup table.

| | $t_1$ | $t_2$ | $t_3$ | $t_4$ | $t_5$ | $t_6$ | $t_7$ | $t_8$ |
|---|---|---|---|---|---|---|---|---|
| $t_1$ | | | | | | | | |
| $t_2$ | $p^2$ | | | | | | | |
| $t_3$ | $p^2$ | 0.5 | | | | | | |
| $t_4$ | $p^2$ | 0.5 | 0.5 | | | | | |
| $t_5$ | $p^2$ | $p$ | $p$ | 0.5 | | | | |
| $t_6$ | $p^2$ | $p$ | $p$ | $p$ | $p$ | | | |
| $t_7$ | $p^2$ | $p^2$ | $p$ | $p$ | $p$ | $p$ | | |
| $t_8$ | $p^2$ | $p^2$ | $p$ | $p$ | $p$ | $p$ | 0.5 | |

For $t_7$ to win the format on the order of probability $p^5$, they must face $t_6$ in the first round, $t_3$ in the second round, $t_1$ in the semifinal, and $t_4$ (crucially not $t_2$) in the final. This happens when $t_4$ survives $t_5$ and defeats $t_2$, which has probability on the order of 0.25. Thus,

$$\mathbb{W}_{\mathcal{A}}(t_7, \mathcal{T}) = 0.25p^5 + O(p^6).$$

Similarly, for $t_8$ to win the format on the order of probability $p^5$, they must face $t_5$ in the first round, $t_3$ in the second round, $t_1$ in the semifinal, and $t_4$ (again, crucially not $t_2$) in the final. However, because $t_4$ will be playing $t_6$ in the second round rather than $t_5$, they will almost certainly win, meaning that $t_4$ advances to the final with probability on the order of 0.5. Thus,

$$\mathbb{W}_{\mathcal{A}}(t_8, \mathcal{T}) = 0.5p^5 + O(p^6).$$

So for small enough $p$, $\mathbb{W}_{\mathcal{A}}(t_7, \mathcal{T}) < \mathbb{W}_{\mathcal{A}}(t_8, \mathcal{T})$, so $\mathcal{A}$ is not monotonic with respect to $\mathcal{T}$ and thus not ordered. □

Recapping,

**Figure 2.5.16: Which Proper Reseeded Brackets are Ordered**

| Ordered | Not Ordered |
|---|---|
| $[[\mathbf{1}]]^R$ | $[[\mathbf{6};\mathbf{1};\mathbf{0};\mathbf{0}]]^R$ |
| $[[\mathbf{2};\mathbf{0}]]^R$ | $[[\mathbf{4};\mathbf{2};\mathbf{2};\mathbf{0};\mathbf{0}]]^R$ |
| $[[\mathbf{4};\mathbf{0};\mathbf{0}]]^R$ | |
| $[[\mathbf{4};\mathbf{2};\mathbf{0};\mathbf{0}]]^R$ | |

Finally, we apply the stapling and containment lemmas to complete the theorem.

**Theorem 2.5.17** **(Fried, 2024)**

The ordered reseeded brackets are exactly those corresponding to signatures that can be generated in the following way.

1. Start with the list $[[\mathbf{0}]]^R$ (note that this is not yet a bracket signature).
2. As many times as desired, prepend the list with $[[\mathbf{1}]]$, $[[\mathbf{3};\mathbf{0}]]$, or $[[\mathbf{3};\mathbf{2};\mathbf{0}]]$.
3. Then, add 1 to the first element in the list, turning it into a bracket signature.

*Proof.* The stapling lemma, combined with the fact that $[[\mathbf{1}]]^R$, $[[\mathbf{2};\mathbf{0}]]^R$, $[[\mathbf{4};\mathbf{0};\mathbf{0}]]^R$, and $[[\mathbf{4};\mathbf{2};\mathbf{0};\mathbf{0}]]^R$ are ordered, ensure that any reseeded brackets generated by the above procedure is indeed ordered. Left is to use the containment lemma to ensure that these are the only ones.

Let $\mathcal{A}$ be a bracket signature that cannot be generated by the procedure. Then, either there is a round in which three or more games are to be played, or there is a round in which exactly two games are played and the next two rounds each have exactly two games played as well.

Let $i$ be the latest such round. If round $i$ is the first of three rounds with two games each, then round $i+3$ must have only one game played (otherwise $i$ would not be the latest such round). But then $\mathcal{A}$ contains $[[\mathbf{4};\mathbf{2};\mathbf{2};\mathbf{0};\mathbf{0}]]^R$, and so is not ordered.

If round $i$ has three or more games, then round $i+1$ must contain exactly two games (any less and not every winner would have a game, any more and $i$ would not be the latest such round.) Then, if round $i+2$ has one game, then $\mathcal{A}$ contains $[[\mathbf{6};\mathbf{1};\mathbf{0};\mathbf{0}]]^R$, and if it has two, then $\mathcal{A}$ contains $[[\mathbf{4};\mathbf{2};\mathbf{2};\mathbf{0};\mathbf{0}]]^R$. In either case, $\mathcal{A}$ is not ordered.

Thus, the ordered reseeded brackets are exactly those generated by the procedure. □

So, the space of ordered reseeded brackets is slightly larger than the space of ordered traditional brackets, although perhaps this is not quite as much of an expansion as we would have liked or expected. Despite this, reseeded brackets definitely *feel* more ordered than traditional brackets of the same signature, even if neither is ordered in the definitional sense.

**Open Question 2.5.18** **(Fried, 2024)**

Is there some sense in which reseeded brackets that are not ordered are closer to being ordered than their traditional bracket analogues?

In the meantime, reseeding remains an important tool in our tournament design toolkit, though it is not without its drawbacks, as discussed by Baumann, Matheson, and Howe [4].

In a reseeded bracket, teams and spectators alike don't know who they will play or where their next game will be until the entire previous round is complete. This can be an especially big issue if parts of the bracket are being played in different locations on short turnarounds: in the NCAA Basketball Tournament, for example, the first two rounds are played over a weekend at various pre-determined locations. It would cause problems if teams had to pack up and travel across the country because they got reseeded and their opponent and thus location changed.

In addition, part of what makes the NCAA Basketball Tournament (affectionately known as "March Madness") such a fun spectator experience is the fact that these matchups are known ahead of time. In "bracket pools," groups of fans each fill out their own brackets, predicting who will win each game and getting points based on how many they get right. If it wasn't clear where in the bracket the winner of a given game was supposed to go, this experience would be diminished.

Finally, reseeding gives the top seed(s) an even greater advantage than they already have: instead of playing against merely the *expected* lowest-seeded team(s) each round, they would get to play against the *actual* lowest-seeded team(s). In March Madness, "Cinderella Stories," that is, deep runs by low seeds, would become much less common.

In many ways, the NFL playoffs is the perfect place to use a reseeded bracket: games are played once a week, giving plenty of time for travel; only seven teams make the playoffs in each, so a huge March Madness-style bracket challenge is unlikely; as a professional league, protecting Cinderella Stories isn't as important; and because the bracket is only three rounds long, reseeding is only required once. Somewhat ironically, the NFL conference playoffs used to employ the format $[[\mathbf{4}; \mathbf{2}; \mathbf{0}; \mathbf{0}]]^R$, which is ordered, but have since allowed a seventh team from each conference into the playoffs and changed to the non-ordered $[[\mathbf{6}; \mathbf{1}; \mathbf{0}; \mathbf{0}]]^R$ [31].

Other leagues with similar structures might consider adopting forms of reseeding to protect their incentives and competitive balance, but in many cases, the traditional bracket structure is too appealing to adopt a reseeded one.

## 2.6 Randomization

Given that reseeding doesn't solve the orderedness problem presented by Edwards's Theorem, we turn to a new approach at generating potentially ordered knockout tournaments: randomization.

**Definition 2.6.1: Totally Randomized Knockout Tournament (Unattributed)**

A *totally randomized knockout tournament* is a bracket, except the teams are randomly placed onto the starting lines instead of being placed according to seed.

Clearly totally randomized knockout tournaments are indeed knockout tournaments.

Chung and Hwang [6] conjectured that all totally randomized knockout tournaments were ordered. After all, the teams are all being treated identically: how could a better team be at a disadvantage relative to a worse one?

**Conjecture 2.6.2 (Chung and Hwang, 1978)**

All totally randomized knockout tournaments are ordered.

Indeed, Lemma 2.6.3, proved by Chen and Hwang [5], seems to provide some evidence for the conjecture.

**Lemma 2.6.3 (Chen and Hwang, 1988)**

Let $\mathcal{A}$ be a totally randomized knockout tournament with signature $[[\mathbf{a_0}; \ldots; \mathbf{a_r}]]$, let $\mathcal{S}$ be a set of teams, and let $\mathcal{T}$ be the set of teams produced by replacing a given team $s \in \mathcal{S}$ with a team $t$ such that for all other teams $u$,

$$\mathbb{P}[t \text{ beats } u] \geq \mathbb{P}[s \text{ beats } u].$$

Then,

$$\mathbb{W}_{\mathcal{A}}(t, \mathcal{T}) \geq \mathbb{W}_{\mathcal{A}}(s, \mathcal{S}).$$

*Proof.* Let $X$ be the power set of $\mathcal{S} \setminus \{s\} = T \setminus \{t\}$, and for each set of teams $Y \in X$, let $P_Y$ be the probability that $s$ or $t$ will have to beat exactly the set of teams $Y$ in order to win the tournament (noting that this probability is the same for $s$ and $t$).

Then,

$$\begin{aligned}\mathbb{W}_{\mathcal{A}}(t, \mathcal{T}) &= \sum_{Y \in X} \left( P_Y \cdot \prod_{u \in Y} \mathbb{P}[t \text{ beats } u] \right) \\ &\geq \sum_{Y \in X} \left( P_Y \cdot \prod_{u \in Y} \mathbb{P}[s \text{ beats } u] \right) \\ &= \mathbb{W}_{\mathcal{A}}(s, \mathcal{S})\end{aligned}$$

□

Unfortunately, despite the lemma, Chung and Hwang's conjecture is false due to a counterexample given by Israel [13]. (The colors are not part of the bracket itself and just used to aid the proof.)

**Figure 2.6.4: Setup of Theorem 2.6.5**

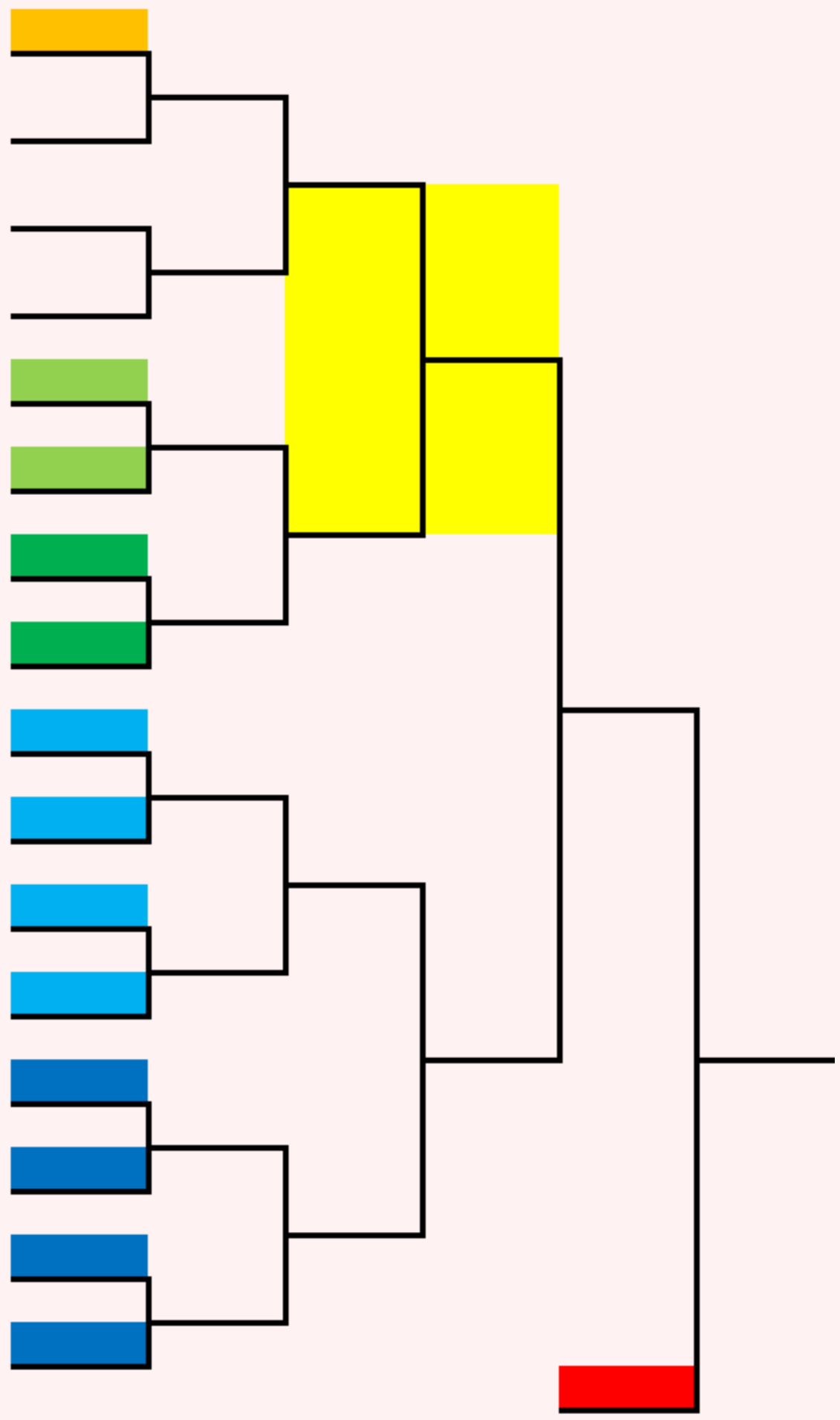

**Theorem 2.6.5** **(Israel, 1981)**

The totally randomized knockout tournament of signature $[[\mathbf{16}; \mathbf{0}; \mathbf{0}; \mathbf{0}; \mathbf{1}; \mathbf{0}]]$ is not ordered.

*Proof.* Let $\mathcal{A}$ be the format in question, and let $\mathcal{T}$ be the list of seventeen teams containing one copy of each of $t_1, t_3, t_4$, and $t_5$, and thirteen copies of $t_2$, with the following matchup table.

| | $t_1$ | $t_2$ | $t_3$ | $t_4$ | $t_5$ |
|---|---|---|---|---|---|
| $t_1$ | | | | | |
| $t_2$ | 0.5 | | | | |
| $t_3$ | 0 | $2p$ | | | |
| $t_4$ | 0 | $p$ | $2p$ | | |
| $t_5$ | 0 | $p$ | $p$ | 0.5 | |

Let $i \in \{4, 5\}$ and let $j = 9 - i$. For $t_i$ to win $\mathcal{A}$ without getting placed on the red starting line, it must win at least four games against teams $t_1$, $t_2$, or $t_3$, which happens with probability $O(p^4)$. Thus we let $\mathcal{B}_i$ be the format identical to $\mathcal{A}$ except we enforce that $t_i$ will be placed on the red starting line and note that

$$\mathbb{W}_{\mathcal{A}}(t_i, \mathcal{T}) = \frac{1}{17}\mathbb{W}_{\mathcal{B}_i}(t_i, \mathcal{T}) + O(p^4).$$

Now, $t_j$ reaches the finals of $\mathcal{B}_i$ with probability $O(p^4)$, $t_3$ reaches the finals of $\mathcal{B}_i$ with probability $O(p^3)$ and so $t_i$ beats them in the finals with probability $O(p^4)$, and of course $t_i$ cannot beat $t_1$ in the finals. Thus,

$$\mathbb{W}_{\mathcal{B}_i}(t_i, \mathcal{T}) = p \cdot \mathbb{P}[t_2 \text{ reaches the finals of } \mathcal{B}_i] + O(p^4).$$

Since $t_3$ and $t_j$ reach the finals of $\mathcal{B}_i$ with probability $O(p^3)$ and $O(p^4)$ respectively,

$$\mathbb{W}_{\mathcal{B}_i}(t_i, \mathcal{T}) = p \cdot \mathbb{P}[t_1 \text{ doesn't reach the finals of } \mathcal{B}_i] + O(p^4).$$

Assume now without loss of generality that $t_1$ gets placed on the orange starting line.

Any difference in $\mathbb{P}[t_1 \text{ doesn't reach the finals of } \mathcal{B}_i]$ between $i \in \{4, 5\}$ will have to come as a result of a game involving $t_j$ (as $t_j$ is the only difference in $t_1$'s route to the finals between $\mathcal{B}_4$ and $\mathcal{B}_5$), and because $t_4$ and $t_5$ have the same probability of beating every team other than $t_3$, it will have to be as a result of a game against $t_3$. However, because neither $t_3$ nor $t_j$ can beat $t_1$, in order to play each other in a game whose winner doesn't immediately play $t_1$, they will have to be placed on two colored starting lines of the same color.

If $t_3$ and $t_j$ are placed on two of the light blue or dark blue starting lines, then any difference in $\mathbb{P}[t_1 \text{ doesn't reach the finals of } \mathcal{B}_i]$ between $i \in \{4, 5\}$ will be induced by $t_j$ winning its first three games, with happens with probability $O(p^3)$.

However, if $t_3$ and $t_j$ are placed on the two dark green or two light green starting lines, then when $i = 4$, $t_1$ will play $t_2$ in the yellow game with probability

$$p_{35}p_{23} + p_{53}p_{25} = ((1-p)(1-2p) + (p)(1-p)) = 1 - 2p + p^2,$$

while when $i = 5$, $t_1$ will play $t_2$ in the yellow game with probability

$$p_{34}p_{23} + p_{43}p_{24} = ((1-2p)(1-2p) + (2p)(1-p)) = 1 - 2p + 2p^2.$$

Thus,

$$\begin{aligned}&\mathbb{P}[t_1 \text{ plays } t_2 \text{ in the yellow game of } \mathcal{B}_5]\\ &- \mathbb{P}[t_1 \text{ plays } t_2 \text{ in the yellow game of } \mathcal{B}_4]\\ &= cp^2 + O(p^3)\end{aligned}$$

for some constant $c$, so

$$\begin{aligned}&\mathbb{P}[t_1 \text{ doesn't reach the finals of } \mathcal{B}_5]\\ &- \mathbb{P}[t_1 \text{ doesn't reach the finals of } \mathcal{B}_4]\\ &= cp^2 + O(p^3)\end{aligned}$$

for some constant $c$, so

$$\mathbb{W}_{\mathcal{B}_5}(t_5, \mathcal{T}) - \mathbb{W}_{\mathcal{B}_4}(t_4, \mathcal{T}) = cp^3 + O(p^4)$$

for some constant $c$, so

$$\mathbb{W}_{\mathcal{A}}(t_5, \mathcal{T}) - \mathbb{W}_{\mathcal{A}}(t_4, \mathcal{T}) = cp^3 + O(p^4)$$

for some constant $c$.

Therefore $\mathcal{A}$ is not monotonic with respect to $\mathcal{T}$ and so $\mathcal{A}$ is not ordered. □

Chung and Hwang's conjecture was rescued by Chen and Hwang [5] who restricted the domain of the claim to balanced formats.

**Theorem 2.6.6** **(Chen and Hwang, 1988)**

All totally randomized balanced knockout tournaments are ordered.

*Proof.* Let $\mathcal{A}_r$ be the totally randomized balanced knockout tournament on $2^r$ teams. We proceed by induction on $r$. Clearly the one-team format $\mathcal{A}_0$ is ordered. For any other $r$, let $\mathcal{T}$ be a list of teams, and let $t_i$ and $t_j$ be teams such that $i < j$.

Let $\mathcal{B}_r$ be the totally randomized balanced knockout tournament on $2^r$ teams except $t_i$ and $t_j$ are forced to play each other in the first round, and let $\mathcal{C}_r$ be the totally randomized balanced knockout tournament on $2^r$ teams except $t_i$ and $t_j$ cannot play each other in the first round. Then,

$$\mathbb{W}_{A_r}(t_i, \mathcal{T}) = \left(\frac{1}{2^r - 1}\right) \mathbb{W}_{B_r}(t_i, \mathcal{T}) + \left(\frac{2^r - 2}{2^r - 1}\right) \mathbb{W}_{C_r}(t_i, \mathcal{T})$$

and likewise for $t_j$.

Because $p_{ij} \geq p_{ji}$, and by Lemma 2.6.3, $\mathbb{W}_{B_r}(t_i, \mathcal{T}) \geq \mathbb{W}_{B_r}(t_j, \mathcal{T})$. Left is to show that $\mathbb{W}_{C_r}(t_i, \mathcal{T}) \geq \mathbb{W}_{C_r}(t_j, \mathcal{T})$.

For two other teams $t_a$ and $t_b$, let $M_{ab}$ be the set of $2^{r-1} - 2$ team subsets of $\mathcal{T} \setminus \{t_i, t_j, t_a, t_b\}$, and for $\mathcal{S} \in M_{ab}$, let $P_\mathcal{S}$ be the probability that the teams in $\mathcal{S}$ all win their first-round games and none of them play any of $t_i, t_j, t_a$, or $t_b$ in the first round. Now,

$$\begin{aligned}
\mathbb{W}_{C_r}(t_i, \mathcal{T}) &= \frac{1}{2} \sum_{t_a, t_b \in \mathcal{T} \setminus \{t_i, t_j\}} \sum_{\mathcal{S} \in M_{ab}} P_\mathcal{S} \cdot ( \quad (p_{ia}p_{jb} + p_{ib}p_{ja}) \cdot \mathbb{W}_{A_{r-1}}(t_i, \mathcal{S} \cup \{t_i, t_j\}) \\
&\qquad + p_{ia}p_{bj} \cdot \mathbb{W}_{A_{r-1}}(t_i, \mathcal{S} \cup \{t_i, t_b\}) \\
&\qquad + p_{ib}p_{aj} \cdot \mathbb{W}_{A_{r-1}}(t_i, \mathcal{S} \cup \{t_i, t_a\})) \\
&\geq \frac{1}{2} \sum_{t_a, t_b \in \mathcal{T} \setminus \{t_i, t_j\}} \sum_{\mathcal{S} \in M_{ab}} P_\mathcal{S} \cdot ( \quad (p_{ja}p_{ib} + p_{jb}p_{ia}) \cdot \mathbb{W}_{A_{r-1}}(t_j, \mathcal{S} \cup \{t_i, t_j\}) \\
&\qquad + p_{ja}p_{bi} \cdot \mathbb{W}_{A_{r-1}}(t_j, \mathcal{S} \cup \{t_j, t_b\}) \\
&\qquad + p_{jb}p_{ai} \cdot \mathbb{W}_{A_{r-1}}(t_j, \mathcal{S} \cup \{t_j, t_a\})) \\
&= \mathbb{W}_{C_r}(t_i, \mathcal{T})
\end{aligned}$$

The inequality follows by comparing each term to its corresponding term: the $\mathbb{W}_{A_{r-1}}(t_i, \mathcal{S} \cup \{t_i, t_j\})$ term inequality is by induction, while the other two terms are by Lemma 2.6.3.

Thus, $\mathcal{A}_r$ is ordered. □

In some ways this is a great revelation: we finally have an ordered balanced knockout tournament for arbitrary numbers of rounds. Of course, this orderedness does not come without drawbacks.

For one, the randomization feels a bit cheap: once the randomization is complete, before any games have even been played, the orderedness is lost. (Compare to the ordered traditional and reseeded brackets, which maintain their orderedness throughout the whole tournament.)

And secondly, totally randomness has the undesirable property that it might make for some very lopsided and anti-climatic knockout tournaments. It could be that top-two teams, whom everyone wants to see face off in the championship game, are set to play each other in the first round! We can extend the notion of *dramatic* from brackets to knockout tournaments, noting that totally randomized knockout tournaments are not dramatic.

**Definition 2.6.7: Dramatic Knockout Tournament** **(Fried, 2024)**

A knockout tournament is *dramatic* if, as long as the knockout tournament goes chalk, in every round, the $m$ remaining teams are guaranteed to be the top $m$ seeds.

To fix this, we define a new class of randomized knockout tournaments: *cohort randomized knockout tournaments*, first defined by Schwenk [16].

**Definition 2.6.8: Cohort Randomized Knockout Tournament** **(Schwenk, 2000)**

The $r$-round *cohort randomized knockout tournament* is the balanced bracket on $2^r$ teams, except, for each $i$, seeds $2^i + 1$ through $2^{i+1}$ are shuffled randomly before play.

Thus the 1- and 2-seeds are locked into their places, the 3- and 4-seeds exchange places half the time, seeds 5-8 are randomly shuffled, and as are 9-16, 17-32, etc.

**Theorem 2.6.9** **(Schwenk, 2000)**

Cohort randomized knockout tournaments are dramatic.

*Proof.* We proceed by induction on $r$. If $r = 0$, then there are no rounds and so the theorem holds. For any other $r$, in the first round, the top $2^{r-1}$ seeds will face the bottom $2^{r-1}$ seeds, and because the format goes chalk, the bottom half of teams will be eliminated. Thus after the first round, the top $2^{r-1}$ seeds will remain. The remaining format is just the $(r-1)$-round cohort randomized knockout tournament, for which the theorem holds by induction. □

In Schwenk's paper, he wrote that in cohort randomized knockout tournaments, "higher-seeded teams are never given a schedule more difficult than that of any lower seed." By this Schwenk meant that if the format went chalk, higher-seeded teams would never play a team from a higher cohort than lower-seeded teams, a property very similar to properness. But what about being ordered? It seems as though cohort randomized knockout tournaments ought to be ordered: being in a higher cohort seems preferable to being a lower cohort, and teams in the same cohort are treated identically.

Unfortunately, like many other formats we've seen thus far, cohort randomized knockout tournaments are not (for more than two rounds) ordered.

**Figure 2.6.10: Setup of Theorem 2.6.11**

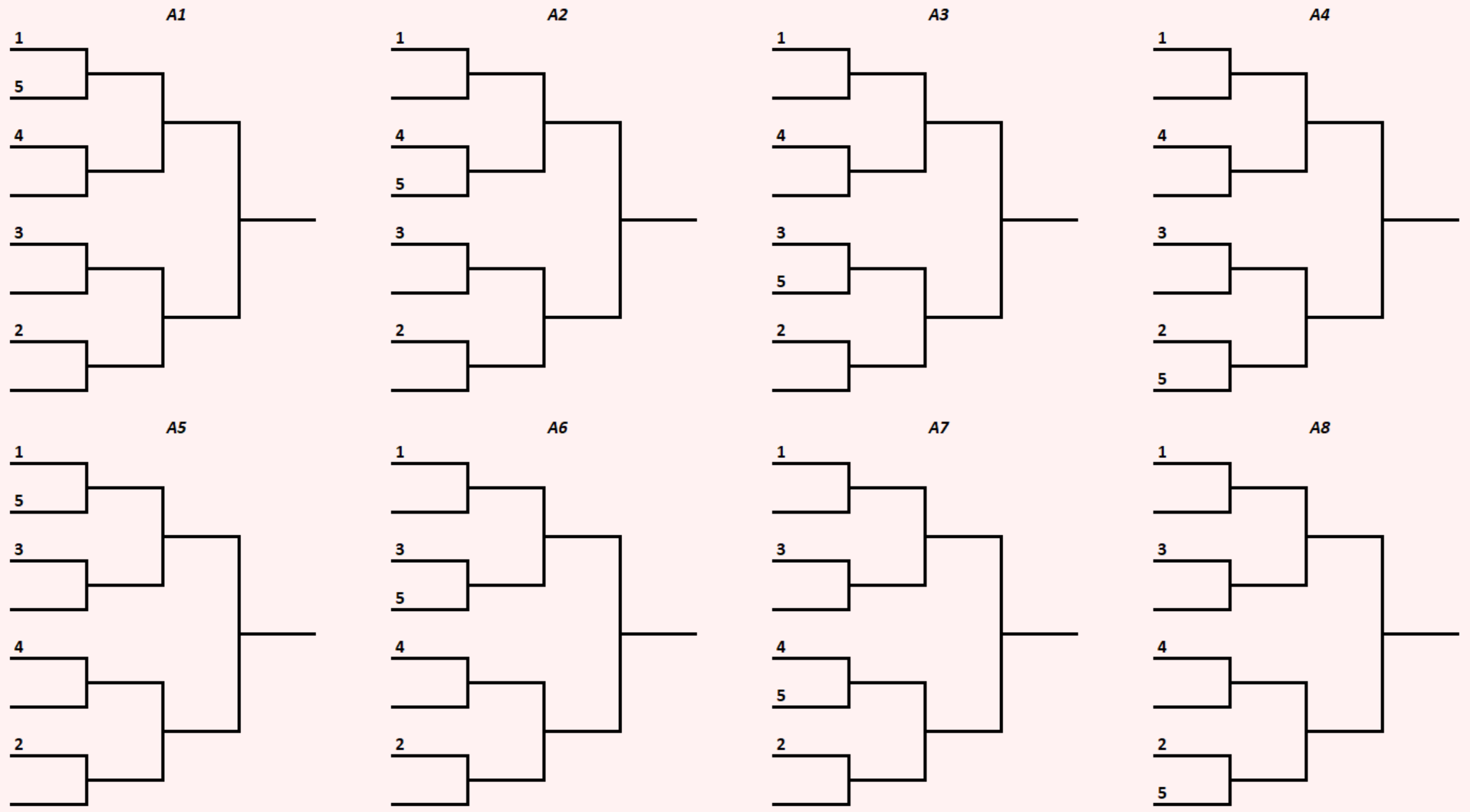

**Theorem 2.6.11** **(Fried, 2024)**

The eight-team cohort randomized knockout tournament is not ordered.

*Proof.* Let $\mathcal{A}$ be the eight-team cohort randomized knockout tournament, and let $\mathcal{T}$ have the following matchup table for $0 < p < 0.5$.

| | $t_1$ | $t_2$ | $t_3$ | $t_4$ | $t_5$ | $t_6$ | $t_7$ | $t_8$ |
|---|---|---|---|---|---|---|---|---|
| $t_1$ | | | | | | | | |
| $t_2$ | 0.5 | | | | | | | |
| $t_3$ | 0.5 | 0.5 | | | | | | |
| $t_4$ | 0.5 | 0.5 | 0.5 | | | | | |
| $t_5$ | $p$ | 0.5 | 0.5 | 0.5 | | | | |
| $t_6$ | $p$ | $p$ | $p$ | $p$ | 0.5 | | | |
| $t_7$ | $p$ | $p$ | $p$ | $p$ | 0.5 | 0.5 | | |
| $t_8$ | $p$ | $p$ | $p$ | $p$ | 0.5 | 0.5 | 0.5 | |

Note that because $t_6, t_7$, and $t_8$ each have identical matchups against every other team,

permutations of those teams don't affect the probability of any other teams winning the tournament. Thus we consider the eight possible randomizations in Figure 2.6.10 noting that for $i \in \{2, 3\}$,

$$\mathbb{W}_{\mathcal{A}}(t_i, \mathcal{T}) = \frac{1}{8} \sum_{j=1}^{8} \mathbb{W}_{\mathcal{A}_j}(t_i, \mathcal{T}).$$

Some calculation finds,

$$\begin{aligned}
\mathbb{W}_{\mathcal{A}_1}(t_3, \mathcal{T}) &= \mathbb{W}_{\mathcal{A}_1}(t_2, \mathcal{T}) \\
\mathbb{W}_{\mathcal{A}_2}(t_3, \mathcal{T}) &= \mathbb{W}_{\mathcal{A}_2}(t_2, \mathcal{T}) \\
\mathbb{W}_{\mathcal{A}_3}(t_3, \mathcal{T}) &= \mathbb{W}_{\mathcal{A}_4}(t_2, \mathcal{T}) \\
\mathbb{W}_{\mathcal{A}_4}(t_3, \mathcal{T}) &= \mathbb{W}_{\mathcal{A}_3}(t_2, \mathcal{T}) \\
\mathbb{W}_{\mathcal{A}_5}(t_3, \mathcal{T}) &= \mathbb{W}_{\mathcal{A}_7}(t_2, \mathcal{T}) \\
\mathbb{W}_{\mathcal{A}_6}(t_3, \mathcal{T}) &= \mathbb{W}_{\mathcal{A}_8}(t_2, \mathcal{T}) \\
\mathbb{W}_{\mathcal{A}_8}(t_3, \mathcal{T}) &= \mathbb{W}_{\mathcal{A}_6}(t_2, \mathcal{T})
\end{aligned}$$

However, letting $q = 1 - p$, $r = \frac{1}{2}q + \frac{1}{4}$, and $s = pq + \frac{1}{2}q$,

$$\begin{aligned}
\mathbb{W}_{\mathcal{A}_5}(t_2, \mathcal{T}) &= qs\left(q\frac{1}{2} + p(pr + qs)\right) \\
&< qs\left(q\frac{1}{2} + p\left(\frac{1}{2}r + \frac{1}{2}s\right)\right) \qquad \text{because } r < s \text{ and } p < \frac{1}{2} < q \\
&= \mathbb{W}_{\mathcal{A}_7}(t_3, \mathcal{T})
\end{aligned}$$

Therefore,

$$\mathbb{W}_{\mathcal{A}}(t_2, \mathcal{T}) < \mathbb{W}_{\mathcal{A}}(t_3, \mathcal{T})$$

so $\mathcal{A}$ is not monotonic with respect to $\mathcal{T}$ and thus not ordered. □

If cohort randomized knockout tournaments don't solve the orderedness problem, why would we use them over traditional proper knockout tournaments? (A close variant of) cohort randomization is most famously found on the ATP Tour [19], a collection of tournaments played by professional tennis players that all use large balanced knockout tournaments. Additionally, the seeding for these tournaments is set by the ATP rankings, which tend to be slow to update. As a result, if every ATP Tour tournament used the proper seeding, the 6-seed and 27-seed would play each other in the first round at every tournament until one of them moved up or moved down. These rematches were deemed undesirable and so this randomization procedure was introduced: The 1-seed's quarterfinals matchup (if everything goes chalk) is now randomly drawn from the 5- through 8-seeds, instead of always being the 8-seed.

But Theorem 2.6.11 tells us that they are not ordered, meaning that the only ordered balanced knockout format we've developed for more than two rounds is the totally randomized one, which is neither deterministic nor dramatic. Unfortunately, we conclude the chapter without a more satisfying design, leaving behind two big open questions.

**Open Question 2.6.12** **(Fried, 2024)**

For all $r$, does there exist an $r$-round deterministic ordered balanced knockout tournament?

**Open Question 2.6.13** **(Fried, 2024)**

For all $r$, does there exist an $r$-round dramatic ordered balanced knockout tournament?

We (pessimistically) conjecture that both answers are no.

# 3 Multibrackets

## 3.1 Consolation Brackets

In the previous chapter, we discussed brackets and knockout tournaments, paying attention only to which team is declared champion. Edwards's Theorem and its analogies make claims only about which teams are most likely to win the tournament: all participants that don't win are grouped together as losers. Real tournaments, however, do not always operate in this way: explicitly or not, the team that lost the championship game is often considered to have earned the second-place finish, for example.

Third-place is often harder to determine. If a team is given a bye all the way to the finals, and thus there is only one semifinal, then the loser of that semifinal can be unambiguously granted third. The 2023 Korean Baseball Organization League Playoffs [30] have this property: they use a bracket of signature $[[\mathbf{2}; \mathbf{1}; \mathbf{1}; \mathbf{1}; \mathbf{0}]]$, and so could easily assign third-place to the NC Dinos, who lost in the sole semifinal. (The LG Twins won the format, and finals' losers KT Wiz came in second.)

**Figure 3.1.1: 2023 KBO League Playoffs**

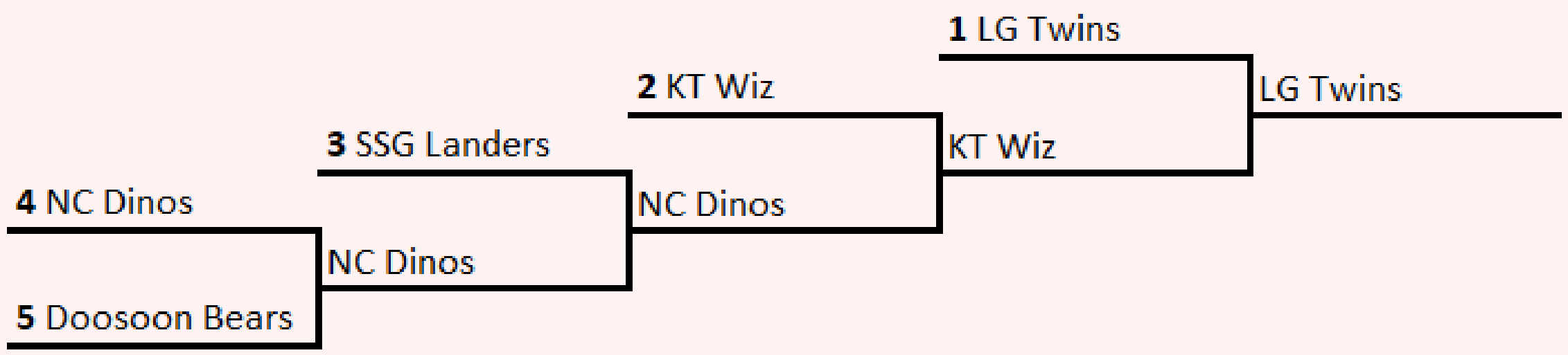

But in most brackets (those brackets whose signature's penultimate digit is a zero), assigning third-place is trickier: there are two teams who lost in the semifinal and have an equal claim to the place. There are a number of strategies that a league might use in the face of this ambiguity.

The first option is to just not assign a third-place at all. In the wise words of Will Ferrell from Talladega Nights [26], "If you ain't first, you're last." Who cares who came in third: you didn't win, you didn't even come in second, so you lost. This approach is taken by all four major American professional sports leagues (the NFL, NBA, NHL, and MLB).

The second option is to declare the two semifinal losers co-third-place finishers. In many ways, this is the same as the first option, but with a single sentence added to the end of a press-release indicating that the teams in question each finished third. (This option also has the unsatisfying property that four teams will be able to claim a top-three finish. This can be easily fixed, however, by just granting both teams fourth-place instead.)

The third option is to use some (relatively) arbitrary tiebreaker to select the third-place team. A few potential such tiebreakers are: whichever team was seeded higher, whichever team lost to the tournament champion (as opposed to the tournament runner-up), or if the teams played each other during the "regular season" portion of a tournament, whichever team won that game.

None of these are particularly satisfying. While they may do alright when giving out third-place isn't important, if we really want to assign third in a fair and equitable way, say because there is a bronze medal or spot in a future tournament up for grabs, these options will not do.

Instead, the best thing to do is play a game: The 2015 Asian Football Confederation Asian Cup [20] did exactly that.

**Figure 3.1.2: 2015 AFC Asian Cup**

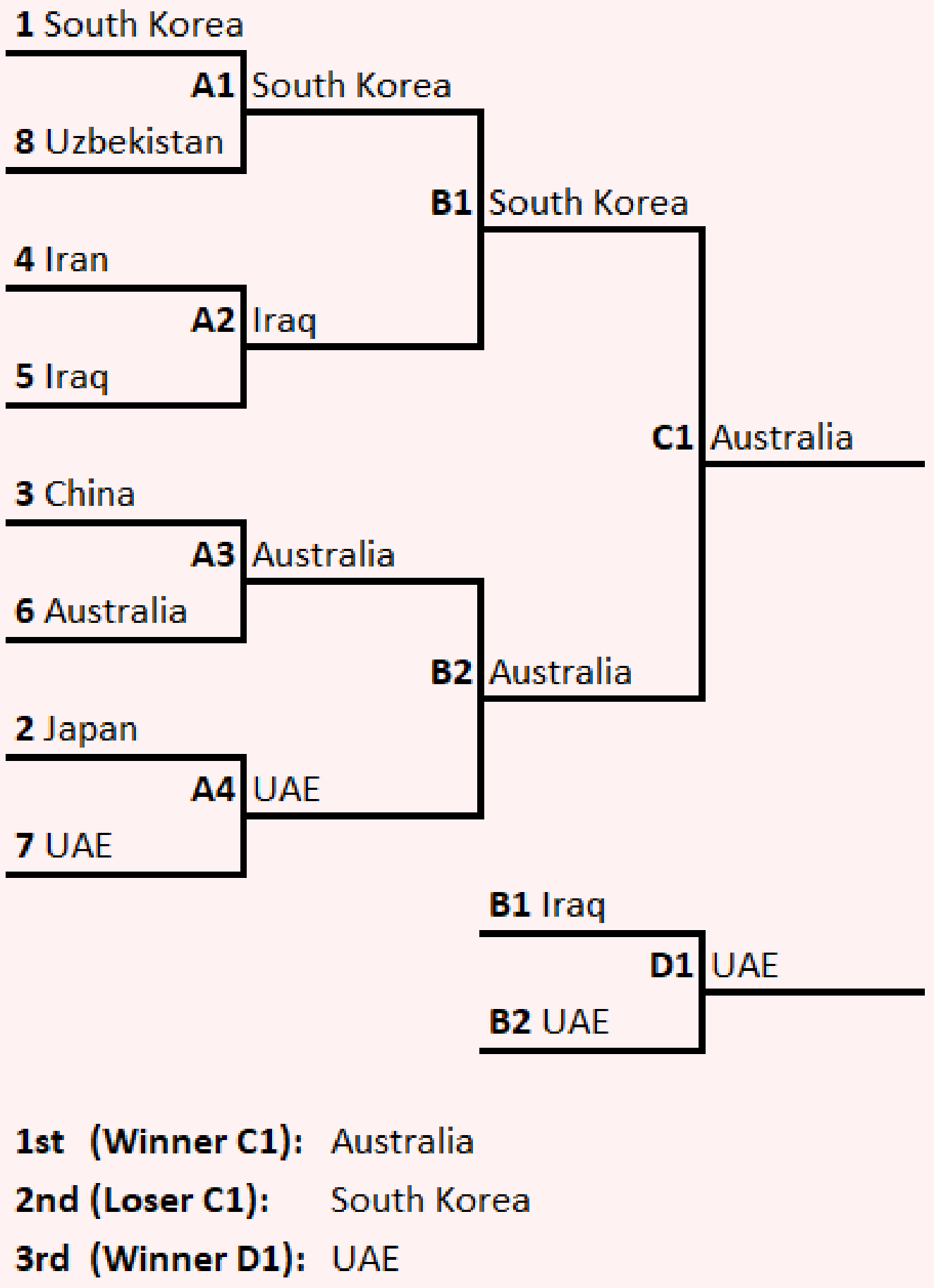

In the 2015 AFC Asian Cup, after the main bracket is complete, with the winner of the final game (Australia) being crowned champion and the loser (South Korea) coming in second, the two semifinal losers (Iraq and the UAE) are matched up in the third-place game.

A quick note about Figure 3.1.2: each game in the figure is labeled. In the primary bracket, first-round games are **A1** through **A4**, while the semifinals are **B1** and **B2**, and the finals is game **C1**. The third-place game is labeled **D1**: even though it could be played concurrently to the championship game, it is part of a different bracket and so we label it as a different round.

We indicate that the third-place game is to be played between the losers of games **B1**

and **B2** by labeling the starting lines in the third-place game with those games. This is not ambiguous because the winners of those games always continue on in the original bracket, so such labels only refer to the losers.

The third-place game, which can also be viewed as a two-team bracket of signature $[[\mathbf{2}; \mathbf{0}]]$, in an example of a *consolation bracket*.

**Definition 3.1.3: Consolation Bracket** **(Unattributed)**

A *consolation bracket* is a bracket in which teams that did not win the tournament compete for an $m$th-place finish for some $m$.

Consolation brackets are as opposed to *primary brackets*.

**Definition 3.1.4: Primary Bracket** **(Unattributed)**

A *primary bracket* is a bracket in a multibracket the winner of which is declared champion.

The formats from the previous chapter, then, consist only of a primary bracket with no consolation ones. The third-place game, however, as used by the 2015 AFC Asian Cup, is a common and well-liked consolation bracket used for selecting the top-three teams. Of course, it is far from the only way that the AFC could have handed out gold, silver, and bronze.

In fact, it's not clear that the loser of **C1**, who comes in second place, is really more deserving than the winner of **D1**, who comes in third. The UAE might argue: South Korea and we both finished with two wins and one loss – a first-round win, a win against Iraq, and a loss against Australia. The only reason that South Korea came in second and we came in third was because South Korea lucked out by having Australia on the other half of the bracket. That's not fair!

If the AFC took this complaint seriously, it could modify the format to add a game **E1** for second-place to be played between the loser of **C1** and the winner of **D1**, with the loser coming in third.

**Figure 3.1.5: 2015 AFC Asian Cup Alternative**

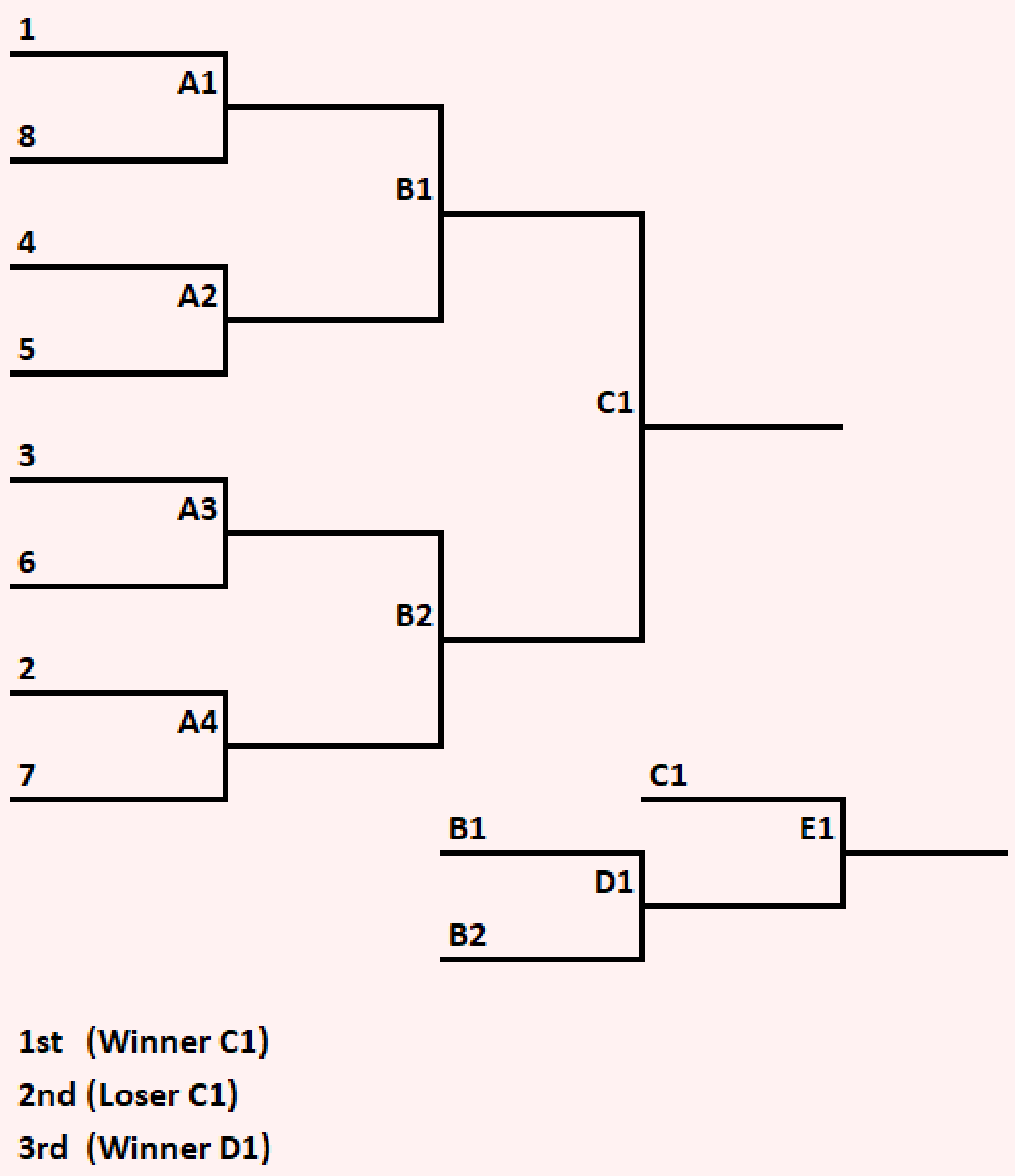

If the AFC used the format in Figure 3.1.5 in 2015, then South Korea and the UAE would have played each other for second place after all of the other games were completed. In some sense, this is a more equitable format than the one used in reality: we have the same data about the UAE and South Korea and so we ought to let them play for second-place instead of deciding almost randomly.

However, swapping formats doesn't come without costs. For one thing, South Korea and the UAE would've had to play a fourth game: if the AFC had only three days to put on the tournament and teams can play at most one game a day, then the format in Figure 3.1.5 isn't feasible.

Another concern: suppose that Iraq had beaten the UAE when they played in game **D1**. Then the two teams with a claim to second-place would have been South Korea and Iraq, except South Korea had already beaten Iraq! It seems a bit unfair to South Korea to make the teams play again, this time for stakes. One option is to say "tough luck, later games being more important than earlier ones is a staple in sports." But another is to designate game **E1** as *contingent*.

**Definition 3.1.6: Contingent Game** **(Fried, 2024)**

A game in a tournament format is *contingent* if, under certain circumstances, (most commonly if the teams have already played earlier in the tournament) the game is skipped and the result of a previous game is used.

Ultimately, whether game **E1** should be included or not depends on the purpose of the tournament. If there is a huge difference between the prizes for coming in second and third – for instance, if the top two finishing teams in the Asian Cup qualified for the World Cup – then **E1** is quite important. If, on the other hand, this is a self-contained format played purely for bragging rights, **E1** could probably be left out. In reality, the 2015 AFC Asian Cup qualified only its winner to another tournament (the 2017 Confederations Cup), and gave medals to its top three, so game **E1**, which distinguishes between second- and third-place, was probably unnecessary.

What if instead of just the champion, the top four teams from the Asian Cup advanced to the Confederations Cup. One could imagine an easy extension of the format used presently, in which the loser of game **D1** is awarded fourth-place, to determine the four teams that qualify. However, this format would be quite poor: whether or not a team qualifies for the Confederations Cup would be solely determined by the result of their first-round game and so the **B**, **C**, and **D**-round games might as well not even be played. A better format for selecting the top-four would allow first-round losers to win their way back onto the podium, as was employed by the 2023 Southern Conference Wrestling Championships [24].

**Figure 3.1.7: 2023 Southern Conference Wrestling Championships**

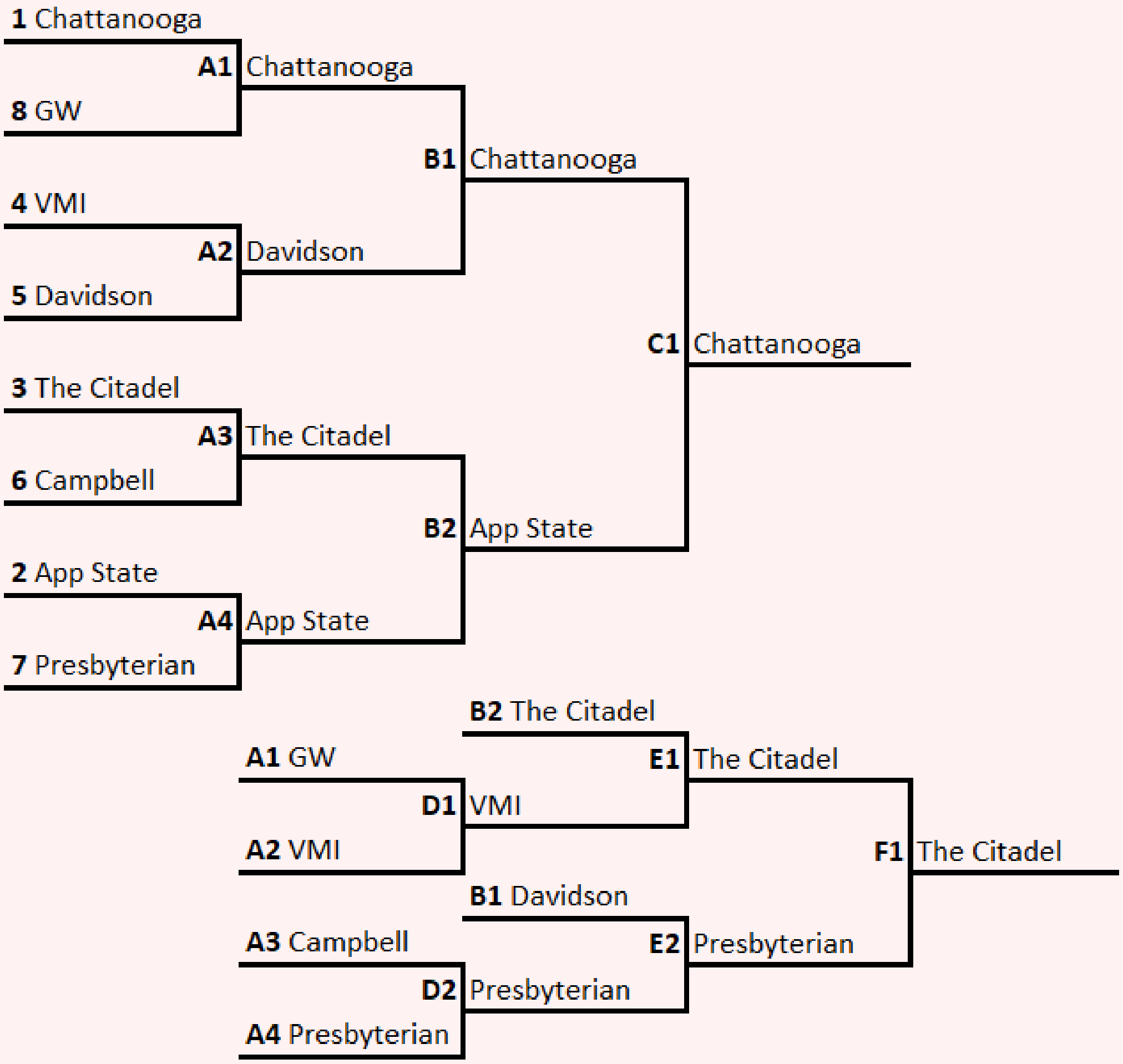

**1st (Winner C1):** Chattanooga
**2nd (Loser C1):** App State
**3rd (Winner F1):** The Citadel
**4th (Loser F1):** Presbyterian

The format in Figure 3.1.7 is a dramatic improvement for selecting a top-four over that in Figure 3.1.2. In the 2023 Southern Conference Wrestling Championships, teams finish in the top-four if and only if they win two games before they lose two, which is a nice property to have. The one downside is that it takes a fourth round: if there is not enough time for a fourth round, or if there is safety risk to teams playing four matches in a row, the format isn't feasible. Although if we only care about the top-four, and not the specifics of which team came in third or in fourth, we could drop game **F1**, ensuring that each team plays at most three games.

**Figure 3.1.8: 2023 SoCon Wrestling Championships Alternative**

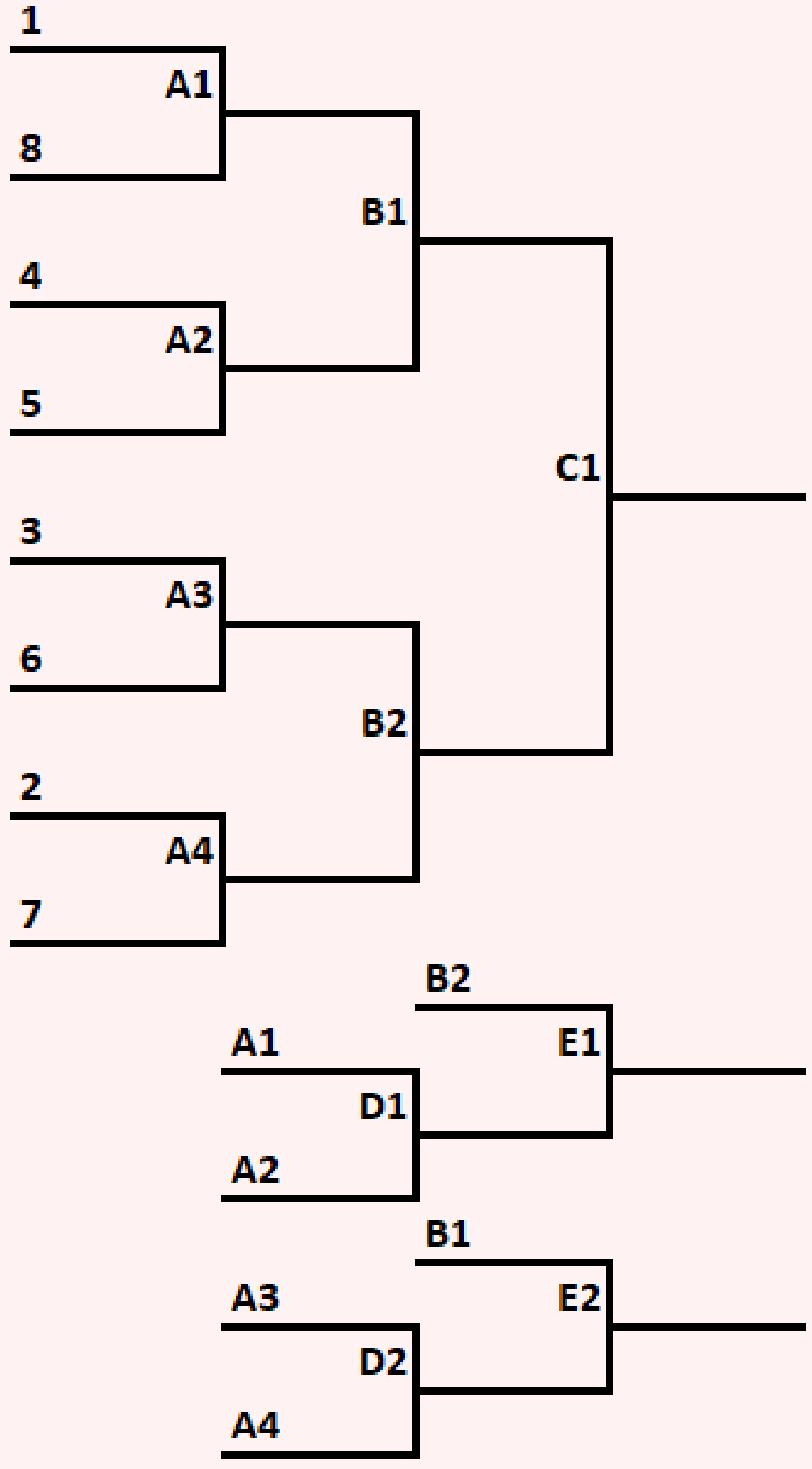

1st (Winner C1)
2nd (Loser C1)
4th (Winner E1, Winner E2)

(As discussed earlier, we opt to rank both the **E**-round winners in fourth, to ensure that no more than $m$ teams can claim a top-$m$ finish for any $m$.)

The four formats with consolation brackets presented thus far are examples of *multibrackets*.

**Definition 3.1.9: Multibracket** **(Fried, 2024)**

A *multibracket* is a collection of one or more brackets coupled with a specification of which winners and losers of which games receive which places. Starting lines in multibrackets can be marked by a seed, or by a game, indicating that the loser of the specified game should be placed there, but no seed or game can be placed on more than one starting line.

Since which game each team plays in next (and which place each team ends up in) can be derived only from which game that team has played in most recently and whether they won or lost that game, this definition is equivalent to saying that the format is networked.

**Definition 3.1.10: Multibracket** **(Fried, 2024)**

A *multibracket* is a networked tournament format.

(Note that this means multibrackets with contingent games are technically not multibrackets at all. However, they are close enough to being multibrackets and are important enough tools for tournament design that we include them in our discussion, in the same way that in the last chapter we discussed reseeded brackets even though they are technically not brackets.)

We will see in the coming sections that many formats used in a variety of settings are actually just examples of multibrackets. Figure 3.1.11 gives an outline of what the space of multibrackets looks like: we will spent the rest of the chapter examining the various categories in more detail.

**Figure 3.1.11: The Space of Multibrackets**

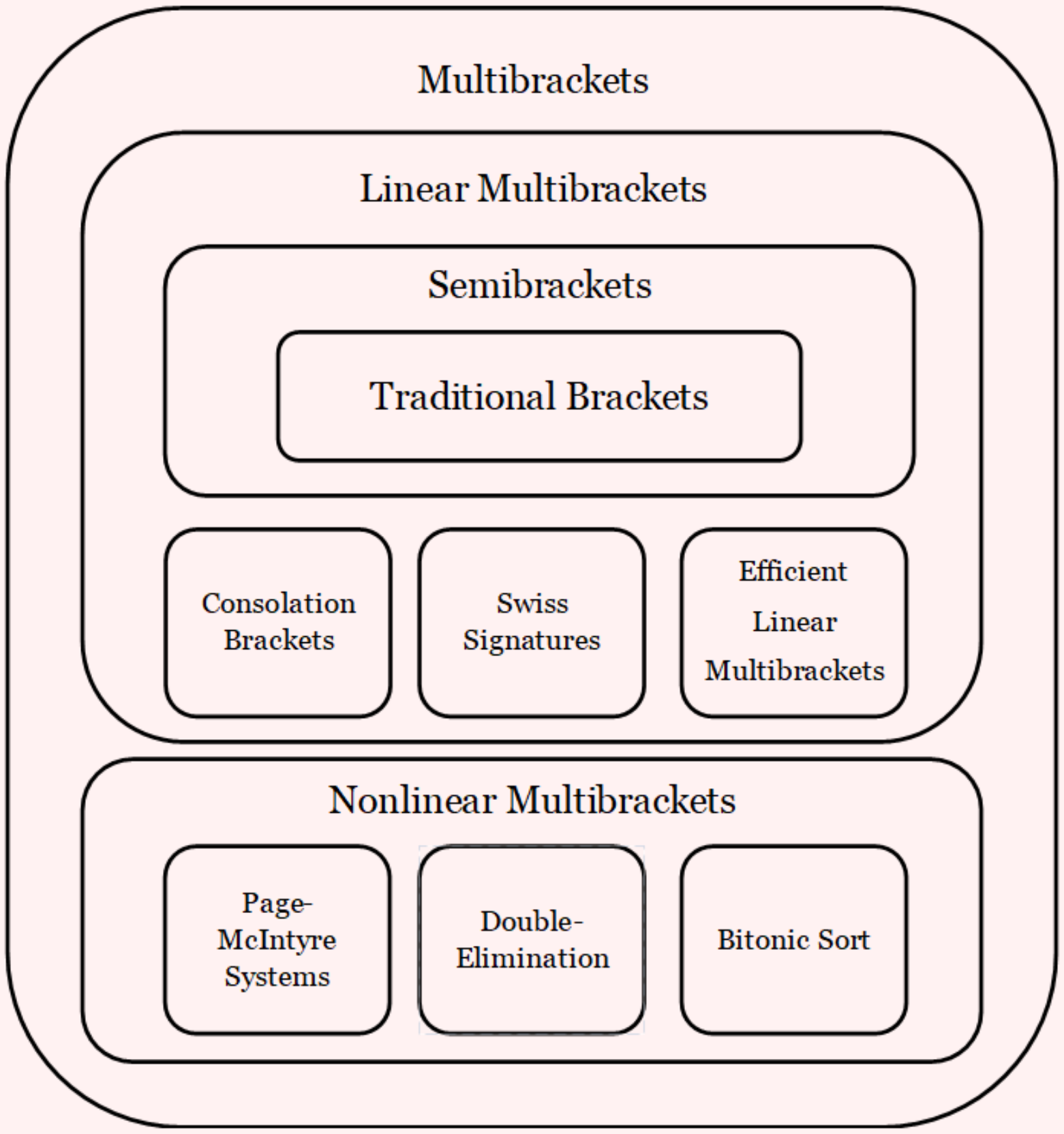

## 3.2 Semibrackets

In this section, we move on from consoloation brackets to focus on semibrackets, which as indicated by Figure 3.1.11 is a generalization of the traditional bracket. We will then use the notion of a semibracket to define linear multibrackets, which we will study for a few sections before addressing nonlinear multibrackets at the end of the chapter.

Consider now the following tournament design problem: we are tasked with designing an eight-team tournament to select the top two teams who will go on to compete in the national tournament. However, there's only enough time for two rounds: perhaps due to field space or team fatigue, each team can only play two games. What design should we use?

The most natural answer to this question is to use a traditional eight-team bracket, but leave the championship game unplayed. This format is displayed in the figure below.

**Figure 3.2.1: $[[8;0;0;0]]$ with no Championship Game**

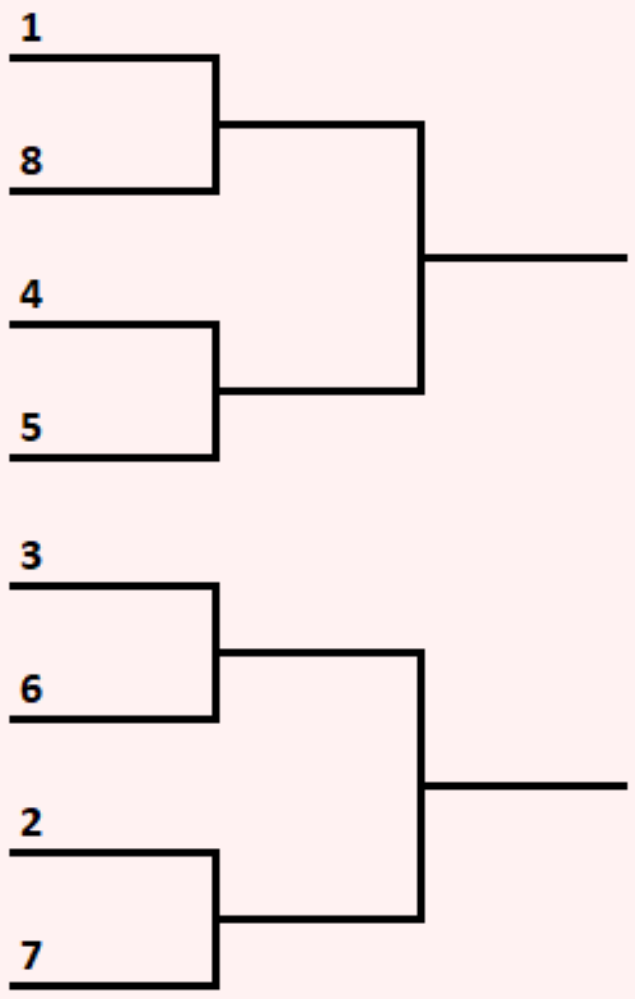

The format in Figure 3.2.1 does exactly what we need. The championship game being left unplayed is not a bug but a feature: each team plays a maximum of two games, and the two teams that advance to the national tournament are clear.

While it would be reasonable to describe the format in Figure 3.2.1 as two brackets that run side-by-side, it would be nice to be able to describe it as a single format: a bracket in which the championship game is left unplayed.

**Definition 3.2.2: Semibracket** **(Fried, 2024)**

A *semibracket* is a networked format in which

(a) Teams don't play any games after their first loss, and

(b) All teams that finish with no losses are declared co-champions.

Thus semibrackets are a generalization of brackets: a bracket is a semibracket in which only one team is left undefeated and declared champion.

Figure 3.2.3 describes which properties various networked formats require.

**Figure 3.2.3: Properties of Networked Formats**

| Format | No Games After First Loss | Only One Team Finishes Undefeated |
|---|---|---|
| Bracket | ✔ | ✔ |
| Semibracket | ✔ | ✘ |
| Multibracket | ✘ | ✘ |

The format in Figure 3.2.1 is not a particularly interesting example of a semibracket: it is just a traditional bracket minus one game. Are there any examples of semibrackets that are not just traditional brackets with some rounds left uncompleted?

Indeed there are. Let's modify the original problem so that we need to pick a top three teams out of twelve. Again, no team can play more then two games. The natural choice is shown below in Figure 3.2.4.

**Figure 3.2.4: A More Interesting Semibracket**

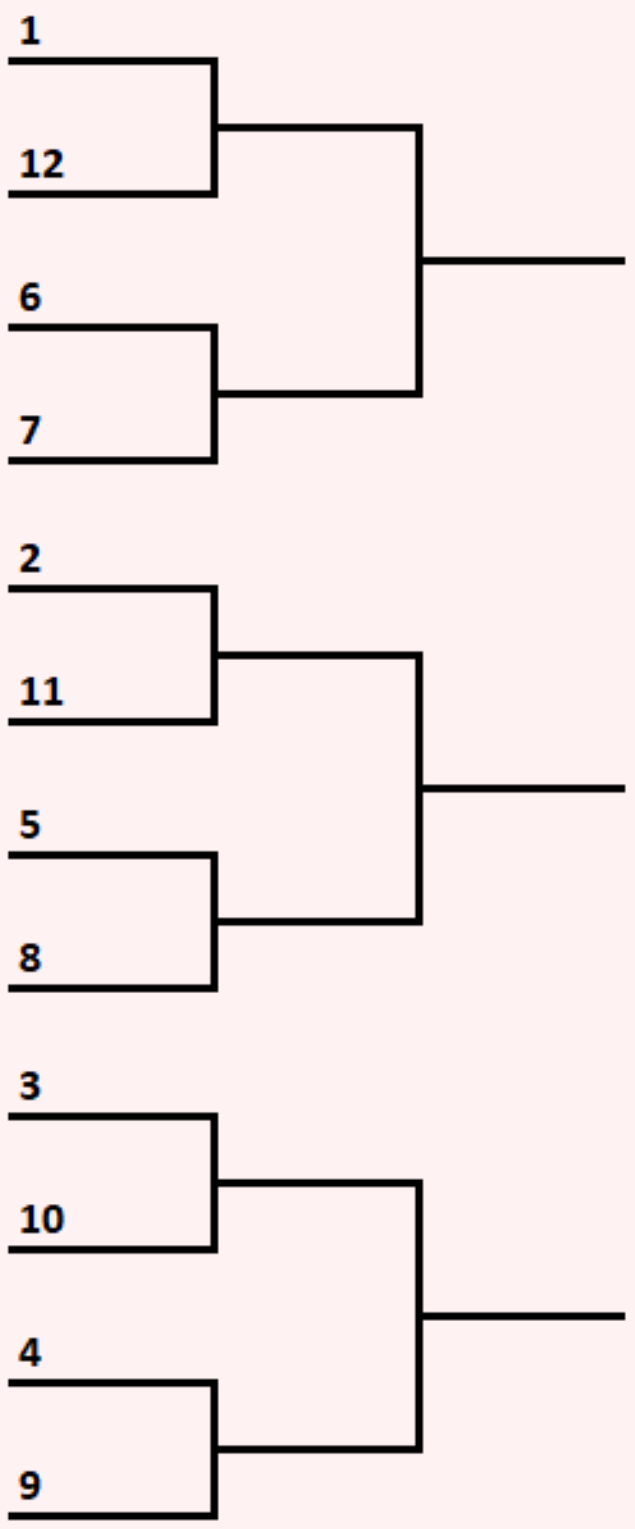

There is no potential for the format in Figure 3.2.4 to be completed into a traditional bracket: the next round would include an odd number of teams. But as a semibracket, this is still a viable format, one that nicely solves the tournament design problem that we were given.

**Definition 3.2.5: Rank of a Semibracket** **(Fried, 2024)**

The *rank* of a semibracket is how many co-champions it crowns. If the semibracket $\mathcal{A}$ has rank $m$, we say $\text{Rank}(\mathcal{A}) = m$ or that $\mathcal{A}$ *ranks* $m$ *teams.*

Traditional brackets are exactly the semibrackets that rank one team. The formats in Figures 3.2.1 and 3.2.4 rank two and three teams, respectively.

We can adapt the concept of a bracket signature to semibrackets.

**Definition 3.2.6: Semibracket Signature** **(Fried, 2024)**

The *signature* of an $r$-round semibracket $\mathcal{A}$ is the list $[[\mathbf{a_0};...;\mathbf{a_r}]]_m$, where $a_i$ is the number of teams that get $i$ byes and $m = \text{Rank}(\mathcal{A})$. (In the case where $m = \text{Rank}(\mathcal{A}) = 1$, it can be omitted.)

Thus the signature of traditional brackets are the same as when they are viewed as semibrackets that rank one team. The signatures of the formats in Figures 3.2.1 and 3.2.4 are $[[\mathbf{8};\mathbf{0};\mathbf{0}]]_2$ and $[[\mathbf{12};\mathbf{0};\mathbf{0}]]_3$, respectively.

In analogy with traditional bracket signature's Theorem 2.1.14, we have Theorem 3.2.7.

**Theorem 3.2.7** **(Fried, 2024)**

Let $\mathcal{A} = [[\mathbf{a_0};...;\mathbf{a_r}]]_m$ be a list of natural numbers. Then $\mathcal{A}$ is a semibracket signature if and only if

$$\sum_{i=0}^{r} a_i \cdot \left(\frac{1}{2}\right)^{r-i} = m.$$

The proof is almost identical to that of Theorem 2.1.14 so we leave it out for brevity. Likewise, properness can be defined in the same way for semibracket, and the fundamental theorem still applies (again with a nearly identical proof that is left out for brevity).

**Theorem 3.2.8** **(Fried, 2024)**

There is exactly one proper semibracket with each semibracket signature.

Semibrackets are used in practice in situations where the excitement of a single elimination tournament is desired, but multiple winners are needed. The 2023 Union of European Football Associations Champions League Qualifying Phase [22], for example, used a (somewhat randomized) semibracket of signature $[[\mathbf{4};\mathbf{0};\mathbf{29};\mathbf{9};\mathbf{8};\mathbf{2};\mathbf{0}]]_6$ to determine the final six teams that would qualify for the group stage.

**Figure 3.2.9: 2023 UEFA Champions League Qualifying Phase**

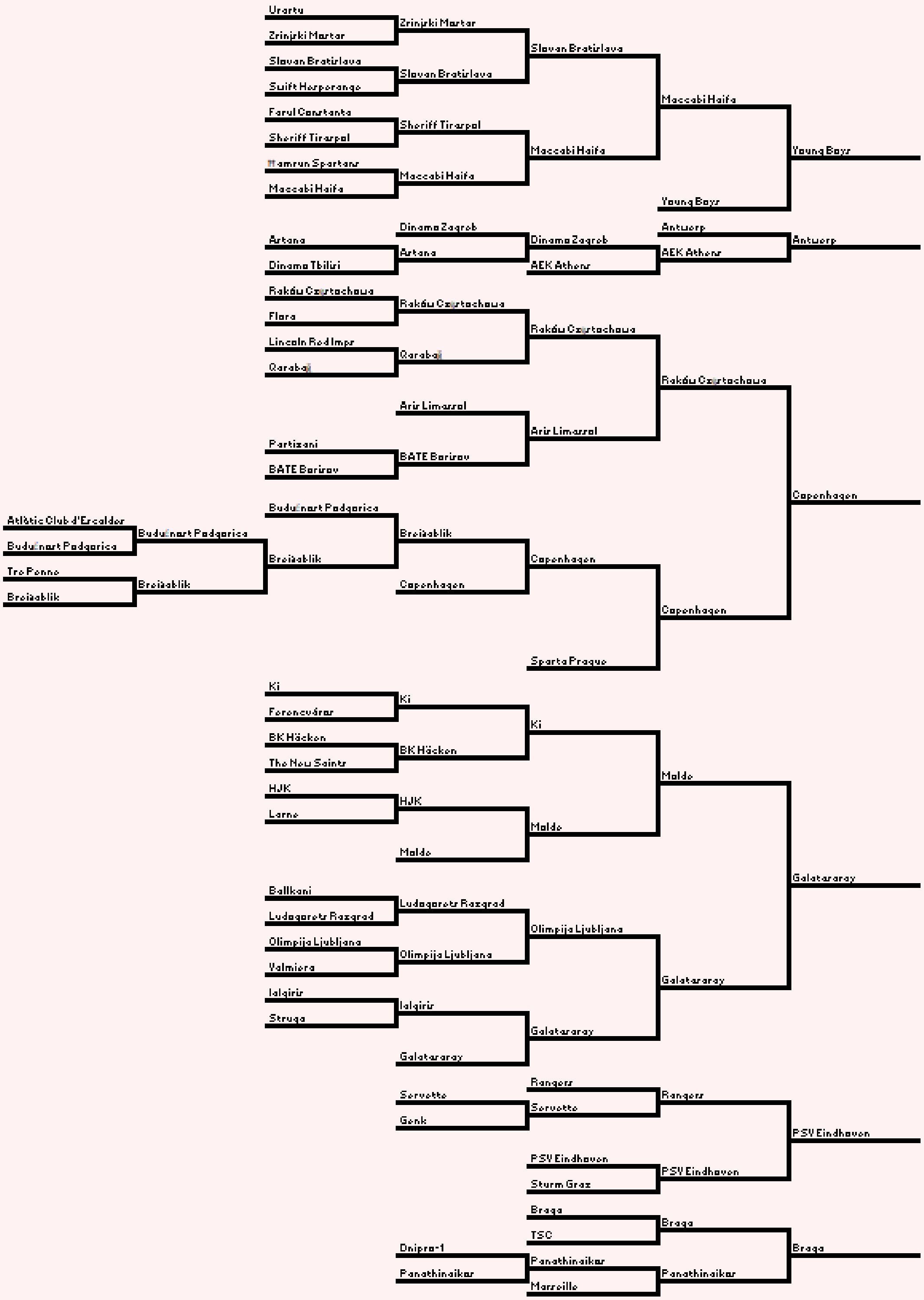

Finally, we give a few descriptors to describe certain semibracket shapes.

**Definition 3.2.10: Trivial Semibracket** **(Fried, 2024)**

A semibracket is *trivial* if every team is declared co-champion without playing any games. Equivalently, a semibracket is trivial if its signature is of the form $[[\mathbf{m}]]_m$.

**Definition 3.2.11: Competitive Semibracket** **(Fried, 2024)**

A semibracket is *competitive* if no teams are declared co-champion without winning at least one game. Equivalently, a semibracket is competitive if its signature ends in a 0.

Clearly the two categories are mutually exclusive. Restricting briefly to the domain of traditional brackets, the two categories are also collectively exhaustive: there is no traditional bracket that is neither competitive nor trivial. (In fact, the only trivial traditional bracket is $[[\mathbf{1}]]$; every other traditional bracket is competitive.) However, this dichotomy does not apply to semibrackets: there are semibrackets that are neither trivial nor competitive. The simplest example is $[[\mathbf{2}; \mathbf{1}]]_2$, where the 1-seed is automatically one co-champion (so it is not competitive), but the 2- and 3-seeds play to be the other co-champion (so it is not trivial).

**Figure 3.2.12:** $[[2; 1]]_2$

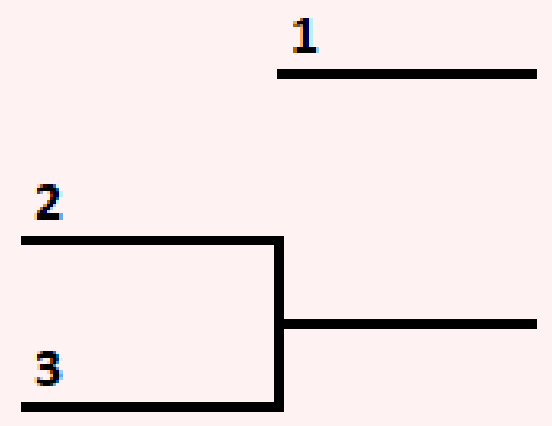

These two properties of semibrackets will sometimes be useful in defining and proving theorems about certain types of multibrackets down the line. In the next section, we will use semibrackets to construct a particularly nice kind of multibracket: the *linear multibracket.*

## 3.3 Linear Multibrackets

In the previous two sections, we have looked at semibrackets, as well as formats with a consolation bracket, as examples of multibrackets. Let's back up a bit from specific examples, however, and ask what information we can learn about arbitrary multibrackets. One potential question to ask is if the fundamental theorem of brackets, which held for traditional brackets and semibrackets, holds for multibrackets as well. But before we can do that, we need to define what a multibracket signature and proper multibracket seeding might look like.

This is trickier than it seems: for arbitrary multibrackets, there isn't a natural generalization of signatures and properness. (See Figure 3.7.1 on page 97 for an example of a multibracket that is difficult to assign a signature to). But there is a subset of multibrackets for which these notions generalize, allowing us to examine the fundamental theorem as it applies to this subset. These multibrackets are called *linear multibrackets*.

**Definition 3.3.1: Linear Multibracket** **(Fried, 2024)**

A *linear multibracket* is a multibracket that can be arranged into a sequence of semibrackets such that

(a) If a team loses in a given semibracket but is not eliminated, they are sent to a later semibracket, and

(b) Each team that wins the $i$th semibracket finishes in $m$th place, where $m$ is the sum of the ranks of the first $i$ semibrackets.

A linear multibracket can then be easily imbued with a signature derived from the signatures of the semibrackets in the sequence.

**Definition 3.3.2: Linear Multibracket Signatures** **(Fried, 2024)**

The *signature* of a linear multibracket that consists of semibrackets with signature $\mathcal{A}_1, ..., \mathcal{A}_k$ is $\mathcal{A}_1 \to ... \to \mathcal{A}_k$.

All four of the multibrackets discussed in the previous section are linear: let's see what their signatures are. First, the 2015 AFC Asian Cup [20].

**Figure 3.3.3: 2015 AFC Asian Cup**

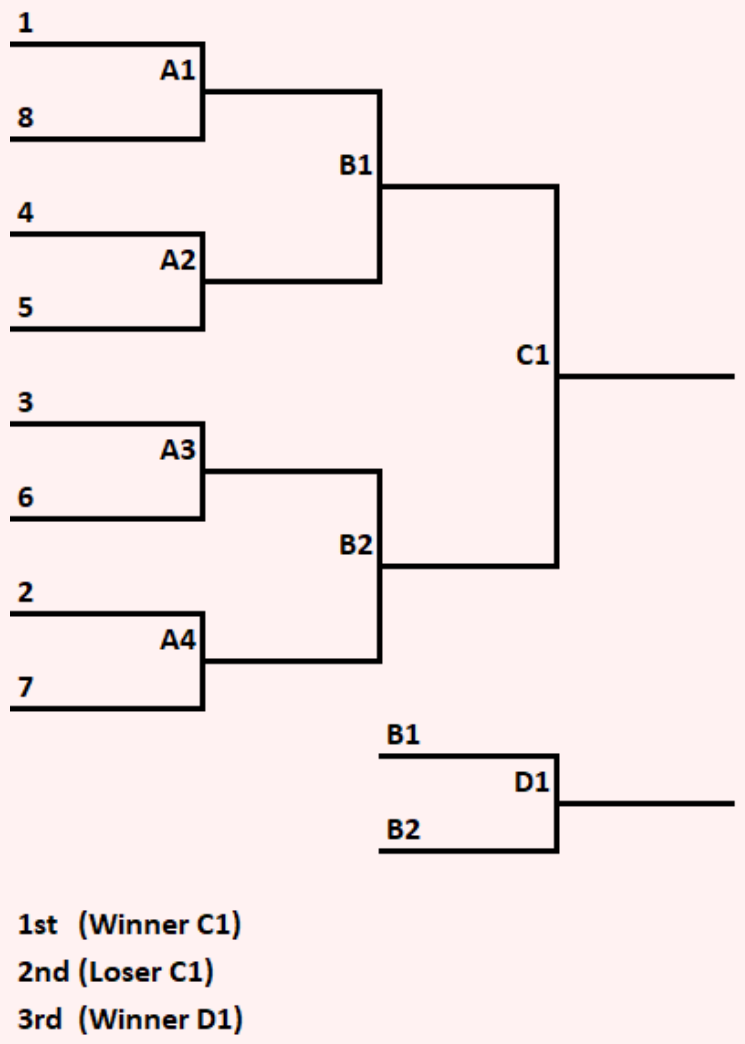

Looking at Figure 3.3.3, it can be tempting to say that the 2015 AFC Asian Cup is a linear multibracket of signature $[[\mathbf{8};\mathbf{0};\mathbf{0};\mathbf{0}]] \to [[\mathbf{2};\mathbf{0}]]$. But this is not quite right: The format with this signature would give second place to the winner of **D1** (as the winner of the second bracket), while outright eliminating the loser of **C1** (as a team that did not win any bracket). But in fact, we want to give second place to the loser of **C1**, and then third place to the winner of the consolation bracket with signature $[[\mathbf{2};\mathbf{0}]]$. We can do this by adding a second bracket with signature $[[\mathbf{1}]]$ while sliding the bracket with signature $[[\mathbf{2};\mathbf{0}]]$ to third.

Thus in total, the 2015 Asian Cup is a linear multibracket with signature $[[\mathbf{8};\mathbf{0};\mathbf{0};\mathbf{0}]] \to [[\mathbf{1}]] \to [[\mathbf{2};\mathbf{0}]]$. To make clear that the middle one-team bracket is included, we include it in the figure. This also allows us to drop the labeling of which teams finish in which place, as they are guaranteed by the linearity.

**Figure 3.3.4:** $[[8;0;0;0]] \to [[1]] \to [[2;0]]$

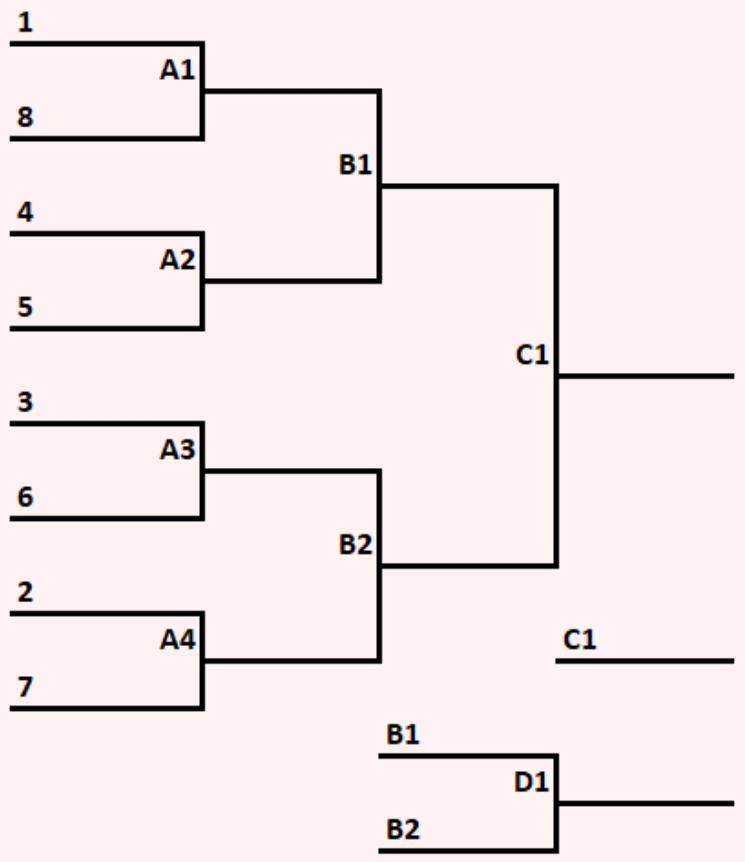

Next, let's examine our alternative to the 2015 AFC Asian Cup.

**Figure 3.3.5: 2015 AFC Asian Cup Alternative**

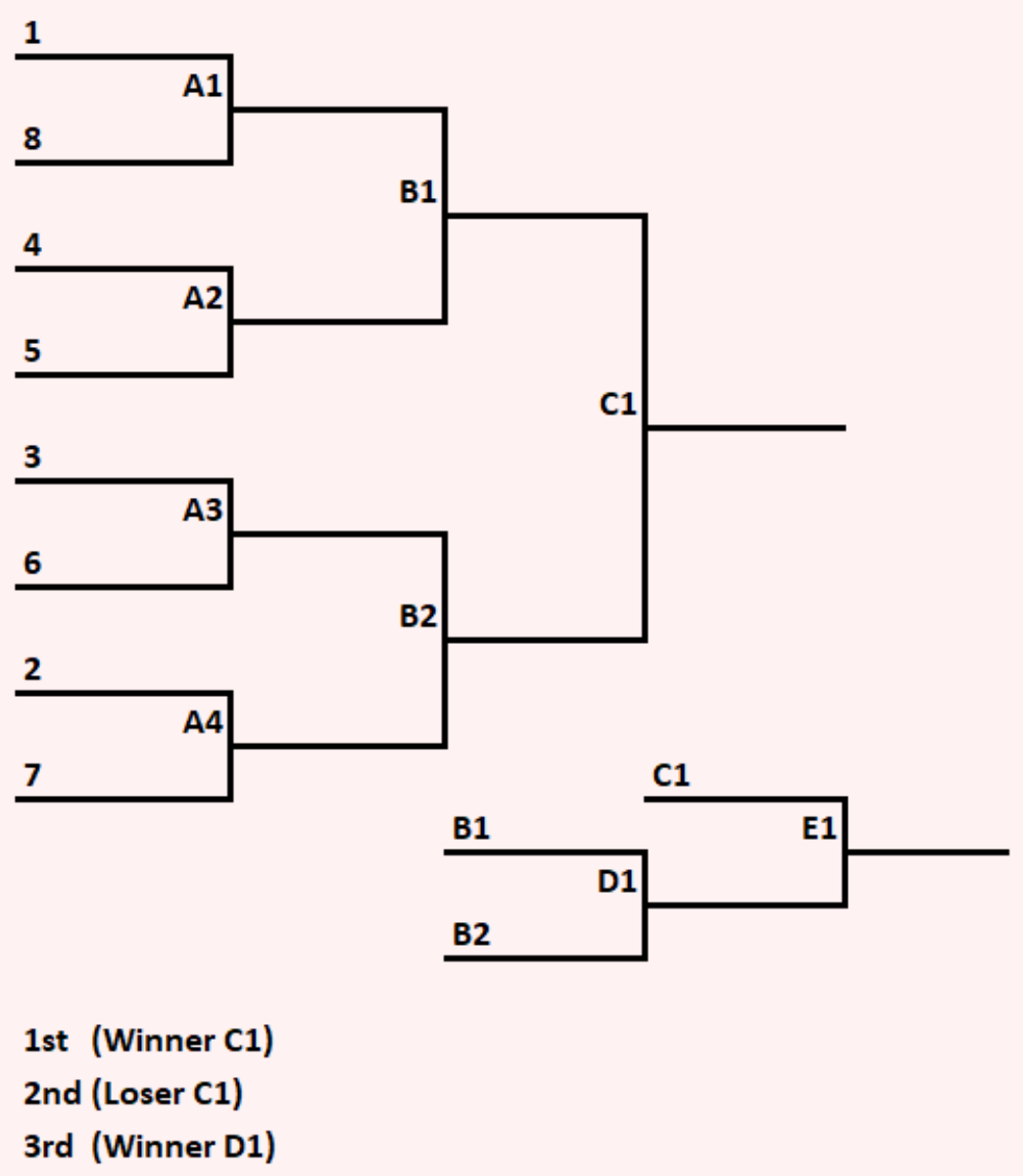

Again, a quick look indicates a signature of $[[\mathbf{8};\mathbf{0};\mathbf{0};\mathbf{0}]] \rightarrow [[\mathbf{2};\mathbf{1};\mathbf{0}]]$. And while this signature would correctly assign a first- and second-place, it doesn't assign third-place. Instead, we need a signature of $[[\mathbf{8};\mathbf{0};\mathbf{0};\mathbf{0}]] \rightarrow [[\mathbf{2};\mathbf{1};\mathbf{0}]] \rightarrow [[\mathbf{1}]]$.

**Figure 3.3.6: $[[8;0;0;0]] \rightarrow [[2;1;0]] \rightarrow [[1]]$**

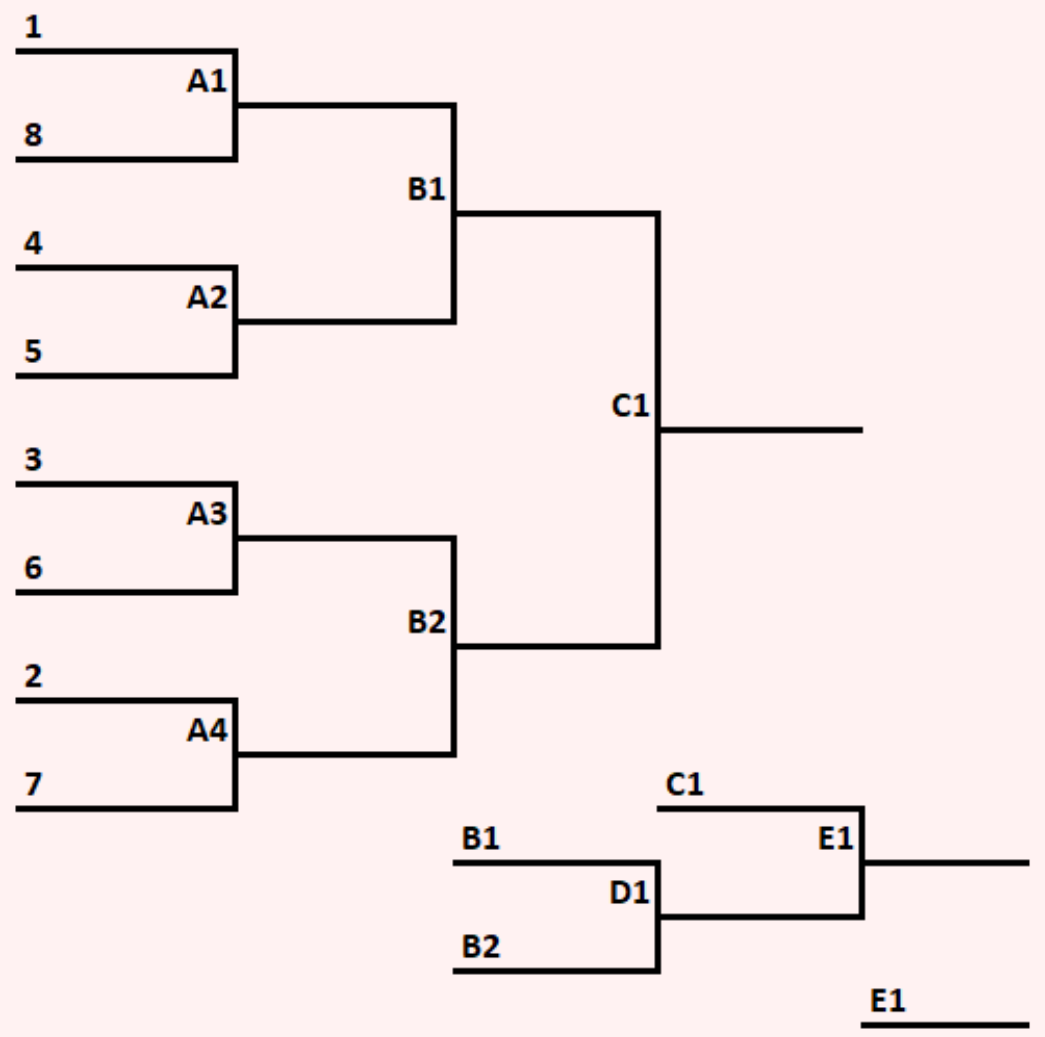

A similar analysis finds that the signature of the 2023 Southern Conference Wrestling Championships [24] is $[[\mathbf{8};\mathbf{0};\mathbf{0};\mathbf{0}]] \rightarrow [[\mathbf{1}]] \rightarrow [[\mathbf{4};\mathbf{2};\mathbf{0};\mathbf{0}]] \rightarrow [[\mathbf{1}]]$.

**Figure 3.3.7:** $[[8; 0; 0; 0]] \to [[1]] \to [[4; 2; 0; 0]] \to [[1]]$.

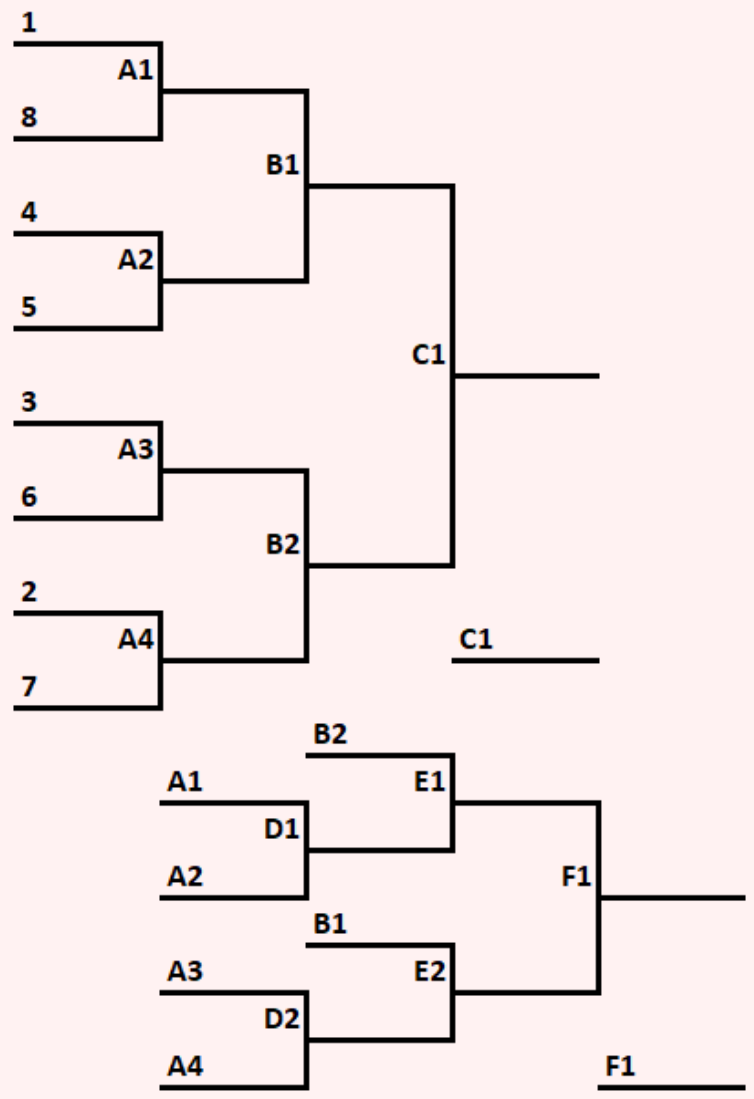

In all three examples so far, every semibracket has had rank one (that is, been a traditional bracket). However our final example, the 2023 Southern Conference Wrestling Championships Alternative, requires a semibracket of greater rank than one. (Recall the motivation for the alternative multibracket: we want to identify the top four teams while not eliminating any team from contention after just a single loss.)

**Figure 3.3.8: 2023 SoCon Wrestling Championships Alternative**

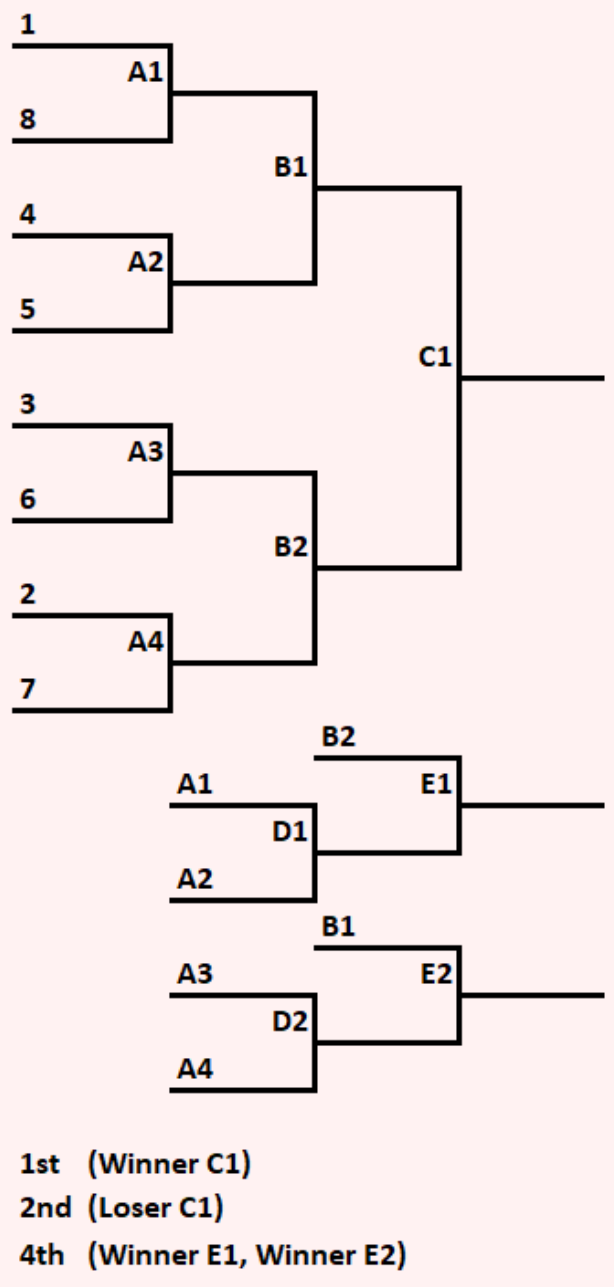

An attempt to give a signature to this format without the use of non-traditional semibrackets might be $[[\mathbf{8};\mathbf{0};\mathbf{0};\mathbf{0}]] \to [[\mathbf{1}]] \to [[\mathbf{2};\mathbf{1};\mathbf{0}]] \to [[\mathbf{2};\mathbf{1};\mathbf{0}]]$. Unfortunately, this isn't quite the same format: it assigns third place to the winner of **E1** and fourth to the loser of **E2**. We want to treat both winners identically: luckily, this is the exact problem that semibrackets were developed to solve. Using semibrackets, we can see that the signature should be $[[\mathbf{8};\mathbf{0};\mathbf{0};\mathbf{0}]] \to [[\mathbf{1}]] \to [[\mathbf{4};\mathbf{2};\mathbf{0}]]_2$.

**Figure 3.3.9:** $[[8;0;0;0]] \to [[1]] \to [[4;2;0]]_2$

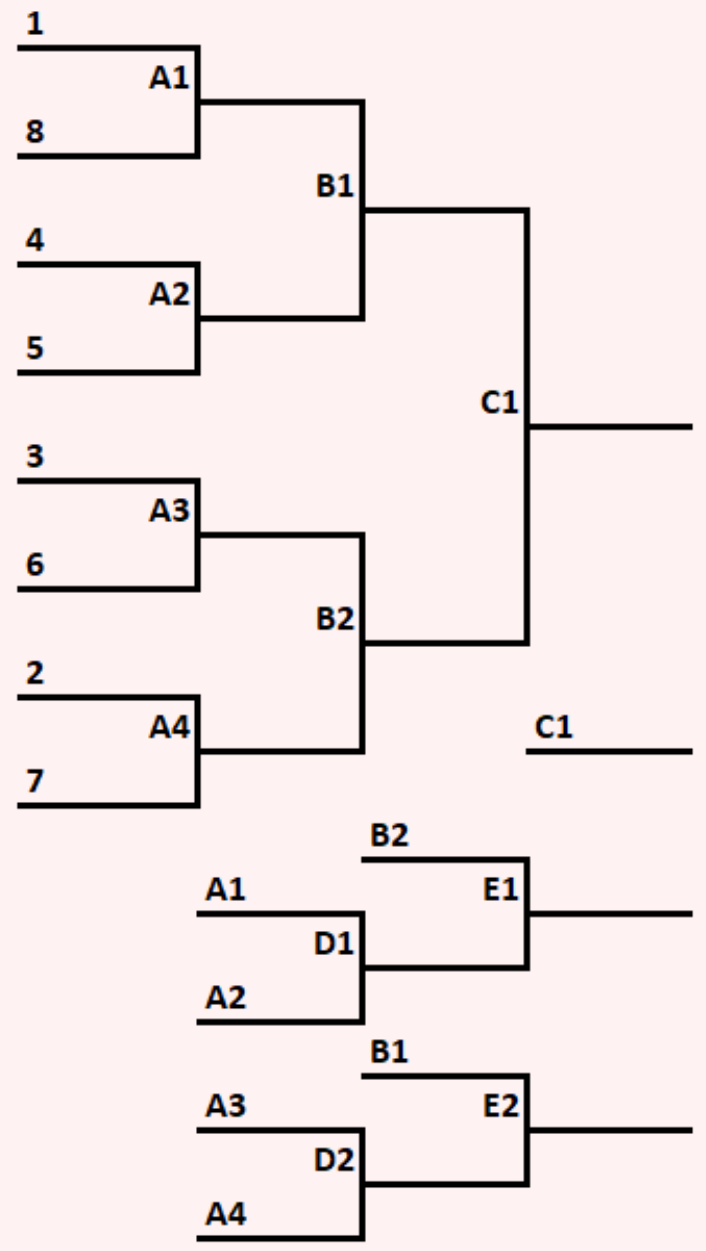

This format is differentiated from the (admittedly a bit strange) format in which the winner of game **E1** comes in third and the winner of game **E2** comes in fourth by the lettering of the games: the fact that games **E1** and **E2** are both **E**-round games means they must come from the same semibracket. If games **D2** and **E2** were instead **F1** and **G1**, respectively, then we would indeed have a linear multibracket of signature $[[\mathbf{8};\mathbf{0};\mathbf{0};\mathbf{0}]] \to [[\mathbf{1}]] \to [[\mathbf{2};\mathbf{1};\mathbf{0}]] \to [[\mathbf{2};\mathbf{1};\mathbf{0}]]$.

Now that we have defined linear multibrackets and developed a notion of signature, we can turn to the other half of the fundamental theorem: properness.

## 3.4 Flowcharts

Recall the definition of properness for a single bracket: a bracket is proper if, as long as the bracket goes chalk, in every round it is better to be a higher seed than a lower seed. We adapt this definition to linear multibrackets.

**Definition 3.4.1: Proper Linear Multibracket** **(Fried, 2024)**

A linear multibracket is *proper* if, as long as the bracket goes chalk, in every round of every semibracket it is better to be a higher-seeded team than a lower-seeded one, where:

(a) It is best to have already won an earlier semibracket.

(b) If you have not yet won an earlier semibracket, it is to better to be competing in the current semibracket than to not.

(c) If you are competing in a semibracket, it is better to have a bye in the current round than to not.

(d) If you are playing a game, it is better to play a lower seed than to play a higher seed.

We state without proof that the fundamental theorem applies to linear multibrackets as well as traditional brackets and semibrackets.

**Theorem 3.4.2** **(Fried, 2024)**

There is exactly one proper linear mulibracket with each linear mulibracket signature.

Consider the 2023 Major League Quadball Championship Play-In Tournament [33], which used a linear multibracket of signature $[[\mathbf{4};\mathbf{2};\mathbf{0};\mathbf{0}]] \to [[\mathbf{4};\mathbf{0};\mathbf{1};\mathbf{0}]]$. Is it proper?

**Figure 3.4.3: 2023 MLQ Championship Play-In Tournament**

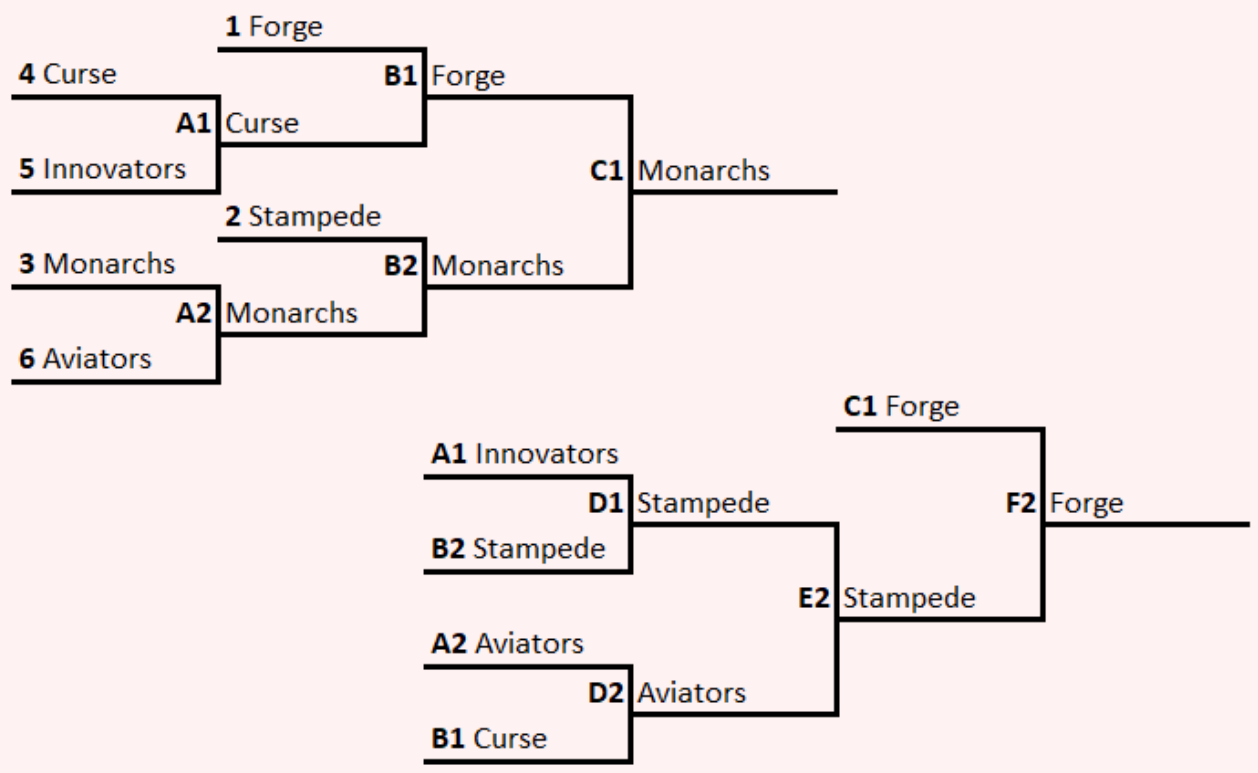

Well the primary semibracket certainly looks good: it is just the proper bracket of signature $[[\mathbf{4};\mathbf{2};\mathbf{0};\mathbf{0}]]$. To analyze the secondary bracket, we begin with Figure 3.4.4 which details which seeds would lose which primary bracket games if it went chalk.

**Figure 3.4.4: Which Seeds Would Lose Which Games if the 2023 MLQ Championship Play-In Tournament Went Chalk**

| Game | Seed |
|---|---|
| **A1** | 5 |
| **A2** | 6 |
| **B1** | 4 |
| **B2** | 3 |
| **C1** | 2 |

Thus for a linear multibracket of signature $[[\mathbf{4};\mathbf{2};\mathbf{0};\mathbf{0}]] \to [[\mathbf{4};\mathbf{0};\mathbf{1};\mathbf{0}]]$ to be proper, the first round of the secondary bracket would have to pair the loser of **A1** with the loser of **B1** and the loser of **A2** with the loser of **B2**. Instead, the loser of **A1** is paired with the loser of **B2** and the loser of **A2** with the loser of **B1**, meaning that the linear multibracket is not proper. Indeed, the 6-seed Aviators had an easier **D**-round matchup than the 5-seed Innovators. And it's not just the MLQ, many leagues that use linear multibrackets use ones that aren't proper. Why? Rematches.

In any tournament format, rematches are far from ideal. From an information theoretical perspective, a rematch is less informative than a new matchup: we already have some data on how those two team compare. From a competitive perspective, they are unsatisfying: without the ability to play a third "rubber" match, if each team wins one game, we are left in a disappointing state of uncertainty. These issues are only exacerbated in a linear multibracket. In the 2023 MLQ Championship Play-In Tournament teams are eliminated after their second loss: it would feel awful for those two losses to come at the hands of the same team.

If the proper linear multibracket of signature $[[\mathbf{4};\mathbf{2};\mathbf{0};\mathbf{0}]] \to [[\mathbf{4};\mathbf{0};\mathbf{1};\mathbf{0}]]$ went chalk, both **D**-round games would be rematches of the **A**-round games, which, as discussed, is disappointing for both the competing teams as well as the audience. In the MLQ's format, by contrast, the **D** rounds games are guaranteed to be new matchups.

The balance of how much to prioritize properness versus dodging rematches is one that every league must answer for itself: while the MLQ's format (if it went chalk) would have only a single rematch, there are formats with non-proper primary brackets that would have none. In any case, unlike for traditional brackets, where non-proper brackets are hardly ever used, in the world of linear multibrackets they are actually quite common.

Is there any notion of properness that we can hang on to as a property that most linear multibrackets ought to have? Indeed there is, but to understand it we first must introduce a new way of looking at linear multibrackets: *flowcharts*.

**Definition 3.4.5: Flowchart** **(Fried, 2024)**

The *flowchart* of a linear multibracket that consists of $k$ semibrackets is a directed graph in which the nodes are arranged into rows, where

(a) There is a node for each team, each round of each semibracket, and each place a team could finish in, plus one additional node representing elimination.

(b) The zeroth row has the nodes representing each team, arranged from lowest seed to highest seed.

(c) The $i$th row for $1 \leq i \leq k$ has the nodes representing the rounds of the $i$th semibracket, arranged in order, plus the node representing the place a team gets for winning the $i$th semibracket.

(d) The final row row has only the node representing elimination.

(e) There is an arrow from each team to the round where that teams plays its first game.

(f) For each round **R**, there is an arrow (or arrows) from **R** to the round(s) where **R**-round losers go.

The flowchart for the 2023 MLQ Championship Play-In Tournament is displayed below, with the team nodes labeled $t_1$ through $t_6$, the round nodes labeled with the letter of the round, the nodes representing finishing in $i$th place labeled as such, and the node representing elimination filled in.

**Figure 3.4.6: 2023 MLQ Championship Play-In Tournament Flowchart**

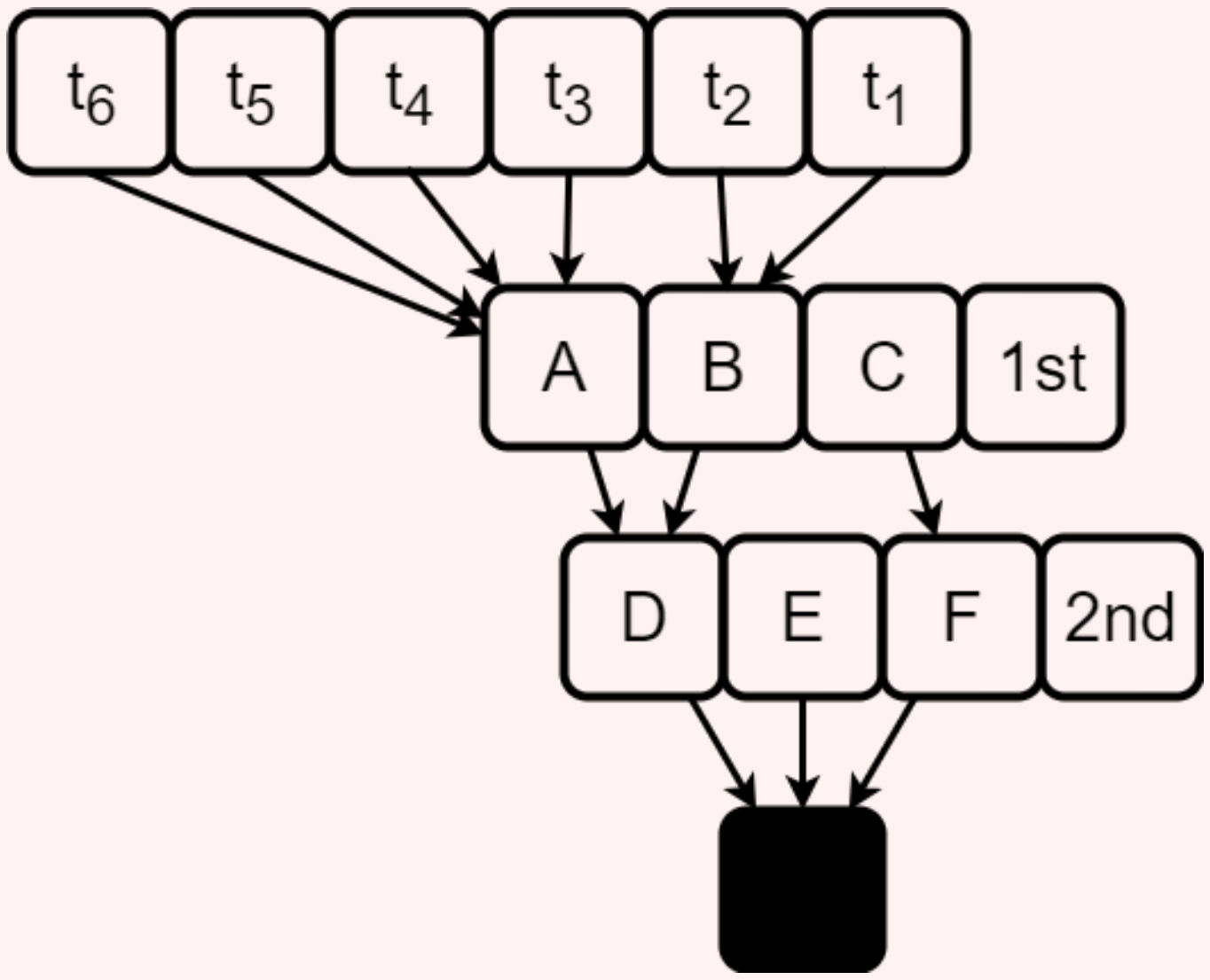

Flowcharts exist in an in-between space between linear multibrackets and linear multibracket signatures: multiple formats can have the same flowchart, and multiple flowcharts can have the same signature.

Imagine, for example, that Major League Quadball was interested in selecting a third and fourth place team using their play-in tournament as well. Figure 3.4.7 shows three formats they could use: the left two have the same flowchart, while the rightmost format has a different flowchart. Both flowcharts are displayed in Figure 3.4.8, and they have the same signature: $[[\mathbf{4}; \mathbf{2}; \mathbf{0}; \mathbf{0}]] \to [[\mathbf{4}; \mathbf{0}; \mathbf{1}; \mathbf{0}]] \to [[\mathbf{4}; \mathbf{0}; \mathbf{0}]] \to [[\mathbf{2}; \mathbf{1}; \mathbf{0}]]$.

**Figure 3.4.7: Three Linear Multibrackets with the Same Signature**

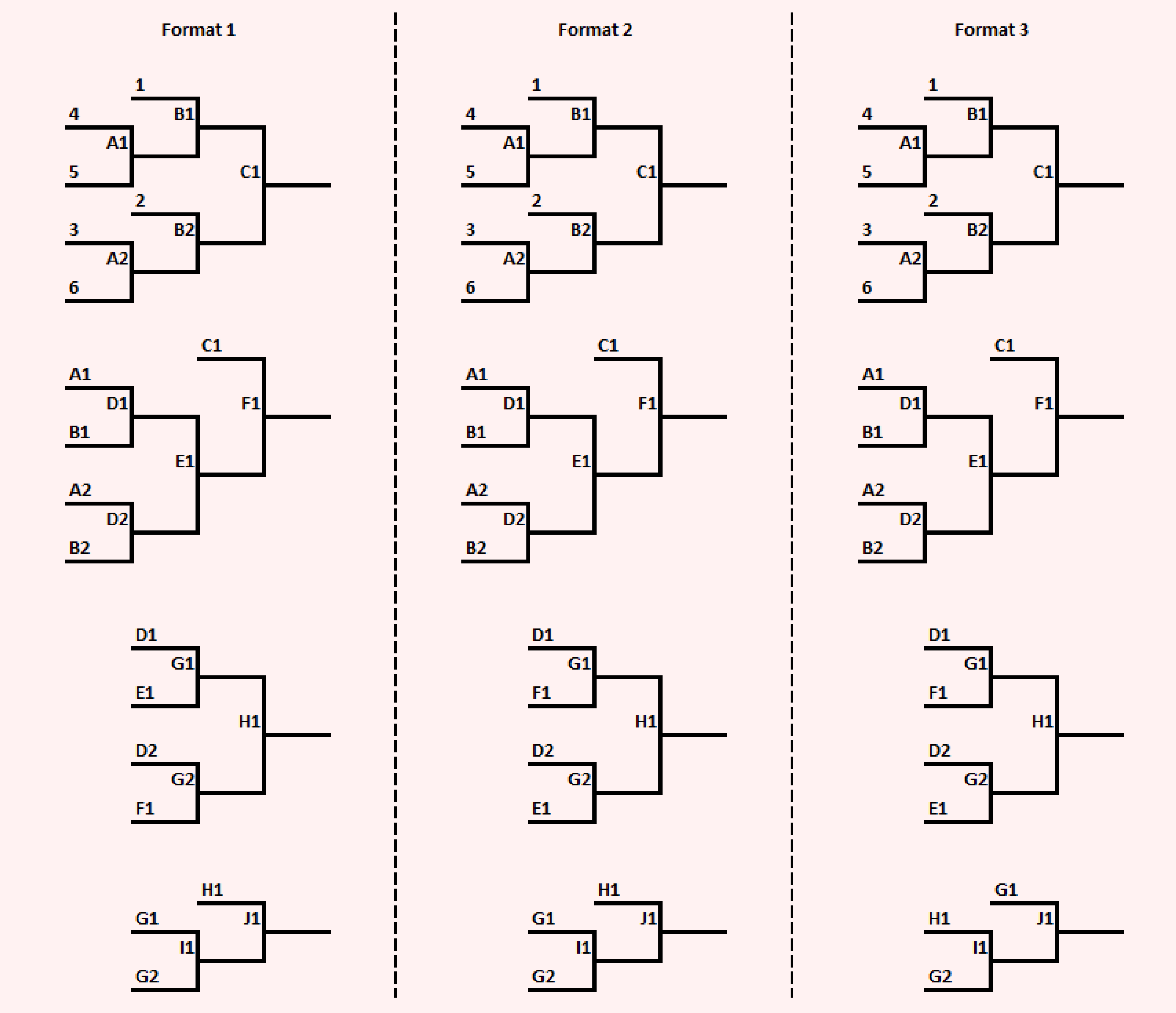

**Figure 3.4.8: Two Flowcharts with the Same Signature**

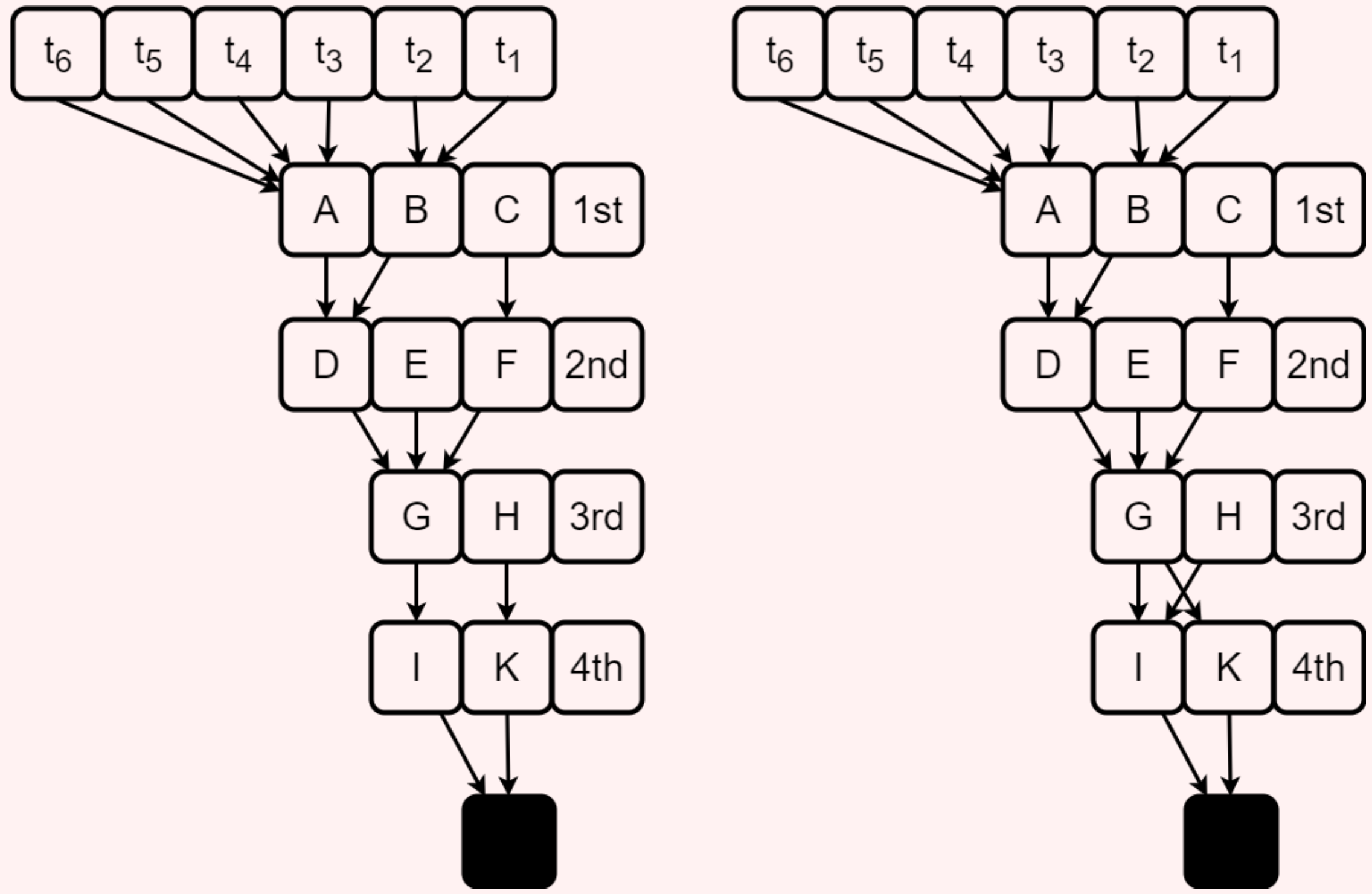

However, these three formats and two flowcharts are not created equal. While the two leftmost formats – and thus the leftmost flowchart – all seem reasonable, the rightmost format and flowchart have some issues. In particular, the fourth semibracket is poorly seeded: the loser of game **H1** did better than both **G**-round losers in the previous semibracket, so they ought to be the one to get the bye. Additionally, the two **G**-round losers ought to be treated the same, instead of one of them getting a bye and one not. Further, none of these issues are resolved by appealing to a decrease in rematches: teams are being sent to wholly the wrong *round*.

These two problems are reflected in the flowchart as well: the fact that the **G**-round losers are sent to different rounds of the 4th-place bracket means that the **G**-node has two arrows coming out of it, and the fact that a less deserving team was treated better than a more deserving team was means that there were two arrows in the flowchart that crossed over. (Note that technically we could wrap the arrow coming out of **H** all the way around the flowchart so it points to **I** from the left side, removing this crossing. However, we can fix hacks like this by imagining arrows coming out of the teams, places, and the elimination node extending infinitely up, to the right, and down respectively.)

Both of these issues fly in the face of our intuitive notion of rewarding better teams that motivated us to define properness in the first place, so we use them to define the notion of *respectfulness*.

**Definition 3.4.9: Respectful Linear Multibracket** **(Fried, 2024)**

A linear multibracket is *respectful* if its flowchart has no arrow crossings and every node in its flowchart has at most one arrow coming out of it.

Note that while respectfulness is like properness in that it requires linear multibrackets to treat more-deserving teams better than less-deserving teams (for some definition of deserving and better), it neither implies nor is implied by properness. The signature $[[\mathbf{2};\mathbf{3};\mathbf{0};\mathbf{0}]] \to [[\mathbf{1}]] \to [[\mathbf{2};\mathbf{1};\mathbf{0}]]$ has a proper instantiation but not a respectful one, while the signature $[[\mathbf{8};\mathbf{0};\mathbf{0};\mathbf{0}]] \to [[\mathbf{3}]] \to [[\mathbf{4};\mathbf{0};\mathbf{0}]]$ has only one proper instantiation but many respectful ones. In any case, due to the issues with properness outlined at the beginning of the section, we will primarily deal with respectful linear multibrackets for the rest of the chapter.

We conclude the section by proving a couple of nice lemmas about respectful linear multibrackets.

**Lemma 3.4.10** **(Fried, 2024)**

If $\mathcal{A}$ is a respectful linear multibracket, and two teams lose in the same round of the same semibracket of $\mathcal{A}$, then they will either both play their next game in the same round of the same semibracket, or both be elimineted.

*Proof.* This holds because every node in the flowchart of $\mathcal{A}$ has at most one arrow coming out of it. □

**Lemma 3.4.11** **(Fried, 2024)**

Let $\mathcal{A} = \mathcal{A}_1 \to \ldots \to \mathcal{A}_k$ be a linear multibracket, and $i \in \mathbb{N}$ be such that $\mathcal{A}_{i+1}$ is noncompetitive. If $\mathcal{A}$ is respectful, then any team that loses in the final round of $\mathcal{A}_i$ will win $\mathcal{A}_{i+1}$ without playing any more games.

*Proof.* We show the contrapositive. Let $\mathbf{R}$ be the final round of $\mathcal{A}_i$, whose losers do not win $\mathcal{A}_{i+1}$ automatically. $\mathcal{A}_{i+1}$ is noncompetitive though, so at least one team wins $\mathcal{A}_{i+1}$ without playing a game in it: let $\mathbf{S}$ be the round that the auto-winner fell from. Then $\mathbf{R}$ is completely boxed in by the arrow coming out of $\mathbf{S}$ on the left and the arrows extending infinitely from the final ranking nodes on the right: the arrow coming out of $\mathbf{R}$ must cross one of those arrows, so $\mathcal{A}$ is not respectful. □

## 3.5 Swiss Signatures

Consider the 1988 Men's College Basketball Maui Invitational [25], which used a respectful linear multibracket of signature $[[\mathbf{8};\mathbf{0};\mathbf{0};\mathbf{0}]] \to [[\mathbf{1}]] \to [[\mathbf{2};\mathbf{0}]] \to [[\mathbf{1}]] \to [[\mathbf{4};\mathbf{0};\mathbf{0}]] \to [[\mathbf{1}]] \to [[\mathbf{2};\mathbf{0}]] \to [[\mathbf{1}]]$.

**Figure 3.5.1: 1988 Men's College Basketball Maui Invitational**

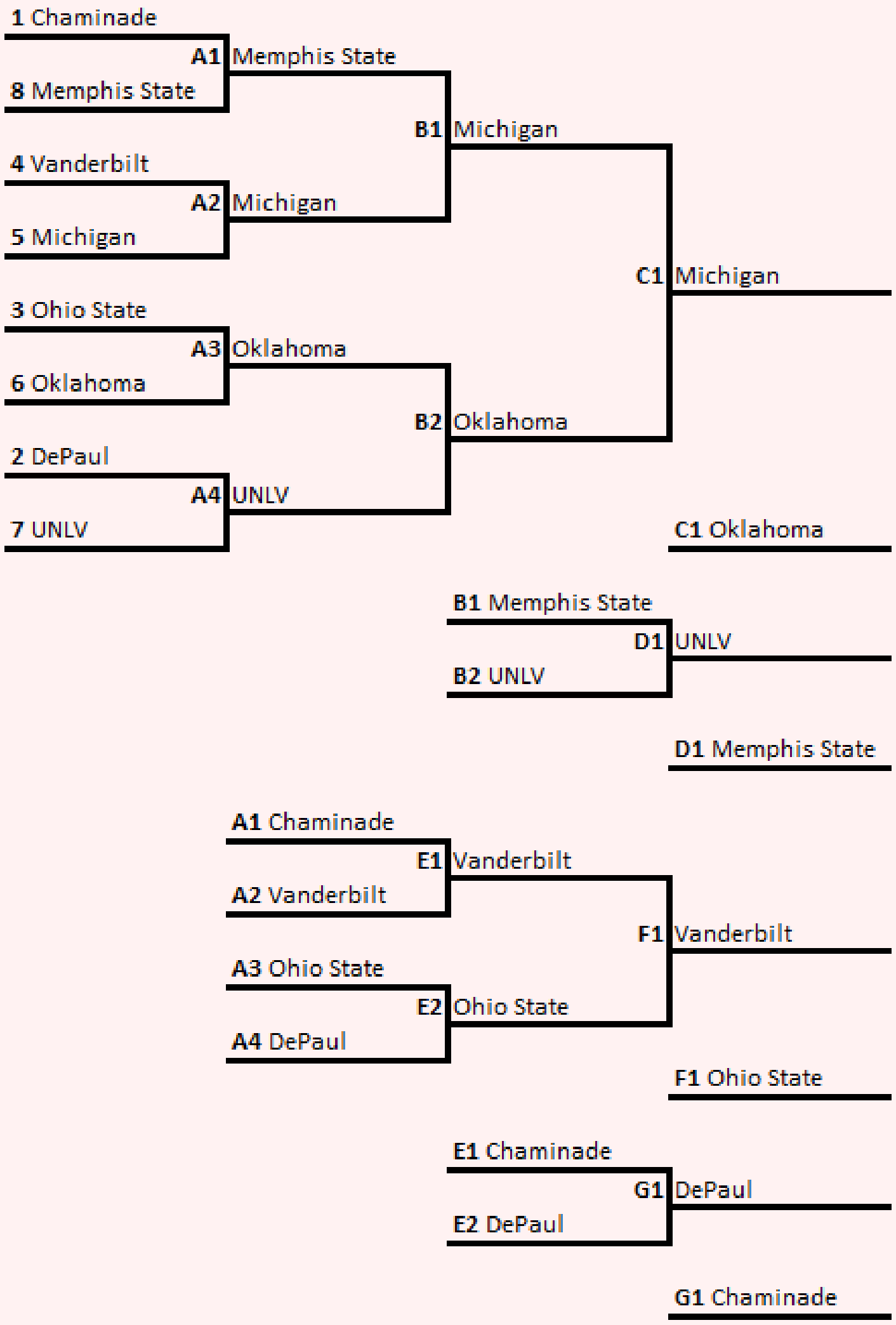

**Figure 3.5.2: 1988 Men's College Basketball Maui Invitational Flowchart**

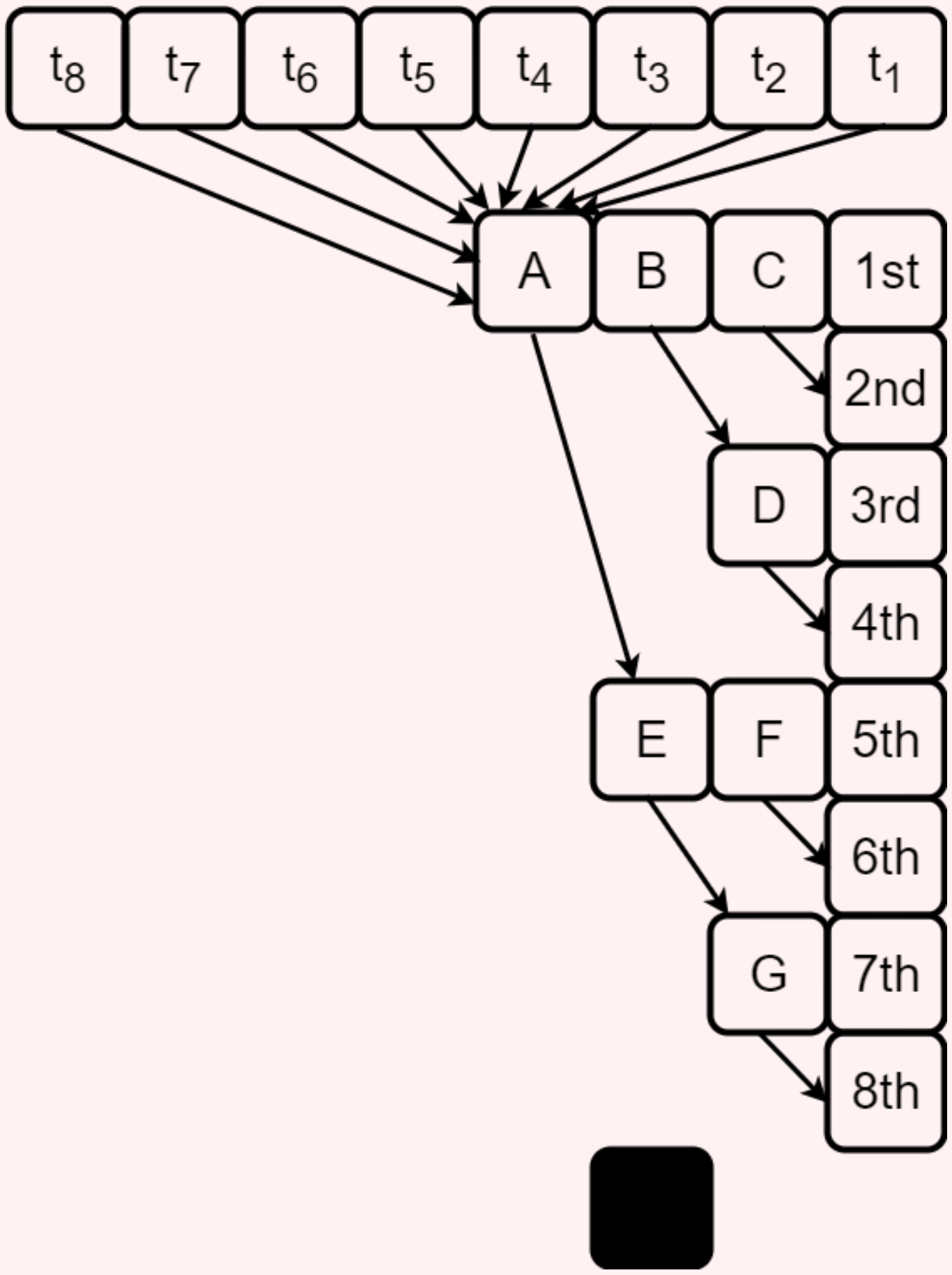

The format used in the Maui Invitational is of course respectful, but its flowchart has several additional properties that even the respectful flowchart in Figure 3.4.6 from last section does not have.

- Every team starts in the same cell, so we can unambiguously drop the top row of nodes.
- Games are always between teams of the same record, so we can unambiguously label each node with a record instead of a letter.
- Every team plays the same number of games, so our flowchart is nicely divided into columns, with each team playing one game in each column. Further, this allows us to unambiguously drop the arrows, as losers always play their next game in the round below the round directly to the right of the round they lost in.
- Every team wins a semibracket and so every team ends in a cell, allowing us to drop the elimination node.

Our newly stylized flowchart, which we will refer to as a *swisschart* for the 1988 Men's College Basketball Maui Invitational is displayed below.

**Definition 3.5.3: Swisschart** **(Fried, 2024)**

A *Swisschart* is the flowchart of a Swiss format except with the arrows, nodes representing the teams, and the node representing elimination removed, and with the labels on the remaining nodes replaced with the record of the teams in that node.

**Figure 3.5.4: 1988 Men's College Basketball Maui Invitational Swisschart**

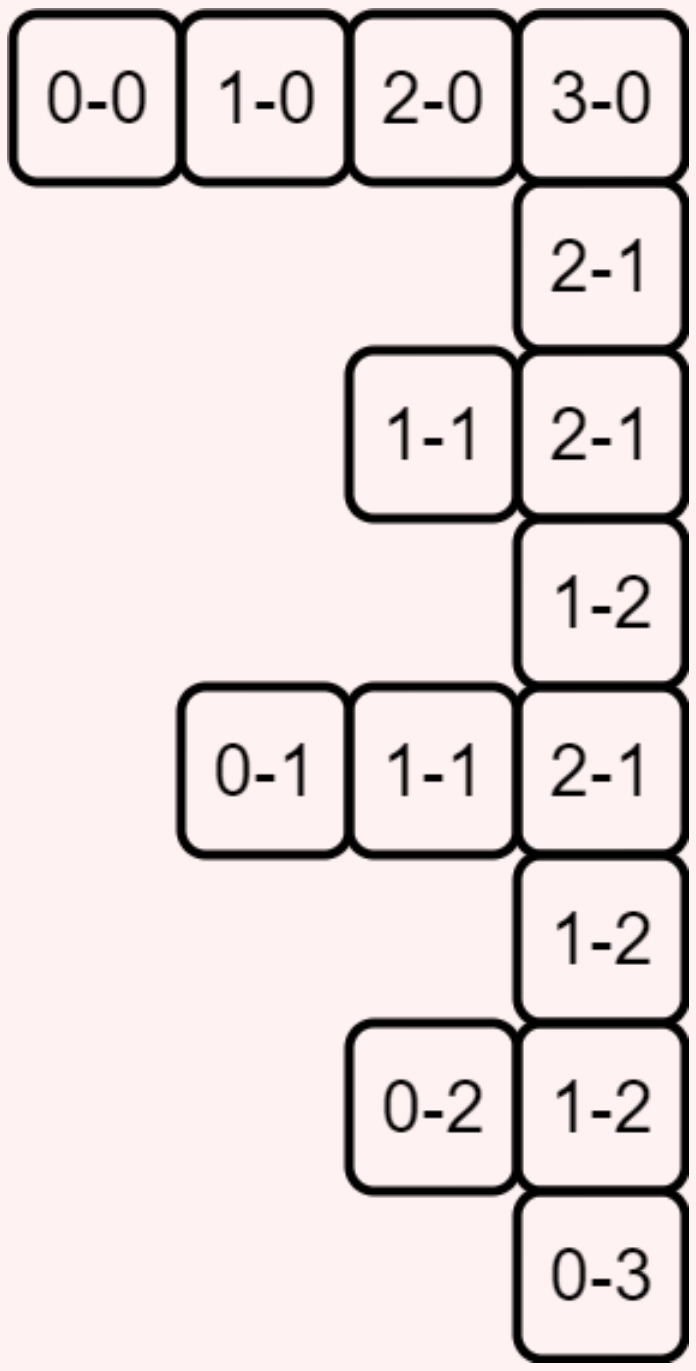

There is one other nice property that the 1988 Men's College Basketball Maui Invitational Flowchart has: all of its semibrackets are either trivial or competitive. A nontrivial noncompetitive semibracket of signature $[[\mathbf{a_1}; ...; \mathbf{a_r}]]_m$ can always be split into the the pair of semibrackets $[[\mathbf{a_1}; ...; \mathbf{a_{r-1}}; \mathbf{0}]]_{m-a_r} \to [[\mathbf{a_r}]]_{a_r}$ without affecting the games played in the tournament or who ends up ranked, and so to avoid double counting formats we require this property as well.

Formats with all of these properties are called *Swiss formats*, and their signatures are called *Swiss signatures*, named because of their first recorded use at a chess tournament in Zürich, Switzerland in 1895 [35].

**Definition 3.5.5: Swiss Format** **(Unattributed)**

A *Swiss format* is a respectful linear multibracket with the following additional properties:

(a) Every team starts in the primary semibracket.

(b) Every game is between two teams with the same record.

(c) Every team plays the same number of games.

(d) Every team wins a semibracket.

(e) Every semibracket is either trivial or competitive.

**Definition 3.5.6: Swiss Signature** **(Fried, 2024)**

A *Swiss signature* is a linear multibracket signature that admits a Swiss format.

**Definition 3.5.7: $r$-Round Swiss Signature** **(Fried, 2024)**

An *$r$-round Swiss signature* is a Swiss signature in which each team plays $r$ games.

Thus $[[\mathbf{8};\mathbf{0};\mathbf{0};\mathbf{0}]] \to [[\mathbf{1}]] \to [[\mathbf{2};\mathbf{0}]] \to [[\mathbf{1}]] \to [[\mathbf{4};\mathbf{0};\mathbf{0}]] \to [[\mathbf{1}]] \to [[\mathbf{2};\mathbf{0}]] \to [[\mathbf{1}]]$ is a 3-round Swiss signature. In fact, it is an example of a particular family of Swiss signatures known as the the *standard Swiss signatures*, which we abbreviate by $\mathcal{S}_r$ for some $r$.

**Definition 3.5.8: Standard Swiss Signature ($\mathcal{S}_r$)** **(Fried, 2024)**

$\mathcal{S}_r$, or the *standard $r$-round Swiss signature*, is the multibracket signature defined recursively by

$$\mathcal{S}_0 = [[\mathbf{1}]],$$

and

$$\mathcal{S}_r = [[\mathbf{2^r};...;\mathbf{0}]] \to \mathcal{S}_0 \to \mathcal{S}_1 \to ... \to \mathcal{S}_i \to ... \to \mathcal{S}_{r-1}.$$

Thus we have

$$\mathcal{S}_0 = [[\mathbf{1}]]$$
$$\mathcal{S}_1 = [[\mathbf{2};\mathbf{0}]] \to [[\mathbf{1}]]$$
$$\mathcal{S}_2 = [[\mathbf{4};\mathbf{0};\mathbf{0}]] \to [[\mathbf{1}]] \to [[\mathbf{2};\mathbf{0}]] \to [[\mathbf{1}]]$$
$$\mathcal{S}_3 = [[\mathbf{8};\mathbf{0};\mathbf{0};\mathbf{0}]] \to [[\mathbf{1}]] \to [[\mathbf{2};\mathbf{0}]] \to [[\mathbf{1}]] \to [[\mathbf{4};\mathbf{0};\mathbf{0}]] \to [[\mathbf{1}]] \to [[\mathbf{2};\mathbf{0}]] \to [[\mathbf{1}]]$$

Figures 3.5.9 and 3.5.10 display the brackets and swisscharts for $\mathcal{S}_0, \mathcal{S}_1$, and $\mathcal{S}_2$, while the Maui Invitational used the standard Swiss signature $\mathcal{S}_3$.

**Figure 3.5.9:** $\mathcal{S}_0, \mathcal{S}_1$, and $\mathcal{S}_2$

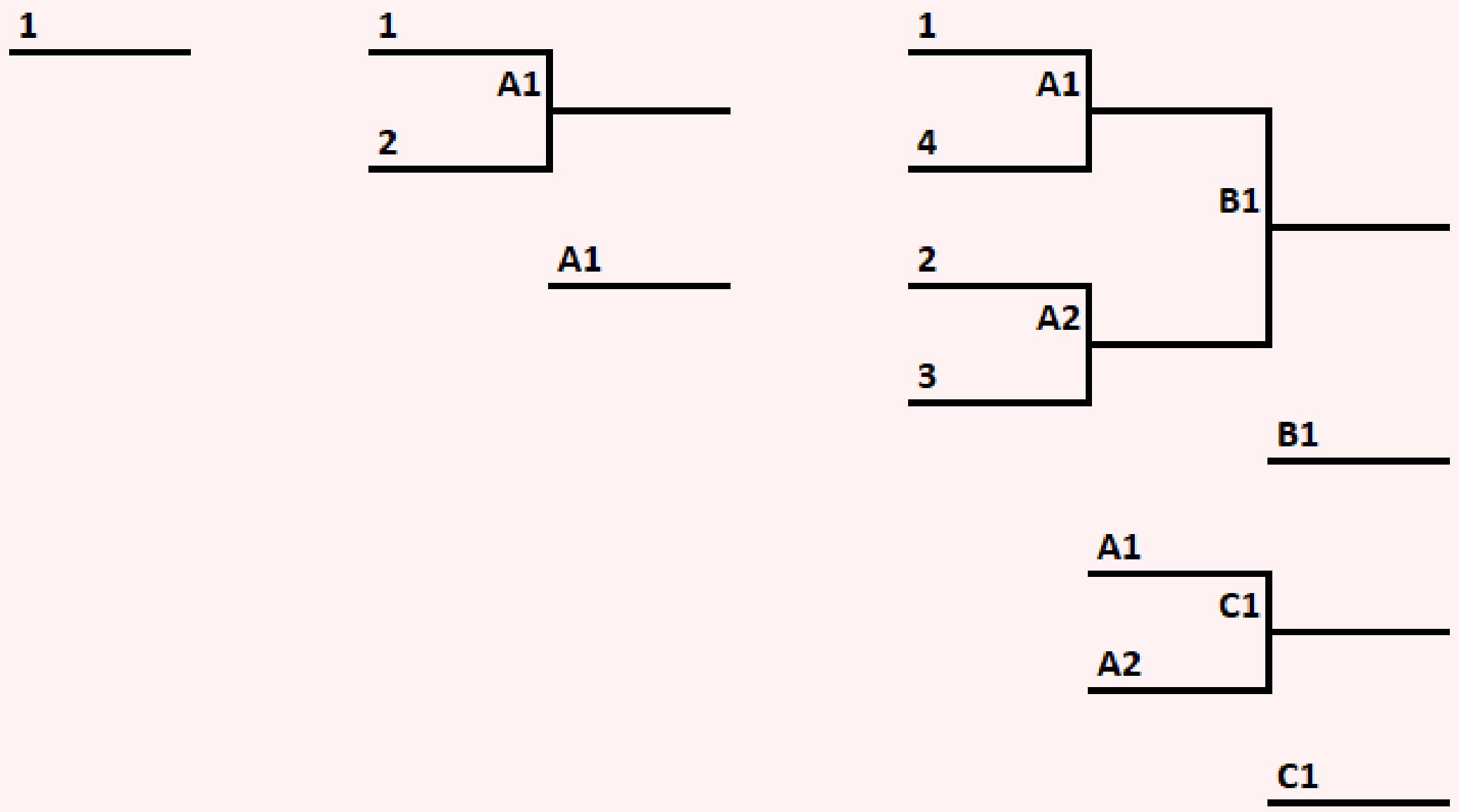

**Figure 3.5.10:** $\mathcal{S}_0, \mathcal{S}_1$, and $\mathcal{S}_2$

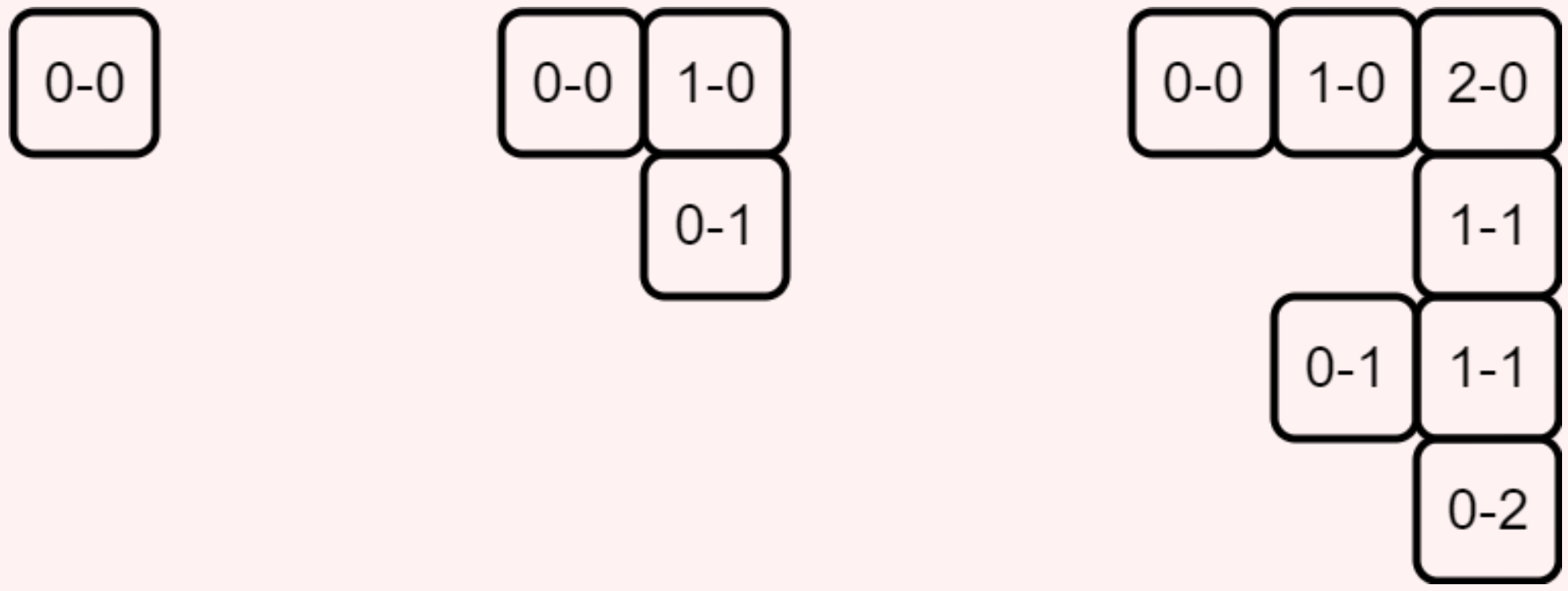

Consider now the following non-standard Swiss format.

**Figure 3.5.11:** $[[8;0;0]]_2 \to [[2]]_2 \to [[4;0]]_2 \to [[2]]_2$

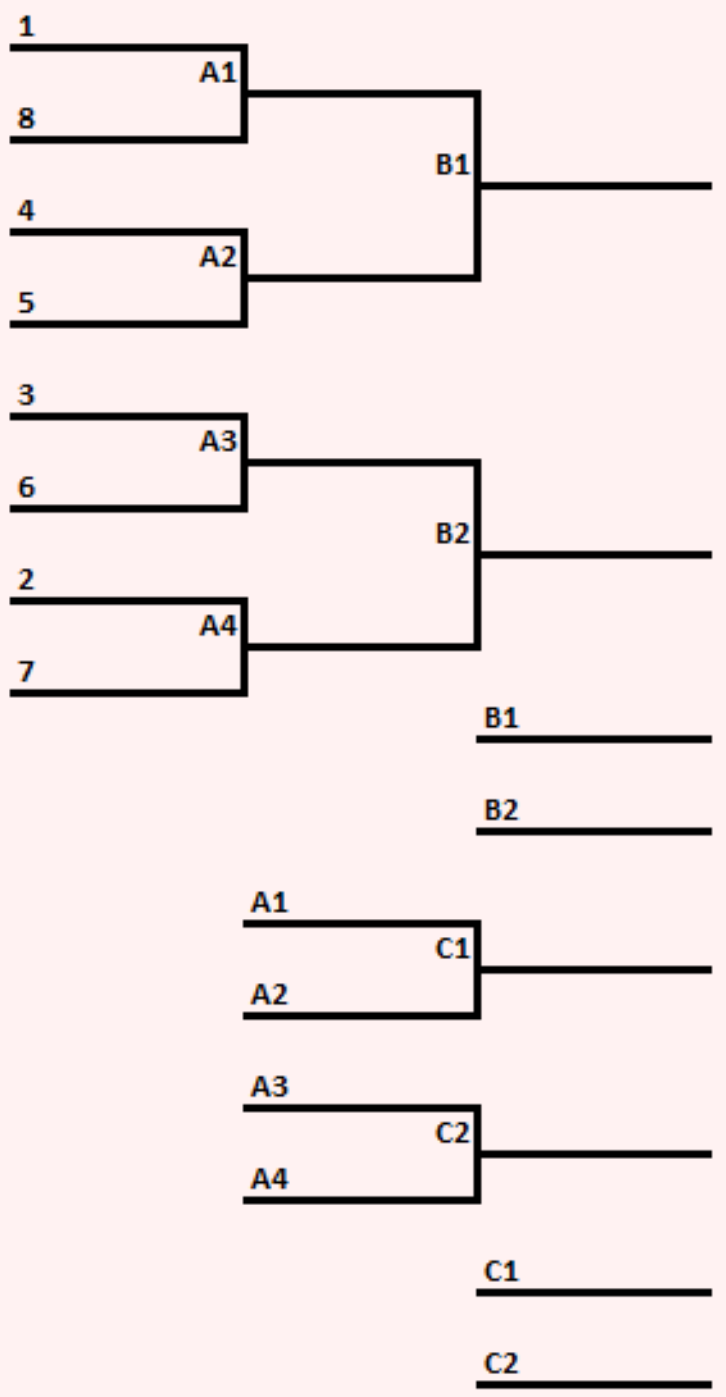

Unlike the standard Swiss signatures, this signature does not crown a single champion: it is not *compact*.

**Definition 3.5.12: Compact Swiss Signature** **(Fried, 2024)**

A Swiss signature is *compact* if only one teams finishes undefeated.

If we attempt to draw the swisschart for the format, we notice something a little strange.

**Figure 3.5.13:** $[[8;0;0]]_2 \to [[2]]_2 \to [[4;0]]_2 \to [[2]]_2$ **Swisschart**

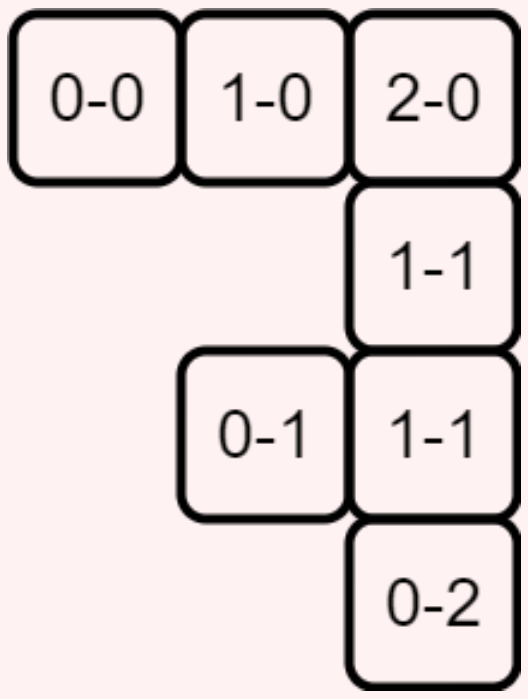

It is identical to $\mathcal{S}_2$! This is not a coincidence: the format is actually just two copies of $\mathcal{S}_2$ being run simultaneously. Taking another look at the signature, it is even the same signature as $\mathcal{S}_2$, just with every number multiplied by 2. We use this to introduce $m$x notation.

**Definition 3.5.14: $m\text{x}\mathcal{A}$** **(Fried, 2024)**

If $m \in \mathbb{N}$ and $\mathcal{A}$ is a multibracket signature, then $m\text{x}\mathcal{A}$ is the multibracket signature formed by multiplying every number in every signature in $\mathcal{A}$ by $m$.

So the format in Figures 3.5.11 and 3.5.13 is $2\text{x}\mathcal{S}_2$. In fact, every noncompact Swiss signature can be represented using $m\text{x}\mathcal{A}$ notation.

**Theorem 3.5.15** **(Fried, 2024)**

Let $\mathcal{A}$ be a noncompact Swiss signature where $m > 1$ teams end undefeated. Then $m$ will divide every number in the signature of $\mathcal{A}$.

*Proof.* We first prove by reverse induction that the number of teams participating in the $i$th-to-last round of every semibracket in $\mathcal{A}$ is divisible by $m \cdot 2^i$.

We begin with the base case of $i = r$. Only the primary semibracket of a Swiss format has an $r$th to last round, and because it is a balanced semibracket of rank $m$, it has $m \cdot 2^r$ teams. For any other $i$, by induction, $m \cdot 2^{i+1}$ divides the number of teams that competed in the $(i+1)$th-to-last round of the semibracket. Half of them won and so are still competing in this semibracket. Meanwhile, if the $i$th round of this semibracket is to take the losers of another round of another semibracket, it must also be the $(i+1)$th-to-last round, so by induction $m \cdot 2^i$ divides the number of teams that will fall into this semibracket. Thus the number of teams playing in the $i$th-to-last round of every semibracket will be divisible by $m \cdot 2^i$.

Thus the $i$th-to-last number of every semibracket signature, which is a collection of $i$th-to-last round losers, is divisible by $m \cdot 2^{i-1}$, proving the theorem. ☐

With the standard Swiss signatures and $m$x notation defined, we are ready for Figure 3.5.16, which details the various Swiss signatures for 1-, 2-, 4-, and 8-teams.

**Figure 3.5.16: The 1-, 2-, 4-, and 8-team Swiss Signatures**

| | 1 Team | 2 Teams | 4 Teams | 8 Teams |
|---|---|---|---|---|
| 0 Rounds | $\mathcal{S}_0$ | $2\text{x}\mathcal{S}_0$ | $4\text{x}\mathcal{S}_0$ | $8\text{x}\mathcal{S}_0$ |
| 1 Round | | $\mathcal{S}_1$ | $2\text{x}\mathcal{S}_1$ | $4\text{x}\mathcal{S}_1$ |
| 2 Rounds | | | $\mathcal{S}_2$ | $2\text{x}\mathcal{S}_2$ |
| 3 Rounds | | | | $\mathcal{S}_3, \mathcal{T}_3$ |

How do we know that this diagram is complete? Well the cells below the diagonal must be empty: only one out of $2^r$ teams can remain undefeated after $r$ rounds and so would have no opponent in a hypothetical $(r+1)$th round. The cells above the diagonal are noncompact, and thus complete by Theorem 3.5.15. The most interesting cells are the compact ones on the diagonal.

**Theorem 3.5.17** **(Fried, 2024)**

The only compact 0-, 1-, and 2-round Swiss signatures are the standard ones.

*Proof.* $[[\mathbf{1}]]$ is the only 1-team linear multibracket, so it is clearly the only compact 0-round Swiss signature.

$[[\mathbf{2};\mathbf{0}]] \to [[\mathbf{1}]]$ is the only 2-team Swiss signature in which a game is played, so it is the only compact 1-round Swiss signature.

Finally, a compact 2-round Swiss signature must start with $[[\mathbf{4};\mathbf{0};\mathbf{0}]]$ to be compact, its secondary bracket must be $[[\mathbf{1}]]$ to ensure every bracket is either trivial or competitive, its tertiary bracket must be $[[\mathbf{2};\mathbf{0}]]$ otherwise the two losers would not have a game to play, and it must end in $[[\mathbf{1}]]$ to ensure every team is ranked. Thus $[[\mathbf{4};\mathbf{0};\mathbf{0}]] \to [[\mathbf{1}]] \to [[\mathbf{2};\mathbf{0}]] \to [[\mathbf{1}]]$ is the only compact 2-round Swiss signature. □

Figure 3.5.16 tells us that there are two compact 3-round Swiss signatures: the standard $\mathcal{S}_3$, and the yet to be defined $\mathcal{T}_3$. It's worth attempting to construct $\mathcal{T}_3$ before reading on.

The key insight is to realize that teams with the same record in vertically adjacent nodes of the swisschart can actually play against each other without violating any of the Swiss format requirements, merging the nodes. Thus the flow chart for $\mathcal{T}_3$ looks like so. (Note that the 1-1 node contains four teams, and the bottommost 2-1 node as well as the topmost 1-2 node each contain two teams.)

**Figure 3.5.18:** $\mathcal{T}_3$

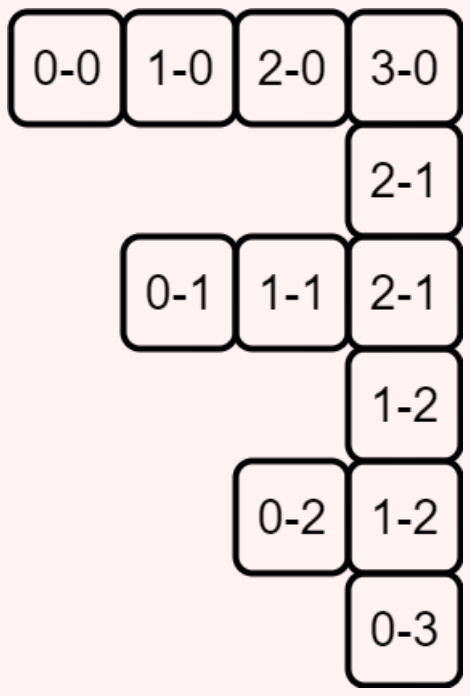

We can use the swisschart to reconstruct the bracket and signature.

**Figure 3.5.19:** $\mathcal{T}_3 = [[8;0;0;0]] \to [[1]] \to [[4;2;0]]_2 \to [[2]]_2 \to [[2;0]] \to [[1]]$.

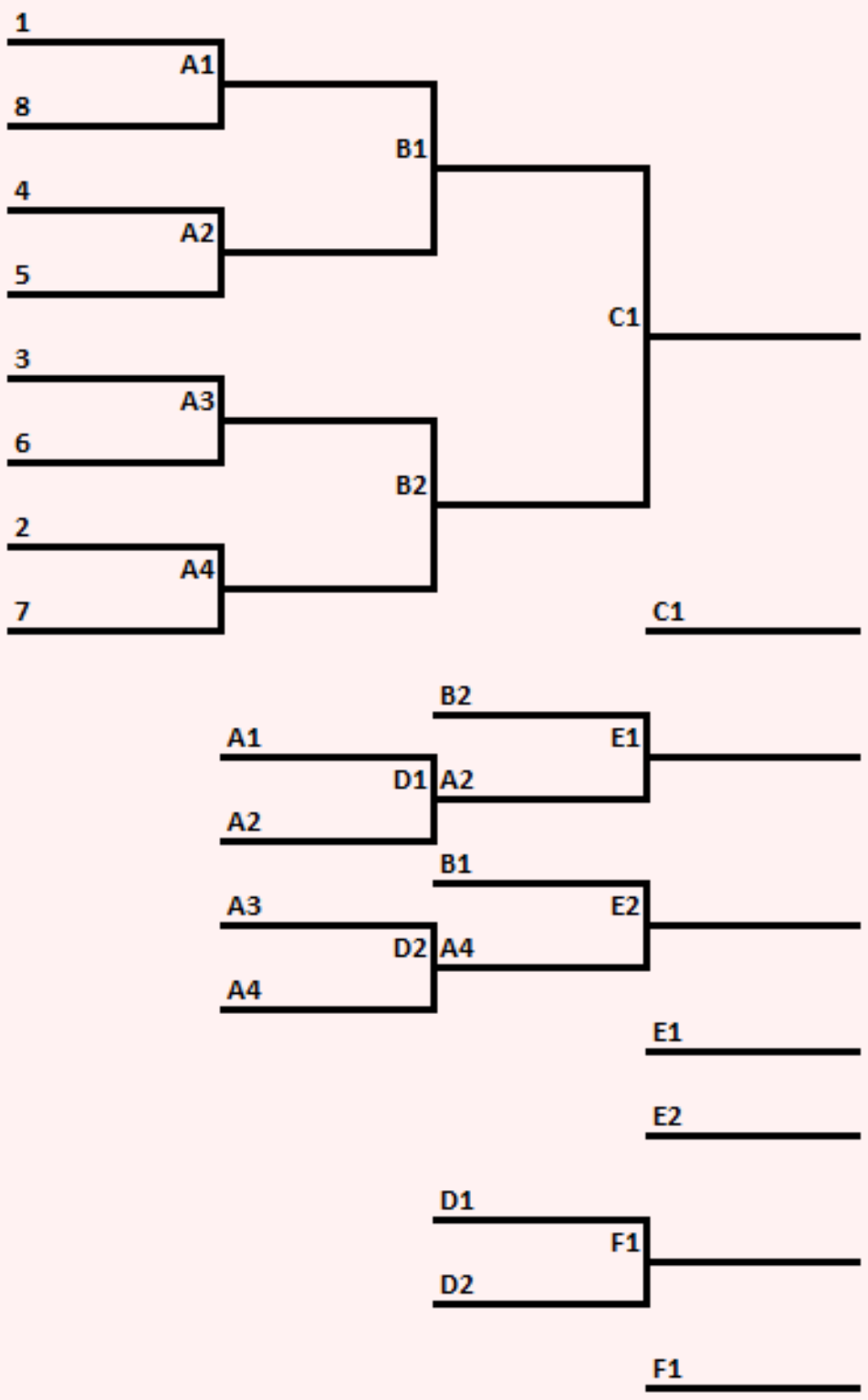

$\mathcal{S}_3$ and $\mathcal{T}_3$ differ in how they treat the teams that went 1-1. While $\mathcal{S}_3$ pairs teams that had their win and loss in the same order in games for either third or fifth place, $\mathcal{T}_3$ pairs teams that had their win and loss in different orders in games for fourth place. $\mathcal{T}_3$ is also very similar to the format used by the 2023 Southern Conference Wrestling Championships in Figure 3.1.7 on page 59: both use a primary eight-team balanced bracket and let their first-round losers fight their way back for a top-half finish.

**Theorem 3.5.20** **(Fried, 2024)**

$\mathcal{S}_3$ and $\mathcal{T}_3$ are the only compact 3-round Swiss signatures.

*Proof.* Any compact 3-round Swiss signature must begin with $[[\mathbf{8};\mathbf{0};\mathbf{0};\mathbf{0}]] \to [[\mathbf{1}]]$. Now let $\mathcal{A}$ be the semibracket that first-round primary brackets losers fall into. $\mathcal{A}$ must have two rounds, and the first-round primary bracket losers must all get no byes (otherwise they would not play the requisite three games). Thus $\mathcal{A} = [[\mathbf{4};\mathbf{a_1};\mathbf{0}]]_{(a_1/2+1)}$ for some $a_1$. As neither of the two semifinal winners can fall into $\mathcal{A}$, $a_1 \leq 2$. Additionally, if $a_1 = 1$, $\mathcal{A}$ would not be a signature. Thus, $a_1 = 0$ or 2.

If $a_1 = 0$, then in between the first two brackets and $\mathcal{A}$, we must have two more brackets for the second-round losers of the primary bracket: $[[\mathbf{2};\mathbf{0}]]$ and $[[\mathbf{1}]]$. Then $\mathcal{A}$ must be followed by $[[\mathbf{1}]]$ for the loser of its championship game, and then $[[\mathbf{2};\mathbf{0}]]$ and $[[\mathbf{1}]]$ so that the last two teams get a third game. This is the Swiss signature $\mathcal{S}_3$.

If $a_1 = 2$, then the losers of the two championship games of $\mathcal{A}$ have already played all three of their games and so need to fall into the bracket $[[\mathbf{2}]]$. Then we need $[[\mathbf{2};\mathbf{0}]]$ and $[[\mathbf{1}]]$ so that the last two teams get a third game. This is the Swiss signature $\mathcal{T}_3$. □

Figure 3.5.16 tells us that there are five 8-team Swiss signatures. How would a tournament designer decide which one to use? Well, it depends on what the prize structure of the format is. If the goal is to identify a top-three, then $\mathcal{S}_3$ is preferable: $\mathcal{T}_3$ doesn't even recognize a third-place, instead assigning fourth-place to two teams. But if the goal is to identify a top-four, $\mathcal{T}_3$ is preferable: the team that comes in fourth in $\mathcal{S}_3$ actually finishes with only one win, while the team that comes in fifth finishes with two. While it is still reasonable to grant the one-win team fourth-place – they had a more difficult slate of opponents – this is a somewhat messy situation that is solved by just using $\mathcal{T}_3$.

(McGarry and Schutz [15] considered outright swapping the positions of the fourth- and fifth-place teams at the conclusion of $\mathcal{S}_3$, but this format is not respectful and provides some incentive for losing in the first round in order to get an easier path to a top-half finish. Simply using $\mathcal{T}_3$ when identifying the top-four teams is preferable.)

For similar reasons, both formats are good for selecting a top-one or top-seven, and $\mathcal{S}_3$ but not $\mathcal{T}_3$ is good for selecting a top-five. Finally, it might seem that $\mathcal{S}_3$ and $\mathcal{T}_3$ are good formats for selecting a top-two or top-six: in both cases, the top two and top six teams are clearly defined, and there are no teams with better records that don't make the cut. However, notice that if we use $\mathcal{S}_3$ or $\mathcal{T}_3$ to select a top-two, the final round of games are meaningless: the two teams that finish in the top-two are the two teams that win their first two games, irrespective of how the third round of games went. Better than using either $\mathcal{S}_3$ or $\mathcal{T}_3$ would be to use the noncompact 2x$\mathcal{S}_2$, shortening the format down to two rounds without losing any important games.

We now count the number of $r$-round signatures.

**Theorem 3.5.21** **(Fried, 2024)**

Let $s_r$ be the number of compact $r$-round Swiss signatures. Then,

$$s_0 = s_1 = 1$$
$$s_r = s_{r-1} \cdot \sum_{i=1}^{r-1} s_i$$

*Proof.* Theorem 3.5.17 shows the cases for $r = 0$ or $r = 1$. For any other case, first name the semibracket that the first-round losers of the primary bracket fall into the "middle semibracket." Now when designing the swisschart for an $r$-round Swiss signature, one can consider the half above the middle semibracket and the half below the middle semibracket separately.

The half of the swisschart below the middle semibracket is straightforward: it looks just like the swisschart of an $(r-1)$ round Swiss signature, and in fact could look like the flow chart of any $(r-1)$ Swiss signature. Thus there are $s_{r-1}$ options.

The half above the middle semibracket is trickier. The first thing to note is that teams from the primary semibracket can continue to fall into the middle semibracket indefinitely, but once one round is skipped, respectfulness precludes any later rounds from falling down. Further, once teams are no longer falling into the middle semibracket, there are $s_i$ different ways the primary bracket could arrange the rest of its losers, where $i$ is the number of remaining rounds. Thus in total, there are $\sum_{i=1}^{r-1} s_i$ options for the half of the swisschart above the middle semibracket.

Therefore, there are $s_{r-1} \cdot \sum_{i=1}^{r-1} s_i$ different Swiss signatures in total. □

The eight compact 4-round signatures are displayed in Figure 3.5.22.

**Figure 3.5.22: The Eight Compact 4-round Swiss Signatures**

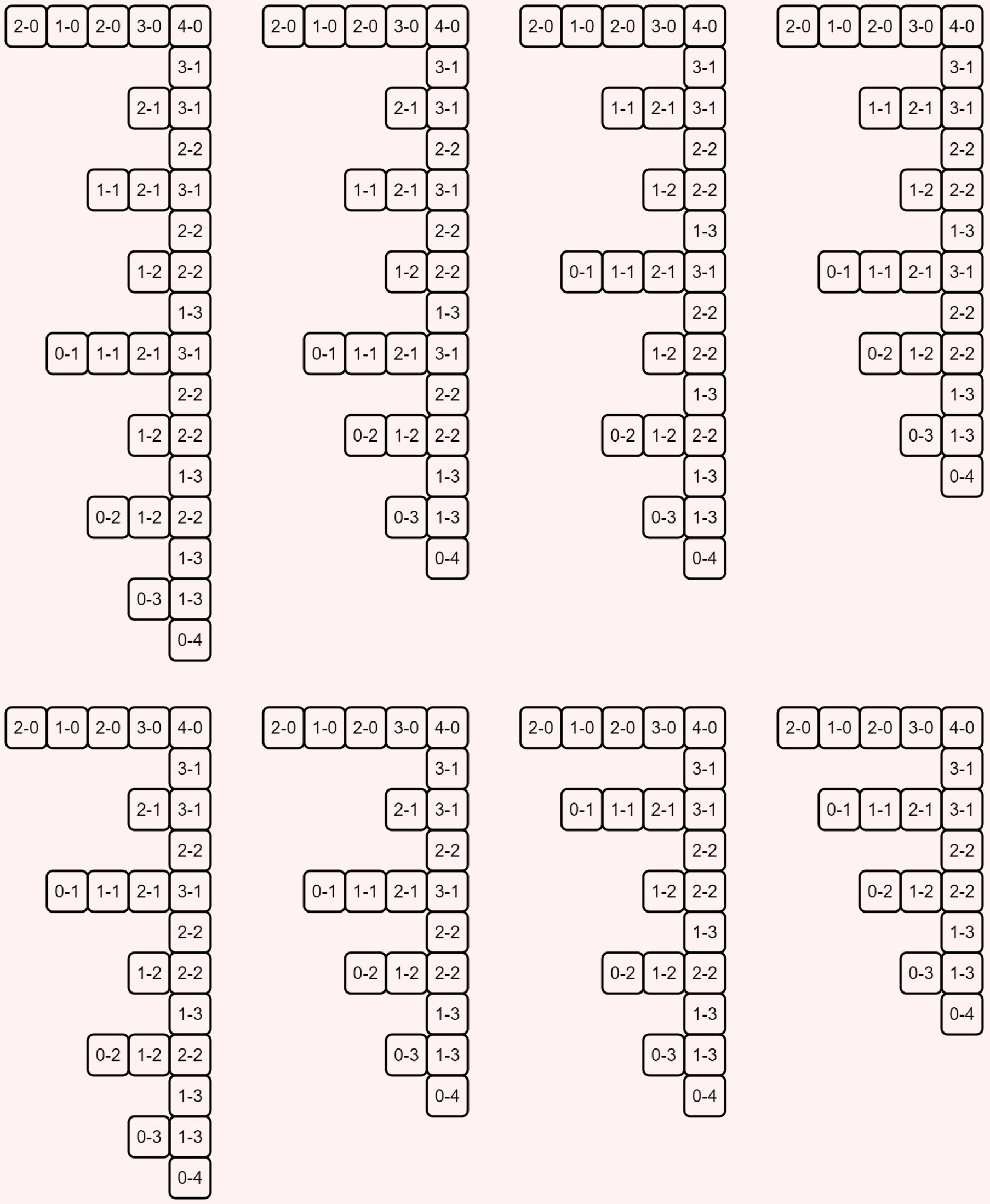

Overall, Swiss formats are very useful and practical tournament designs: they give each team the same number of games, they ensure that games are being played between teams that have the same record and thus, hopefully, similar skill levels, and, for many values of $m$, they efficiently identify a top-$m$ in a fair and satisfying way.

Further, Swiss or near-Swiss formats are great when the number of teams is exceedingly large. Even if not every requirement in Definition 3.5.5 is met, or the number of teams isn't a power of two, or the signature is not compact, or there is a round at the end that doesn't affect placement for important places, formats that are Swiss in spirit tend to do a great job of gathering a lot of meaningful data about a large number of teams in a small number of rounds. For this reason, they are often used in large tournaments for board or cards games, such as chess or Magic: The Gathering [40].

## 3.6 Efficient Multibrackets

In the past few sections, we have looked at multibrackets (and in particular linear multibrackets) as a solution to the tournament design question of how to crown a champion as well as give out certain consolation places.

We now consider a slightly different tournament design problem: we no longer care about which teams finish in first or any other specific place, only about which teams finish in the top-$m$ for a particular $m$. This is a problem commonly faced at regional tournaments in which the top-$m$ teams qualify for a national tournament: the ranking of the teams within the region aren't relevant, only which teams are above and below the cutoff.

Recall the format used in the 2023 Southern Conference Wrestling Championships [24].

**Figure 3.6.1: 2023 Southern Conference Wrestling Championships**

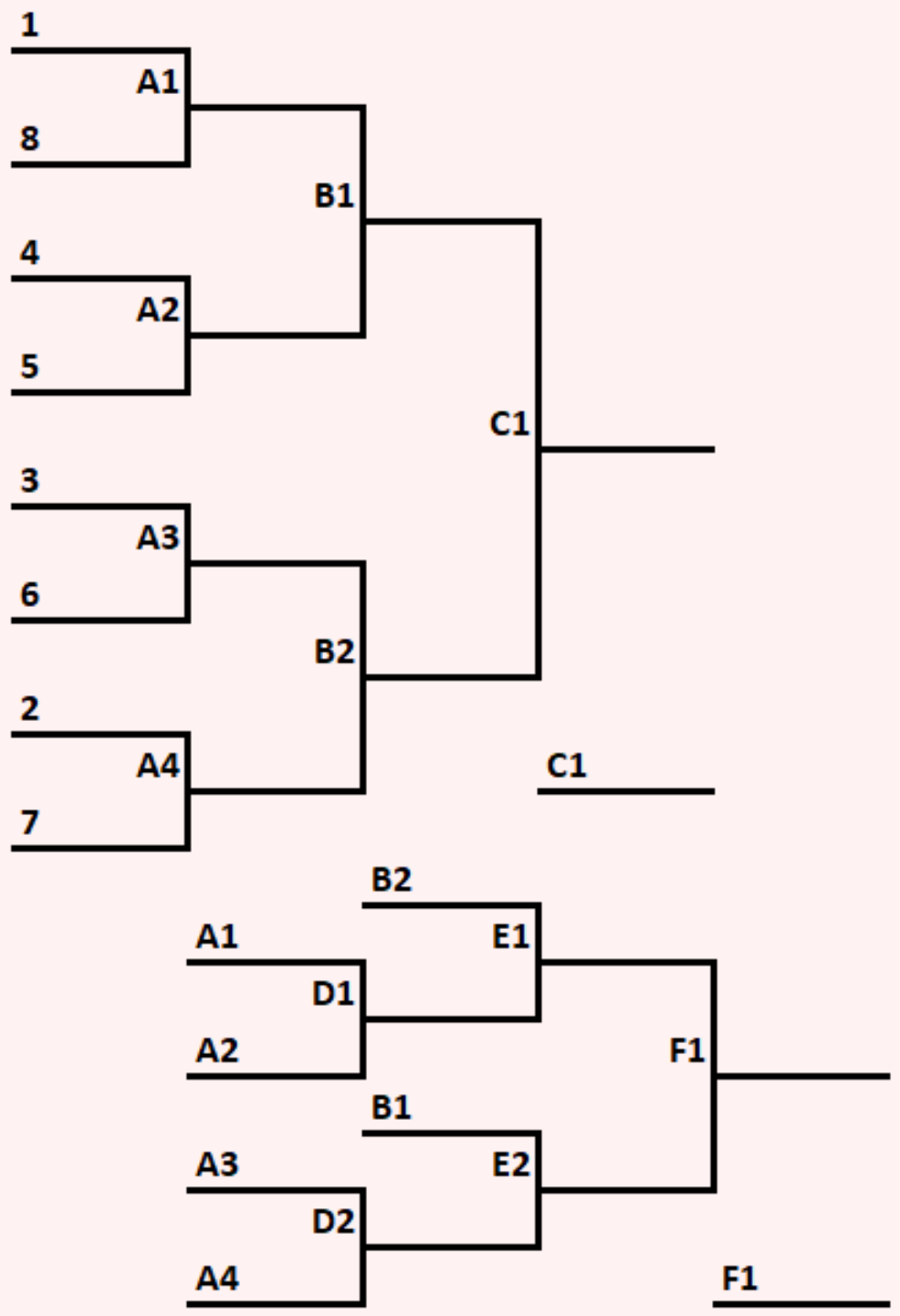

If we were only interested in the top four teams, rather than the rank of the team within those top four slots, games **C1** and **F1** become unnecessary: no matter what the results of those games are, the top four teams are the same. A more efficient format would leave those games unplayed, resulting in the following format.

**Figure 3.6.2: An Efficient Format for Selecting a Top Four**

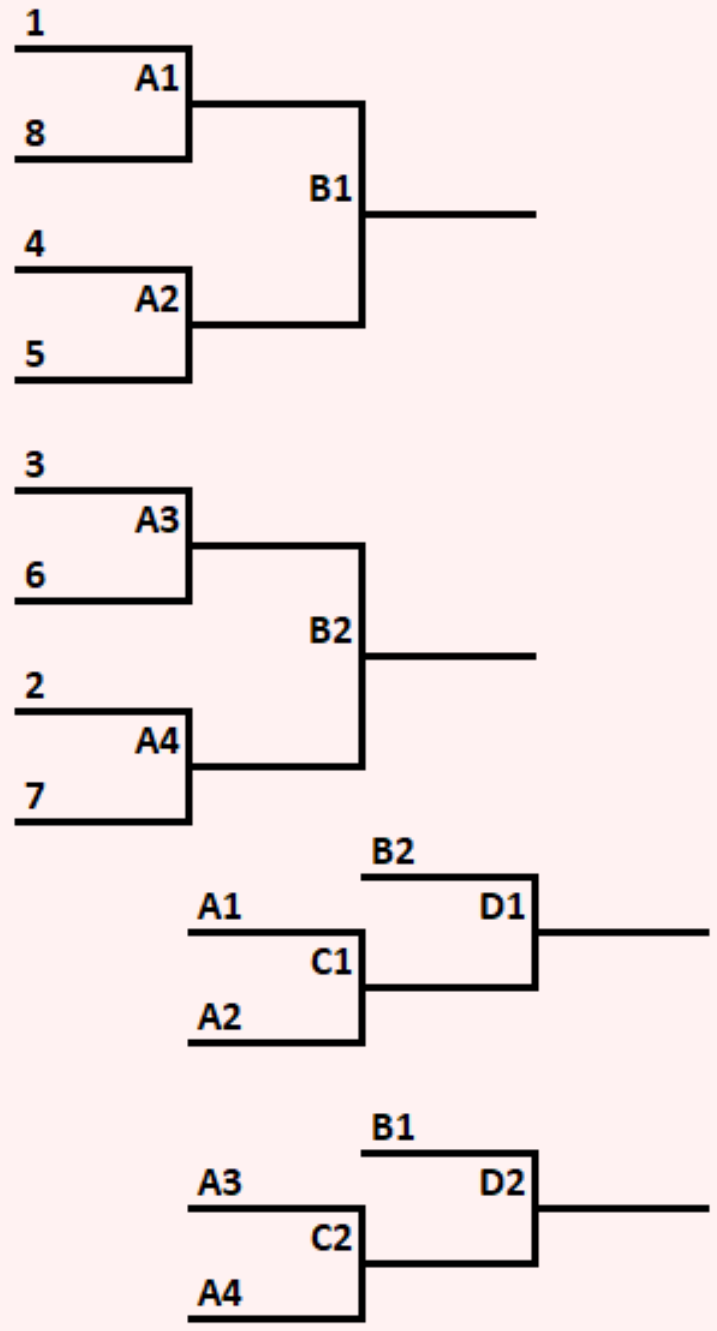

Instead of being composed of four traditional brackets, the format in Figure 3.6.2 is composed of two semibrackets each of which have rank two: one with the **A** and **B** round games, and one with the **C** and **D** round games. As desired, no games are played between two teams where both the winner and loser of each of those games are guaranteed to finish in the top four. (Nor are there any games where both the winner and loser are guaranteed to finish in the bottom four.)

This format has signature $[[\mathbf{8};\mathbf{0};\mathbf{0}]]_2 \to [[\mathbf{4};\mathbf{2}]]_2$, and it is *weakly efficient*.

**Definition 3.6.3: Weakly Efficient** (Fried, 2024)

A multibracket is *weakly efficient* if, once a team is guaranteed to be ranked by the format or guaranteed to finish unranked, they stop playing games.

Checking whether an arbitrary multibracket is weakly efficient requires fully examining the format, a process that can be quite arduous. But for respectful linear multibrackets, you can tell just from the signature.

**Theorem 3.6.4** (Fried, 2024)

A respectful linear multibracket with signature $\mathcal{A} = \mathcal{A}_1 \to \ldots \to \mathcal{A}_k$ is weakly efficient if and only if there is some integer $j$ with $1 \le j \le k$ such that every semibracket $\mathcal{A}_i$ with $i < j$ is trivial and every semibracket $\mathcal{A}_i$ with $i > j$ is competitive.

*Proof.* Assume first that such a $j$ exists. Let $\mathbf{G}$ be a game. Because all semibrackets $\mathcal{A}_i$ with $i < j$ are trivial, $\mathbf{G}$ must be in a semibracket $\mathcal{A}_i$ for $i \geq j$. Thus the loser of $\mathbf{G}$ is either eliminated outright, or falls into a competitive semibracket $\mathcal{A}_i$ for $i > j$, in which case they will play another game. If they continue losing, they will continue falling into competitive semibrackets, until they are eliminated outright and do not get ranked. Thus if a team competing in $\mathbf{G}$ loses the rest of their games, they will finish unranked. But of course if they win the rest of their games they will finish ranked, so $\mathcal{A}$ is weakly efficient.

Assume now that no such $j$ exists, so there exists some some $i$ such that $\mathcal{A}_i$ is nontrivial and $\mathcal{A}_{i+1}$ is noncompetitive. Thus by Lemma 3.4.11, any team that loses in the championship game of $\mathcal{A}_i$ will win $\mathcal{A}_{i+1}$. $\mathcal{A}_i$ is nontrivial so it has at least one championship game: the winner of that game wins $\mathcal{A}_i$, and the loser of that game wins $\mathcal{A}_{i+1}$, so $\mathcal{A}$ is not weakly efficient. □

USA Ultimate, the governing body for the sport of ultimate frisbee in America, runs a series of sectional and regional tournaments to determine which $m$ out of $n$ teams should advance to the regional or national tournament, respectively. Unsurprisingly, the USA Ultimate Manual of Championship Series Tournament Formats [36], contains a host of weakly efficient linear multibrackets for various values of $1 \leq m \leq 12$ and $3 \leq n \leq 24$ after a "regular season" portion of the tournament has been played to establish seeds.

A couple of examples are Figure 3.6.5, which selects a top six out of seven, and Figure 3.6.6, which selects a top five out of sixteen. (In reality, sometimes additional games are played to determine placements within the top-$m$, but we display only the weakly efficient part of the format here.)

**Figure 3.6.5:** $[[\mathbf{1}]] \to [[\mathbf{1}]] \to [[\mathbf{1}]] \to [[\mathbf{4};\mathbf{0}]]_2 \to [[\mathbf{2};\mathbf{0}]]$

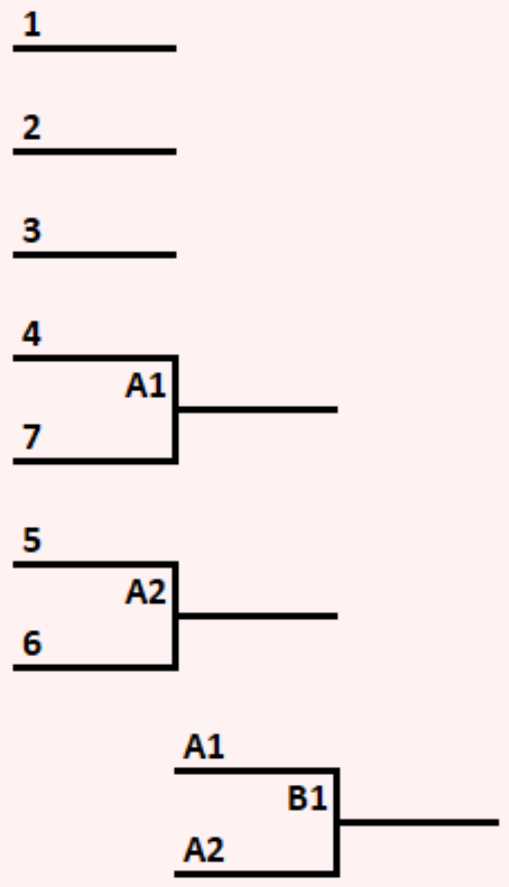

**Figure 3.6.6:** $[[8; 0; 0]]_2 \to [[4; 2; 0]]_2 \to [[8; 0; 2; 2; 0; 0]]$

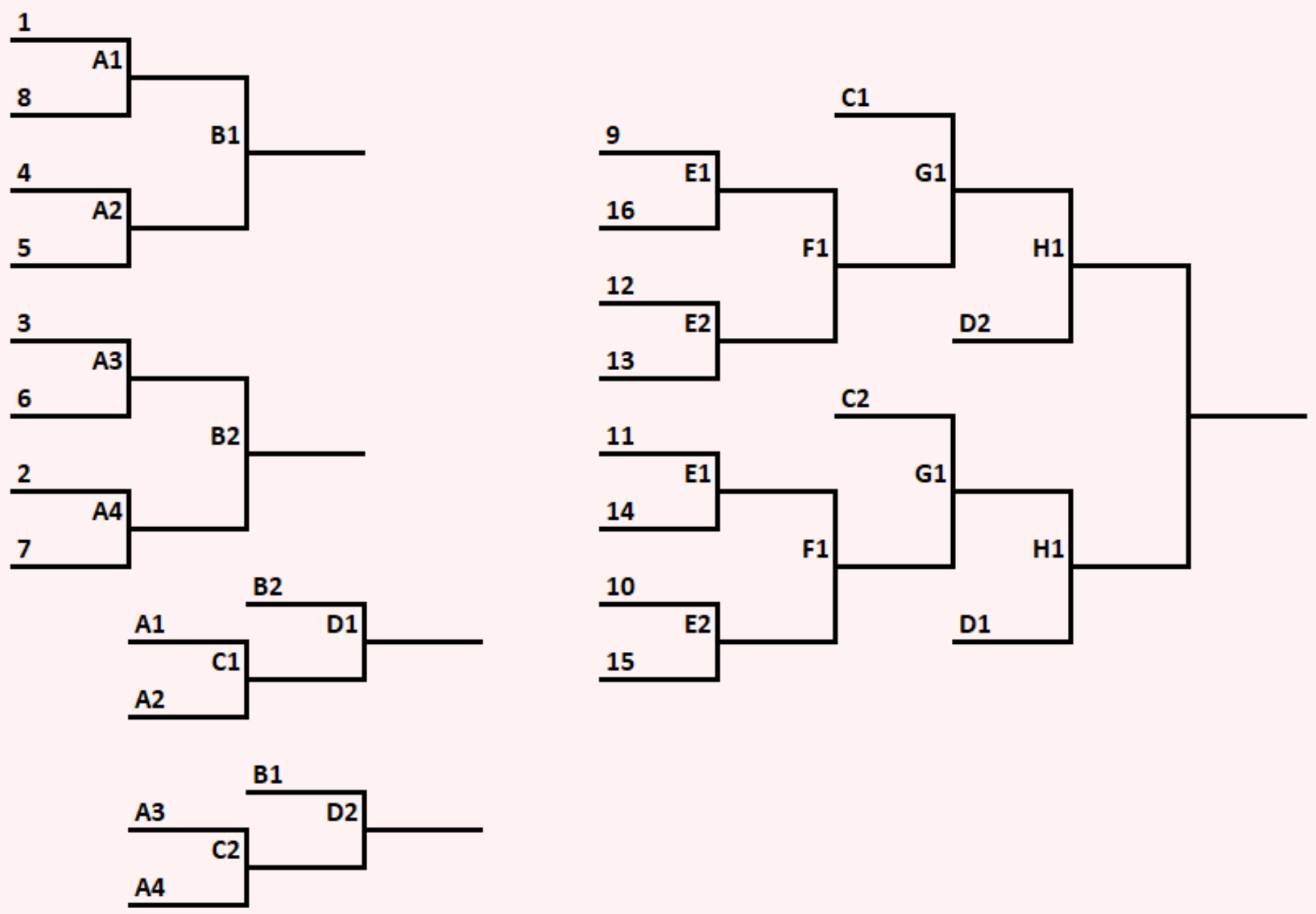

We note two things about the notion of weak efficiency presented above. First, Theorem 3.6.4 implies that a weakly efficient multibracket can begin with a long string of trivial semibrackets before the nontrivial ones begin. While this is sufficient for avoiding playing unnecessary games, it does not completely remove unnecessary semibrackets: the set of leading trivial semibrackets

$$[[\mathbf{m_1}]]_{m_1} \to \ldots \to [[\mathbf{m_j}]]_{m_j}$$

of a weakly efficient multibracket can be combined into a single trivial semibracket

$$[[\mathbf{m_1} + \ldots + \mathbf{m_j}]]_{(m_1 + \ldots + m_j)}$$

without affecting which teams end up ranked. Applying this to the format in Figure 3.6.5 yields a signature of

$$[[\mathbf{3}]]_3 \to [[\mathbf{4}; \mathbf{0}]]_2 \to [[\mathbf{2}; \mathbf{0}]].$$

In fact, if there is at least one game played in a weakly efficient multibracket, trivial semibrackets can removed entirely, converting a multibracket of signature

$$[[\mathbf{m_1}]]_{m_1} \to [[\mathbf{a_1}; \ldots; \mathbf{a_r}]]_{m_2} \to \ldots \to \mathcal{A}_k$$

into one of signature

$$[[\mathbf{a_1}; \ldots; \mathbf{a_r} + \mathbf{m_1}]]_{m_1 + m_2} \to \ldots \to \mathcal{A}_k.$$

Applying this to the format in Figure 3.6.5 yields a signature of

$$[[\mathbf{4}; \mathbf{3}]]_5 \to [[\mathbf{2}; \mathbf{0}]].$$

To patch this, we strengthen the notion of weak efficiency into just *efficiency.*

**Definition 3.6.7: Efficient** **(Fried, 2024)**

A respectful linear multibracket is *efficient* if one of three conditions hold:

(a) It is a single trivial semibracket.

(b) It is a sequence of competitive semibrackets.

(c) It is a single nontrivial noncompetitive semibracket followed by a sequence of competitive semibrackets.

Theorem 3.6.4 says that in each of these three cases no games are played between teams guaranteed to be ranked, and the process detailed above can reduce any weakly efficient signature into a signature that takes one of those three forms.

The second thing to note is that efficiency makes a lot of sense if we are only interested in the top-$m$ teams (where $m$ is the sum of the ranks of the semibrackets in our format) and not in the rankings of the teams within them. But sometimes we might be interested in the intermediate rankings as well. For example, let's say we want to design an eight-team tournament format in which the top team receives the grand prize, second-place receives a second-place prize, while the third- through seventh-place each get equivalent consolation prizes, and last place gets nothing. While not efficient (or even weakly efficient), the following format assigns the desired places without playing any games between teams that are guaranteed to receive the same prize.

**Figure 3.6.8:** $[[8; 0; 0; 0]] \to [[1]] \to [[4; 2]]_4 \to [[2; 0]]$

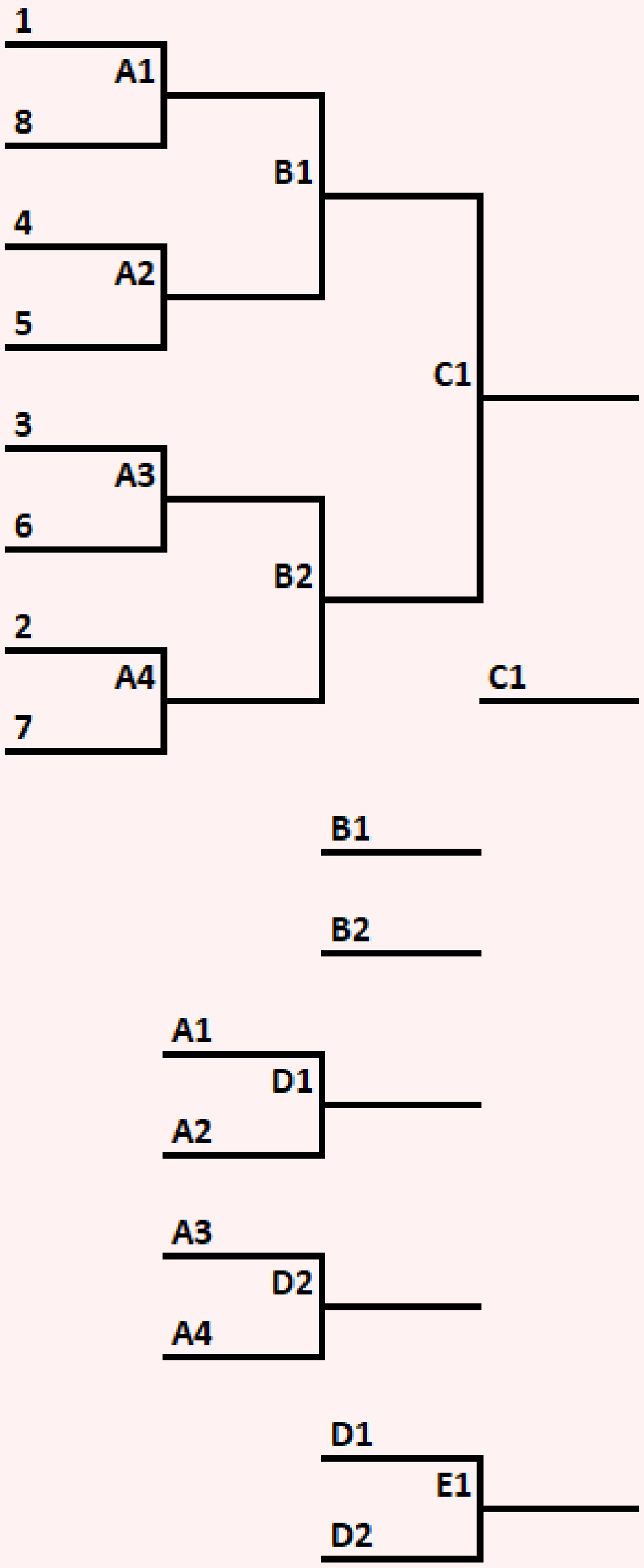

To account for this, we introduce the notion of a *prize structure*.

**Definition 3.6.9: Prize Structure** (Fried, 2024)

A *prize structure* $\mathcal{P}$ is a sequence $(p_1, ..., p_m)$ indicating that the top $p_1$ teams in a format receive some prize, the next $p_2$ receive some smaller prize, etc. Any teams finishing in place $1 + \sum_{i=1}^{m} p_i$ or worse receive no prize.

Then,

**Definition 3.6.10: Efficient with Respect to a Prize Struture (Fried, 2024)**

A respectful linear multibracket $\mathcal{A} = \mathcal{A}_1 \to \dots \to \mathcal{A}_k$ is *efficient with respect to a prize structure* $\mathcal{P} = (p_1, \dots, p_m)$ if

(a)

$$\sum_{i=1}^{k} \text{Rank}(\mathcal{A}_i) = \sum_{i=1}^{m} p_i,$$

(b) $\mathcal{A}_j$ being noncompetitive implies that for some $\ell < m$,

$$\sum_{i=1}^{j-1} \text{Rank}(\mathcal{A}_i) = \sum_{i=1}^{\ell} p_i,$$

and

(c) $\mathcal{A}_j$ being trivial implies that for some $\ell \leq m$,

$$\sum_{i=1}^{j} \text{Rank}(\mathcal{A}_i) = \sum_{i=1}^{\ell} p_i.$$

(The first condition ensures that teams stop playing games once they have eliminated from prize contention, the second condition ensures that teams stop playing games once their prize can no longer change, and the last condition ensures that there are no trivial semibrackets that could be combined with another semibracket as per the process detailed before Definition 3.6.7.)

So the respectful linear multibracket $[[\mathbf{8}; \mathbf{0}; \mathbf{0}; \mathbf{0}]] \to [[\mathbf{1}]] \to [[\mathbf{4}; \mathbf{2}]]_4 \to [[\mathbf{2}; \mathbf{0}]]$ is efficient with respect to the prize structure $(1, 1, 5)$. A linear multibracket being efficient is the same as it being efficient with respect the prize structure $(m)$, where $m$ is the sum of the ranks of its semibrackets.

Efficient formats are great for tournaments whose primary goal is to select the top $m$ teams to move on to the next stage of the competitions, as discussed in the beginning of this section. They do so excitingly, with each spot in the top-$m$ being awarded as the winner of a particular game; efficiently, with no games being played between teams who will receive the same prize; and fairly, as respectfulness ensures that winning is always better than losing. It is not surprising that many sports with regional tournaments that qualify teams for a national one use such formats.

## 3.7 Nonlinear Multibrackets

In the last four sections we have focused our study on *linear* multibrackets, which have the property that when a team loses in a given semibracket they drop into a different, lower semibracket. But many leagues use nonlinear multibrackets as well, and so while the difficulty in assigning signatures makes us less equipped to study them, we will still take a look at the space of such formats.

An simple example of a nonlinear multibracket was the format used by the 2019 Suncorp Super Netball Playoffs [29], sometimes called the Page-McIntyre system [23].

**Figure 3.7.1: 2019 Suncorp Super Netball Playoffs**

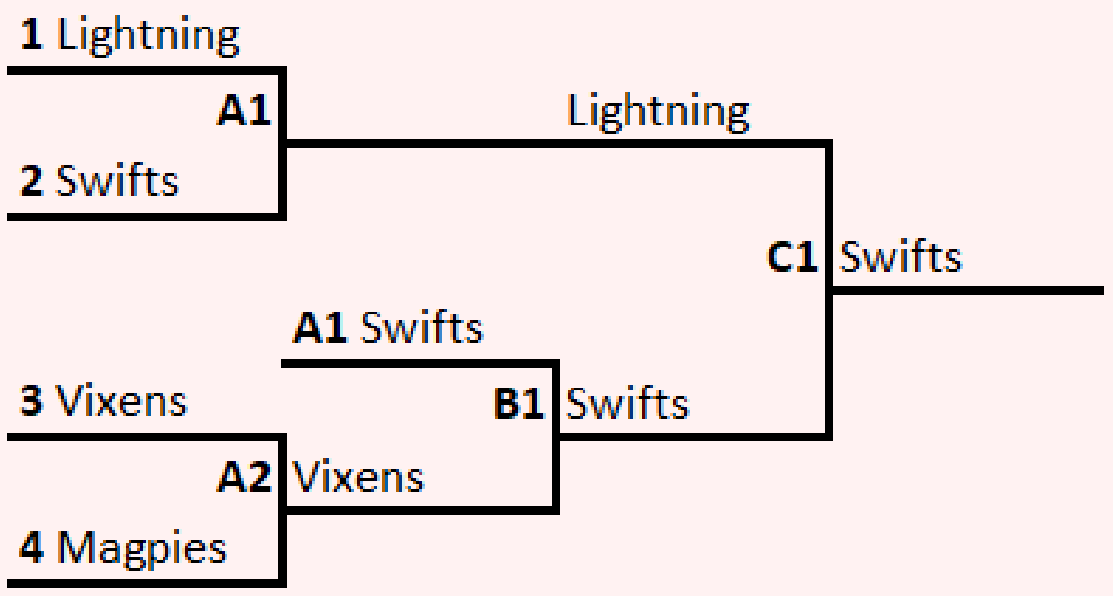

Nonlinear multibrackets are a bit strange: while the winner of game **A1** goes directly to the final, the loser falls into the semifinal of the *same bracket*. This poses problems in out attempts to define both a notion of signature as well as a notion of properness.

Beginning with signature, the shape of the multibracket is a bit strange: the winner of game **A1** gets a bye *after winning a game*, something that never happens in a traditional bracket. Attempts to give this multibracket a signature might lead to $[[\mathbf{4};\mathbf{1};\mathbf{0};\mathbf{0}]]$ or even $[[\mathbf{4};\mathbf{0};\mathbf{0};\mathbf{0}]]$, neither of which are actually bracket signatures (they both violate Theorem 2.1.14). The issue here is that game **A1** is actually a semifinal, and so "should" (if it didn't deliver its loser to the other semifinal) live in the second round, producing a signature of $[[\mathbf{2};\mathbf{3};\mathbf{0};\mathbf{0}]]$. But then of course this format is quite different from traditional brackets with that same signature. Signatures on nonlinear multibrackets are in general not well-defined.

To make matters worse, the first round appears to have an "improper" set of matchups: the games are 1v2 and 3v4 rather than then "proper" 1v4 and 2v3. However, properness is a much trickier concept for nonlinear multibrackets. While the 1- and 2-seeds to have tougher first round matchups than the 3- and 4-seeds, this is compensated by them getting an extra life: if they lose, they play the winner of the 3v4 matchup, while the 3v4 loser is just eliminated. In total, no team would prefer to be seeded lower than they are. One could imagine developing this intuition of extra lives into a formal notion of properness, but we leave that question untreated.

**Open Question 3.7.2** **(Fried, 2024)**

How can we define signatures and properness for nonlinear multibrackets?

One thing of note about the Page-McIntyre system in particular is that **C1** can be a rematch of game **A1**. In fact, this is pretty likely: if the format goes chalk, the 1- and 2-seed will find themselves replaying the game they played just two rounds ago. In Figure 3.7.1, the tournament did not go chalk, but game **C1** was still a rematch. As we discussed in Section 3.4, this can be pretty unsatisfying: indeed, in the 2019 Suncorp Super Netball Playoffs, the Swifts and the Lightning each beat each other once, but the Swifts won the game that mattered and so was declared champion.

One option would be to be to make game **C1** contingent on it not being a rematch: if it is a rematch, then the game is skipped and whichever team won the previous game is declared champion. While this solution is effective for the second-place game in our alternative AFL Asian Cup format (Figure 3.1.5 on page 57), it doesn't work here. Making the game contingent would mean that the loser of **A1** is actually eliminated upon their loss: even if they win **B1**, they wouldn't have the ability to play in the championship game.

A better solution might be a *double-elimination tournament*, as employed by the 2016 NCAA Softball Ann Arbor Regional [28].

**Figure 3.7.3: 2016 NCAA Softball Ann Arbor Regional**

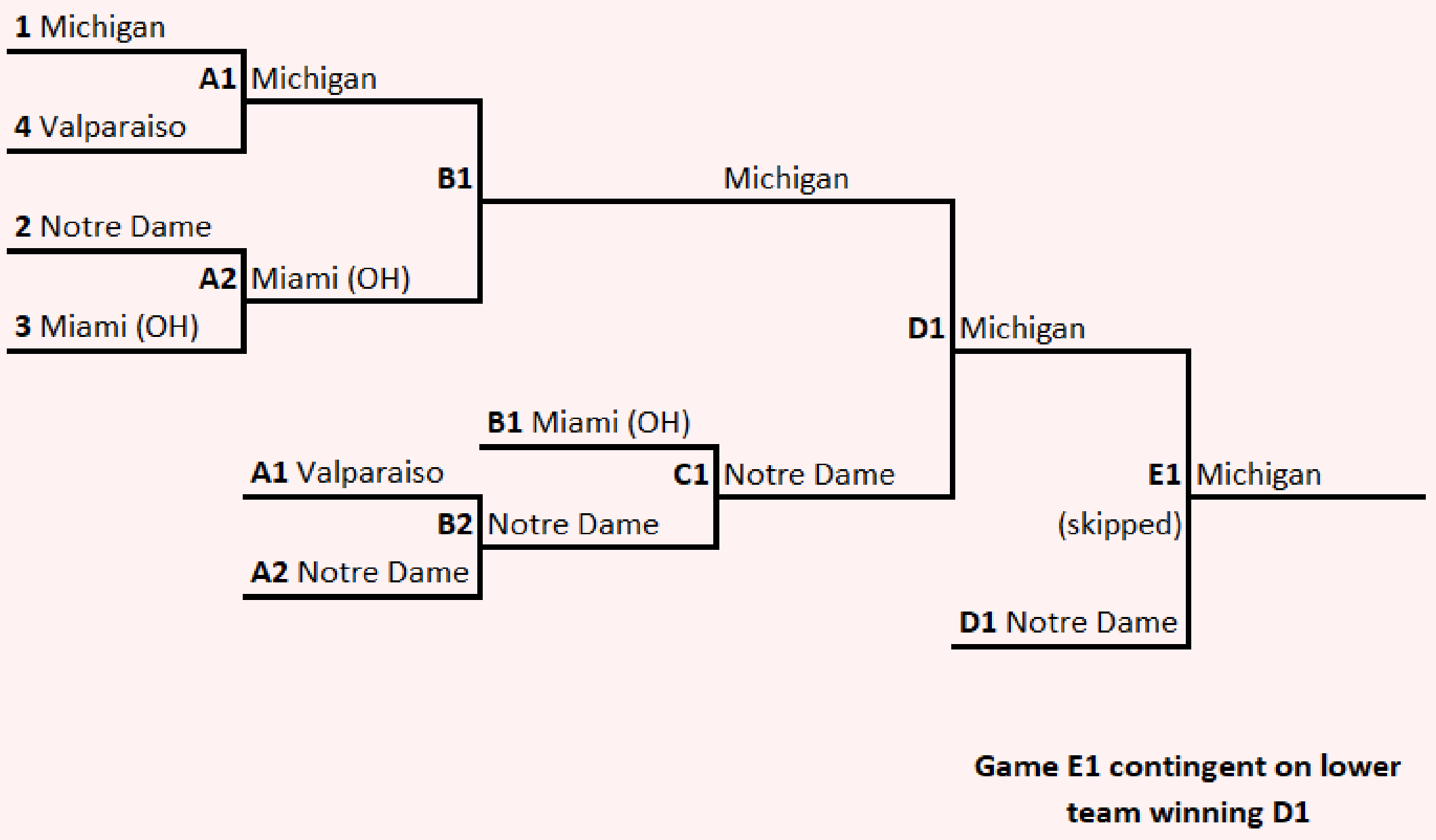

**Definition 3.7.4: Double-Elimination Tournament** **(Unattributed)**

A *double-elimination* tournament is a multibracket (plus one contingent game) consisting of a *winners' bracket*, where every team starts, a *losers' bracket*, that every winners' bracket loser falls into, and a *grand finals*, in which the winner of the winners' bracket plays the winner of the losers' bracket for the championship, with the losers' bracket winner needing to win twice, while the winners' bracket winner only needs to win once.

A double-elimination tournament guarantees that the winner will finish undefeated or with only one loss, while every other team finishes with two.

The 2016 NCAA Softball Ann Arbor Regional is an example of a double-elimination format: the winners' bracket consists of games **A1**, **A2**, and **B1**; the losers' bracket consists of games **B2** and **C1**; and the grand finals of game **D1** and then, if necessary, **E1**. Michigan finished undefeated while Valparaiso, Notre Dame, and Miami (OH) each finished with two losses.

Because double-elimination tournaments are so common, and all use a contingent game that is played only if the lower team wins (**E1** in the case of Figure 3.7.3), that contingent game has a name.

**Definition 3.7.5: Recharge Game** **(Dabney, 2017)**

A *recharge game* is a contingent game in a multibracket that is a rematch of a previous game and played only if the lower team won the first game.

Recharge games are so common that we introduce a special notation: if the name of a game has a star after it, then that game is followed by a recharge game (if necessary). This allows us to condense the format in the Figure 3.7.3 a little bit, as displayed in Figure 3.7.6.

**Figure 3.7.6: 2016 NCAA Softball Ann Arbor Regional**

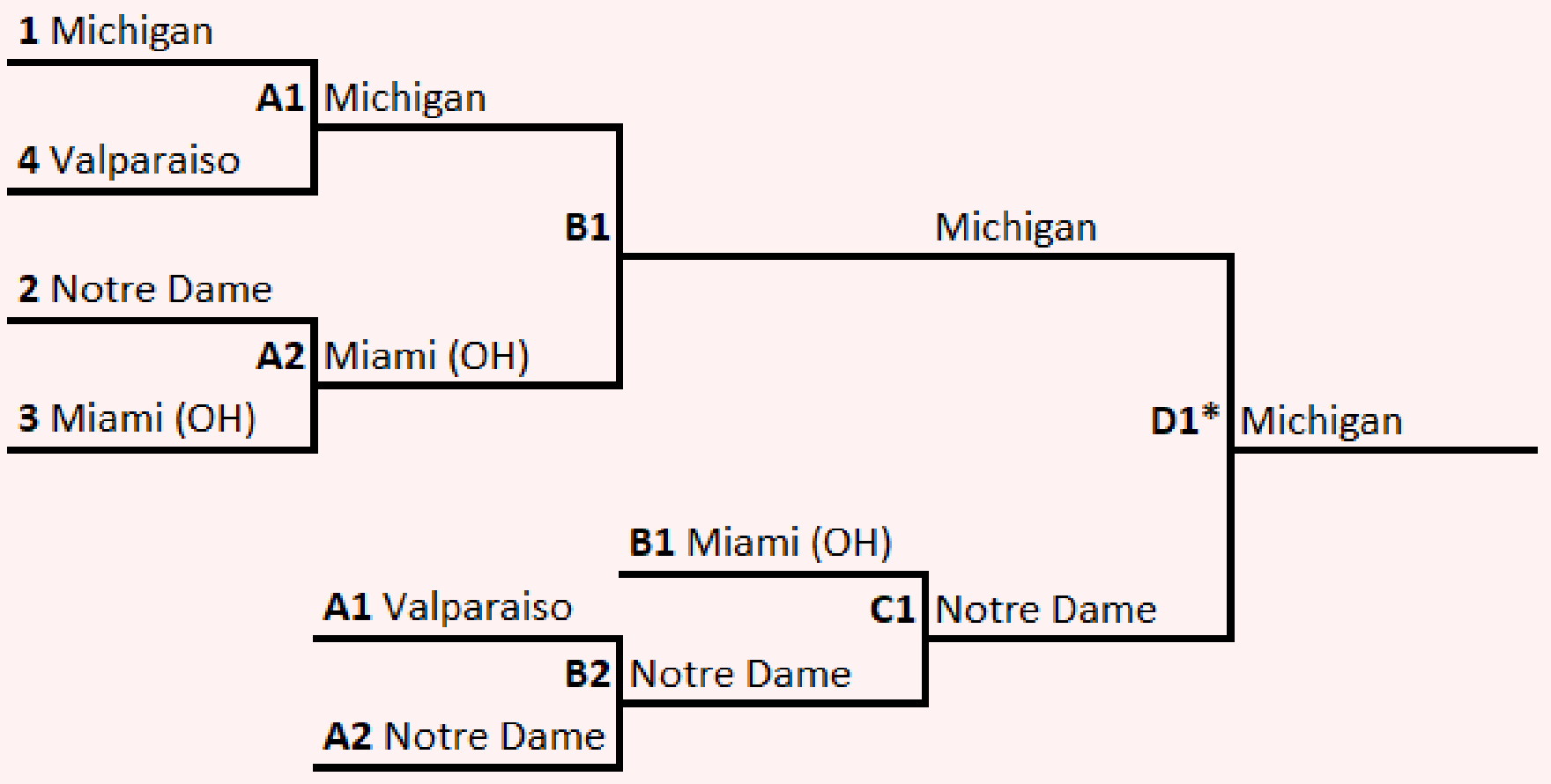

The only issue with this notation is that, if the recharge game was triggered but won by the upper team, there is no natural place to denote that the recharge game was played. We adopt the convention of writing the lower team *under* the line that the winner of the recharge game is placed over in this case. This is depicted in Figure 3.7.7.

**Figure 3.7.7: Figure 3.7.6 if Notre Dame Beat Michigan Once**

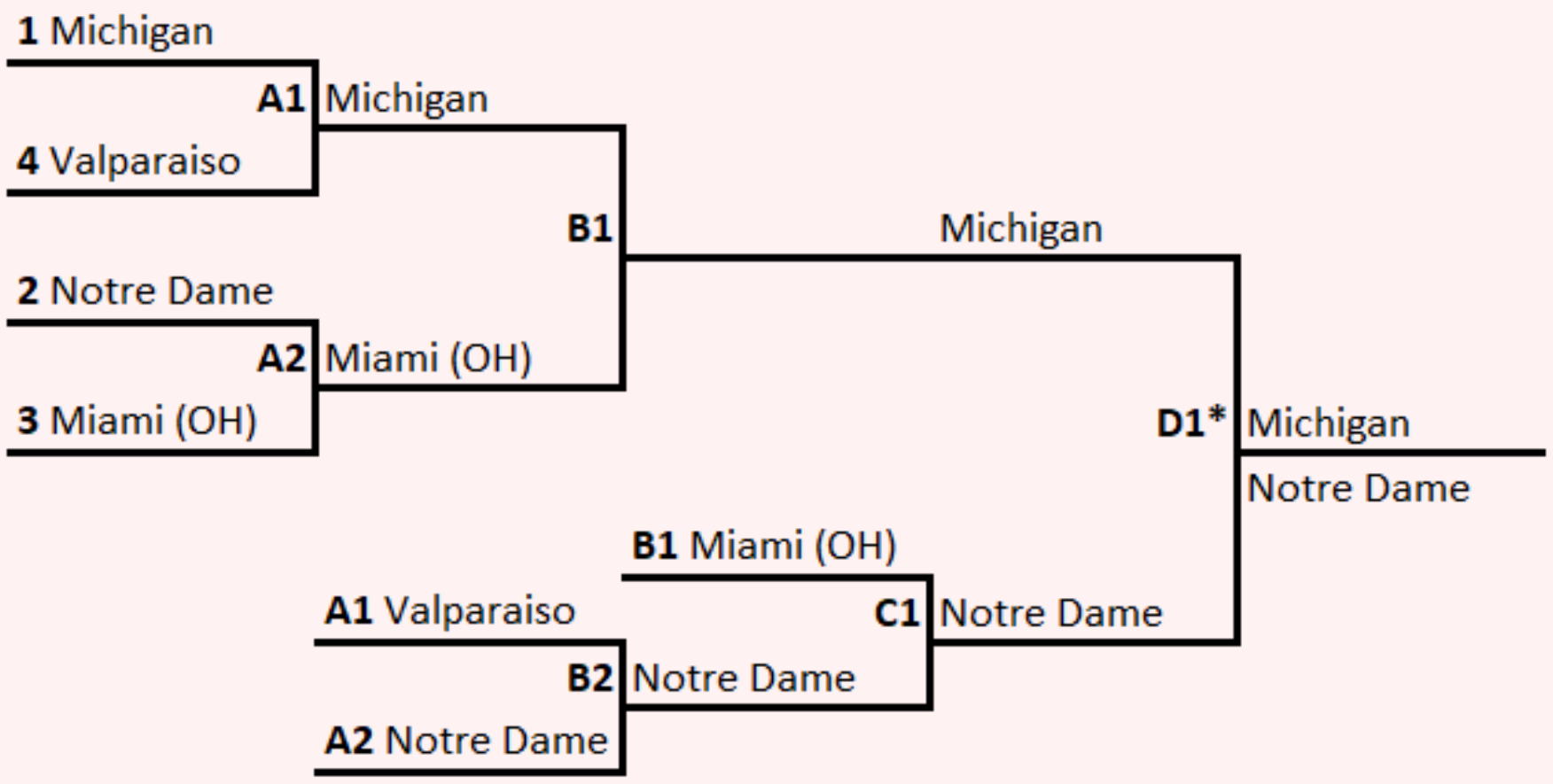

While the recharge game is necessary to ensure that the format is a truly a double-elimination tournament, as well as preventing the problem in the Page-McIntyre System where the champion and runner-up each finish with one-loss, it's not all upside. For one thing, Dabney [8] found some evidence that a tournament with no recharge game actually does a better job of crowning the best team as champion than the double-elimination with the recharge game included. Additionally, formats with recharge games tend to be less exciting, as they risk not playing a true championship game (a game in which either team wins the format if they win that game).

In any case, double-elimination tournaments are a powerful tool in a tournament designer's arsenal whether or not a recharge game is used, as they are in some sense more accurate than their single-elimination counterparts. We prove this fact for a simplified case where the winners' and losers' bracket are relatively nice, and where there is a single best team that is favored against every other team with a constant probability $1/2 < p < 1$.

**Theorem 3.7.8** **(Fried, 2024)**

Let $r$ be a positive integer, $p$ be a probability such that $1/2 < p < 1$, and $\mathcal{T}$ be a list of $2^r$ teams with a team $t \in \mathcal{T}$ such that for every other team $s$,

$$\mathbb{P}[t \text{ beats } s] = p.$$

Let $\mathcal{A}$ be the balanced bracket on $2^r$ teams, let $\mathcal{B}$ be a bracket on $2^r - 1$ teams such that the linear multibracket $\mathcal{A} \to \mathcal{B}$ is respectful, and let $\mathcal{C}$ be the double-elimination format with winners' bracket $\mathcal{A}$ and losers' bracket $\mathcal{B}$. Then,

$$\mathbb{W}_{\mathcal{C}}(t, \mathcal{T}) \geq \mathbb{W}_{\mathcal{A}}(t, \mathcal{T})$$

with equality only when $r = 1$ and there is no recharge round.

*Proof.* To win $\mathcal{A}$, $t$ simply has to win $r$ games. Thus

$$\mathbb{W}_{\mathcal{A}}(t, \mathcal{T}) = p^n.$$

Now consider $\mathcal{C}$. Let $s$ be the number of rounds in $\mathcal{B}$, let $s_i$ be the round of $\mathcal{B}$ that teams that lose in the $i$th round of $\mathcal{A}$ fall into, and let $c_i = s - s_i + 1$, so teams that lose in the $i$th round of $\mathcal{A}$ need to win $c_i$ games in $\mathcal{B}$ in order to make the grand finals.

Since there are $2^{r-i}$ $i$-round losers, by Theorem 2.1.14,

$$\sum_{i=1}^{r} 2^{r-i} \cdot \left(\frac{1}{2}\right)^{c_i} = 1,$$

so,

$$\sum_{i=1}^{r} \left(\frac{1}{2}\right)^{c_i+i-1} = \left(\frac{1}{2}\right)^{r-1}. \tag{$*$}$$

Letting $q = 1 - p$, note that $t$ wins the winners' bracket with probability $p^r$, and the losers' bracket with probability

$$\sum_{i=1}^{r} p^{i-1} \cdot q \cdot p^{c_i} = q \cdot \sum_{i=1}^{r} p^{c_i+i-1} \geq q \cdot p^{r-1},$$

with the inequality coming by equation $(*)$, and with equality only when $r = 1$.

Now, if there is a recharge round, then

$$\begin{aligned}
\mathbb{W}_{\mathcal{C}}(t, \mathcal{T}) &= \mathbb{W}_{\mathcal{A}}(t, \mathcal{T}) \cdot (p + qp) + \mathbb{W}_{\mathcal{B}}(t, \mathcal{T}) \cdot p^2 \\
&\geq p^r(p + qp) + (q \cdot p^{r-1}) \cdot p^2 && \text{with equality only when } r = 1 \\
&= p^r(p + 2qp) \\
&> p^r \\
&= \mathbb{W}_{\mathcal{A}}(t, \mathcal{T}).
\end{aligned}$$

If there is no recharge round, then

$$\begin{aligned}
\mathbb{W}_{\mathcal{C}}(t, \mathcal{T}) &= \mathbb{W}_{\mathcal{A}}(t, \mathcal{T}) \cdot p + \mathbb{W}_{\mathcal{B}}(t, \mathcal{T}) \cdot p \\
&\geq p^r \cdot p + (q \cdot p^{r-1}) \cdot p && \text{with equality only when } r = 1 \\
&= p^r \\
&= \mathbb{W}_{\mathcal{A}}(t, \mathcal{T}).
\end{aligned}$$

Thus,

$$\mathbb{W}_{\mathcal{C}}(t, \mathcal{T}) \geq \mathbb{W}_{\mathcal{A}}(t, \mathcal{T})$$

with equality only when $r = 1$ and there is no recharge round. ∎

We conclude our discussion of nonlinear multibrackets with a few more interesting examples. The first is the 2022 Big Ten Baseball Tournament [39].

**Figure 3.7.9: 2022 Big Ten Baseball Tournament**

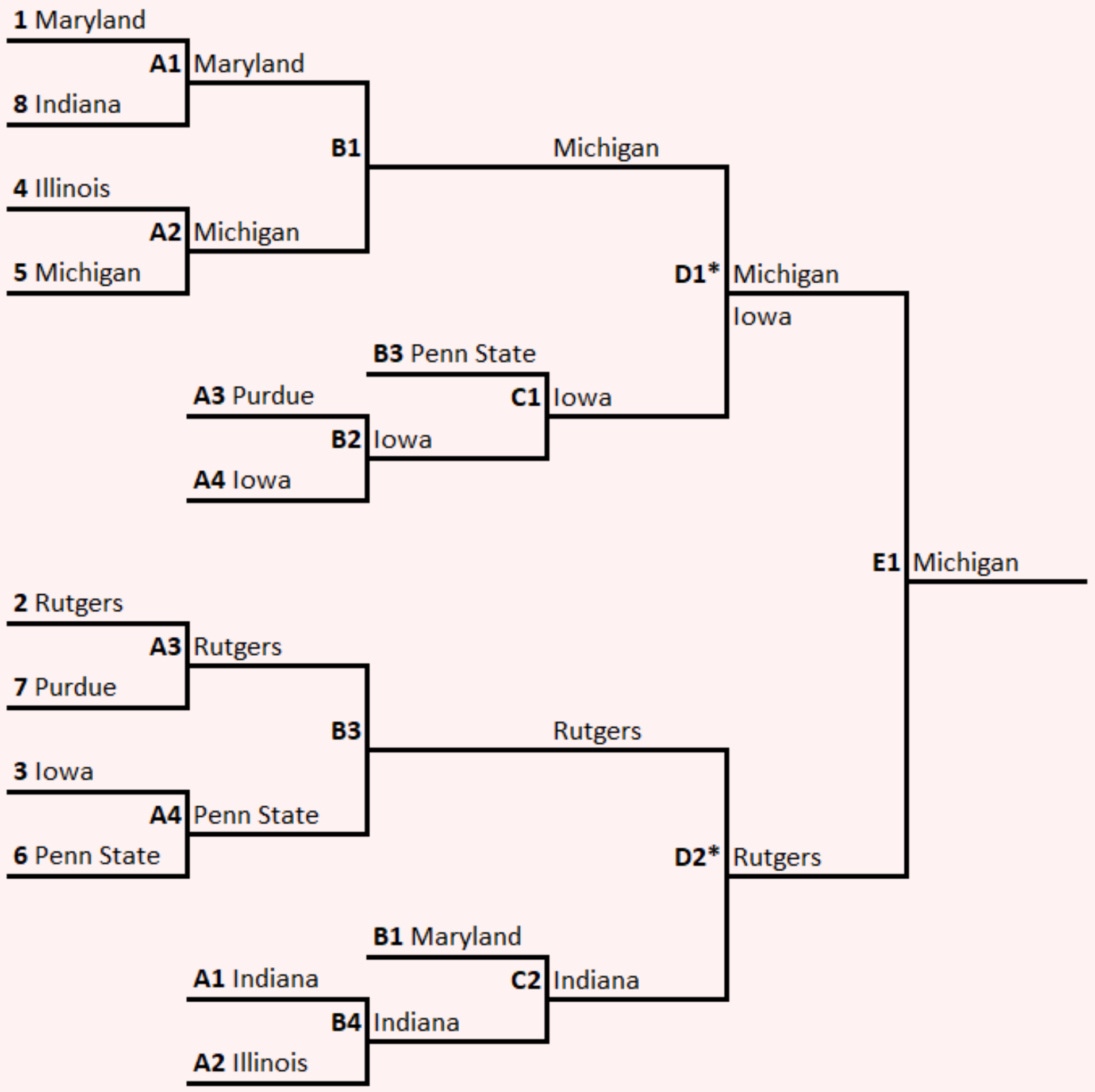

The 2022 Big Ten Baseball Tournament wanted to balance two effects: first, that double-elimination formats lead to more accurate results, but second, that championship games are exciting and double-elimination games risk not including one. 2022 Big Ten Baseball Tournament innovates to solve the latter issue by including recharge games in the *semifinals*, and then having the championship game be single winner-take-all game.

Note that this format does not fully solve all the problems it is attempting to tackle: for one thing, it is not a true double-elimination, as Rutgers gets eliminated with only a single loss. That said, Michigan is unambiguously the most deserving winner: every team other than Michigan and Rutgers lost once, and Michigan defeated Rutgers in their one matchup.

However, this property was not guaranteed: had Penn State beaten Iowa in game **C1**, Michigan twice in game **D1** and the recharge game, and then Rutgers in the final, we would be back to the issue with the Page-McIntyre System. Penn State and Rutgers would have each finished with only one loss to the other team, with the champion being determined somewhat arbitrarily by who won the most recent game. This illustrates an important point: the desire for an unambiguous champion and the desire for an unambiguous championship *game* are fundamentally in conflict in the world of nonlinear multibrackets.

Another interesting nonlinear multibracket of note is the NBA Playoffs. You may recall from Figure 2.2.4 that in 2004, the NBA Eastern Conference Playoffs used a simple bracket of signature $[[\mathbf{8}; \mathbf{0}; \mathbf{0}; \mathbf{0}]]$ to determine its champion (the Western Conference did the same, and then the two conference champions played each other in the NBA finals). However, in 2020, after a much of the NBA regular season was cut short due to Covid, there was a feeling that the regular season wasn't as accurate a measure as it usually is. So the playoffs were expanded slightly: if the 8th and 9th place teams were close enough in record, the playoff for that conference expanded to $[[\mathbf{2}; \mathbf{7}; \mathbf{0}; \mathbf{0}; \mathbf{0}]]$, allowing both teams in [41]. After the success of that system, the playoffs were expanded further starting in 2021 to the following nonlinear multibracket [34].

**Figure 3.7.10: 2022 NBA Eastern Conference Playoffs**

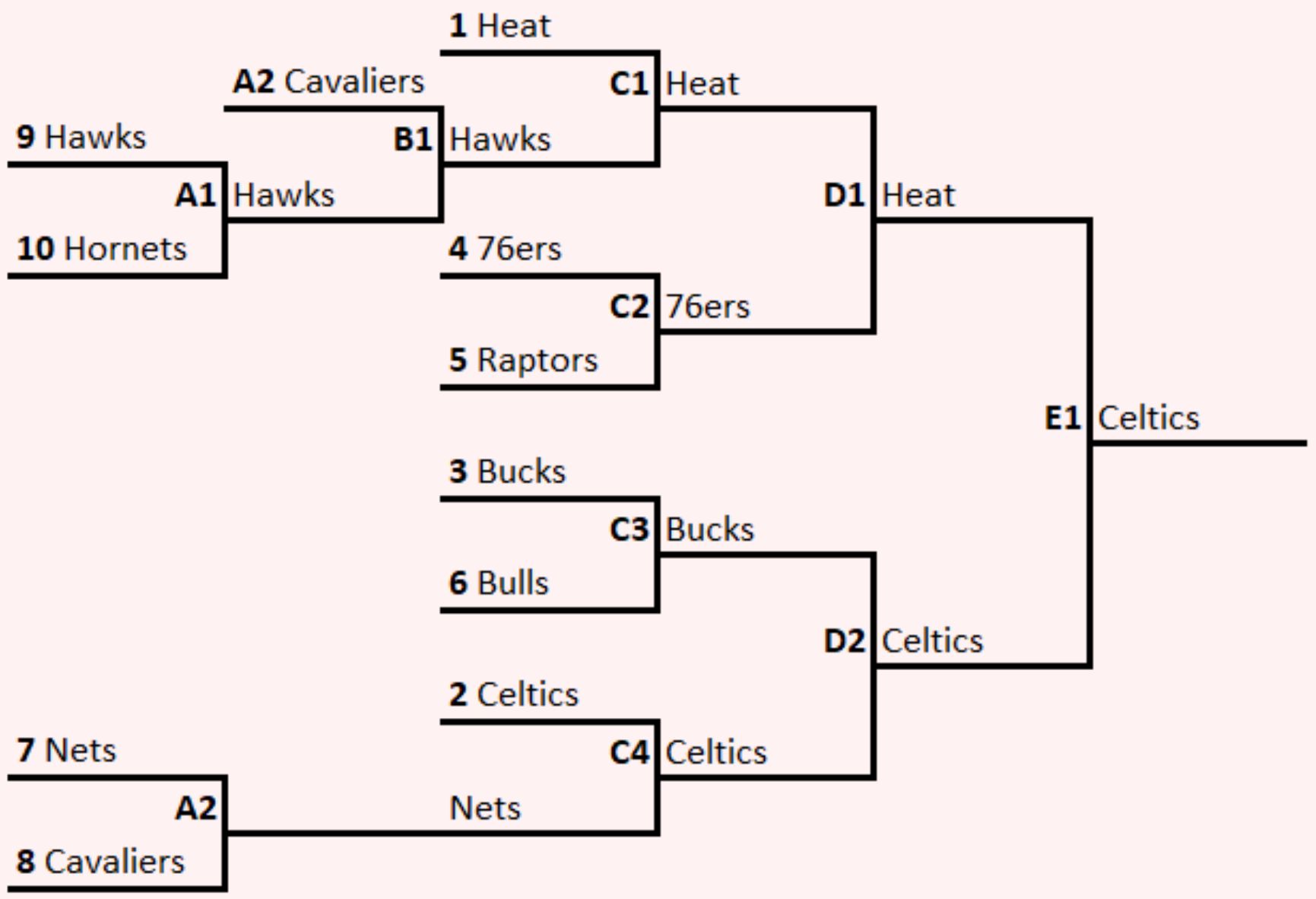

The first two rounds of the new NBA playoffs are similar in structure to the Page-McIntyre system: two lower-seeded teams play each other and two higher-seeded teams play each other, and then the winner of the first game plays the loser of the second. But because the two qualifying teams get dumped into a larger eight-team bracket, rather than facing off immediately, the issues of the original Page-McIntyre system are avoided.

A final nice example of nonlinear multibrackets is bitonic sort. Bitonic sort was developed by Batcher [3] as a networked sorting algorithm with low delay (the sorting theory equivalent to a low number of rounds). As every sorting algorithm can be transformed into a tournament format, and every networked sorting algorithm can be transformed into a multibracket, we can construct an nonlinear multibracket that executes Batcher's bitonic sort.

**Definition 3.7.11: Bitonic Sort** **(Batcher, 1968)**

The *bitonic sort* on $2^r$ teams proceeds by diving the teams into two groups of $2^{r-1}$ teams, recursively running the bitonic sort on $2^{r-1}$ teams on each group, and then running the standard swiss format $\mathcal{S}_r$ on the full group of $2^r$ teams, with one of the groups getting the odd seeds in $\mathcal{S}_r$ and the other group getting the even seeds.

The 8-team bitonic sort is displayed in Figure 3.7.12. The **A**-, **B**-, and **C**- round games facilitate the running of two parallel instantiations of the 4-team bitonic sort, while the **D**-, **E**-, and **F**−round games carry out $\mathcal{S}_3$.

**Figure 3.7.12: 8-Team Bitonic Sort**

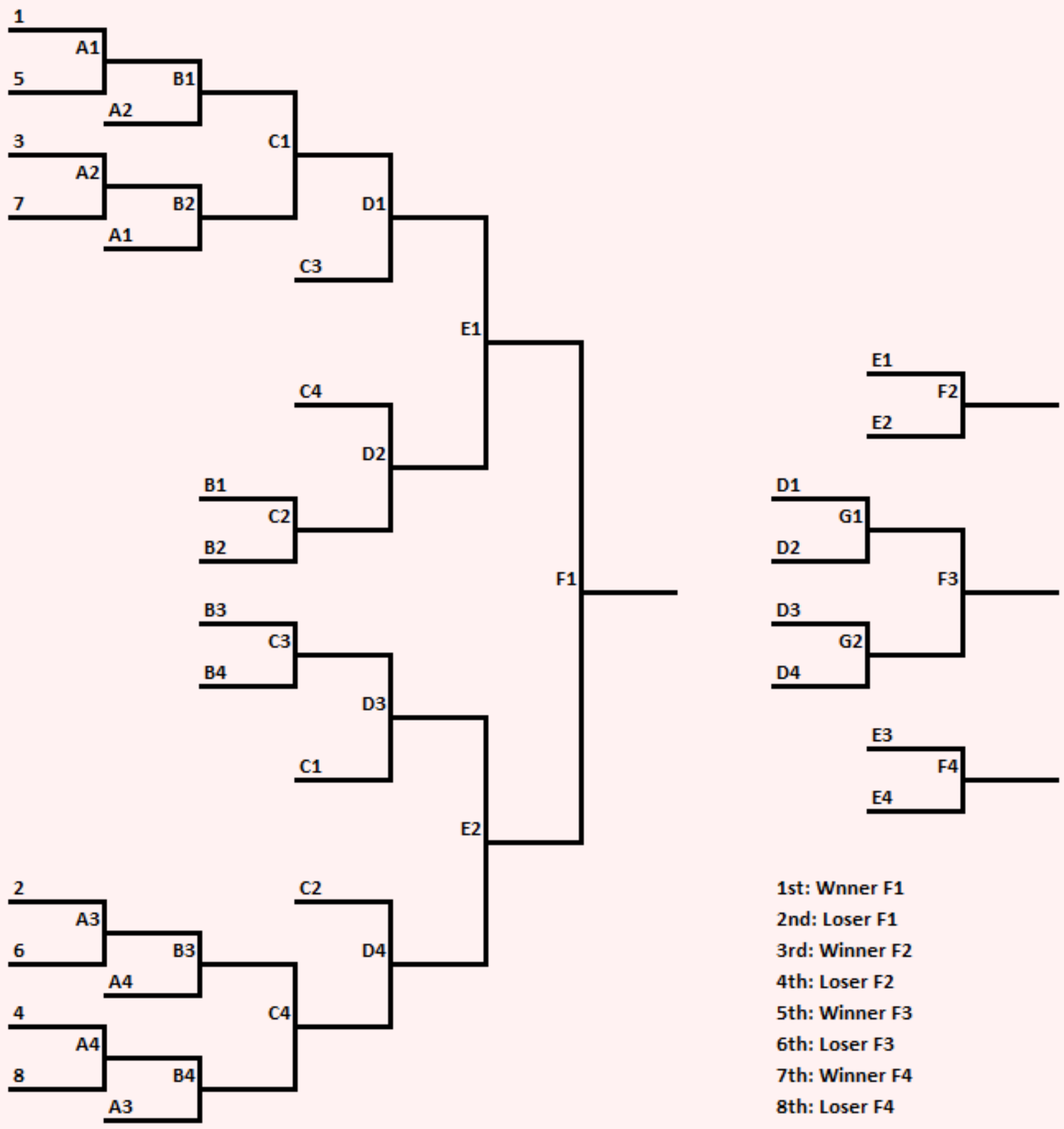

We leave it to the reader to verify that bitonic sort is in fact a sorting algorithm: that is, if the matchup table is SST with all win probabilities being 0 or 1 (even if the teams are not seeded in the correct order initially!), bitonic sort will correctly sort the teams. Impressively, the 8-team bitonic sort does this in only six rounds: no team needs to play every other team in order to complete the sort.

# 4 Postmatter

## 4.1 Future Work

There are four areas of research in the field of tournament design that we hope to pursue in the future.

The first are answers to the two open questions presented at the end of Chapter 2: for all $r$, does there exist an $r$-round deterministic ordered balanced knockout tournament? And for all $r$, does there exist an $r$-round dramatic ordered balanced knockout tournament? Of course, finding such formats would be ideal, but we think an impossibility theorem is more likely.

The second is a continuation of the study of multibrackets. Some areas in particular include: defining a full range of the degrees of properness and respectfulness that linear multibrackets can exhibit, counting the number of multibrackets that are efficient with respect to a given prize structure, determining which Swiss signature to use to select a top-$m$ out of $n$ teams for arbitrary $m$ and $n$, and extending the notions of signatures and properness to nonlinear multibrackets.

The third is an expansion of our analysis to formats that are not networked: round robins and pool-based formats in particular. Additionally, we are interested in a treatment of the space of formats as a whole, working from the top down by defining universal properties and observing which formats uphold them, rather than continuing to define specific formats for analysis.

And the fourth is a statistical model to measure how fair and accurate tournament formats are, allowing us to distinguish between multiple ordered (or unordered) formats. In particular, we are interested in an algorithm that returns the most accurate format that is efficient with respect to a given prize structure for use at club sport regional tournaments.

## 4.2 Glossary of Terms

**Definition 2.1.15: Balanced Bracket** **(Unattributed, p. 15)**

A bracket is *balanced* if none of the participating teams have byes.

**Definition 3.7.11: Bitonic Sort** **(Batcher, 1968, p. 104)**

The *bitonic sort* on $2^r$ teams proceeds by diving the teams into two groups of $2^{r-1}$ teams, recursively running the bitonic sort on $2^{r-1}$ teams on each group, and then running the standard swiss format $\mathcal{S}_r$ on the full group of $2^r$ teams, with one of the groups getting the odd seeds in $\mathcal{S}_r$ and the other group getting the even seeds.

**Definition 2.1.1: Bracket** **(Unattributed, p. 11)**

A *bracket* is a networked format in which

(a) Teams don't play any games after their first loss, and

(b) Games are played until only one team has no losses, and that team is crowned champion.

**Definition 2.1.11: Bracket Signature** **(Fried, 2024, p. 13)**

The *signature* of an $r$-round bracket $\mathcal{A}$ is the list $[[\mathbf{a_0}; \ldots; \mathbf{a_r}]]$ where $a_i$ is the number of teams that get $i$ byes.

**Definition 2.1.8: Bye** **(Unattributed, p. 12)**

A team has a *bye* in round $r$ if it plays no games in round $r$ or before.

**Definition 2.2.6: Chalk** **(Unattributed, p. 18)**

A tournament went *chalk* if the higher-seeded team won every game during the tournament.

**Definition 2.6.8: Cohort Randomized Knockout Tournament** **(Schwenk, 2000, p. 49)**

The $r$-round *cohort randomized knockout tournament* is the balanced bracket on $2^r$ teams, except, for each $i$, seeds $2^i + 1$ through $2^{i+1}$ are shuffled randomly before play.

**Definition 3.5.12: Compact Swiss Signature** **(Fried, 2024, p. 83)**

A Swiss signature is *compact* if only one teams finishes undefeated.

**Definition 3.2.11: Competitive Semibracket** **(Fried, 2024, p. 66)**

A semibracket is *competitive* if no teams are declared co-champion without winning at least one game. Equivalently, a semibracket is competitive if its signature ends in a 0.

**Definition 3.1.3: Consolation Bracket** **(Unattributed, p. 56)**

A *consolation bracket* is a bracket in which teams that did not win the tournament compete for an $m$th-place finish for some $m$.

**Definition 2.4.4: Containment** **(Fried, 2024, p. 32)**

Let $\mathcal{A}$ and $\mathcal{B}$ be bracket signatures. $\mathcal{A}$ *contains* $\mathcal{B}$ if there exists some $i$ such that

(a) At least as many games are played in the $(i+1)$th round of $\mathcal{A}$ as in the first round of $\mathcal{B}$, and

(b) For $1 < j \leq r$ where $r$ is the number of rounds in $\mathcal{B}$, there are exactly as many games played in the $(i+j)$th round of $\mathcal{A}$ as in the $j$th round of $\mathcal{B}$.

**Lemma 2.4.6: The Containment Lemma** **(Fried, 2024, p. 33)**

If $\mathcal{A}$ contains $\mathcal{B}$, and $\mathcal{B}$ is not ordered, then neither is $\mathcal{A}$.

**Definition 3.1.6: Contingent Game** **(Fried, 2024, p. 58)**

A game in a tournament format is *contingent* if, under certain circumstances, (most commonly if the teams have already played earlier in the tournament) the game is skipped and the result of a previous game is used.

**Definition 1.3.10: Deterministic Tournament Format** **(Unattributed, p. 9)**

A tournament format is *deterministic* if it employs no randomness other than the randomness inherent in the gameplay function $g$.

**Definition 3.7.4: Double-Elimination Tournament** **(Unattributed, p. 99)**

A *double-elimination* tournament is a multibracket (plus one contingent game) consisting of a *winners' bracket*, where every team starts, a *losers' bracket*, that every winners' bracket loser falls into, and a *grand finals*, in which the winner of the winners' bracket plays the winner of the losers' bracket for the championship, with the losers' bracket winner needing to win twice, while the winners' bracket winner only needs to win once.

**Definition 2.2.10: Dramatic Bracket** **(Fried, 2024, p. 19)**

A bracket is *dramatic* if, as long as the bracket goes chalk, in every round, the $m$ remaining teams are the top $m$ seeds.

**Definition 2.6.7: Dramatic Knockout Tournament** **(Fried, 2024, p. 49)**

A knockout tournament is *dramatic* if, as long as the knockout tournament goes chalk, in every round, the $m$ remaining teams are guaranteed to be the top $m$ seeds.

**Theorem 2.4.7: Edwards's Theorem** **(Edwards, 1991, p. 33)**

The set of ordered brackets is exactly the set of proper brackets whose signature is formed by the following process:

1. Start with the list $[[\mathbf{0}]]$ (note that this not yet a bracket signature).
2. As many times as desired, prepend the list with $[[\mathbf{1}]]$ or $[[\mathbf{3}; \mathbf{0}]]$.
3. Then, add 1 to the first element in the list, turning it into a bracket signature.

**Definition 3.6.7: Efficient** **(Fried, 2024, p. 94)**

A respectful linear multibracket is *efficient* if one of three conditions hold:

(a) It is a single trivial semibracket.

(b) It is a sequence of competitive semibrackets.

(c) It is a single nontrivial noncompetitive semibracket followed by a sequence of competitive semibrackets.

**Definition 3.6.10: Efficient with Respect to a Prize Struture** **(Fried, 2024, p. 96)**

A respectful linear multibracket $\mathcal{A} = \mathcal{A}_1 \to \ldots \to \mathcal{A}_k$ is *efficient with respect to a prize structure* $\mathcal{P} = (p_1, \ldots, p_m)$ if

(a)

$$\sum_{i=1}^{k} \operatorname{Rank}(\mathcal{A}_i) = \sum_{i=1}^{m} p_i,$$

(b) $\mathcal{A}_j$ being noncompetitive implies that for some $\ell < m$,

$$\sum_{i=1}^{j-1} \operatorname{Rank}(\mathcal{A}_i) = \sum_{i=1}^{\ell} p_i,$$

and

(c) $\mathcal{A}_j$ being trivial implies that for some $\ell \leq m$,

$$\sum_{i=1}^{j} \operatorname{Rank}(\mathcal{A}_i) = \sum_{i=1}^{\ell} p_i.$$

**Definition 3.4.5: Flowchart** **(Fried, 2024, p. 74)**

The *flowchart* of a linear multibracket that consists of $k$ semibrackets is a directed graph in which the nodes are arranged into rows, where

(a) There is a node for each team, each round of each semibracket, and each place a team could finish in, plus one additional node representing elimination.

(b) The zeroth row has the nodes representing each team, arranged from lowest seed to highest seed.

(c) The $i$th row for $1 \leq i \leq k$ has the nodes representing the rounds of the $i$th semibracket, arranged in order, plus the node representing the place a team gets for winning the $i$th semibracket.

(d) The final row row has only the node representing elimination.

(e) There is an arrow from each team to the round where that teams plays its first game.

(f) For each round $\mathbf{R}$, there is an arrow (or arrows) from $\mathbf{R}$ to the round(s) where $\mathbf{R}$-round losers go.

**Theorem 2.2.16: The Fundamental Theorem of Brackets** **(Fried, 2024, p. 21)**

There is exactly one proper bracket with each bracket signature.

**Definition 1.3.1: Gameplay Function** **(Unattributed, p. 7)**

A *gameplay function* $g$ on $\mathcal{T}$ is a nondeterministic function $g : \mathcal{T} \times \mathcal{T} \to \mathcal{T}$ with the following properties:

(a) $\mathbb{P}[g(t_i, t_j) = t_i] + \mathbb{P}[g(t_i, t_j) = t_j] = 1$.

(b) $\mathbb{P}[g(t_i, t_j) = t_i] = \mathbb{P}[g(t_j, t_i) = t_i]$.

**Definition 2.2.3: Higher and Lower Seeds** **(Unattributed, p. 17)**

Somewhat confusingly, convention is that smaller numbers are the *higher seeds*, and bigger numbers are the *lower seeds*.

**Definition 2.2.2: $i$-seed** **(Unattributed, p. 17)**

In a list of teams $\mathcal{T} = [t_1, ..., t_n]$, we refer to $t_i$ as the $i$-seed.

**Definition 2.5.1: Knockout Tournament** **(Unattributed, p. 35)**

A *knockout tournament* is a tournament in that is played over a series of rounds subject to the following constraints:

(a) Each team plays at most one game in each round.

(b) If a team loses in a round, they don't play any games in later rounds.

(c) If a team wins in a round, they play a game in the next round.

(d) Exactly one team finishes undefeated, and that team is crowned champion.

**Definition 2.5.2: Knockout Tournament Signature** **(Fried, 2024, p. 35)**

The *signature* of an $r$-round knockout tournament $\mathcal{A}$ is the list $[[\mathbf{a_0}; ...; \mathbf{a_r}]]$ where $a_i$ is the number of teams that get $i$ byes.

**Definition 3.3.1: Linear Multibracket** **(Fried, 2024, p. 67)**

A *linear multibracket* is a multibracket that can be arranged into a sequence of semibrackets such that

(a) If a team loses in a given semibracket but is not eliminated, they are sent to a later semibracket, and

(b) Each team that wins the $i$th semibracket finishes in $m$th place, where $m$ is the

sum of the ranks of the first $i$ semibrackets.

**Definition 3.3.2: Linear Multibracket Signatures** **(Fried, 2024, p. 67)**

The *signature* of a linear multibracket that consists of semibrackets with signature $\mathcal{A}_1, ..., \mathcal{A}_k$ is $\mathcal{A}_1 \to ... \to \mathcal{A}_k$.

**Definition 1.3.4: Matchup Table** **(Unattributed, p. 7)**

The *matchup table* implied by a gameplay function $g$ on a list of teams $\mathcal{T}$ of length $n$ is an $n$-by-$n$ matrix $\mathcal{M}$ such that $\mathcal{M}_{ij} = p_{ij}$.

**Definition 2.3.5: Monotonic** **(Unattributed, p. 25)**

A tournament format $\mathcal{A}$ is monotonic with respect to a list of teams $\mathcal{T}$ if, for all $i < j$, $\mathbb{W}_{\mathcal{A}}(t_i, \mathcal{T}) \geq \mathbb{W}_{\mathcal{A}}(t_j, \mathcal{T})$.

**Definition 3.1.10: Multibracket** **(Fried, 2024, p. 61)**

A *multibracket* is a networked tournament format.

**Definition 3.5.14: $m$x$\mathcal{A}$** **(Fried, 2024, p. 84)**

If $m \in \mathbb{N}$ and $\mathcal{A}$ is a multibracket signature, then $m$x$\mathcal{A}$ is the multibracket signature formed by multiplying every number in every signature in $\mathcal{A}$ by $m$.

**Definition 1.3.11: Networked Tournament Format** **(Armstrong, Nelson, and O'Connor, 1957, p. 9)**

A tournament format is *networked* if it is deterministic, and after each game between $t_i$ and $t_j$, the rest of the format is identical no matter which team won, except that $t_i$ and $t_j$ are swapped.

**Definition 2.3.6: Ordered** **(Edwards, 1991, p. 25)**

An $n$-team tournament format $\mathcal{A}$ is *ordered* if it is monotonic with respect to every SST list of $n$ teams.

**Definition 1.3.3: $p_{ij}$** **(Unattributed, p. 7)**

$p_{ij} = \mathbb{P}[t_i \text{ beats } t_j]$.

**Definition 1.3.2: Playing, Winning, Losing, and Beating (Unattributed, p. 7)**

When $g$ is queried on input $(t_i, t_j)$ we say that $t_i$ and $t_j$ *played a game*. We say that the team that got output by $g$ *won*, that the team that did not *lost*, and that the winning team *beat* the losing team.

**Definition 3.1.4: Primary Bracket (Unattributed, p. 56)**

A *primary bracket* is a bracket in a multibracket the winner of which is declared champion.

**Definition 3.6.9: Prize Structure (Fried, 2024, p. 95)**

A *prize structure* $\mathcal{P}$ is a sequence $(p_1, ..., p_m)$ indicating that the top $p_1$ teams in a format receive some prize, the next $p_2$ receive some smaller prize, etc. Any teams finishing in place $1 + \sum_{i=1}^{m} p_i$ or worse receive no prize.

**Definition 2.2.8: Proper Bracket (Fried, 2024, p. 19)**

A bracket is *proper* if its seeding is proper.

**Definition 2.5.3: Proper Knockout Tournament (Fried, 2024, p. 35)**

A knockout tournament is *proper* if, as long as the tournament goes chalk, in every round it is better to be a higher-seeded team than a lower-seeded one, where:

(a) It is better to have a bye than to play a game.

(b) It is better to play a lower seed than to play a higher seed.

**Definition 3.4.1: Proper Linear Multibracket (Fried, 2024, p. 72)**

A linear multibracket is *proper* if, as long as the bracket goes chalk, in every round of every semibracket it is better to be a higher-seeded team than a lower-seeded one, where:

(a) It is best to have already won an earlier semibracket.

(b) If you have not yet won an earlier semibracket, it is to better to be competing in the current semibracket than to not.

(c) If you are competing in a semibracket, it is better to have a bye in the current round than to not.

(d) If you are playing a game, it is better to play a lower seed than to play a higher seed.

**Definition 2.2.7: Proper Seeding** **(Fried, 2024, p. 19)**

A seeding of a bracket is *proper* if, as long as the bracket goes chalk, in every round it is better to be a higher-seeded team than a lower-seeded one, where:

(a) It is better to have a bye than to play a game.

(b) It is better to play a lower seed than to play a higher seed.

**Definition 3.2.5: Rank of a Semibracket** **(Fried, 2024, p. 64)**

The *rank* of a semibracket is how many co-champions it crowns. If the semibracket $\mathcal{A}$ has rank $m$, we say $\text{Rank}(\mathcal{A}) = m$ or that $\mathcal{A}$ *ranks* $m$ *teams*.

**Definition 3.7.5: Recharge Game** **(Dabney, 2017, p. 99)**

A *recharge game* is a contingent game in a multibracket that is a rematch of a previous game and played only if the lower team won the first game.

**Definition 2.5.4: Reseeded Bracket** **(Hwang, 1982, p. 36)**

A *reseededed bracket* is a knockout tournament in which, after each round, the highest-seeded team playing that round is matched up with the lowest-seeded team playing that round, second-highest vs second-lowest, etc.

**Definition 3.4.9: Respectful Linear Multibracket** **(Fried, 2024, p. 77)**

A linear multibracket is *respectful* if its flowchart has no arrow crossings and every node in its flowchart has at most one arrow coming out of it.

**Definition 2.1.7: Round** **(Unattributed, p. 12)**

A *round* is a set of games such that the winners of each of those games have the same number of games remaining to win the tournament.

**Definition 3.5.7: $r$-Round Swiss Signature** **(Fried, 2024, p. 81)**

An $r$*-round Swiss signature* is a Swiss signature in which each team plays $r$ games.

**Definition 2.2.1: Seeding** **(Unattributed, p. 17)**

The *seeding* of an $n$-team bracket is the arrangement of the numbers 1 through $n$ on the starting lines of a bracket.

**Definition 3.2.2: Semibracket** **(Fried, 2024, p. 62)**

A *semibracket* is a networked format in which

(a) Teams don't play any games after their first loss, and

(b) All teams that finish with no losses are declared co-champions.

**Definition 3.2.6: Semibracket Signature** **(Fried, 2024, p. 64)**

The *signature* of an $r$-round semibracket $\mathcal{A}$ is the list $[[\mathbf{a_0}; \ldots; \mathbf{a_r}]]_m$, where $a_i$ is the number of teams that get $i$ byes and $m = \text{Rank}(\mathcal{A})$. (In the case where $m = \text{Rank}(\mathcal{A}) = 1$, it can be omitted.)

**Definition 2.1.9: Shape** **(Unattributed, p. 13)**

The *shape* of a bracket is the tree that underlies it.

**Definition 3.5.8: Standard Swiss Signature ($\mathcal{S}_r$)** **(Fried, 2024, p. 81)**

$\mathcal{S}_r$, or the *standard $r$-round Swiss signature*, is the multibracket signature defined recursively by

$$\mathcal{S}_0 = [[\mathbf{1}]],$$

and

$$\mathcal{S}_r = [[\mathbf{2^r}; \ldots; \mathbf{0}]] \to \mathcal{S}_0 \to \mathcal{S}_1 \to \ldots \to \mathcal{S}_i \to \ldots \to \mathcal{S}_{r-1}.$$

**Lemma 2.4.2: The Stapling Lemma** **(Fried, 2024, p. 30)**

If $\mathcal{A} = [[\mathbf{a_0}; \ldots; \mathbf{a_r}]]$ and $\mathcal{B} = [[\mathbf{b_0}; \ldots; \mathbf{b_s}]]$ are ordered brackets, then $\mathcal{C} = [[\mathbf{a_0}; \ldots; \mathbf{a_r} + \mathbf{b_0} - \mathbf{1}; \ldots; \mathbf{b_s}]]$ is an ordered bracket as well.

**Definition 2.1.6: Starting Line** **(Unattributed, p. 12)**

A *starting line* is a line in a bracket where a team is placed before it has played any games.

**Definition 2.3.3: Strongly Stochastically Transitive** **(David, 1963, p. 24)**

A list of teams $\mathcal{T}$ is *strongly stochastically transitive* if for each $i, j, k$ such that $j < k$,

$$\mathbb{P}[t_i \text{ beats } t_j] \leq \mathbb{P}[t_i \text{ beats } t_k].$$

**Definition 3.5.5: Swiss Format** **(Unattributed, p. 81)**

A *Swiss format* is a respectful linear multibracket with the following additional properties:

(a) Every team starts in the primary semibracket.

(b) Every game is between two teams with the same record.

(c) Every team plays the same number of games.

(d) Every team wins a semibracket.

(e) Every semibracket is either trivial or competitive.

**Definition 3.5.3: Swisschart** **(Fried, 2024, p. 80)**

A *Swisschart* is the flowchart of a Swiss format except with the arrows, nodes representing the teams, and the node representing elimination removed, and with the labels on the remaining nodes replaced with the record of the teams in that node.

**Definition 3.5.6: Swiss Signature** **(Fried, 2024, p. 81)**

A *Swiss signature* is a linear multibracket signature that admits a Swiss format.

**Definition 2.6.1: Totally Randomized Knockout Tournament** **(Unattributed, p. 44)**

A *totally randomized knockout tournament* is a bracket, except the teams are randomly placed onto the starting lines instead of being placed according to seed.

**Definition 1.3.8: Tournament Format** **(Unattributed, p. 8)**

A *tournament format* is an algorithm that takes as input a list of teams $\mathcal{T}$ and a gameplay function $g$ and outputs a ranking (potentially including ties) on $\mathcal{T}$.

**Definition 3.2.10: Trivial Semibracket** **(Fried, 2024, p. 66)**

A semibracket is *trivial* if every team is declared co-champion without playing any games. Equivalently, a semibracket is trivial if its signature is of the form $[[\mathbf{m}]]_m$.

**Definition 1.3.9: $\mathbb{W}_{\mathcal{A}}(t, \mathcal{T})$** **(Unattributed, p. 8)**

$\mathbb{W}_{\mathcal{A}}(t, \mathcal{T})$ is the probability that team $t \in \mathcal{T}$ wins tournament format $\mathcal{A}$ when it is run on the list of teams $\mathcal{T}$.

**Definition 3.6.3: Weakly Efficient** **(Fried, 2024, p. 91)**

A multibracket is *weakly efficient* if, once a team is guaranteed to be ranked by the format or guaranteed to finish unranked, they stop playing games.

## 4.3 Glossary of Formats

**Figure 2.1.4: 2024 College Football Playoff** (p. 12)

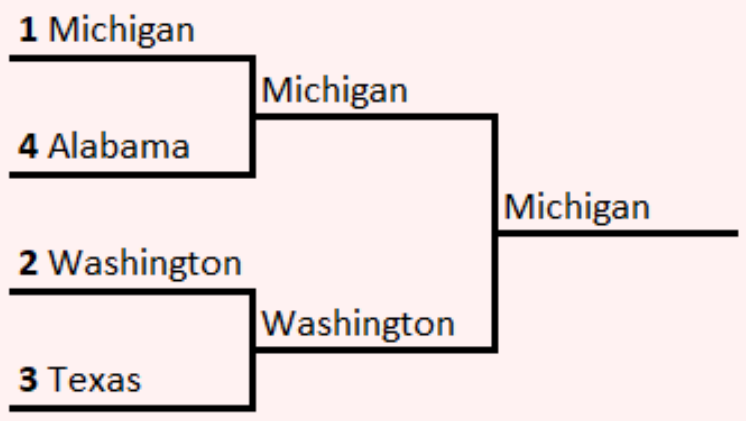

**Figure 2.1.16: 2023 WCC Women's Basketball Tournament** (p. 15)

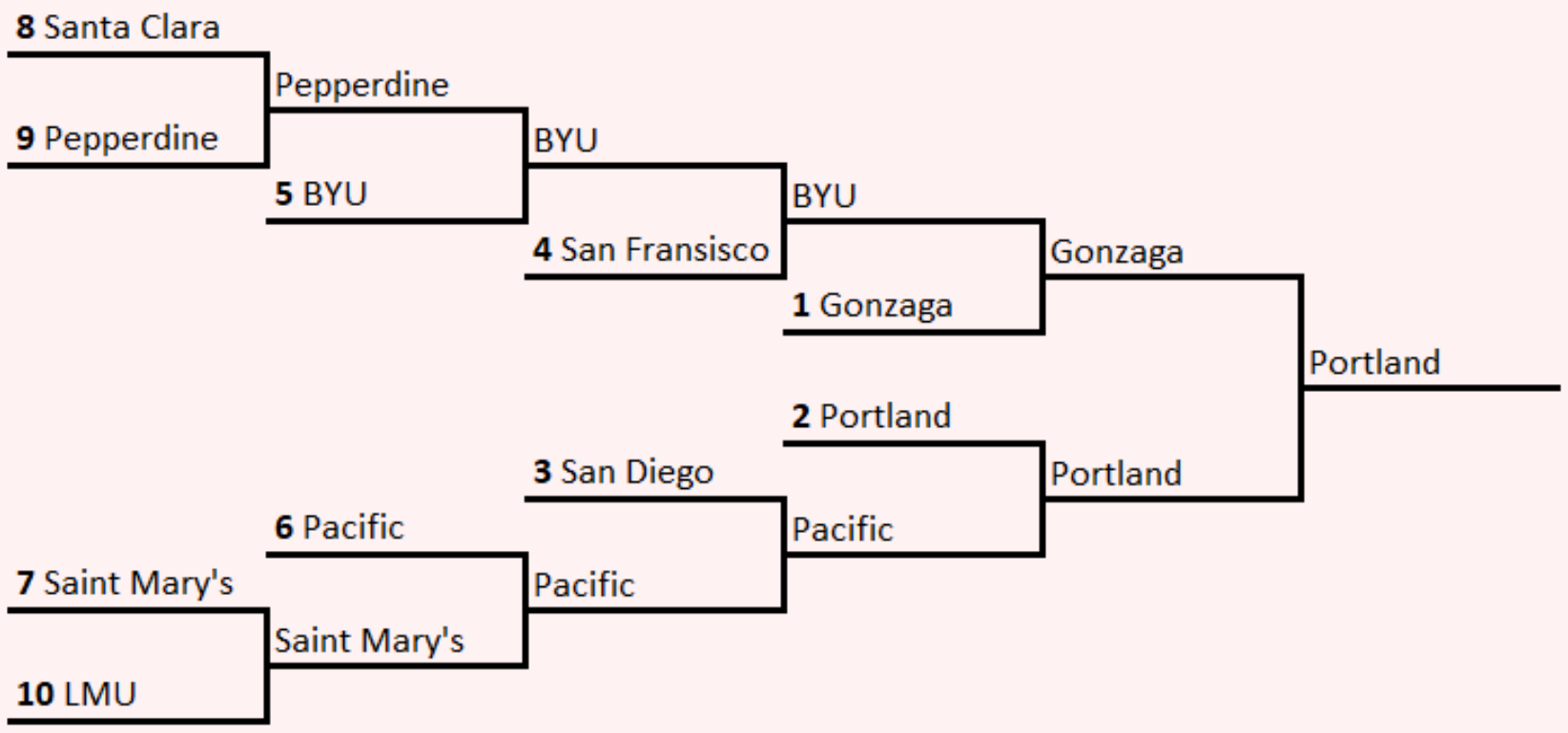

**Figure 2.2.4: 2005 NBA Eastern Conference Playoffs** (p. 17)

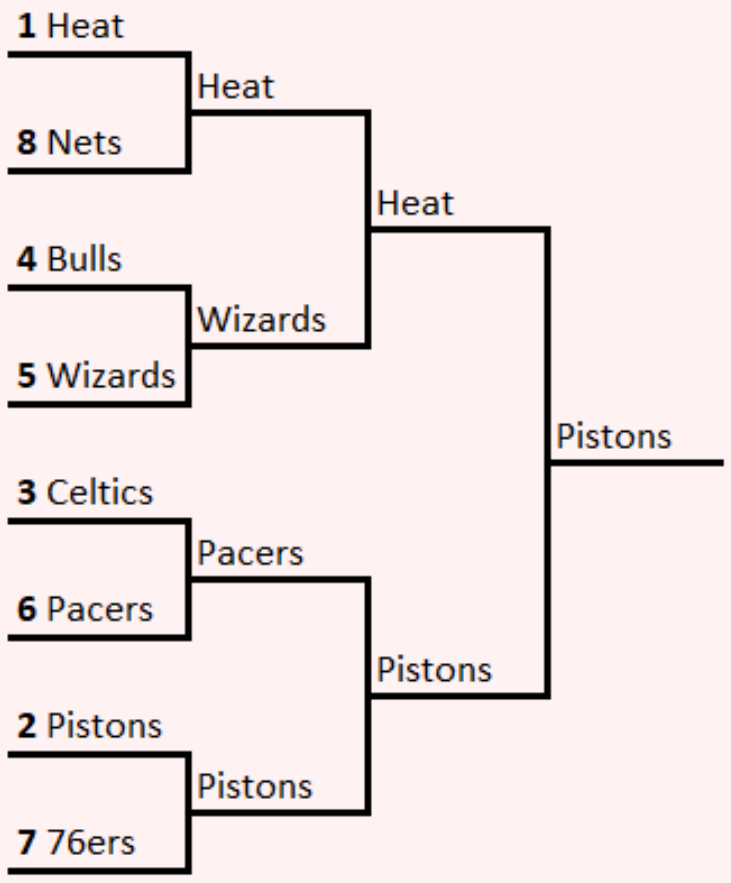

**Figure 2.3.1: 2013 NCAA Men's Basketball Tournament South Region (p. 23)**

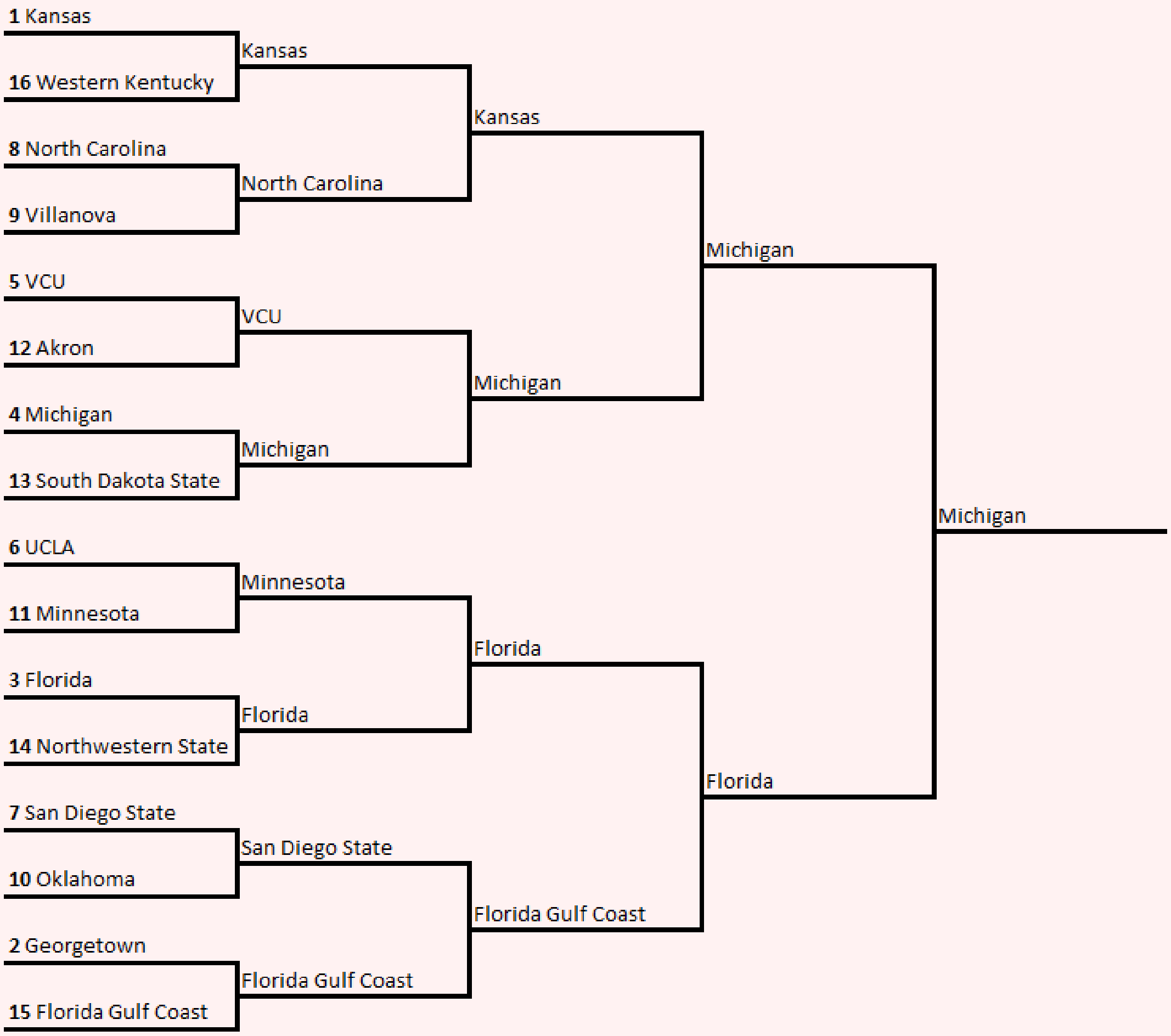

**Figure 2.5.5: 2024 National Football League AFC Playoffs** (p. 36)

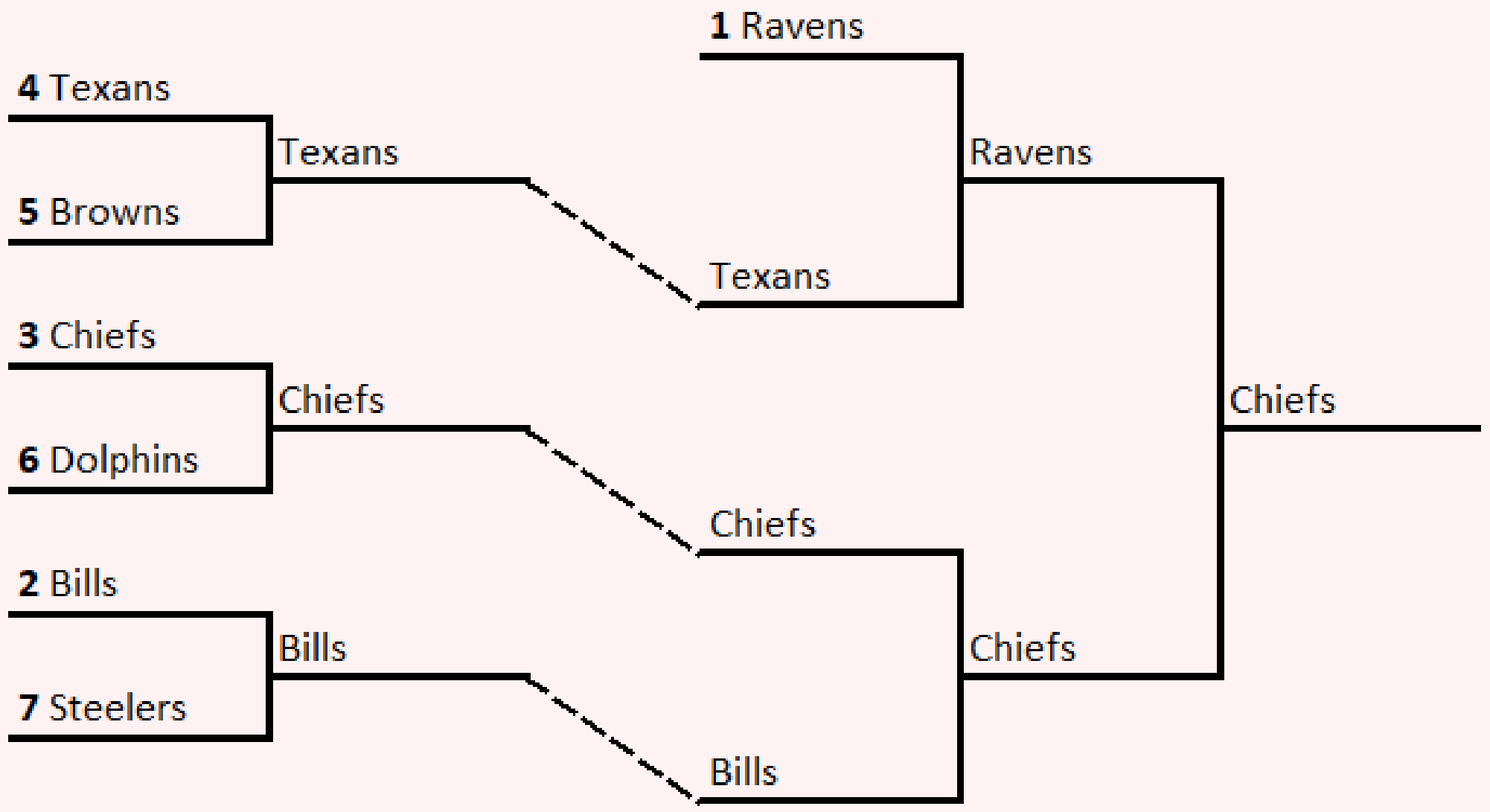

**Figure 2.5.6: 2024 National Football League NFC Playoffs** (p. 37)

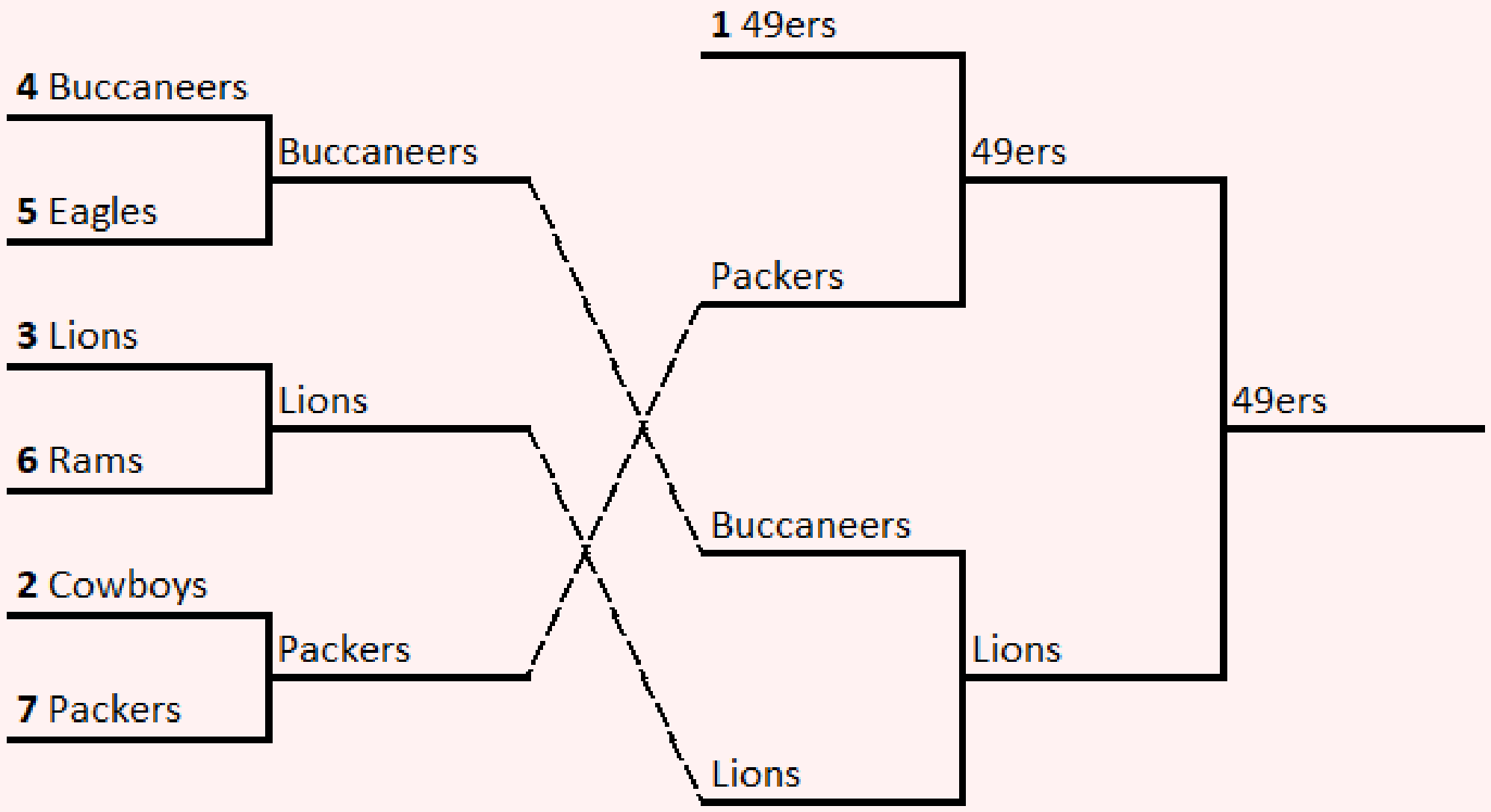

**Figure 3.1.1: 2023 KBO League Playoffs** (p. 54)

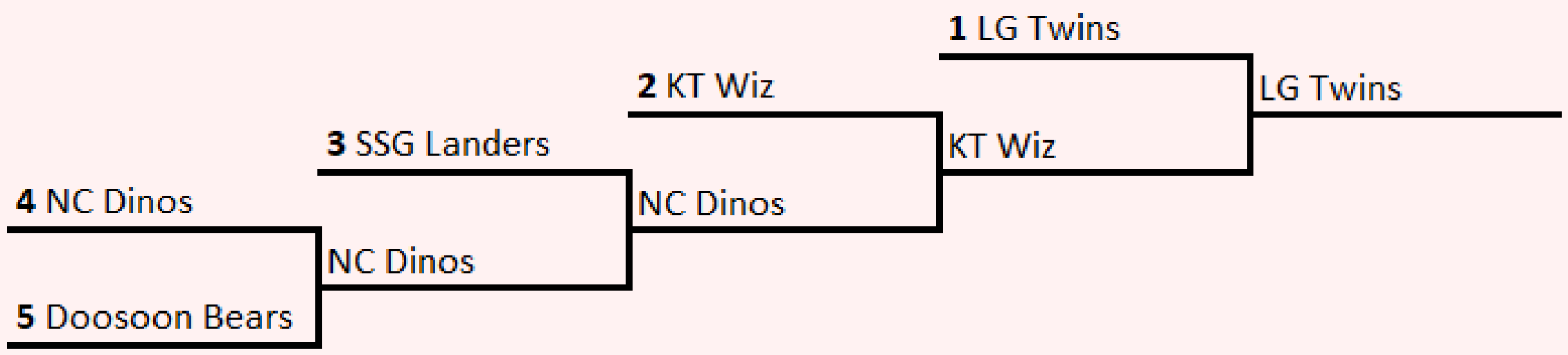

**Figure 3.1.2: 2015 AFC Asian Cup** (p. 55)

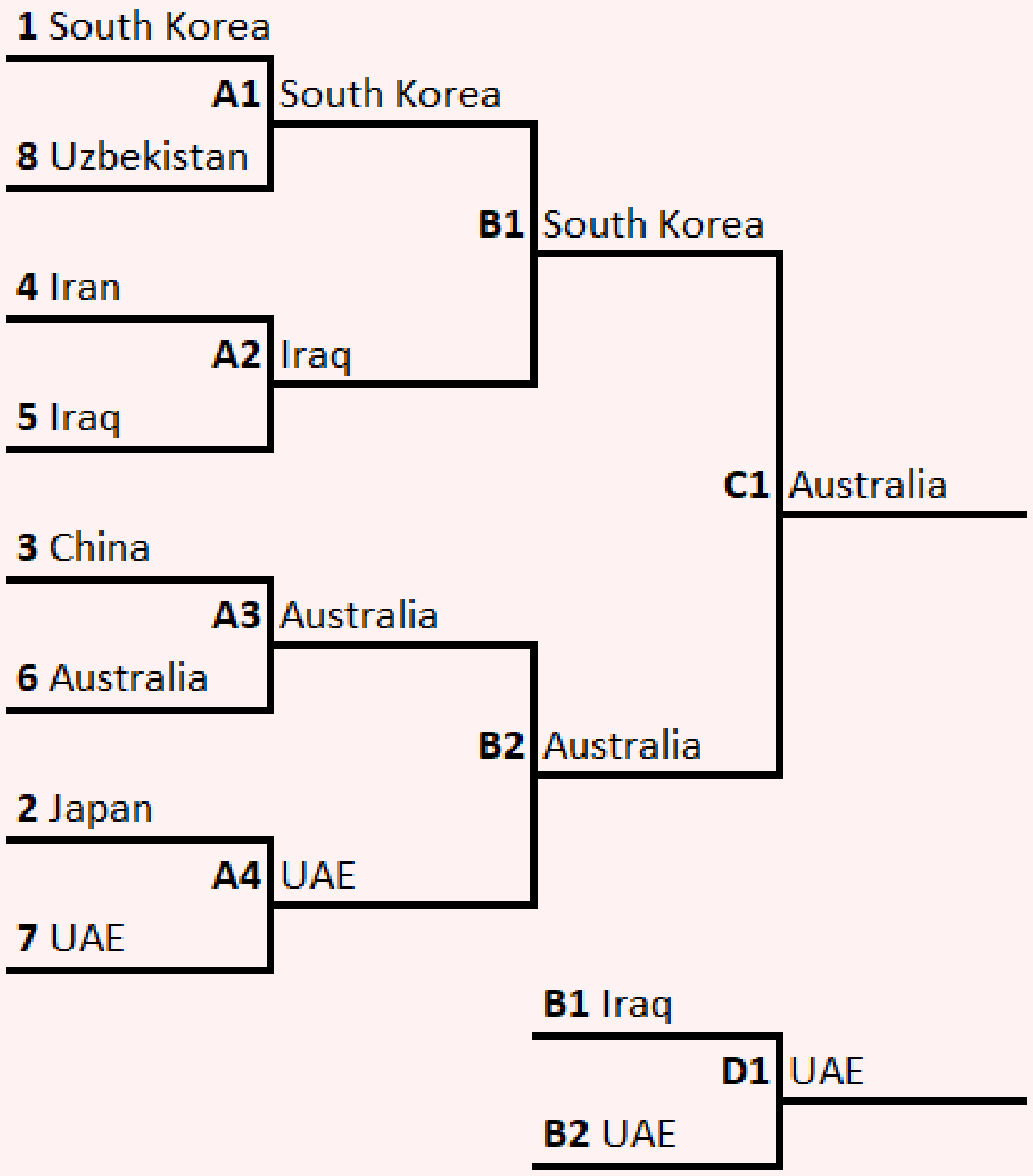

**1st (Winner C1):** Australia
**2nd (Loser C1):** South Korea
**3rd (Winner D1):** UAE

**Figure 3.1.7: 2023 Southern Conference Wrestling Championships** (p. 59)

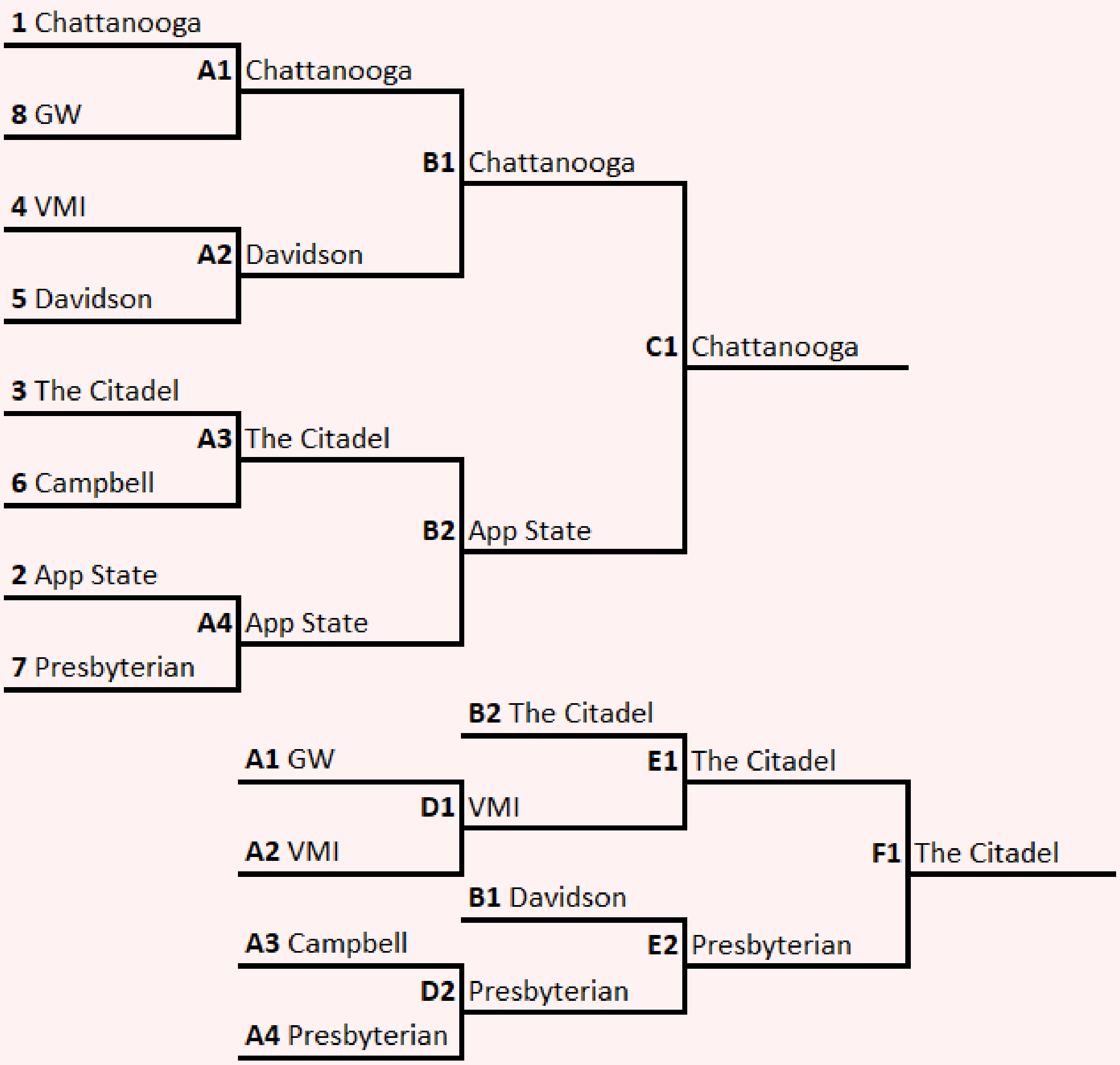

**1st (Winner C1):** Chattanooga
**2nd (Loser C1):** App State
**3rd (Winner F1):** The Citadel
**4th (Loser F1):** Presbyterian

**Figure 3.2.9: 2023 UEFA Champions League Qualifying Phase (p. 65)**

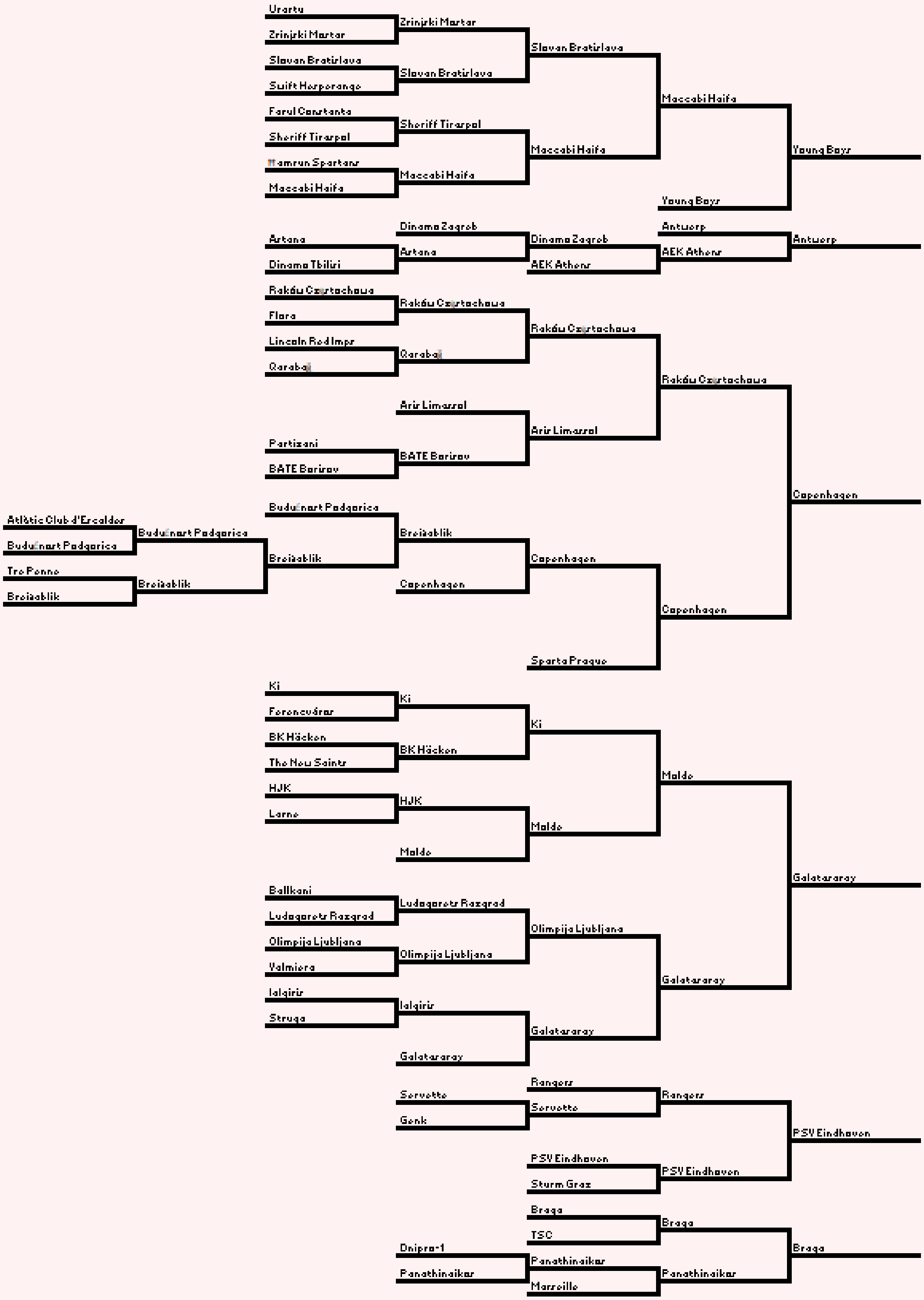

## Figure 3.4.3: 2023 MLQ Championship Play-In Tournament (p. 72)

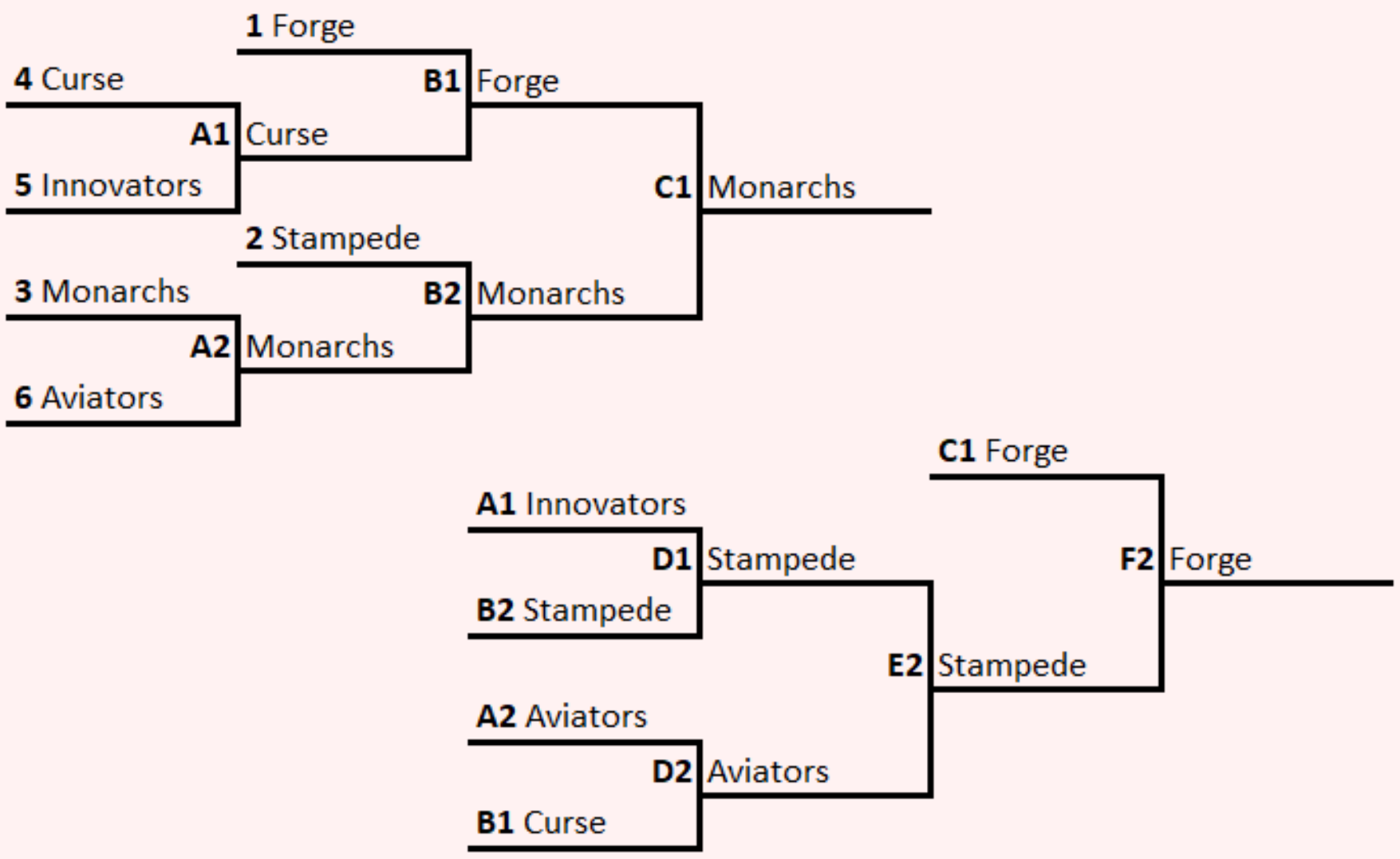

## Figure 3.5.1: 1988 Men's College Basketball Maui Invitational (p. 78)

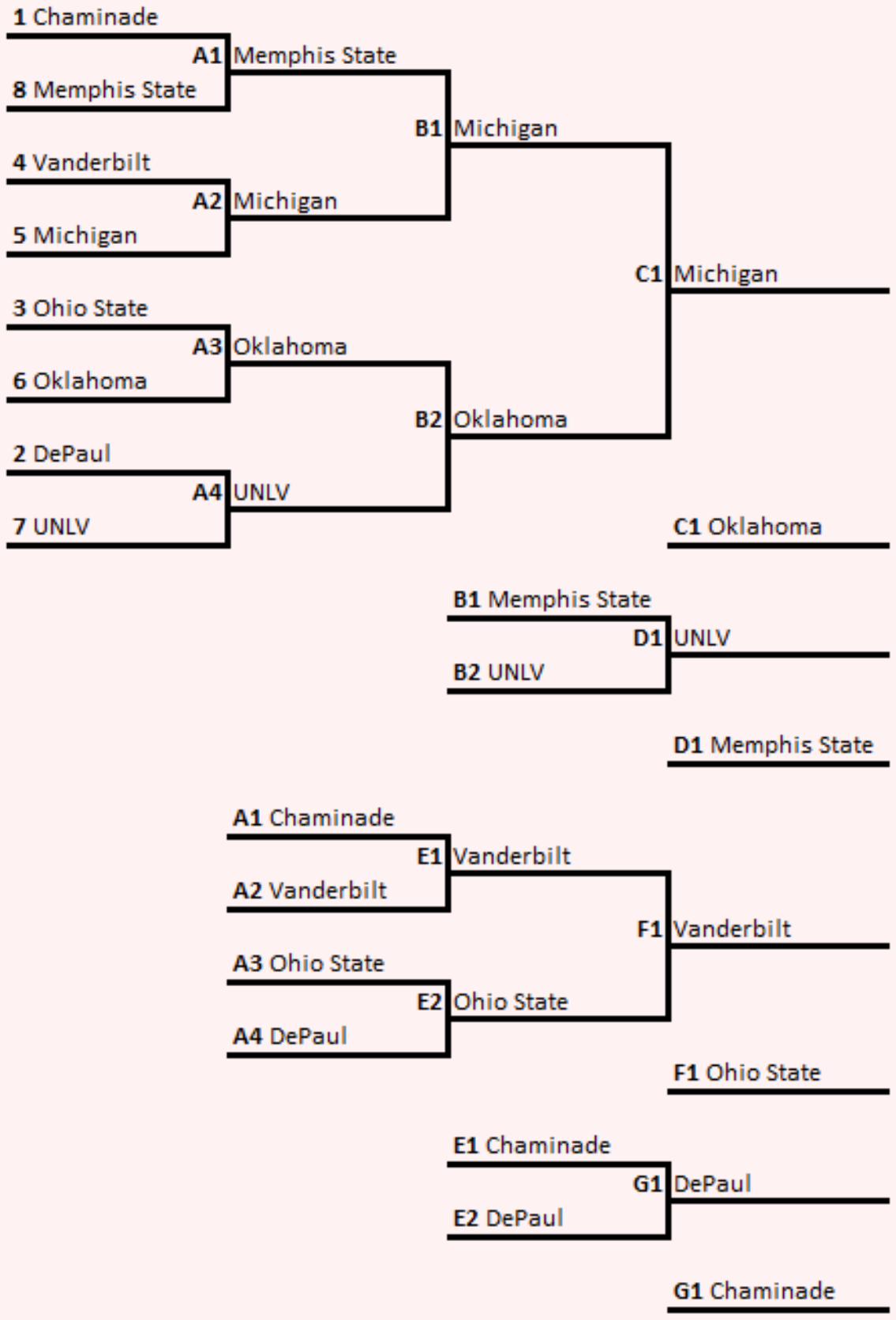

**Figure 3.7.1: 2019 Suncorp Super Netball Playoffs** (p. 97)

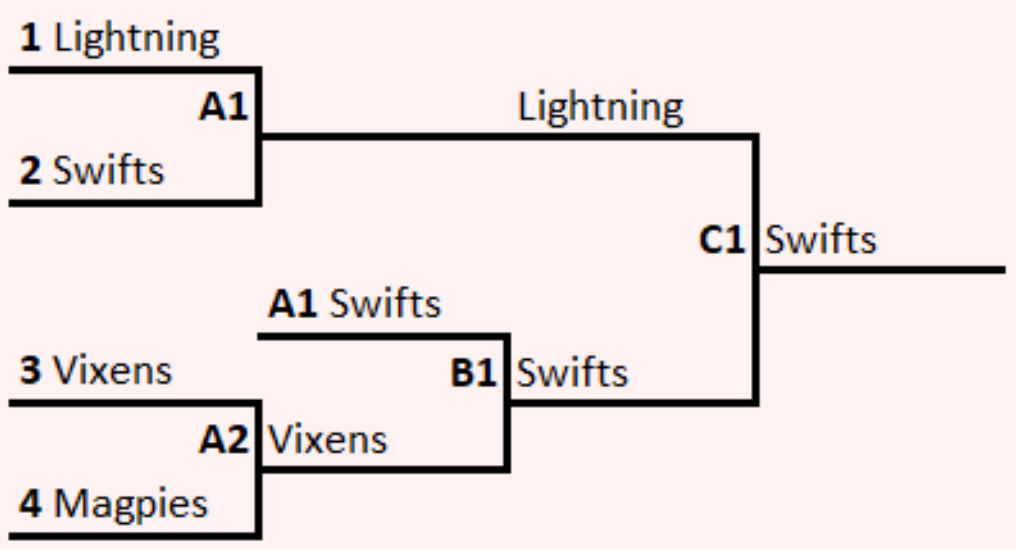

**Figure 3.7.6: 2016 NCAA Softball Ann Arbor Regional** (p. 99)

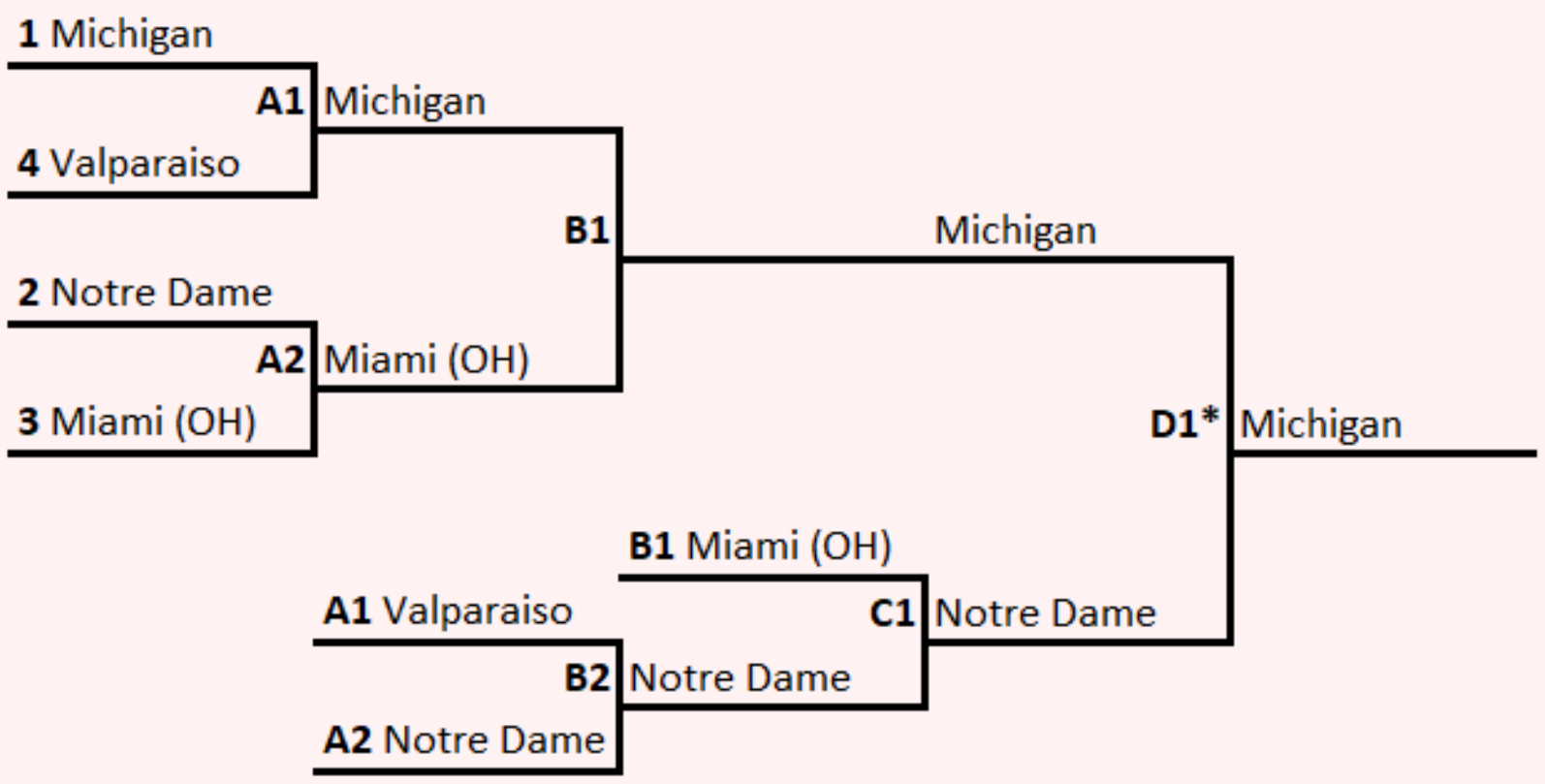

Figure 3.7.9: 2022 Big Ten Baseball Tournament (p. 102)

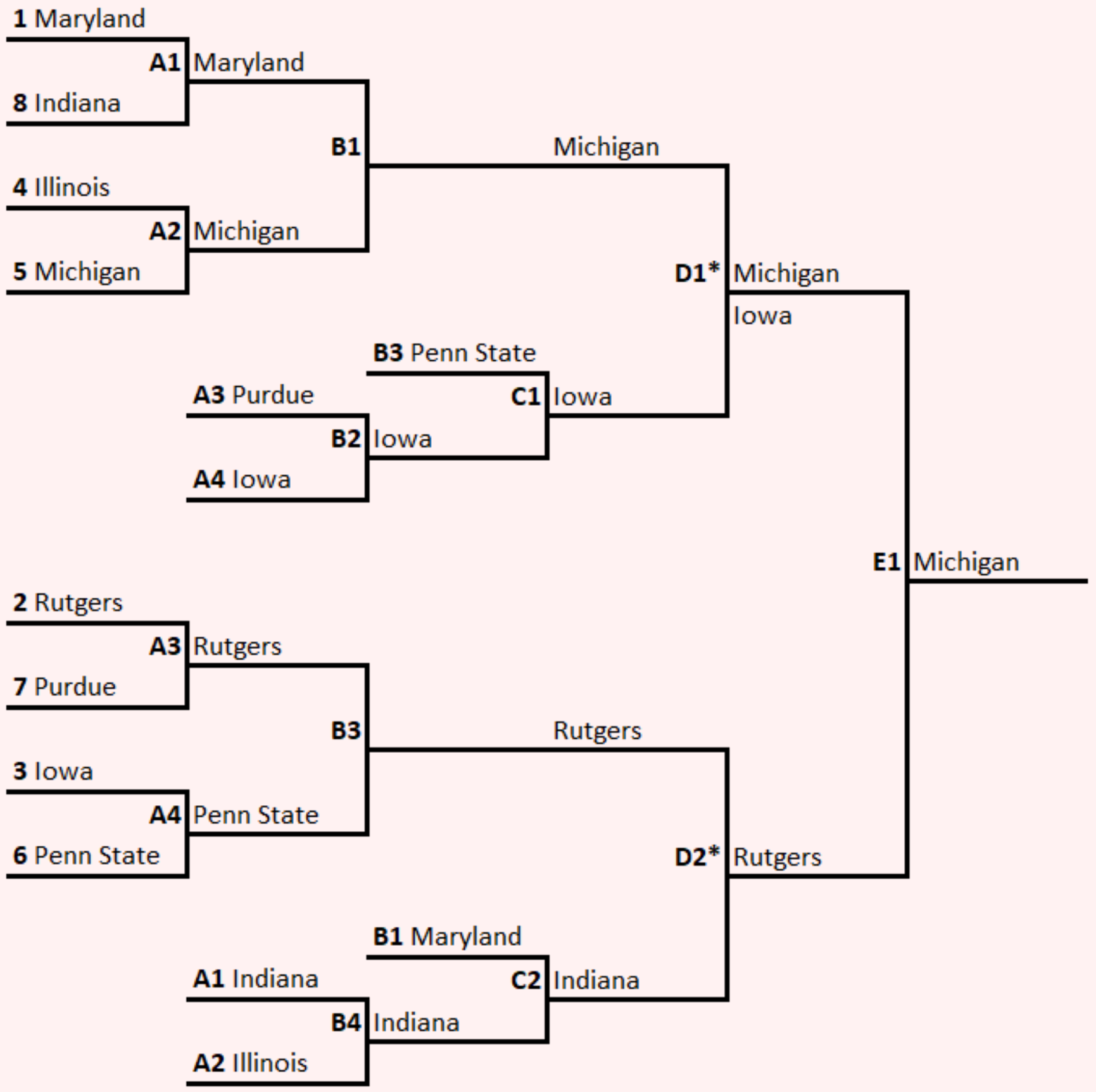

Figure 3.7.10: 2022 NBA Eastern Conference Playoffs (p. 103)

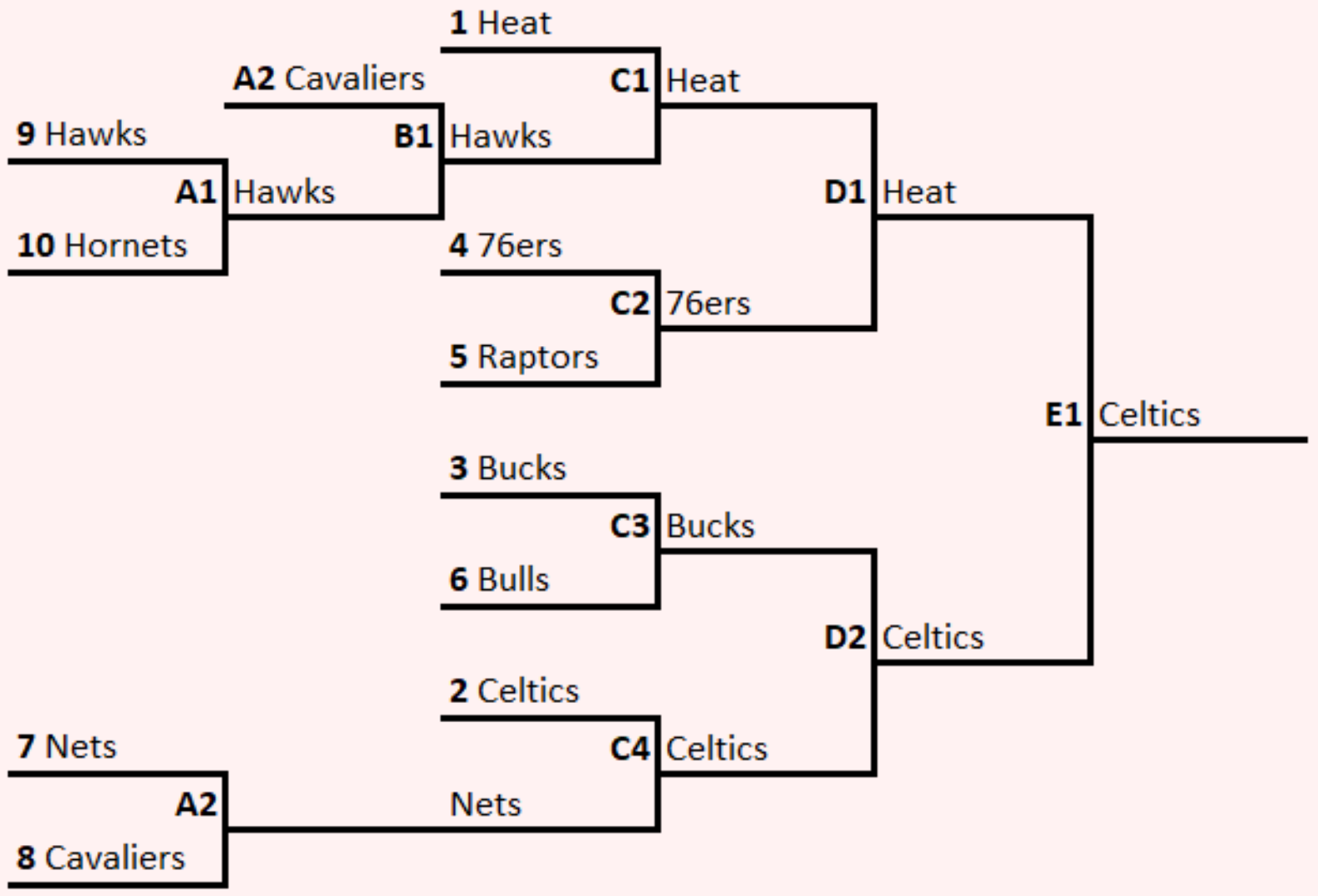

## 4.4 Analysis References

## 4.5 Format References